\documentclass[reqno,oneside]{amsart} 
\pdfoutput=1
\usepackage{amsmath}
\usepackage{amsthm,
	amssymb,amsmath,enumerate,stackrel,
	verbatim,
}
\usepackage{xr-hyper}
\usepackage[all,cmtip]{xy}
\usepackage{enumitem}
\usepackage{soul}
\usepackage{mathrsfs}
\usepackage[T1]{fontenc} 
\usepackage[pagebackref,hyperindex,linktocpage=true]{hyperref}
\hypersetup{
	colorlinks,
	linkcolor={red!50!black},
	citecolor={red!50!black}, 
	urlcolor={red!50!black}, 
	filecolor={red!50!black} 
}

\usepackage{theoremref}

\usepackage{tikz}
\usetikzlibrary{matrix,arrows}
\usepackage{tikz-cd}

\usepackage{leftidx} 

\makeatletter
\@addtoreset{equation}{section}
\makeatother

\numberwithin{equation}{section}

\theoremstyle{definition}
\newtheorem{Defi}{Definition}[section] \newcommand{\defi}{\begin{Defi}} \newcommand{\xdefi}{\end{Defi}} \newtheorem{DefiLemm}[Defi]{Definition and Lemma} \newcommand{\defilemm}{\begin{DefiLemm}} \newcommand{\xdefilemm}{\end{DefiLemm}} 
\newtheorem{Bsp}[Defi]{Example} \newcommand{\exam}{\begin{Bsp}} \newcommand{\xexam}{\end{Bsp}} 
\newtheorem{ToyBsp}[Defi]{Toy Example} \newcommand{\toyexam}{\begin{ToyBsp}} \newcommand{\xtoyexam}{\end{ToyBsp}} 
\newtheorem{Syno}[Defi]{Synopsis} \newcommand{\syno}{\begin{Syno}} \newcommand{\xsyno}{\end{Syno}} 
\newtheorem{Bem}[Defi]{Remark} \newcommand{\rema}{\begin{Bem}} \newcommand{\xrema}{\end{Bem}} 
\newtheorem{Notation}[Defi]{Notation} \newcommand{\nota}{\begin{Notation}} \newcommand{\xnota}{\end{Notation}} 
\newtheorem{Convention}[Defi]{Convention} \newcommand{\conv}{\begin{Convention}} \newcommand{\xconv}{\end{Convention}} 
\newtheorem{Warning}[Defi]{Warning} \newcommand{\warn}{\begin{Warning}} \newcommand{\xwarn}{\end{Warning}} 
\newtheorem{Situation}[Defi]{Situation} \newcommand{\situ}{\begin{Situation}} \newcommand{\xsitu}{\end{Situation}}
\newtheorem{Assumption}[Defi]{Assumption} \newcommand{\assu}{\begin{Assumption}} \newcommand{\xassu}{\end{Assumption}} 

\theoremstyle{plain}
\newtheorem{Theo}[Defi]{Theorem} \newcommand{\theo}{\begin{Theo}} \newcommand{\xtheo}{\end{Theo}} 
\newtheorem{Satz}[Defi]{Proposition} \newcommand{\prop}{\begin{Satz}} \newcommand{\xprop}{\end{Satz}} 
\newtheorem{Lemm}[Defi]{Lemma} \newcommand{\lemm}{\begin{Lemm}} \newcommand{\xlemm}{\end{Lemm}} 
\newtheorem{Coro}[Defi]{Corollary} \newcommand{\coro}{\begin{Coro}} \newcommand{\xcoro}{\end{Coro}}
\newtheorem{Ques}[Defi]{Question} \newcommand{\ques}{\begin{Ques}} \newcommand{\xques}{\end{Ques}}
\newtheorem{Conj}[Defi]{Conjecture} \newcommand{\conj}{\begin{Conj}} \newcommand{\xconj}{\end{Conj}}

\newcommand{\eqn}{\begin{equation}} \newcommand{\xeqn}{\end{equation}}
\newcommand{\eqnarr}{\begin{eqnarray*}} \newcommand{\xeqnarr}{\end{eqnarray*}}
\newcommand{\eqnarra}{\begin{eqnarray}} \newcommand{\xeqnarra}{\end{eqnarray}}

\newcommand{\pf}{\begin{proof}} \newcommand{\xpf}{\end{proof}}

\newcounter{heyheyCounter}[section]

\newif \ifHideSomeComments
\HideSomeCommentstrue

\ifHideSomeComments
\newcommand{\heyheyX}[1]{}
\else
\newcommand{\heyheyX}[1]{\heyhey{#1}}
\fi

\newcommand{\heyhey}[1]{\refstepcounter{heyheyCounter}$\bigstar^\text{\theheyheyCounter}$\marginpar{\footnotesize $\bigstar^\text{\theheyheyCounter}$ #1}}

\newif \ifDraft 
\Draftfalse

\ifDraft
\usepackage[top=1.5cm, bottom=1.5cm, left=.5cm, right=7.5cm, heightrounded, marginparwidth=7cm, marginparsep=.2cm]{geometry}
\else
\usepackage[top=1.5cm, bottom=1.5cm, left=2cm, right=2cm, heightrounded]{geometry}
\usepackage{letltxmacro}
\LetLtxMacro\Oldfootnote\footnote
\renewcommand{\footnote}[2][]{\relax}
\renewcommand{\heyhey}[1]{}
\fi

\newcommand{\nc}{\newcommand}
\nc{\StP}[1]{\cite[Tag~\href{http://stacks.math.columbia.edu/tag/#1}{#1}]{StacksProject}}
\nc{\StPd}[2]{\cite[Tags~\href{http://stacks.math.columbia.edu/tag/#1}{#1}, \href{http://stacks.math.columbia.edu/tag/#2}{#2}]{StacksProject}} 
\nc{\Kero}[1]{\cite[Tag~\href{https://kerodon.net/tag/#1}{#1}]{Kerodon}}
\nc{\Kerod}[2]{\cite[Tags~\href{https://kerodon.net/tag/#1}{#1}, \href{https://kerodon.net/tag/#2}{#2}]{Kerodon}}

\nc{\on}{\operatorname}
\nc{\aff}{{\on{aff}}}
\nc{\modi}{{\on{mod}}} 
\nc{\even}{{\on{even}}}
\nc{\odd}{{\on{odd}}}
\nc{\naive}{{\on{naive}}}
\nc{\hofib}{\on{hofib}}
\nc{\Bun}{\on{Bun}}
\nc{\ad}{{\on{ad}}}
\nc{\lft}{{\on{lft}}}
\nc{\modulo}{\on{mod}} 
\nc{\FinSet}{\on{FinSet}} 
\nc{\surj}{\on{surj}} 
\nc{\threerightarrows}{  \mathrel{
    \substack{
      \xrightarrow{} \\[-0.4ex]
      \xrightarrow{} \\[-0.4ex]
      \xrightarrow{}
    }
  }}
\nc{\Weil}{{\on{Weil}}} 
\nc{\FWeil}{{\on{FWeil}}} 
\nc{\cons}{{\on{cons}}} 
\nc{\tot}{{\on{Tot}}} 

\nc{\str}{\on{-}}
\nc{\perf}{{\on{perf}}}
\nc{\Rel}{{\on{Pos}}}
\nc{\lan}{\langle}
\nc{\ran}{\rangle}

\nc{\bba}{{\mathbf a}}
\nc{\bbb}{{\mathbf b}}
\nc{\bbc}{{\mathbf c}}
\nc{\bbd}{{\mathbf d}}
\nc{\bbe}{{\mathbf e}}
\nc{\bbf}{{\mathbf f}}
\nc{\bbg}{{\mathbf g}}
\nc{\bbh}{{\mathbf h}}
\nc{\bbi}{{\mathbf i}}
\nc{\bbj}{{\mathbf j}}
\nc{\bbk}{{\mathbf k}}
\nc{\bbl}{{\mathbf l}}
\nc{\bbm}{{\mathbf m}}
\nc{\bbn}{{\mathbf n}}
\nc{\bbo}{{\mathbf o}}
\nc{\bbp}{{\mathbf p}}
\nc{\bbq}{{\mathbf q}}
\nc{\bbr}{{\mathbf r}}
\nc{\bbs}{{\mathbf s}}
\nc{\bbt}{{\mathbf t}}
\nc{\bbu}{{\mathbf u}}
\nc{\bbv}{{\mathbf v}}
\nc{\bbw}{{\mathbf w}}
\nc{\bbx}{{\mathbf x}}
\nc{\bby}{{\mathbf y}}
\nc{\bbz}{{\mathbf z}}
\nc{\calA}{{\mathcal A}}
\nc{\calB}{{\mathcal B}}
\nc{\calC}{{\mathcal C}}
\nc{\calD}{{\mathcal D}}
\nc{\calE}{{\mathcal E}}
\nc{\calF}{{\mathcal F}}
\nc{\calG}{{\mathcal G}}
\nc{\calH}{{\mathcal H}}
\nc{\calI}{{\mathcal I}}
\nc{\calJ}{{\mathcal J}}
\nc{\calK}{{\mathcal K}}
\nc{\calL}{{\mathcal L}}
\nc{\calM}{{\mathcal M}}
\nc{\calN}{{\mathcal N}}
\nc{\calO}{{\mathcal O}}
\nc{\calP}{{\mathcal P}}
\nc{\calQ}{{\mathcal Q}}
\nc{\calR}{{\mathcal R}}
\nc{\calS}{{\mathcal S}}
\nc{\calT}{{\mathcal T}}
\nc{\calU}{{\mathcal U}}
\nc{\calV}{{\mathcal V}}
\nc{\calW}{{\mathcal W}}
\nc{\calX}{{\mathcal X}}
\nc{\calY}{{\mathcal Y}}
\nc{\calZ}{{\mathcal Z}}
\nc{\co}{\colon}

\newcommand{\category}[1]{\mathrm{#1}}

\newcommand\restr[2]{{ 
  \left.\kern-\nulldelimiterspace 
  #1 
  \vphantom{\big|} 
  \right|_{#2} 
  }}

\newcommand{\Mod}{\category{Mod}} 
\newcommand{\Ani}{\category{Ani}}

\newcommand{\Perf}{\category{Perf}}

\renewcommand{\L}{\category{L}}
\renewcommand{\Pr}{\category{Pr}}

\newcommand{\PrL}{\Pr^\L}

\newcommand{\PrLSt}{\PrL_{\mathrm{St}}} 
\newcommand{\PrLst}{\PrLSt}

\newcommand{\Cat}{\category{Cat}} 
\newcommand{\Fun}{\category{Fun}} 
\newcommand{\FunL}{\Fun^{\mathrm L}} 
\newcommand{\Corr}{\category{Corr}} 
\newcommand{\Sch}{\category{Sch}} 
\newcommand{\Vbdl}{\category{Vbdl}} 
\newcommand{\can}{\mathrm{can}} 

\newcommand{\IndSch}{\category{IndSch}}

\def\ft{\mathrm{ft}} 
 
\newcommand{\dual}{\vee} 
\newcommand{\Zar}{\mathrm{Zar}} 

\newcommand{\Arr}{\category{Arr}} 

\newcommand{\Fl}{Fl}
\newcommand{\D}{\mathrm{D}}
\newcommand{\et}{\mathrm{et}} 
\newcommand{\opp}{\mathrm{op}}

\def\Gm{\mathbf {G}_\mathrm m} 
\def\Ga{\mathbf {G}_\mathrm a} 
\def\add{\mathrm{add}} 
\def\sub{\mathrm{sub}} 
				
\newcommand{\colim}{\operatornamewithlimits{colim}}

\def\id{{\rm id}}

\def\pr{{\rm pr}} 
\def\opp{{\rm op}} 
 
\def\To#1#2{\mathop{\count0=#1 \loop\ifnum\count0>0 \smash-\mkern-7mu \advance\count0 -1 \repeat \mathord\rightarrow}\limits^{#2}}

\def\Maps{\mathop{\rm Maps}\nolimits}

\def\av{\mathop{\rm av}\nolimits} 
 
\def\cofib{\mathop{\rm cofib}\nolimits} 
 
\def\Hom{\mathop{\rm Hom}\nolimits} 
\def\bx{\boxtimes} 

\def\Du{\mathrm D} 
\def\PGL{\mathrm {PGL}}

\def\fib{\mathop{\rm fib}\nolimits}

\def\IC{\mathrm{IC}} 
				
\def\Ab{\mathop{\rm Ab}\nolimits}
\def\Ind{\mathop{\category{Ind}}} 
\def\RHom{\mathop{\rm RHom}\nolimits} 
\def\Gr{\mathop{\rm Gr}\nolimits} 
\def\Fl{\mathop{\rm Fl}\nolimits} 
\def\Ext{\mathop{\rm Ext}\nolimits} 
\def\IHom{\underline{\Hom}}

\def\Map{\Maps} 
\def\Rep{\category{Rep}} 
 
\def\et{\mathrm{\acute et}} 
 
\def\coker{\operatorname{coker}} 
 
\def\adj{\mathrm{adj}}

\def\sph{\mathrm{sph}}

\def\act{\mathrm{act}}

\def\Z{{\mathbf Z}} 
\def\Fp{{{\mathbf F}_p}} 
\def\Fq{{{\mathbf F}_q}} 
\def\Fqq{\overline {\mathbf F}_q} 
\def\Fpq{\overline {\mathbf F}_p}

\def\Qlq{{\overline \Q_\ell}} 
 
\def\Q{{\mathbf Q}}

\def\Ql{{\Q_\ell}}

\def\C{{\mathbf {C}}} 
\def\A{{\bf A}} 
\renewcommand{\P}[1][1]{\mathbf P^{#1}}

\def\H{{\rm H}} 
\def\pe{{{}^\mathrm{p}} \! }

\def\red{\mathrm{r}} 
\def\redx{(\red)} 
\def\DM{\category{DM}} 
\def\DA{\category{DA}} 
\def\DAW{\DA_{\mathrm W}} 
\def\Ker{\category{Ker}} 
\def\DMr{\DM_\red} 
\def\DMrx{\DM_{\redx}} 			
\def\DTM{\category{DTM}} 
\def\DTMr{\DTM_\red} 
\def\DTMrx{\DTM_{\redx}} 
\def\anti{\mathrm{anti}} 
\def\xanti{{(\anti)}} 
 
\def\bff{\mathbf{f}}
\def\DMexp{\DM_{\exp}} 
\def\DAe{\DA_{\exp}} 
\def\DAest{\DA_\est} 
\def\DAesh{\DA_\esh} 
\def\DMe{\DMexp} 
\def\est{{\exp*}} 
\def\esh{{\exp!}} 
\def\DMest{\DM_\est} 
\def\DMesh{\DM_\esh} 
\def\DTMe{\DTM_{\exp}} 
\def\DTMea{\DTMe^\anti} 
\def\DTMexa{\DTMe^\xanti} 
\def\DTMa{\DTM^\anti} 
\def\MTMe{\MTM_{\exp}} 
\def\DTMexp{\DTMe} 
\def\MTMexp{\MTMe} 
\def\DTMest{\DTM_{\est}} 
\def\DTMesh{\DTM_{\esh}} 
\newcommand{\Obj}{\category{Obj}} 
\def\Gm{\mathbf {G}_\mathrm m} 
\def\Ga{\mathbf {G}_\mathrm a} 
\newcommand{\GaX}[  1]{\mathbf {G}_{\mathrm {a}, #1}} 
\newcommand{\GmX}[  1]{\mathbf {G}_{\mathrm {m}, #1}} 
\def\Perv{\mathrm{Perv}} 

\def\MTM{\category{MTM}} 
\def\MTMr{\MTM_\red} 
\def\MTMrx{\MTM_{\redx}} 
 
\def\ii{$\infty$}

\renewcommand{\r}{\rightarrow}			
\newcommand{\x}{\times}			
\renewcommand{\t}{\otimes}			
\newcommand{\ol}{\overline}	
\newcommand{\refsect}[1]{§\ref{sect--#1}}	

\def\Spec{\mathop{\rm Spec}}

\newcommand{\comp}{\omega} 
\newcommand{\locc}{\mathrm{lc}}

\def\bd{\mathrm{b}}

\DeclareMathOperator{\im}{im}
\newcommand{\rhoexp}{\hat{\rho}\text{-}\!\exp}
\newcommand{\rhocirc}{\hat{\rho}\text{-}\circ}
\DeclareMathOperator{\Stab}{Stab}
\DeclareMathOperator{\IndCoh}{IndCoh}
\DeclareMathOperator{\Coh}{Coh}
\DeclareMathOperator{\QCoh}{QCoh}
\DeclareMathOperator{\CS}{CS} 
\newcommand{\CSr}{\CS_{\red}}
\newcommand{\B}{\mathbf{B}} 

\def\Delexp{\underline{\Delta}^{\exp}}
\def\nabexp{\overline{\nabla}^{\exp}}
\def\Delrhoexp{\underline{\Delta}^{\rhoexp}}
\def\nabrhoexp{\overline{\nabla}^{\rhoexp}}

\newcommand{\pot}[1]{ [\hspace{-0,5mm}[ {#1} ]\hspace{-0,5mm}] }
\newcommand{\rpot}[1]{ (\hspace{-0,7mm}( {#1} )\hspace{-0,7mm}) }
				
\mathchardef\mhyphen="2D

\pagecolor{white}

\begin{document}

\title{Exponential motives on the affine Grassmannian}
\author{Robert Cass, Thibaud van den Hove, Jakob Scholbach}

\begin{abstract}
We develop a notion of exponential motives on general prestacks equipped with a $\Ga$-action, and compare them with Whittaker motives via Gaitsgory's Kirillov model.
We then establish foundational results for exponential motives on affine flag varieties concerning Tate motives and t-structures.
We use this to prove a motivic Casselman--Shalika equivalence, relating exponential Tate motives on the affine Grassmannian to ind-coherent sheaves on the classifying stack of the Langlands dual group.
The decategorification of this equivalence provides a new construction of the Whittaker module for the spherical Hecke algebra which works for arbitrary coefficients, including a generic version.
\end{abstract}

\maketitle

\setcounter{tocdepth}{1}
\tableofcontents

\section{Introduction}
Whittaker functions are an important tool in the Langlands program, and they are used to study both smooth representations of reductive groups over local fields, as well as automorphic representations.
In the spirit of geometrizing the Langlands program, various authors have defined notions of Whittaker sheaves.
Both the global approach from \cite{FrenkelGaitsgoryVilonen:Whittaker}, and the local \emph{Iwahori--Whittaker} (also known as \emph{baby Whittaker}) approach from \cite{ArkhipovBezrukavnikov:Perverse} categorify the classical Whittaker functions via Grothendieck's sheaf-function dictionary. 

In this paper we introduce a motivic version of Iwahori--Whittaker sheaves, and we use them to prove a motivic Casselman--Shalika equivalence. In upcoming work, we will study in detail the category of Iwahori--Whittaker motives on the full affine flag variety and prove a motivic refinement of the Arkhipov--Bezrukavnikov equivalence \cite{ArkhipovBezrukavnikov:Perverse}.
In particular, this paper is a next step in our project \cite{CassvdHScholbach:MotivicSatake,CassvdHScholbach:Central} of enhancing Bezrukavnikov's work \cite{Bezrukavnikov:Two} to prove a tamely ramified motivic local Langlands correspondence.
Moreover, it fits into the recent trend of using motivic sheaves to prove independence of \(\ell\) in the Langlands program at the categorical level \cite{RicharzScholbach:Intersection,RicharzScholbach:Motivic,Scholze:GeometrizationMotivically}.

\subsection{Exponential motives and Whittaker motives}
The statement of our motivic Casselman--Shalika equivalence is in terms of \emph{exponential motives} (also known as the \emph{Kirillov model} in the context of our main applications).
Compared to classical Whittaker sheaves, exponential motives bring a number of benefits.
Classical Whittaker sheaves are restricted to schemes over $\Spec \Fp$ and depend on a choice of a generic character \(\Fq \to \Lambda^\times\) (where \(\Lambda\) is the ring of coefficients used).
This excludes, say, \(\Lambda =\Fpq\).
Exponential motives, by contrast, exist for integral coefficients and over arbitrary base schemes. 
We denote by \(\DM(-)\) the category of (Nisnevich) motivic sheaves with \(\Z\)-coefficients, as constructed in \cite{Spitzweck:Commutative}.

\defi
For a scheme or, more generally, a prestack $X$ endowed with an action of the additive group $\Ga$, the category of \emph{exponential motives} is defined as
$$\DMe(X) := \DM(X) / \DM(X / \Ga),$$ 
i.e., motives on $X$ modulo those that are $\Ga$-equivariant. 
\xdefi

In the prototypical case $X = \Ga \x Y$ (for some scheme $Y$) this is the category considered by Gallauer--Pepin Lehalleur in \cite{GallauerPepinLehalleur:Exponentiation}, whose work was in turn inspired by Fresán--Jossen's \cite{FresanJossen:Exponential} paper on exponential Nori motives.
Rudiments of this construction have also appeared in \cite{Gaitsgory:LocalGlobal,GaitsgoryLysenko:Metaplectic}.

As we show in \refsect{exponential motives}, the assignment $X \mapsto \DMe(X)$ is a four-functor formalism in the sense that $f_*, f_!, f^*, f^!$ exist (for $\Ga$-equivariant maps $f$) and satisfy the properties familiar from (regular) motives.

Applications in geometric representation theory will often use the following natural extension of the above definition, where $\Ga \rtimes \Gm$ acts on $X$:
\begin{equation}\label{intro:eqexp}
	\DMe(X/\Gm) := \DM(X/\Gm) / \DM(X / \Ga \rtimes \Gm).
\end{equation}
Recall also that when using motives in geometric representation theory, one often restricts to \emph{stratified Tate motives} \cite{SoergelWendt:Perverse,RicharzScholbach:Intersection}, which form a class of motives which is more controllable than the full category of motives.
To transfer this to the setting of exponential motives, we consider the situation of an (ind-)scheme $X$ equipped with a Whitney--Tate stratification $X^\dagger$, satisfying an additional geometric condition concerning the behavior of the $\Ga$-averaging functor called \emph{\(\Ga\)-averageability} (\thref{DTMexp}).
Then \(\DMexp(X)\) admits a subcategory of \emph{exponential stratified Tate motives} $\DTMe(X, X^\dagger)$, which inherits the nice properties that usual Tate motives satisfy.
We also consider a variant where only anti-effective Tate motives are allowed (\refsect{exponential Tate motives}).
A basic example, which is the seed of our later considerations of exponential motives on (affine) flag varieties, is the category $\DTMe(\A^1) := \DTM(\A^1_S, \GmX S \sqcup 0) / \DTM(S)$, which is equivalent to $\DTM(\GmX S)$.

To facilitate a comparison of exponential Tate motives on the affine Grassmannian with the classical approach of Whittaker sheaves, we provide a comparison between exponential motives and Whittaker motives in \refsect{Whittaker motives}. For a scheme or prestack $X$  over $\Spec \Fp$ with a $\Ga$-action, we consider a category $\DAW(X)$ of \emph{(étale) Whittaker motives}, which can be thought of as a  categorification of Whittaker functions, i.e., satisfying $f(\lambda \cdot x)=\psi(\lambda) f(x)$, where $\psi : \Fp \r \Lambda^\x$ is a non-trivial character, and $\Lambda$ a coefficient ring containing $1/p$ and a $p$-th root of unity. This category admits various equivalent descriptions, one of them being the full subcategory of $\DA(X)$ (étale motives with $\Lambda$-coefficients) 
given by the essential image of convolution by a motivic Artin--Schreier sheaf in $\DA(\Ga)$ associated to $\psi$.

\theo\thlabel{intro--Kirillov model}
For any prestack $X / \Fp$ acted upon by $\Ga \rtimes \Gm$, there is an equivalence
$$\DAe(X/\Gm) \cong \DAW(X)$$
between the category  of $\Gm$-equivariant exponential étale motives (defined in \eqref{intro:eqexp}) and the category of Whittaker motives.
\xtheo 
This assertion, proved in \thref{Whit vs Kir} below, confirms that the argument outlined by Gaitsgory--Lysenko in \cite[§A]{GaitsgoryLysenko:Metaplectic} holds in any homotopy-invariant six functor formalism.

\subsection{Exponential geometry of affine flag varieties}
Let $G$ be a pinned split reductive group over a base scheme $S$, with maximal torus and Borel $T \subset B \subset G$. We assume that $S$ is of finite type over a Dedekind domain or field, and that $G$ is not a torus. For a facet $\mathbf{f}$ contained in the closure of the standard alcove $\mathbf{a}_0$ associated to $B$, let $\Fl_{\mathbf{f}}$ be the corresponding partial affine flag variety over $S$. 

To construct a suitable category $\DTM_{\exp}(\Fl_{\mathbf{f}})$, we require a stratification that is finer than the usual one by orbits under the Iwahori loop group $\calI = L^+G_{\mathbf{a}_0}$. To define it, let $\calU \subset \calI$ be the pro-unipotent radical, and let $\calU_0 \subset \calU$ be the kernel of the homomorphism $\calU \r \Ga$ given by the sum of the simple root groups. The half-sum of positive coroots $\hat{\rho}$ induces by conjugation an action of $\Gm$ on $\calU_0$, and we may form the semi-direct product $\calU^{\exp} := \calU_0 \rtimes \Gm$. 

\theo \thlabel{WTstratTheorem} Let $\mathbf{f}$ be a facet contained in the closure of the standard alcove.
Then the $\calU^{\exp}$-orbits on $\Fl_{\mathbf{f}}$ form a stratification. Furthermore, this stratification is anti-effective universally Whitney--Tate.
\xtheo

Here {\em{Whitney--Tateness}} means that we can glue Tate motives on the $\calU^{\exp}$-orbits, which is what is needed to consider stratified Tate motives in the first place.
The {\em{universality}} of the stratification means that push-pull functors between strata commute with base change in $S$, and {\em{anti-effectiveness}} means that nonpositive Tate twists are also preserved.
The corresponding assertion for the coarser stratification by orbits under $\calI$ and more general parahorics was proved in \cite{RicharzScholbach:Intersection} and further utilized in \cite{CassvdHScholbach:MotivicSatake, CassvdHScholbach:Central}.

This statement appears in \thref{exponential orbits WT} below. Its proof starts with the finite flag variety $G/B \subset \Fl_{\mathbf{a}_0}$. Here, the stratification of $G/B$ by $\calI$-orbits is the usual Bruhat stratification by $U$-orbits, where $U \subset B$ is the unipotent radical. Each Bruhat cell is also a $\calU^{\exp}$-orbit, except for the largest Bruhat cell, which is acted upon simply transitively by $U$ and splits into two $\calU^{\exp}$-orbits, namely one being isomorphic to $\ker (U \r \Ga)$, and its complement in $U$. The fact that the closure of the former $\calU^{\exp}$-orbit is a union of $\calU^{\exp}$-orbits, independently of the base $S$, is \thref{theo--strat}.  The proof makes use of the Demazure resolution as well as Deodhar's decomposition of Richardson varieties. The proof of Whitney--Tateness then involves a study of push-pull functors along the projection $G/B \r G/P$ for a minimal parabolic $P \supset B$, making use of the fact that $\calU^{\exp}$-orbits in $G/P$ are all Bruhat cells.

By \thref{WTstratTheorem}, we have the category $\DTM(\Fl_{\mathbf{f}})$ of stratified Tate motives with respect to the stratification by $\calU^{\exp}$-orbits. Noting that $(\calU \rtimes \Gm)/\calU_0 \cong \Ga \rtimes \Gm$, we apply the above formalism of exponential motives and define the category of exponential stratified Tate motives as
\begin{equation}\label{intro-defi Flexp}
	\DTM_{\exp}(\Fl_{\mathbf{f}}) = \DTM(\calU^{\exp} \backslash \Fl_{\mathbf{f}})/ \DTM((\calU \rtimes \Gm) \backslash \Fl_{\mathbf{f}}).
\end{equation}

A refined geometric analysis is necessary in order to compute the adjoints of the forgetful functor $\DTM((\calU \rtimes \Gm) \backslash \Fl_{\mathbf{f}}) \r \DTM(\calU^{\exp} \backslash \Fl_{\mathbf{f}})$ (Propositions \ref{GaAverageable} and \ref{nunu}).
This is used to show that the above category $\DTM_{\exp}(\Fl_{\mathbf{f}})$ is well-behaved, and that it inherits from $\DTM(\calU^{\exp} \backslash \Fl_{\mathbf{f}})$ a t-structure with heart $\MTM_{\exp}(\Fl_{\mathbf{f}})$ consisting of {\em {exponential stratified mixed Tate motives}}, as long as \(S\) satisfies the Beilinson--Soulé vanishing conjecture (which is known for \(\Spec \Z\), \(\Spec \Q\), and \(\Spec \Fp\)). 

With the category $\DTM_{\exp}(\Fl_{\mathbf{f}})$ at hand, we proceed to study its module structure. The category \(\DM(L^+G_{\bbf} \backslash \Fl_{\bbf})\) is monoidal for the convolution product, which also restricts to a convolution product on stratified Tate motives as is shown in \cite{RicharzScholbach:Motivic,CassvdHScholbach:MotivicSatake}.
For general motives, the convolution product induces a right action of \(\DM(L^+G_{\bbf} \backslash \Fl_{\bbf})\) on \(\DM_{\exp}(\Fl_{\bbf})\). The following cellular description of the fibers of exponential convolution morphisms, as in \thref{fibers of exponential convolution}, implies that Tate motives are preserved, so that $\DTM_{\exp}(\Fl_{\mathbf{f}})$ is a module over $\DTM(L^+G_{\bbf} \backslash \Fl_{\bbf})$.

\theo\thlabel{intro--cellularity}
Let \(\bbf\) be a facet in the closure of the standard alcove, let \(\Fl_z^{\exp}(\bbf)\) be any \(\calU^{\exp}\)-orbit in \(\Fl_{\bbf}\), and let \(\Fl_w^{\circ}(\bbf)\) be any \(L^+G_{\bbf}\)-orbit in \(\Fl_{\bbf}\).
Consider the convolution morphism
\[m_{z,w}\colon \Fl_z^{\exp}(\bbf) \widetilde{\times} \Fl_w^\circ(\bbf) \to \Fl_{\bbf}.\]
Then any fiber of \(m_{z,w}\) admits a filtrable decomposition into cells, i.e., into products of schemes of the form \(\A^1\), \(\Gm\), or \(\Gm\setminus \{1\}\).
\xtheo

This theorem refines the main result of \cite{Haines:Pavings}, and we follow its proof.
In contrast to \cite{Haines:Pavings}, we have to allow cells of the form \(\Gm\setminus \{1\}\), as can already be seen in the case of the finite flag variety \(\P\) for \(\mathrm{SL}_2\) in \thref{example punctured Gm}.

In the original non-motivic works which we refine, a fundamental technical input comes in the form of t-exactness results on convolution by (co)standard objects. In classical sheaf theories such as $\ell$-adic or Betti sheaves, these results can mostly be deduced from Artin vanishing or the t-exactness of nearby cycles, but these are not available in the motivic setting. Similar obstacles also occur when constructing the t-structure on $\DTM_{\exp}(\Fl_{\mathbf{f}})$, as discussed above. In the context of exponential Tate motives, we prove the following result on (co)standard objects.

\theo \thlabel{mainthm--conv} Let $\mathbf{f}, \mathbf{f}'$ be facets in the closure of the standard alcove and let $\iota_w : \Fl_{w}^{\exp}{(\mathbf{f}')} \r \Fl_{\mathbf{f}'}$ be the inclusion of a $\calU^{\exp}$-orbit.
Then convolution by $\Delta_w^{\exp}(
\Z) := \iota_{w!} \Z[\dim \Fl_{w}^{\exp}{(\mathbf{f}')}]$ defines a left $t$-exact functor $\DTM(
L^+G_{\mathbf{f}'} \backslash \Fl_{\mathbf{f}}) \r \DTM(\calU^{\exp} \backslash \Fl_{\mathbf{f}})$.
Dually, convolution with $\nabla_w^{\exp}(\Z) := \iota_{w*}\Z[\dim \Fl_{w}^{\exp}{(\mathbf{f}')}]$ is right $t$-exact.
\xtheo

As a special case, this implies that the (co)standard functors $\Delta_w^{\exp}, \nabla_w^{\exp} \colon \DTM(\calU^{\exp} \backslash S) \r \DTM(\calU^{\exp} \backslash \Fl_{\mathbf{f}})$ are t-exact.
By passing to quotient categories, we immediately get analogous results for $\DTM_{\exp}(\Fl_{\mathbf{f}})$ in \thref{exponential convolution}. The proof of \thref{mainthm--conv} is closely related to the proof of Whitney--Tateness in \thref{WTstratTheorem}. The proof also leverages similar results obtained by the authors in \cite{CassvdHScholbach:Central} for (co)standard objects on Iwahori-orbits, combined with a direct study on the finite flag variety $G/B$.

\subsection{A motivic Casselman--Shalika equivalence}

We now specialize to the case of the affine Grassmannian $\Gr=\Fl_{\bbf_0}$, and the positive loop group \(L^+G=L^+G_{\bbf_0}\), associated to the standard facet \(\bbf_0\) containing the origin.
The monoidal category of \(L^+G\)-equivariant sheaves on \(\Gr\) is a categorification of the spherical Hecke algebra. By the Casselman--Shalika formula \cite{CasselmanShalika:Unramified}, 
the Whittaker module of the spherical Hecke algebra is free of rank one.
Combining this with the Satake isomorphism \cite{Satake:Theory} shows that the Whittaker module is isomorphic to the representation ring of the Langlands dual group.
The Casselman--Shalika formula was geometrized in \cite{FrenkelGaitsgoryVilonen:Whittaker,NgoPolo:Resolutions}, and \cite{BGMRR:IwahoriWhittaker} categorified it by establishing an equivalence between the categories of perverse \(L^+G\)-equivariant étale sheaves on \(\Gr\), and perverse Iwahori--Whittaker sheaves on \(\Gr\) (cf.~\cite{ABBGM:Modules} for a similar equivalence involving D-modules).
The following theorem, proved in \thref{Equivalence in the reduced case} and \thref{Motivic CasselmanShalika}, provides a motivic enhancement of the Casselman--Shalika formula.
It uses the category of reduced motives from \cite{EberhardtScholbach:Integral}, which are denoted by a subscript $\red$, such as $\MTMr$ and $\DTMr$. The process of reducing motives amounts to modding out the higher motivic cohomology of the base $S$, so that $\DTMr(S)$ is just the derived category of graded abelian groups.

\theo\thlabel{intro--main theorem}
Let \(G\) be a split reductive group over \(S\), which is not a torus.
For reduced motives, there is a natural equivalence 
\begin{equation}\label{intro-reduced equivalence}
	\Phi^{\heartsuit}\colon \MTM_{\red,L^+G}(\Gr) \cong \MTM_{\red,\rhoexp}(\Gr),
\end{equation}
given by convolving with a fixed exponential motive \(\Delrhoexp_0(\Z) \cong \nabrhoexp_0(\Z)\in \MTM_{\red,\rhoexp}(\Gr)\).

For regular motives, there is the following commutative diagram, where the horizontal arrows are equivalences and the vertical arrows are reduction functors: 
\begin{equation}\label{intro-regular equivalence}
\begin{tikzcd}
	\IndCoh\left(\Spec \Z/(\hat{G}\rtimes \Gm)\right) \otimes_{\DTMr(S)} \DTM(S) \arrow[d] \arrow[r, "\cong"] & \DTM_{\rhoexp}(\Gr)\arrow[d]\\
	\IndCoh\left(\Spec \Z/(\hat{G}\rtimes \Gm)\right) \arrow[r, "\cong"'] & \DTM_{\red,\rhoexp}(\Gr).
\end{tikzcd}
\end{equation}
Here the lower horizontal equivalence is (after first passing to compact objects) induced by deriving the precomposition of \eqref{intro-reduced equivalence} with the motivic Satake equivalence $$\Rep_{\hat{G} \rtimes \Gm}(\Ab) \cong \MTM_{\red,L^+G}(\Gr)$$ from \cite{CassvdHScholbach:MotivicSatake}.
\xtheo

The notation \(\rhoexp\) indicates that, in \eqref{intro-defi Flexp}, we have replaced \(\calU\) and \(\calU^{\exp}\) by their \(\hat{\rho}(t^{-1})\)-conjugates, where \(\hat{\rho}\) is the half-sum of the positive coroots of \(G\); cf.~\thref{defi twisted exponentials} for details.
This allows us to remove an assumption on \(G\) appearing in \cite{BGMRR:IwahoriWhittaker}, and for e.g.~adjoint \(G\) we have \(\DTMexp(\Gr)\cong \DTM_{\rhoexp}(\Gr)\).
However, we emphasize that \(\DTM_{\rhoexp}(\Gr)\) is more natural than \(\DTM_{\exp}(\Gr)\); cf.~\cite{Raskin:Walgebras}, or the discussion after \thref{intro-Hecke} below.

Let us highlight some of the motivic difficulties that appear when proving \thref{intro--main theorem}.
The main ingredient in the proof of \cite{BGMRR:IwahoriWhittaker} is the geometric Casselman--Shalika formula from \cite{FrenkelGaitsgoryVilonen:Whittaker,NgoPolo:Resolutions}.
This requires the choice of an Artin--Schreier local system, and hence does not admit a formulation in the setting of motivic sheaves with integral coefficients over general bases.

To get around this issue, recall that the main use of the geometric Casselman--Shalika formula is to prove that \(\Phi^\heartsuit\) identifies the standard objects in \(\MTM_{L^+G}(\Gr)\) and \(\MTM_{\rhoexp}(\Gr)\).
We show in \thref{Equivalent conditions} that this property of \(\Phi^\heartsuit\) is equivalent to a certain dimension bound for the fibers of exponential convolution morphisms.
Although we do not prove this bound directly, by \thref{intro--cellularity} these fibers admit decompositions into cells.
Since this holds over \(\Spec \Z\), the dimensions of these fibers is independent of the base, so that the dimension estimate can be checked over \(\Fpq\).
There, we deduce the bound from the results of \cite{BGMRR:IwahoriWhittaker}, by using \thref{intro--Kirillov model} on the level of \(\ell\)-adic sheaves. 
In particular, \thref{intro--main theorem} implicitly relies on the geometric Casselman--Shalika formula.

Another key input in \cite{BGMRR:IwahoriWhittaker}, especially when working with torsion coefficients, is the fact that the category of perverse Iwahori--Whittaker sheaves on \(\Gr\) with coefficients in a field is a highest weight category.
This also holds for reduced motives, which allows us to prove the equivalence \eqref{intro-reduced equivalence}.
However, this usually fails when working with regular motives: in general \(\MTM(S)\) does not have any nontrivial projective objects, even with rational coefficients.
Instead, we pass to stable \(\infty\)-categories, and we show that \(\DTM_{\rhoexp}(\Gr)\) can be obtained by \emph{unreducing} \(\DTM_{\red,\rhoexp}(\Gr)\), answering a question from \cite[§1.6.1]{EberhardtScholbach:Integral} in a special case.
This procedure of unreducing reduced motives is not formal, and it crucially uses that the Langlands dual group controlling the motivic Satake category is reduced \cite[§6.2.1]{CassvdHScholbach:MotivicSatake}.

\subsection{The generic Whittaker module}

Using our framework of exponential motives on affine flag varieties, we can construct a natural module under the spherical Hecke algebra which exists for general coefficients, and agrees with the Whittaker module after the choice of an additive character of \(\Fq\).
This construction does not seem to have appeared in the literature before.

In fact, this construction works in the setting of generic Hecke algebras, as in \cite{Vigneras:Algebres,PepinSchmidt:Generic}.
Similarly to \cite[Definition 6.1]{CassvdHScholbach:Central}, we define the \emph{generic exponential module} \(M_{\exp}(\mathbf{q})\) as the Grothendieck group of compact objects in the category \(\DTM_{\rhoexp}^{\anti}(\Gr)\) of anti-effective (twisted) exponential motives on the affine Grassmannian.
This is naturally a module under the generic spherical Hecke algebra \(\mathcal{H}^{\sph}(\mathbf{q}) = K_0(\MTM^{\anti}_{L^+G}(\Gr)^{\omega}) \) from \cite{PepinSchmidt:Generic,CassvdHScholbach:MotivicSatake}, which is itself free as a \(\Z[\mathbf{q}]\)-module, where \(\mathbf{q}\) is a free parameter.
Here we use anti-effective motives, since using all stratified Tate motives would yield a \(\Z[\mathbf{q}^{\pm 1}]\)-module instead.

\theo\thlabel{intro-Hecke}
The generic exponential module \(M_{\exp}(\mathbf{q})\) is a free rank one \(\mathcal{H}^{\sph}(\mathbf{q})\)-module.
For a prime power \(q\), its specialization along \(\mathbf{q}\mapsto q\) can be described more explicitly as
\[M_{\exp}(\mathbf{q}) \otimes_{\Z[\mathbf{q}], \mathbf{q}\mapsto q} \Z\cong \frac{C_c\left(\calU_{\hat{\rho}}^{\exp}(\Fq) \backslash G(\Fq\rpot{t})/G(\Fq\pot{t}),\Z\right)}{C_c\left(\calU_{\hat{\rho}}(\Fq)\backslash G(\Fq\rpot{t})/G(\Fq\pot{t}),\Z\right)},\]
where the right hand side is a quotient of two modules of equivariant functions \(G(\Fq\rpot{t}) \to \Z\). 
Finally, after the choice of a nontrivial character \(\psi\colon \Fq\to \Lambda^\times\), for \(\Lambda\) a \(\Z[\frac{1}{q}]\)-algebra, there is a natural isomorphism
\begin{equation}\label{equation intro} M_{\exp}(\mathbf{q}) \otimes_{\Z[\mathbf{q}], \mathbf{q}\mapsto q} \Lambda \cong C_c\left((U(\Fq\rpot{t}),\Psi_U )\backslash G(\Fq\rpot{t})/G(\Fq\pot{t}),\Lambda\right)\end{equation}
of the \(q\)-specialization of \(M_{\exp}(\mathbf{q})\) with the Whittaker module.
\xtheo
The precise definition of \(\Psi_U\) can be found in \eqref{whittaker map}.
The groups \(\calU_{\hat{\rho}}\) and \(\calU_{\hat{\rho}}^{\exp}\) are the \(\hat{\rho}(t^{-1})\)-conjugates of \(\calU\) and \(\calU^{\exp}\), cf.~§\ref{subsec--twisted exponential motives}.

The fact that \(M_{\exp}(\mathbf{q})\) is free of rank one as a \(\calH^{\sph}(\mathbf{q})\)-module follows by decategorifying \thref{intro--main theorem}.
Moreover, the proof of \eqref{equation intro} essentially follows from the same argument as \thref{intro--Kirillov model}, but on the level of functions, together with the observation that the so-called Whittaker and baby Whittaker modules coincide, i.e., that
\[C_c\left((U(F),\Psi_U) \backslash G(F)/G(\mathcal{O}),\Lambda\right) \cong C_c\left((\calU_{\hat{\rho}}(\Fq),\Psi_{\calU_{\hat{\rho}}}) \backslash G(F)/G(F),\Lambda\right);\]
cf.~\thref{Whittaker vs baby}.
This in particular explains why using \(\calU_{\hat{\rho}}\) is more natural than \(\calU\) in this context (even for adjoint groups).
Namely, the proofs of \thref{combinatorics of Whittaker} and \thref{Whittaker vs baby} show that it is essential that the maps \(\Psi_U\colon U(F) \to \Lambda^\times\) and \(\Psi_{\calU_{\hat{\rho}}}\colon \calU_{\hat{\rho}}(\Fq)\to \Lambda^\times\) agree on the affine root groups \(U_{\alpha,-1}(\Fq)\), where \(\alpha\) ranges over the simple roots.
The analogous statement for \(\Psi_{\calU}\colon \calU(\Fq)\to \Lambda^\times\) does not even make sense, since \(U_{\alpha,-1}\nsubseteq \calU\), but it holds after conjugating \(\calU\) by \(\hat{\rho}(t^{-1})\). 

\subsection*{Acknowledgments} We thank Sean Cotner, Arnaud Eteve, Dennis Gaitsgory, Simon Pepin Lehalleur, Timo Richarz, Simon Riche, Peter Scholze, George Seelinger, and Calvin Yost-Wolff for helpful discussions.

R.C.~was supported by the National Science Foundation grant DMS 1840234. R.C.~also thanks Claremont McKenna College, the University of Michigan, and the Sydney Mathematical Research Institute for excellent working conditions. 

T.v.d.H. and J.S. thank the Max Planck Institute for Mathematics in Bonn for its hospitality during the completion of this work. J.S.~acknowledges support by the European Union, project number 101168795-REMOLD.
This project was started while T.v.d.H. was a PhD student at the TU Darmstadt, where he was supported by the European research council (ERC) under the European Union’s Horizon 2020 research and innovation programme (grant agreement No 101002592), and the LOEWE professorship in Algebra (through Timo Richarz), project number LOEWE/4b//519/05/01.002(0004)/87.

\section{Exponential motives}
\label{sect--exponential motives}

In this section, we develop the rudiments of a theory of exponential motives. These were previously approached from slightly different angles by Fresán--Jossen \cite{FresanJossen:Exponential}, Gallauer--Pepin Lehalleur \cite{GallauerPepinLehalleur:Exponentiation} and Gaitsgory--Lysenko \cite{GaitsgoryLysenko:Metaplectic}.
Their purpose is to provide a base-independent replacement for Artin--Schreier sheaves over $\Spec \Fp$ and the exponential D-module over $\C$.

\subsection{Preliminaries and recollections}

\subsubsection{Motives}

Throughout this paper, we will use the following notations and conventions: $S$ denotes a scheme smooth and of finite type over a field or a Dedekind ring.
The category of finite-type $S$ schemes is denoted $\Sch_S^\ft$. For $X \in \Sch_S^\ft$, we write $\DM(X)$ for the category of motivic sheaves on $X$ with integral coefficients and defined using the Nisnevich topology, as constructed by Spitzweck \cite{Spitzweck:Commutative} (building on the works of Ayoub, Bloch, Cisinski, Déglise, Geisser, Levine and Morel--Voevodsky; we refer to op.~cit.~for further references).
It is known (e.g., \cite[Appendix]{RicharzScholbach:Motivic}), that motives form a six functor formalism in the sense that there is a lax symmetric monoidal functor
$$\DM: \Corr(\Sch_S^\ft) \r \PrL.$$
This includes, in particular, the existence of a functor
$$\bx : \DM(X) \t \DM(Y) \r \DM(X \x_S Y)$$
which is compatible with *-pullbacks and !-pushforwards.

\subsubsection{Motives on prestacks; equivariant motives}

We will also use the extension of the formalism of motives to prestacks as developed in \cite[§2]{RicharzScholbach:Intersection} and briefly reviewed in \cite[§2.1.5, §2.5]{CassvdHScholbach:MotivicSatake}.

For a prestack $X$ equipped with an action of a smooth group scheme $G/S$, the category $\DM(G \setminus X)$ of $G$-equivariant motives can be computed as 
\begin{equation}
\DM(G \setminus X) = \lim_{\Delta_{inj}} \DM(G^\bullet \x X), \label{DM G X}
\end{equation}
where the transition functors in the diagram are !-pullbacks along the various action and projection maps, beginning with $\act, \pr : G \x X \r X$ in low degrees.
By relative purity, we may equivalently use *-pullback functors.

We write $u : X \r G \setminus X$ for the canonical quotient map and $u^! : \DM(G \setminus X) \r \DM(X)$ for the forgetful functor.
If $X$ is placid (in the sense of \cite[Appendix A]{RicharzScholbach:Motivic}), this functor is $\DM(S)$-linear if we let $\DM(S)$ act on $\DM(X)$ and $\DM(G \setminus X)$ by $\boxtimes$ (along the placid structural maps $f : X \r S$ and $\ol f: G \setminus X \r S$.)

\lemm
\thlabel{av}
The functor $u^!$ admits both a right adjoint $\av_*$ and a left adjoint $\av_!$. There are natural isomorphisms
\begin{align}
u^! \av_* & = \act_* \pr^* \label{u! av*}\\
u^! \av_! & = \act_! \pr^!. \nonumber
\end{align}
If $X$ is placid, the functor $\av_!$ is naturally $\DM(S)$-linear. The same is true for $\av_*$ if $S = \Spec k$ (a field), and we work with $\Z[1/e]$-coefficients, where $e$ is the exponential characteristic of $k$.
\xlemm

\pf
This is \cite[Lemma~2.22]{CassvdHScholbach:MotivicSatake};
indeed the restriction to schemes (as opposed to prestacks) in loc.~cit. is unnecessary for the above statement. Regarding the $\DM(S)$-linearity of $\av_*$, the characteristic 0 assumption in point (6) in loc.~cit. can be replaced by working with $\Z[1/e]$ coefficients:
in the proof of loc.~cit., one can still use \cite[Theorem~2.4.6]{JinYang:Kuenneth}, cf.~condition (RS2) there.
\xpf

\exam \thlabel{example-avg}
If $G = \Ga$, the additive group or, more generally any group that is isomorphic (as a scheme) to some $\A^n$, then $u^!$ is fully faithful; equivalently $\av_! u^! \stackrel \cong \r \id \stackrel \cong \r u^! \av_*$ by homotopy invariance and \thref{av}.
\xexam

\rema
\thlabel{functoriality equivariant}
Let $f : X \r X'$ be a $G$-equivariant map of prestacks.
Consider the induced map $\ol f :  G \setminus X \r G \setminus X'$.
As is explained in \cite[Remark~2.21]{CassvdHScholbach:MotivicSatake}, we then have a natural adjunction
$$\ol f_! : \DM(G \setminus X) \rightleftarrows \DM(G \setminus X') : \ol f^!.$$
These functors are compatible with the usual !-functors after applying $u^!$ (i.e., $u_{X'}^! \ol f_! = f_! u_X^!$ etc.)

If $f$ is schematic and of finite type, then there is similarly an adjunction $(\ol f^*, \ol f_*)$, and again these functors are compatible with $(f^*, f_*)$. 
This holds since the !-pullback along the \emph{smooth} action and projection maps in the semi-simplicial diagram in \eqref{DM G X} commute with all four functors.
\xrema

A key computational tool is the following isomorphism known as \emph{hyperbolic localization}.

\prop
\thlabel{hyperbolic localization}
Let $S \stackrel i \r \A^1_S \stackrel q \r S$ be the zero section and the structural maps. Denoting their $\Gm$-prestack quotients by an overline, we have an isomorphism of functors $\DM(\A^1 / \Gm) \r \DM(S / \Gm)$:
\begin{equation}
	\ol q_! \stackrel \cong \r \ol i^!.\label{HL}
\end{equation}
\xprop

\pf
It suffices to show this after applying $u^!$, in which case it is the salient special case of \cite[Proposition~2.5]{CassvdHScholbach:MotivicSatake}. (It can also be proven more directly as a consequence of the vanishing $q_! j_* p^* = 0$, where $\A^1 \stackrel j \gets \GmX S \stackrel p \r S$.)
\xpf

\subsubsection{Complements on mixed Tate motives}
\label{sect--complements MTM}

Whenever we want to consider a motivic t-structure, we will need to assume that $S$ satisfies the Beilinson--Soulé vanishing condition, as discussed in \cite[§2.3]{CassvdHScholbach:MotivicSatake}.
For example, this is the case for $S= \Spec \Z$, $\Spec \Q$ or $\Spec \Fp$.

\lemm
\thlabel{MTM Grothendieck}
(\cite[Corollary~3.2.6]{RicharzScholbach:Intersection} and \cite[Lemma~2.15]{CassvdHScholbach:MotivicSatake})
Let $S$ be a scheme satisfying the Beilinson--Soulé vanishing condition.
Let $(X, X^\dagger)$ be an ind-scheme of ind-finite type over $S$ endowed with an admissible Whitney--Tate stratification (cf.~\cite[Definition~2.9]{CassvdHScholbach:MotivicSatake}, if each stratum $X_w$ is isomorphic to some $\A^m_S \x_S \GmX S^n$, the stratification is admissible).
Then $\DTM(X, X^\dagger)$ carries the so-called \emph{motivic t-structure} whose aisle $\DTM(X, X^\dagger)^{\le 0}$ (in cohomological notation) is compactly generated (under colimits) by $\iota^w_! \Z(k)[\dim X_w]$, where $\iota_w : X_w \to X$ are the inclusions of the strata $X_w$ (which are by assumption smooth over $S$).
We denote its heart by $\MTM(X, X^\dagger)$ or just $\MTM(X)$.
\xlemm

\exam
\thlabel{kernel n-multiplication in MTM}
For $n \ne 0$, $K := \ker (\Z \stackrel n \r \Z)$ (kernel in $\MTM(S)$) is trivial. Indeed, it lies a priori in (cohomological) degrees $[-1,0]$, and it suffices to see $H := \Hom(\Z(k), K)=0$. We have an exact sequence $0 = \Ext^{-1}(\Z(k), K) \r H \r \Hom(\Z(k), \Z) \stackrel n \r \Hom(\Z(k), Z)$, and the latter group vanishes if $k > 0$ for any base scheme $S$, and for $k < 0$ as a consequence of Beilinson--Soulé vanishing \cite[Lemma~3.4]{Spitzweck:Mixed}. For $k = 0$, we have $\Hom(\Z, \Z) = \Z^{\pi_0(S)}$ on which multiplication by $n$ is injective.
\xexam

The following statement essentially appears in \cite[Corollary~3.2.24]{RicharzScholbach:Intersection}. We include a more detailed proof for the convenience of the reader.

\lemm
\thlabel{MTM equivariant}
Let $1 \r U \r H \r G \r 1$ be an extension of $S$-group schemes with $U$ being split pro-unipotent  \cite[Definition~A.4.5]{RicharzScholbach:Intersection}, and $G$ split reductive.
We consider the trivial action of $H$ on $S$ and the resulting category $\DTM(H \setminus S) := \DTM_H(S) \subset \DM(H \setminus S)$ of $H$-equivariant motives whose underlying non-equivariant motive on $S$ is a Tate motive.
\begin{enumerate}
\item \label{item--recap}
The category $\DTM_H(S)$ carries a compactly generated t-structure such that the forgetful functor $\DTM_H(S) \r \DTM(S)$ is t-exact. 
\item \label{item--recap more}
The category $\DTM_H(S)$ (resp.~its subcategory $\DTM_H(S)^{\le 0}$, resp.~the heart, denoted $\MTM_H(S)$) is generated under colimits and shifts (resp.~under colimits, resp.~under extensions, direct sums and cokernels) by $\Z(k)$, $k \in \Z$ (regarded as an $H$-equivariant motive via the trivial action).
\end{enumerate}
\xlemm

\pf
By homotopy invariance for $\DM$, forgetting $U$-equivariance induces an equivalence $\DTM_G(S) \cong \DTM_H(S)$, 
so we will assume $H$ is split reductive.
Statement \eqref{item--recap} is then a recap of 
\cite[Lemma 2.26]{CassvdHScholbach:MotivicSatake}.
For \eqref{item--recap more} we note that the adjunction $\av_H : \DM(S) \rightleftarrows \DM(H \setminus S) : u^!$ restricts to an adjunction between $\DTM_H(S)$ and $\DTM(S)$.
Indeed, $u^! \av_H \Z = h_! h^! \Z$ (where $h : H \r S$ is the structural map).
This object lies in finitely many cohomological degrees in $\DTM(S)$: if $H = \bigsqcup_w H_w$ is the Bruhat decomposition by cells of the form $H_w = \A^{n_w} \x \Gm^{m_w}$, this object is an extension of objects of the form $(\Z \oplus \Z(1)[1])^{\otimes m_w}$.
Since $u^!$ preserves colimits, $\av_H \Z \in \DTM_H(S)$ is a compact object. The Tate twists of this object generate $\DTM_H(S)^{\le 0}$ (and also $\DTM_H(S)$ if we additionally include shifts).
We therefore have a compactly generated t-structure, whose heart $\MTM_H(S)$ is again Grothendieck abelian.
The forgetful functor induces an equivalence $\MTM_H(S) \cong \MTM(S)$ \cite[Corollary~3.2.21]{RicharzScholbach:Intersection} (applicable since $H$ is fiberwise connected over $S$), and the latter category is generated by $\Z(k)$, $k \in \Z$.
\xpf

The following lemma, proved in \cite{RicharzScholbach:IntersectionCorrigendum}, adds a necessary boundedness condition that was missing in \cite[Lemma~3.2.8]{RicharzScholbach:Intersection}.

\lemm\thlabel{conservativity}
Consider the $\ell$-adic realization functor (for a prime $\ell$ invertible on $S$) and the reduction functor \cite[(3.3)]{EberhardtScholbach:Integral}
$$\rho_\ell : \DM(S, \Q) \r \D_\et(S, \Ql),$$
$$\rho_\red : \DM(S, \Q) \r \DMr(S, \Q).$$
\begin{enumerate}
	\item \label{conservative BS} Their  restriction to $\DTM(S, \Q)^\bd$ (Tate motives that are bounded in the t-structure) is t-exact and conservative for any $S$ satisfying the Beilinson--Soulé vanishing condition.
	\item \label{conservative Fq} For $S = \Spec \Fq$ or $\Spec \Fpq$, their restriction to $\DTM(S, \Q)$ is t-exact and conservative.
\end{enumerate}
\xlemm

\subsection{Local and colocal objects; quotient categories}
\label{sect--localizations}

Throughout §\ref{sect--localizations}, we consider the following abstract situation.

\situ
\thlabel{notation localization}
Let $\iota : C \subset D$ be a (fully faithful exact) inclusion of a stables \ii-categories.
We let 
${}^\bot C \subset D$ be the full subcategory of $C$-colocal objects (i.e., $d \in D$ such that $\Map_D(d, c)=0$ for all $c \in C$). We write $\nu : {}^\bot C \r D$ for the inclusion.
Dually, $C^\bot \subset D$  denotes the full subcategory of $C$-local objects.
\xsitu

\lemm
\thlabel{(co)localization nonsense}
(\cite[Proposition~4.9.1]{Krause:Localization}, see also \cite[§4]{GallauerPepinLehalleur:Exponentiation})
Suppose that the inclusion $\iota$ admits a left adjoint $\iota^L$.
Then the inclusion $\nu$ admits a right adjoint 
\begin{align}\nu^R & = \fib(\id \r \iota \iota^L).\label{nu R}
\end{align}
The adjunction $(\nu, \nu^R)$ exhibits ${}^\bot C$ as a localization $D/C$.
We will refer to this by saying that the following composite is an equivalence:
\begin{equation}
  {}^\bot C \subset D \r D/C.\label{bot C D/C}
\end{equation}
\xlemm

\pf
The formula \eqref{nu R} is readily confirmed. 
Let $W$ be the morphisms $f$ in $D$ such that $\fib f \in C$. By \Kero{02G3}, we need to check this is colocalizing in the sense of \Kero{02G0}. Properties (1) and (2) there are clear. 
Property (3) there holds since the left term in the cofiber sequence $\nu^R(d) \r d \r \iota \iota^L(d)$ is $W$-colocal, and the right one lies in $C$.
\xpf

In the sequel, we use this description of the localization to propagate various categorical properties from $D$ to $D/C$.

\lemm
\thlabel{compact generation}
Suppose $C$ and $D$ are compactly generated and $\iota$ and $\iota^L$ preserve compact objects.
Then ${}^\bot C$ and $D/C$ are compactly generated and both functors in \eqref{bot C D/C} preserve compact objects.
\xlemm

\pf
In this situation, \thref{(co)localization nonsense} applies to $C^\comp \subset D^\comp$. Then use the equivalence 
$\Ind(D^\comp / C^\comp) = \Ind(D^\comp) / \Ind(C^\comp) = D/C$.
\xpf

\lemm
\thlabel{(co)local preserved}
Suppose we are given a diagram at the left below, where $\iota$ and $\iota'$ admit left adjoints, and $f$ admits a right adjoint $f^R$. 
We suppose 
\begin{align}f^R(C') \subset C.\label{condition (co)local}\end{align} 
Then $f$ restricts to a functor $f : {^\bot C} \r {^\bot C'}$ and the diagram at the right below commutes:
$$
\xymatrix{
C \ar@{^{(}->}[r]^\iota \ar[d]_{g := f|_C} & D \ar[d]_f \\
C' \ar@{^{(}->}[r]^{\iota'} & D' \ar@/_1pc/[u]_{f^R}
} 
\ 
\leadsto 
\
\xymatrix{
^\bot C \ar@/^1pc/[rr]^\cong \ar[d]^f \ar[r] & D \ar[d]^f \ar[r] & D / C \ar[d]^{\ol f} \\
^\bot C' \ar@/_1pc/[rr]_\cong  \ar[r] & D' \ar[r] & D'/C'.
}$$
In particular:
\begin{itemize}
   \item if $f$ is fully faithful, then so is $\ol f$.
   \item the functor $\ol {f^R}$ induced by $f^R$ is a right adjoint of $\ol f$.
\end{itemize} 
\xlemm

\rema
\thlabel{colocal condition explained}
\begin{enumerate}
\item 
The condition \eqref{condition (co)local} is satisfied if the left hand diagram above is vertically right adjointable, i.e., if $g := f|_C : C \r C'$ admits a right adjoint $g^R$ and if the natural map $\iota g^R \r f^R \iota'$ is an isomorphism. It is also satisfied if the diagram is horizontally left adjointable, i.e., we have an isomorphism
\begin{equation}\iota'^L f \stackrel \cong \r g \iota^L. \label{condition left adjointable}
\end{equation}
\item \label{dual colocal preserved}
\thref{(co)local preserved} admits a dual version, in which ``left'' and ``right'' are interchanged throughout, and ${}^\bot C$ is replaced by $C^\bot$.
\end{enumerate}

\xrema

\lemm
\thlabel{t-structure colocal}
In \thref{notation localization}, suppose $D$ is equipped with a t-structure such that $\nu \circ \nu^R : D \r D$ is t-exact. Let $({^\bot C})^{\le 0} \subset {^\bot C}$ be the full subcategory spanned by the objects in $\nu^R(D^{\le 0})$, and likewise with ``$\ge 0$.'' 
\begin{enumerate}
   \item This defines a t-structure on ${^\bot C}$.
   \item Both $\nu$ and $\nu^R$ are t-exact.
   \item The equivalence ${^\bot C} \cong D/C$ endows $D/C$ with a t-structure such that the canonical functor $D \r D/C$ is t-exact.
 \end{enumerate}
\xlemm

\pf
This is immediate from the definitions. 
\xpf

\rema
All of the above statements admit a dual version (since $(D/C)^\opp = D^\opp / C^\opp$).
Thus, if $\iota$ admits a \emph{right} adjoint, the inclusion $\nu': C^\bot \r D$ admits a left adjoint $\nu'^\L = \cofib(\iota \iota^R \r \id)$, which yields an equivalence
$$C^\bot \subset D \r D/C.$$

In the case where $\iota$ has both a right \emph{and} a left adjoint, we have two equivalences
\begin{equation}
	C^\bot \stackrel \cong \gets D/C \stackrel \cong \r {}^\bot C.\label{C bot blah}
\end{equation}
This composite is \emph{not} usually compatible with the inclusions into $D$.
Instead, it sends a $C$-local object $d$ to the $C$-colocal object $\fib(d \r \iota \iota^L (d))$.

Regarding t-structures, \thref{t-structure colocal} admits a dual version, provided that $\nu' \circ \nu'^\L$ is t-exact. In the case of both adjoints, both equivalences in \eqref{C bot blah} are t-exact.
\xrema

\subsection{Definitions of exponential motives and basic examples}

\defi
\thlabel{definition DMexp}
Let $X$ be a prestack with an action of the additive group $\Ga$.
The category of \emph{exponential motives} is defined as 
$$\DMe(X) := \DM(X) / \DM(X / \Ga).$$
We also consider the full subcategory of $\DM(X / \Ga)$-colocal objects and local objects, respectively:
\begin{align}
\DMesh (X) & := {}^\bot \DM(X/\Ga) = \ker \left ( \DM(X) \stackrel{\av^{\Ga}_!} \r \DM(X/\Ga) \right ) \label{DMesh} \\
\DMest (X) & := \DM(X/\Ga)^\bot = \ker \left ( \DM(X) \stackrel{\av^{\Ga}_*} \r \DM(X/\Ga) \right ). \nonumber
\end{align}
\xdefi

To justify this definition, note that the forgetful functor
\begin{align} u^! : \DM(X / \Ga) \r \DM(X)\label{u! Ga G}
\end{align}
is fully faithful by homotopy invariance for $\DM$. It admits a left adjoint $\av_!$ and a right adjoint $\av_*$ (\thref{av}).
We can therefore apply the general discussion from \refsect{localizations}. \thref{(co)localization nonsense} supplies equivalences of stable \ii-categories
\begin{align}
\DMesh(X) \stackrel \cong \r \DMe(X) \stackrel \cong \gets \DMest(X).\label{Kir* Kir!}
\end{align}

\rema
\thlabel{Kir Ga}
This definition is the natural generalization of the work of Gallauer--Pepin Lehalleur \cite{GallauerPepinLehalleur:Exponentiation}, which amounts to the special case $X = \GaX S$ (with usual $\Ga$-action, with $S$ some scheme). In this case $\DMe(X) = \DM(\GaX S) / \DM(S)$.
(This relates to earlier work by Fresán and Jossen \cite{FresanJossen:Exponential} who introduced a category of exponential Nori-type motives.)
In this case, the equivalence \eqref{Kir* Kir!} appears in \cite[Lemma~4.4]{GallauerPepinLehalleur:Exponentiation}.
Writing $q : \GaX S \r S$ for the structural map, we have in this case $u^! = q^!$, and $\DMesh(X) = \ker q_!$, $\DMest(X) = \ker q_*$.
In particular, this shows that composite of the two equivalences in \eqref{Kir* Kir!} is \emph{not} compatible with the inclusions into $\DM(X)$.
We think of $\DMexp$ as being the choice-free incarnation of exponential motives, while $\DMesh$ and $\DMest$ serve as computational tools, see also \thref{DMexp better}.
\xrema

\exam
\thlabel{j! j* DMexp}
Consider the usual complementary inclusions $\GmX S \stackrel j \r \A^1_S \stackrel i \gets S$.
In $\DMe(\A^1)$, the boundary maps of the two localization cofiber sequences below are isomorphisms:
$$j_* \Z \stackrel[\text{in }\DMe]\cong \r i_* i^! \Z[1] \cong i_* \Z(-1)[-1] \stackrel[\text{in }\DMe]\cong \r j_! \Z(-1).$$
\xexam

Before adapting the above definition to the case of prestacks with a $\Ga \rtimes \Gm$-action, we recall a generality concerning quotients by semidirect products.
For a monoid object $M$ in a complete and cocomplete \ii-category $C$ (such as the category of prestacks), we write $\Obj_H C$ for the \ii-category of $H$-objects in $C$.
There are adjunctions (with left adjoints on top)
$$C \mathrel{\mathop{\rightleftarrows}^{H \x -}_{U}} \Obj_H C \stackrel[\text{triv}]{(-) / H} \rightleftarrows C.$$
If $G$ acts by monoid homomorphisms on another monoid object $H \in C$, then $H$ is a monoid object in $\Obj_G C$.
There is an equivalence $\Obj_{H \rtimes G}(C) = \Obj_{H}(\Obj_G C)$. 
Using this and the right hand adjunction above (for $H \rtimes G, H$, and $G$), we obtain an isomorphism, for any object $X \in C$ with an action of $H \rtimes G$,
\begin{equation}
(X/H) /G = X / (H \rtimes G). \label{quotient semidirect}
\end{equation}

We apply this to $H=\Ga$ acted upon by $G= \Gm$ in the usual way. For a prestack $X$ with a $\Ga \rtimes \Gm$-action, the counit map of the left hand adjunction gives a $\Ga \rtimes G$-equivariant map $p := \act_{\Ga} : \Ga \x X \r X$.
We write $\ol p := p / \Ga$ and $\ol {\ol p} := p / (\Ga \rtimes G)$.
By homotopy invariance for $\DM$ and \thref{functoriality equivariant}, $\ol {\ol p}^!$ is fully faithful, and it admits both a right and a left adjoint which we again denote by $\av_*^{\Ga}$ and $\av_!^{\Ga}$, respectively.
By construction, $\av_!^{\Ga} = \ol {\ol p}_!$.
This entitles us to the following definition.

\defi
\thlabel{DMexp G-equivariant}
Suppose $\Ga \rtimes \Gm$ acts on a prestack $X$.
We consider the induced $\Gm (\subset \Ga \rtimes \Gm)$-action on $X$ and the associated category of $\Gm$-equivariant motives, $\DM(X/\Gm)$.
We then consider the category of \emph{$\Gm$-equivariant exponential motives}, and the full subcategories of (co)local objects
\begin{align}
\DMe(X/\Gm) & := \DM(X/\Gm) / \DM(X / \Ga \rtimes \Gm), \label{eqn DMe X/G} \\
\DMesh (X/\Gm) & = \ker \left ( \DM(X/\Gm) \stackrel{\av^{\Ga}_!} \r \DM(X/\Ga \rtimes \Gm) \right ) \label{DMesh G} \\
\DMest (X/\Gm) & = \ker \left ( \DM(X/\Gm) \stackrel{\av^{\Ga}_*} \r \DM(X/\Ga \rtimes \Gm) \right ). \nonumber
\end{align}
\xdefi

\rema
\thref{DMexp G-equivariant} is a transfer of the notion of the Kirillov category due to Gaitsgory--Lysenko  \cite[§A]{GaitsgoryLysenko:Metaplectic} to motivic sheaves.

\thref{DMexp G-equivariant} is not a special case, but rather a variant of \thref{definition DMexp}: $\Ga$ does not usually act on the prestack quotient $X / \Gm$ (while $\Gm$ does act on $X / \Ga$), so that it is not possible to apply \thref{definition DMexp} to $X / \Gm$. We therefore trust that no confusion will arise.

As before, we have natural equivalences
$$\DMesh(X/\Gm) \stackrel \cong \r \DMe(X/\Gm) \stackrel \cong \gets \DMest(X/\Gm).$$
We have $\DMesh(X/\Gm) = \{M \in \DM(X/\Gm), u^! M \in \DMesh(X)\}$ and likewise for $\DMest$.
Indeed, by \thref{functoriality equivariant} and \eqref{u! av*}, there is a commutative diagram where the two vertical forgetful functors are conservative:
$$\xymatrix{
\DM(X/\Gm) \ar[d]^{u^!}  \ar[rr]^(.5){\av_!^{\Ga} = \ol p_!} & & \DM(X / \Ga \rtimes \Gm) \ar[d]^{u^!} \\
\DM(X) \ar[rr]^(.5){\av_!^{\Ga} = \ol {\ol p}_!} & & \DM(X/\Ga)
}$$
\xrema

\toyexam
\thlabel{DMexp A^1}
For the usual action of $\Ga \rtimes \Gm$ on $\A^1$, we compute the category of $\Gm$-equivariant exponential motives $\DMe(\A^1 / \Gm)$  using hyperbolic localization (\thref{hyperbolic localization}): 
$$\av_!^{\Ga} = \ol q_! = \ol i^!: \DM(\A^1 / \Gm) \r \DM(S / \Gm).$$
Here, in addition to the notation from \thref{j! j* DMexp}, $q : \A^1 \r S$ is the structural map and overlines denote the respective $\Gm$-quotients.
The kernel of this functor, $\DMesh(\A^1 / \Gm)$, is thus equivalent to the full subcategory $\ol j_* \DM(\Gm / \Gm) \subset \DM(\A^1 / \Gm)$.
Therefore as an abstract category it is equivalent to $\DM(S)$, which in turn is abstractly equivalent also to $\ol j_! \DM(\Gm / \Gm) = \DMest(\A^1 / \Gm)$.
Note that this assertion contains in particular the fact that 
\begin{equation}
q_! j_* \Z = 0 \text{ and } q_* j_! \Z = 0.
\label{some exponential motives on A1}
\end{equation}  
\xtoyexam

\exam
\thlabel{usage  exponential motives later on}
In \refsect{exponential motives Fl}, we use the above definition as follows:
given an exact sequence $1 \r K \r U \r \Ga \r 1$ of group schemes, with $\Gm$ acting on $U$ by group automorphisms lifting the usual action on $\Ga$, and a prestack $Y$ acted upon by $U \rtimes \Gm$, the prestack $X := Y / K $ is acted upon by $\Ga \rtimes \Gm = (U \rtimes \Gm) / K$. 
(Computing everything pointwise, we can replace prestacks by anima. The surjection $U \r \Ga$ induces an isofibration and therefore a cocartesian fibration after passing to classifying spaces \Kerod{01ER}{023K}. Thus the left adjoint of the restriction functor $\Fun(B\Ga, \Ani) \r \Fun(BU, \Ani)$, which is given by left Kan extension, can be computed by taking homotopy orbits of the action of $BK = \fib(BU \r B\Ga)$ \Kero{02ZM}.)
We can therefore consider $\DMe (X / \Gm) = \DM(Y / (K \rtimes \Gm)) / \DM(Y / (U \rtimes \Gm))$.
\xexam

\subsection{Functoriality for exponential motives}
\lemm
\thlabel{functoriality DMexp}
Let $f : X \r X'$ be a $\Ga$-equivariant map of finite type between two prestacks. (The finite type condition is not needed for the !-functors below.)
\begin{enumerate}
  \item 
  \label{functor Kir}
  All four functors $f^*, f_!, f_*, f^!$ preserve the full subcategories of $\Ga$-equivariant sheaves and therefore descend to $\DMe$ (compatibly with the canonical projection functor $\DM \to \DMe$).
  The functors between the $\DMe$-categories are adjoint functors in the usual way.
  \item 
  \label{functor Kir!} 
  The functors $f_!$ and $f^*$ preserve the $\DMesh$-subcategories. Thus there is a commutative diagram (and likewise with $f^*$)
  $$\xymatrix{
  \DMesh(X) \ar[d]^{f_!} \ar@/^1pc/[rr]^\cong \ar@{^{(}->}[r] & \DM(X) \ar[d]^{f_!} \ar[r] & \DMe(X) \ar[d]^{f_!} \\
  \DMesh(Y) \ar@/_1pc/[rr]_\cong \ar@{^{(}->}[r] & \DM(Y) \ar[r] & \DMe(Y).
  }$$
  Dually, the same statement holds for $f_*$ and $f^!$ and the subcategories $\DMest$ of $\Ga$-local objects.
\end{enumerate}
\xlemm

\pf
\eqref{functor Kir} holds since the forgetful functors $u^!$ as in \eqref{u! Ga G} commute with these four functors by \thref{functoriality equivariant}.

\eqref{functor Kir!} holds by \thref{colocal condition explained} (respectively its dual statement), for the right adjoints $f^!$ and $f_*$ commute with the functors \eqref{u! Ga G} forgetting the $\Ga$-action (\thref{functoriality equivariant}).
\xpf

\rema
\thlabel{DMexp better}
The lemma clarifies why it is more convenient to work with $\DMe$ than with $\DMesh$ or $\DMest$. 
For example, consider the usual ($\Ga$-equivariant) embedding $f : \A^1 \r \P[1]$.
The functor $f_! : \DM(\A^1) \r \DM(\P)$ does not preserve the subcategories $\DMest$ of $\Ga$-local objects: the constant sheaf $\Z_{\P[1]}$ is $\Ga$-equivariant, and for $j : \Gm \subset \A^1$, $j_! \Z$ is $\Ga$-local, but
$\Map_{\DM(\P[1])}(\Z, f_! j_! \Z) \ne 0$.
\xrema

\rema
\thlabel{DM acts on DMexp}
If $X$ is a placid prestack and $Y$ a placid prestack with a $\Ga$-action, we consider the induced $\Ga$-action on $X \x_S Y$. There is a natural functor
$$\bx : \DM(X) \t \DMe(Y) \r \DMe(X \x_S Y)$$
since $F \bx G$ is $\Ga$-equivariant if $G$ is.
This turns $\DMe$ into a $\DM$-module (with respect to the usual lax symmetric monoidal structure on $\DM$ on placid prestacks given by $\bx$).
(However, the exterior product $\bx$ does \emph{not} induce a functor $\DMexp(X) \t \DMexp(Y) \r \DMexp(X \x Y)$.)
\xrema

\rema
\thlabel{Verdier duality DMexp}
For any prestack $X$, we have the Verdier duality functor 
\begin{equation}
\Du : \DM(X)^\opp \r \DM(X), \Du (M) = \IHom(M, \omega_X), \label{Verdier duality}
\end{equation}
where $\omega_X := x^! \Z$, $x : X \r S$.
Since $\Du$ preserves the full subcategory $\DM(X / \Ga) \subset \DM(X)$,
it induces a Verdier duality functor $\Du : \DMexp(X)^\opp \r \DMexp(X)$ on the categories of exponential motives.
For a $\Ga$-equivariant map $f: X \r Y$, we have $\Du_Y f_! = f_* \Du_X$, just as for $\DM$ \cite[Theorem~2.4.50(5)]{CisinskiDeglise:Triangulated}.
For étale motives, Verdier duality is an involution on constructible objects \cite[Theorem~6.2.17]{CisinskiDeglise:Etale}, hence the same is true on the category $\DAe$ of étale exponential motives considered in \refsect{Whittaker motives}.
\xrema

\lemm
\thlabel{DMexp glueing}
Suppose $\Ga$ acts on some scheme $X$. Let $Z \stackrel i \r X \stackrel j \gets U$ be complementary closed and open $\Ga$-equivariant embeddings.
The functors $i_*$ and $j^*$ on the categories of exponential motives (\thref{functoriality DMexp}) form a recollement in the sense of \cite[(1.4.3.1)--(1.4.3.5)]{BeilinsonBernsteinDeligne:Faisceaux}
$$\DMe(Z) \stackrel{i_*} \r \DMe(X) \stackrel{j^*} \r \DMe(U).$$
\xlemm

\pf
We apply the characterization in \cite[§2.3]{CisinskiDeglise:Triangulated}, cf.~also \cite[Remarque~1.4.8]{BeilinsonBernsteinDeligne:Faisceaux}:
the functor $i_* : \DMe(Z) \r \DMe(X)$ is fully faithful by \thref{(co)local preserved}. Moreover, the pair $(j^*, i^*)$ is conservative on the categories of exponential motives: applying \thref{functoriality DMexp}, it suffices to note the same is true for regular motives, and thus for $\Ga$-colocal objects.
\xpf

The following descent statement for exponential motives will be used to prove the comparison with Whittaker sheaves.

\prop
\thlabel{descent DMexp}
Let $f : Y \r X$ be a smooth $\Ga \rtimes G$-equivariant map such that $f^* : \DM(X) \r \DM(Y)$ is conservative.
(This holds whenever $f$ has a not necessarily $\Ga \rtimes G$-equivariant section; see also the discussion in the proof of \thref{descent Whit} below.)
Then there is an equivalence
\begin{equation}
\DMe(X / G) \stackrel[f^*]{\cong} \r \lim \left (
	\DMe(Y / G) \stackrel[p_2^*]{p_1^*} \rightrightarrows 
	\DMe(Y \x_X Y / G) \threerightarrows \dots \right ). \label{descent DMexp eq}
\end{equation}
\xprop

\pf
By assumption and \thref{functoriality DMexp}\eqref{functor Kir!}, we have a conservative functor $f^* : \DMesh(X) \r \DMesh(Y)$.
The functors $f_!$ also preserve the $\DMesh$-subcategories. They commute with *-pullback by smooth base change.
Applying the Beck--Chevalley descent condition, we then get a descent equivalence as in  \eqref{descent DMexp eq}, but for $\DMesh(-)$.
We then conclude using the equivalence in \eqref{Kir* Kir!} (and its functoriality).
\xpf

\rema
\thlabel{not specific to motives}
The preceding considerations and, appropriately formulated, also the results in the sequel, are not at all specific to motives. Instead, the above may be summarized by saying that for a 6-functor formalism, i.e., a lax symmetric monoidal functor $D: \Corr(\Sch^\ft_S) \r \PrLst$ that is homotopy invariant, there is an induced functor $D_{\exp} : \Corr(\Sch^{\ft, \Ga}_S) \r \PrLSt$ defined on the category of correspondences of $\Ga$-equivariant schemes.
For example, this applies to $\ell$-adic or Betti sheaves as well.
\xrema

\subsection{Exponential Tate motives}
\label{sect--exponential Tate motives}
In this section, consider an ind-scheme $X$ with an action of $\Ga$. We also assume $X$ carries a Whitney--Tate stratification $X^\dagger$ (not necessarily stable under the $\Ga$-action).
Recall the category $\DTM(X) := \DTM(X, X^\dagger)$ of stratified Tate motives \cite[Definition and Lemma~3.1.11]{RicharzScholbach:Intersection} and the category $\DTM(X/\Ga)$ of $\Ga$-equivariant stratified Tate motives \cite[Definition~3.1.21]{RicharzScholbach:Intersection}, i.e., the full subcategory of $\DM(X)$ consisting of motives that are $\Ga$-equivariant and lie in $\DTM(X)$.

\defi
\thlabel{DTMexp}
A Whitney--Tate stratification $(X, X^\dagger)$ is called \emph{$\Ga$-averageable} if the two adjoints $\av^{\Ga}_!$ and $\av^{\Ga}_*$ of the forgetful functor $u^! : \DM(X / \Ga) \subset \DM(X)$ map $\DTM(X)$ to $\DTM(X / \Ga)$, i.e., $\Ga$-(co)averaging preserves the property of being stratified Tate.
In this situation we define the category of \emph{exponential stratified Tate motives}
$$\DTMe(X) := \DTM(X) / \DTM(X / \Ga).$$
We also denote it by $\DTMe(X, X^\dagger)$ in order to highlight the choice of the stratification.
We also define $\DTMesh(X)$ and $\DTMest(X)$ to be the colocal objects, resp.~local objects in $\DTM(X)$.
\xdefi

A Whitney--Tate stratification by $\Ga$-stable strata is $\Ga$-averageable. However, we will rather be interested in cases where the strata are not $\Ga$-stable, such as in the \thref{DTMexp A1} and the more substantial \thref{GaAverageable} below.

In \thref{Motivic CasselmanShalika}, we will also use the following variant. 
The reason for only considering the right adjoint in the following definition will become clear in the \thref{DTMexp A1}.
\defi
\thlabel{anti-effective DTMexp}
We say that an \emph{anti-effective Whitney--Tate} \cite[Definition~2.6]{CassvdHScholbach:MotivicSatake} stratification is \emph{anti-effectively $\Ga$-averageable} if it is $\Ga$-averageable and if in addition the right adjoint $\av_*^{\Ga}$ preserves anti-effective Tate motives. We then define
$$\DTMea(X) := \DTMa(X) / \DTMa(X/\Ga).$$
\xdefi

\lemm \thlabel{full exp subcat}
In the situation of \thref{DTMexp} (resp.~\thref{anti-effective DTMexp}), the natural functor $\DTMe(X, X^\dagger) \r \DMe(X)$ is fully faithful (resp.~for the anti-effective categories).
Its image is the presentable stable subcategory generated by $\iota_{w*} \Z(n)$, where $\iota_w: X_w \r X$ is the inclusion of a stratum and $n \in \Z$ (resp.~$n \le 0$).
\xlemm

\pf
For $\DTM$, this holds by \thref{(co)local preserved} and \thref{colocal condition explained} (applied through \eqref{condition left adjointable}), using that the inclusions $\DTM(-) \subset \DM(-)$ admit right adjoints by the adjoint functor theorem.
To prove it for anti-effective Tate motives, we use the dual of \thref{(co)local preserved}, as mentioned in 
\thref{colocal condition explained}\eqref{dual colocal preserved}: the inclusion $\DTMa(X) \subset \DTM(X)$ admits a \emph{left} adjoint since this full subcategory is stable under products by virtue of \cite[Lemma~2.18]{CassvdHScholbach:MotivicSatake}.
These left adjoints commute with the forgetful functor $u^!$ by the assumption of being anti-effectively $\Ga$-averageable.
\xpf

\coro
\thlabel{functoriality DTMexp anti}
If a $\Ga$-equivariant map $f: (X, X^\dagger) \r (Y, Y^\dagger)$ of $\Ga$-averageable Whitney--Tate stratified ind-schemes is such that one of the four functors $f^*$, $f_*$, $f_!$ or $f^!$ preserves $\DTM(-)^{\xanti}$, then it induces a functor on $\DTMexa(-)$.
\xcoro

Recall from \refsect{complements MTM} the basic definitions pertaining to motivic t-structures on stratified Tate motives.
Recall from \cite[Definition~2.9]{CassvdHScholbach:MotivicSatake} the definition of an admissible Whitney--Tate stratification. For example, any cellular stratification is admissible. That being said, we introduce a condition that ensures a well-behaved motivic t-structure on exponential stratified Tate motives.

\defilemm
\thlabel{t-structure DTMexp}
Suppose $S$ satisfies the Beilinson--Soulé vanishing, and the Whitney--Tate stratification on $X$ is admissible. 
Suppose the stratification is $\Ga$-averageable and such that the endofunctor
\begin{eqnarray} \label{endofunctor check t-exactness}
	\fib_!^{\Ga}: \DTM(X, X^\dagger) & \r & \DTM(X, X^\dagger),\\ 
	F & \mapsto & \fib (F \r u^! \av^{\Ga}_! F) = \fib (\act_! 0_! F \r \act_! \pr^! F) \nonumber
\end{eqnarray}
is t-exact. (Here $0 : X \r \Ga \x X$ is the zero section and the map arises by adjunction from $\id = 0^! \pr^! F$.)
Then there is a unique t-structure on the full subcategory (at the left) and the quotient category (at the right) such that the canonical functors
$$\DTMesh(X) \subset \DTM(X, X^\dagger) \r \DTMe(X)$$
are t-exact (with respect to the motivic t-structure on $\DTM(X, X^\dagger)$, cf.~\thref{MTM Grothendieck}).
If the stratification is in addition anti-effectively $\Ga$-averageable, then the t-structure restricts to $\DTMea(X) \subset \DTMe(X)$.
The heart of this t-structure, which we will denote by $\MTMexp^\xanti(X)$ or $\MTMexp^\xanti(X, X^\dagger)$, is a Grothendieck abelian category.
\xdefilemm

\pf
In the notation of \thref{notation localization}, the endofunctor in \eqref{endofunctor check t-exactness} is the composite $\nu \nu^R$. Hence this follows at once from \thref{t-structure colocal}; the final claim holds again by \cite[Remark~1.3.5.21]{Lurie:HA}.
\xpf

\toyexam
\thlabel{DTMexp(A^1)}
\thlabel{DTMexp A1}
Consider $X = \A^1$ with the usual stratification $X^\dagger = \Gm \sqcup S$, so $\DTM(X)$ is generated by Tate twists of $\Z$ and $i_* \Z$.
The coaveraging functor $u^! \av^{\Ga}_! = \act_! p^!$ sends these to $\Z$ and $\Z(1)[2]$, while the averaging functor $u^! \av^{\Ga}_* = \act_* p^*$ gives $\Z$ and $\Z$, respectively.
In particular, the stratification is anti-effectively $\Ga$-averageable (but the coaveraging does not preserve anti-effective Tate motives).
The argument in \thref{DMexp A^1} still applies to (non-equivariant) stratified Tate motives on $\A^1$, i.e., the natural map $q_! F\r i^! F$ is an isomorphism for any $F \in \DTM(\A^1, \Gm \sqcup 0)$.
Therefore we have
$$\DTMesh(\A^1) = \ker q_! = \ker i^! = j_* \DTM(\Gm) \cong \DTMe(\A^1) \cong j_! \DTM (\Gm) = \DTMest (\A^1).$$
The condition in \thref{t-structure DTMexp} is also satisfied. Indeed the endofunctor \eqref{endofunctor check t-exactness} maps the generators 
of the t-structure $i_* \Z$ to $j_* j^* \Z(1)[1]$, and $\Z[1]$ to $0$.
Thus, the functor is right t-exact.
On the other hand, if $M$ lies in (cohomological) degrees $\ge 0$, then $q^! i^! M [-1] \stackrel{\ref{hyperbolic localization}} = q^! q_! M [-1]$ lies in degrees $\ge 1$ since  $i^!$ is (by definition of the t-structure) left t-exact, as is $q^![-1]$.
\xtoyexam

In \refsect{Ga}, we will apply the above paradigm in the following slightly more general situation (cf.~\thref{usage  exponential motives later on}):
given a short exact sequence of pro-smooth group schemes $1 \r K \r U \rtimes \Gm \r \Ga \rtimes \Gm \r 1$, where the right hand map is induced by a $\Gm$-equivariant map $U \r \Ga$,
we consider an ind-scheme $Y$ acted upon by $U \rtimes \Gm$.
The prestack quotient $X := Y/K$ then carries an action of $\Ga \rtimes \Gm$.
In addition, we assume that $Y$ carries an admissible Whitney--Tate stratification $Y^\dagger$ which is stable under the action of $K \rtimes \Gm$ (but not necessarily under the action of the full $U \rtimes \Gm$).
Finally, we assume that $Y = \colim Y_i$ and the $K \rtimes \Gm$-action on $Y_i$ factors through $K_i \rtimes \Gm$, where $K = \lim K_i$ is a presentation as a cofiltered limit of smooth group schemes, as in~\cite[Definition~3.1.26, Proposition~3.1.27]{RicharzScholbach:Intersection}.
We can then consider the category $\DTM(Y / (K \rtimes \Gm))$ of $K \rtimes \Gm$-equivariant stratified Tate motives on $Y$, i.e., by definition the full subcategory of $\DM(X / \Gm)$ whose underlying non-equivariant motive (on $Y$) is stratified Tate.

\defi \thlabel{def--exp}
With the above notation, assume that the stratification on $Y$ is $\Ga$-averageable in the sense that the two adjoints $\av_!^{\Ga}$, $\av_*^{\Ga}$ of the forgetful functor $\DM(Y/(U \rtimes \Gm)) \r \DM(Y / (K \rtimes \Gm))$ preserve stratified Tate motives.
Then we define 
$$\DTMexp(X / \Gm) := \DTMexp(Y / K \rtimes \Gm) := \DTM(Y / (K \rtimes \Gm)) / \DTM(Y / (U \rtimes \Gm)).$$
\xdefi 

\rema\thlabel{DTMexp full subcategory}
Similarly to \thref{full exp subcat}, this is a full subcategory of \(\DM(X/\Gm)\).
\xrema

\toyexam \thlabel{DTMexp equiv A1}
Continuing the discussion of \thref{DTMexp(A^1)}, we have the following description of $\Gm$-equivariant exponential stratified Tate motives on $\A^1$:
$$\DTM(S) = \DTM(\Gm / \Gm) \stackrel [\cong]{j_!} \r \DTMe(\A^1 / \Gm) \stackrel [\cong]{j_*} \gets \DTM(\Gm / \Gm) = \DTM(S).$$
\xtoyexam 

\section{Whittaker motives}
\label{sect--Whittaker motives}
\nota
\thlabel{notation Artin--Schreier}
Let $S$ be a scheme over $\Spec \Fp$.
In this section, all schemes will be over $S$; we denote $\x$ for the fiber product over $S$ and $\Ga := \GaX S$ etc.
For $X/S$, we denote by $\DA(X) := \DA(X, \Lambda)$ the category of étale motives over $X$ with coefficients in a ring $\Lambda$ \cite{Ayoub:Realisation}.
Throughout, we assume that $\Lambda$ contains a primitive $p$-th root of unity $\xi \in \Lambda$, for example $\Lambda = \ol \Q$.
We also assume that $p \in \Lambda^\x$; by \cite[Proposition~A.3.1]{CisinskiDeglise:Etale} this latter condition imposes no loss of generality.
We finally fix a group homomorphism 
\begin{equation}
\chi : \Fp \r \Lambda^\x. \label{chi}
\end{equation}
\xnota 

In this section, we introduce the category of \emph{Whittaker motives} on any prestack $X / S$ with a $\Ga$-action, in analogy to the classical theory of Whittaker sheaves. 
We then show that Whittaker motives agree with the exponential motives whenever the $\Ga$-action extends to a $\Ga \rtimes \Gm$-action (\thref{Whit vs Kir}).

\subsection{The Artin--Schreier motive}

The Artin--Schreier covering
$$f : Y := \Ga \r X := \Ga, x \mapsto x^p - x$$
is a Galois cover with Galois group $G = \Fp$.

\defi
\label{Artin-Schreier motive}
For a character $\chi : G \r \Lambda^\x$ we also denote by $\chi \in \DA(X)$ the object that corresponds, 
under the étale descent equivalence
\begin{eqnarray*}
\DA(X) &=& \lim \left (\DA(Y) \rightrightarrows \DA(Y \x_X Y) \threerightarrows \dots \right )  \\
& = & \lim \left (\DA(Y) \rightrightarrows \DA (G \x Y) \threerightarrows \dots \right ),
\end{eqnarray*}
to $\Lambda \in \DA(Y)$ and the isomorphism (in $\DA(G \x Y)$)
$$\Lambda_{G \x Y} \cong \Lambda_{G \x Y}$$ 
whose restriction to the component $\{g\} \x Y$ is multiplication with $\chi(g)$ etc.
We will refer to $\chi$ as the \emph{Artin--Schreier motive}.
\xdefi

\lemm 
$\chi$ is an additive sheaf, i.e., there are natural isomorphisms (where $\add : \Ga \x \Ga \r \Ga$ is the addition map, $\sub$ the subtraction, and $i$ is the zero section)
\begin{align}
\add^* \chi & = \chi \boxtimes \chi, \label{add chi} \\
\sub^* \chi & = \chi \boxtimes \chi^{-1}, \label{sub* chi} \\
i^* \chi & = \Lambda.\label{i^* chi} \\
i^! \chi & = \Lambda(-1)[-2].\label{i^! chi}
\end{align}
\xlemm

\pf
There is a commutative diagram
$$\xymatrix{
\dots \ar@<1ex>[r] \ar[r] \ar@<-1ex>[r]& 
(G \x Y)^2 \ar[d]_{\add_G \x \add_Y} \ar@<.5ex>[r] \ar@<-.5ex>[r] &
Y^2 \ar[d]_{\add_Y} \ar[r] &
X^2 \ar[d]_{\add_X} \\
\dots  \ar@<1ex>[r] \ar[r] \ar@<-1ex>[r] &
G \x Y \ar@<.5ex>[r] \ar@<-.5ex>[r]  &
Y \ar[r] &
X,
}$$
where $\add$ denotes the addition map (of the abstract group $G$, and the group schemes $X = Y = \Ga$).
Inspecting the descent data of $\add^*\chi$ and $\chi \boxtimes \chi$ then gives \eqref{add chi}.
The same argument works for the subtraction map (where at the right hand side of \eqref{sub* chi}, $\chi^{-1}$ denotes the sheaf associated to the character $\chi^{-1}$).
The composite $\Ga \stackrel q \r S \stackrel i \r \Ga$ equals $\sub \circ \Delta$, so that $q^* i^* \chi = \Delta^* \sub^* \chi = \chi \t \chi^{-1} = \Lambda = q^* \Lambda$. 
The homotopy invariance for $\DA$ now implies the isomorphism \eqref{i^* chi}.

The isomorphism \eqref{i^! chi} follows similarly, using that $\add$ is smooth, so that $\add^! \chi = \add^* \chi (1)[2]$.
\xpf

\lemm
\thlabel{q chi}
The (co)averaging with respect to the $\Ga$-action on $\A^1$ is given by 
\begin{align}
\av^{\Ga}_! \chi = q_! \chi & = \left \{ \begin{tabular}{ll} 0 & $\chi \ne 1,$ \\ $\Lambda(-1)[-2]$ & $\chi = 1$, \end{tabular} \right . \label{p! chi} \\
\av^{\Ga}_* \chi = q_* \chi & = \left \{ \begin{tabular}{ll} 0 & $\chi \ne 1,$ \\ $\Lambda$ & $\chi = 1.$ \end{tabular} \right . \label{p* chi}
\end{align}
\xlemm

\pf
Since $\Lambda$ is by assumption a $\Z[1/p]$-algebra that contains a primitive $p$-th root of unity, we have an isomorphism $\bigoplus_\chi \chi = f_! \Lambda$. Applying $q_!$ and using $f q = q$ and homotopy invariance for $\DA$ gives 
$$q_! \bigoplus_{\chi} \chi = (q \circ f)_! \Lambda = q_! \Lambda = \Lambda(-1)[-2] = \Lambda(-1)[-2].$$
Thus, all summands except for $\chi = 1$ vanish in this expression.
The claim for $q_*$ is shown the same way.
\xpf

\coro
\thlabel{p chi}
For non-trivial characters $\chi$, we have isomorphisms as follows, where $S \stackrel p \gets \Gm \stackrel j \r \A^1$:
$$p_* j^* \chi = \Lambda(-1)[-1] \ \text{ and } \ p_! j^! \chi = \Lambda[-1].$$
\xcoro

\pf
The localization fiber sequence collapses as follows $\Lambda(-1)[-2] \stackrel{\text{\eqref{i^! chi}}} = i^! \chi \r 0 
 \stackrel {\text{\eqref{p! chi}}} =  q_* \chi \r p_* j^* \chi,$
 and dually for $p_! j^! \chi$.
\xpf

\lemm
\thlabel{av chi}
The (co)averaging of a non-trivial character $\chi$ with respect to the $\Gm$-action (cf.~\thref{av}) can be computed as follows, where $\ol j$ is the $\Gm$-quotient of $j : \Gm \subset \A^1$:
\begin{align}
\av^{\Gm}_* \chi & = \ol j_! \Lambda(-2)[-3] \ (\in \DA(\A^1 / \Gm)), \label{av* chi} \\
\av^{\Gm}_! \chi & = \ol j_* \Lambda[-1]. \label{av! chi}
\end{align}
\xlemm

\pf
Throughout this proof, (co)averaging is done with respect to $\Gm$.
We have
$$\ol i^! \av_! \chi \stackrel{\eqref{HL}} = \ol q_! \av_! \chi = \av_! q_! \chi \stackrel{\eqref{p! chi}} = 0.$$
Similarly, using \eqref{i^* chi} instead, we get $\ol i^* \av_* \chi = 0$.
Thus we have isomorphisms
$\av_! \chi \stackrel \cong \r \ol j_* \ol j^* \av_! \chi$ and 
$\ol j_! \ol j^! \av_* \chi \stackrel \cong  \r \av_* \chi$.

We now compute $\ol j^* \av_* \chi = \av_* j^* \chi$ where the $\av_*$ at the right denotes the $\Gm$-averaging of sheaves on $\Gm$, which is given by $p_*(-1)[-2]$. Dually, $\av_! = p_!$.
Hence \thref{p chi} implies
\begin{align*}
\ol j^* \av_* \chi & = \av_* j^* \chi = p_* j^* \chi(-1)[-2] = \Lambda(-2)[-3], \\
\ol j^* \av_! \chi & = p_! j^* \chi = \Lambda[-1].
\end{align*}
\xpf

\subsection{The Fourier transform}

As a preparation for the discussion of Whittaker motives, we briefly review the Fourier transform with respect to the Artin--Schreier motive $\chi$.

\subsubsection{Definition and first properties of the Fourier transform}

The Fourier transform can be concisely set up using the 2-category of kernels $\Ker_\DA$. Recall that its objects are finite type $S$-schemes $X$, and the category of morphisms from $X$ to $Y$ is given by $\DA(Y \x_S X)$. See \cite[§4.1]{HeyerMann:6-functor} for the full definition.
We will use the canonical lax symmetric monoidal functor $\Psi : \Ker_\DA \r \Cat$ \cite[Proposition~4.1.5(ii)]{HeyerMann:6-functor} which on objects is $\Psi(X) = \DA(X)$ and a morphism, i.e., a ``kernel'' $F \in \DA(Y \x_S X)$ is mapped to the functor $\pi_{Y!}(\pi_X^* - \otimes F) : \DA(X) \r \DA(Y)$.

For a vector bundle $E$ on $S$ with dual $E^\vee$ we use the following notation, where $\mu : E \x_S E^\vee \r \A^1$ the natural evaluation map:
$$\xymatrix{
S & E \ar[l]_q & E \x E^\vee \ar[l]_(.6)\pr \ar[r]^(.6){\pr^\vee}  \ar[d]^{\mu_E}  & E^\vee \ar[r]^{q^\vee} & S \\
& & \A^1.
}$$
For a map $f : E \r F$ of vector bundles, there is a commutative diagram 
\begin{equation}
\xymatrix{
	 F^\vee \x E \ar[r]^{\id \x f} \ar[d]^{f^\vee \x \id} & 
F^\vee \x F \ar[d]^{\mu_F} \\ 
E^\vee \x E \ar[r]^{\mu_E} &
\A^1.}
\label{mu f E F}    
\end{equation}
Under the canonical functor $\Corr(\Sch_S) \r \Ker_\DA$ (denoted $\Phi_{D, S}$ in \cite[Proposition~4.1.5]{HeyerMann:6-functor}), the correspondence $\text{``}f_!\text{''} : E = E \stackrel f \r F$ gives rise to a morphism (from $E$ to $F$) in $\Ker_\DA$, again denoted by $\text{``}f_!\text{''}$. Concretely, it is given by $\Gamma_{f!} \Z$, where $\Gamma_f : E \r F \x E$ is the graph of $f$.
Similarly, the correspondence $\text{``}(f^\vee)^*\text{''} : E^\vee \stackrel{f^\vee} \gets F^\vee = F^\vee$ gives a kernel given by $\Gamma_{f^\vee!} \Z \in \DA(F^\vee \x E^\vee)$ which we also denote by $\text{``}(f^\vee)^*\text{''}$.
The commutativity of the diagram \eqref{mu f E F} yields $(\id \x f)^* \mu_F^* \chi = (f^\vee \x \id)^* \mu_E^* \chi$.
Unwinding the definition of the composition in $\Ker_\DA$, this implies
$$\mu_F^* \chi \circ \text{``}f_!\text{''} = \text{``}(f^\vee)^*\text{''} \circ \mu_E^* \chi.$$
This yields a functor from the category of vector bundles on $S$ to the arrow category of the category of kernels
$$c := c_\chi : \Vbdl(S) \r \Arr(\Ker_\DA)$$
that sends a vector bundle $E / S$ to the morphism (in $\Ker_{\DA}$) from $E$ to its dual bundle $E^\dual$ given by $\mu_E^* \chi \in \DA(E \x E^\vee)$ and that sends a morphism $f : E \r F$ to the commutative diagram (in $\Ker_\DA$)
\begin{equation}
\xymatrix{
E \ar[r]^{c(E)} \ar[d]_{\text{``}f_!\text{''}} & E^\vee \ar[d]^{\text{``}(f^\vee)^*\text{''}} \\
F \ar[r]^{c(F)} & F^\vee.}
\label{Fourier fun}
\end{equation}

\defi
The \emph{Fourier transform} (for the character $\chi$) is the composition 
$$T := T_\chi : \Vbdl(S) \stackrel c \r \Arr(\Ker_\DA) \stackrel \Psi \r \Arr(\Cat).$$
\xdefi

\rema
\thlabel{Fourier functoriality}
Evaluated on objects, $T$ sends a vector bundle $E$ to the functor given by convolution with the kernel $\mu^* \chi$, i.e.,
$$T_E := T_{\chi, E} : \DA(E, \Lambda) \r \DA(E^\vee, \Lambda), F \mapsto T(F) := \pr^\vee_!(\pr^* F \otimes \mu^* \chi).$$
For a morphism $f$ as above, the commutativity of \eqref{Fourier fun} implies the formula
\begin{equation}
(f^\vee)^* T_{E, \chi}(M) \stackrel \cong \r T_{F, \chi} (f_! M). 
\label{eqn Fourier functoriality}
\end{equation}

\thlabel{Fourier boxtimes}
The functor $c$ is symmetric monoidal (with respect to the cartesian product of vector bundles). Indeed, the assertion $c(E_1) \t c(E_2) = c(E_1 \x E_2)$ follows from the identity $\mu_{E_1 \x E_2} = \add_{\A^1} \circ (\mu_{E_1} \x \mu_{E_2})$, and the additivity of $\chi$, i.e., \eqref{add chi}.
This implies that $T$ is lax symmetric monoidal. In particular, this includes a natural isomorphism (for any two vector bundles $E_1, E_2$):
\begin{equation}
T_{E_1} (-) \boxtimes T_{E_2} (-) \stackrel \cong \r T_{E_1 \x E_2} (- \boxtimes -).
\label{Fourier boxtimes}
\end{equation}
\xrema

\exam
\thlabel{T rk 0}
If $E = S = E^\dual$ is of rank 0, we have $T_\chi = \id_{\DA(S)}$ by \eqref{i^* chi}.
\xexam

\coro
For any vector bundle $E$, $T_E$ is $\DA(S)$-linear.
In particular, for $M \in \DA(S)$, there is a natural isomorphism
\begin{equation}
	T_E(N) \t q^{\vee*} M = T_E(N \t q^* M).\label{T DM linear}
\end{equation}
\xcoro

In order to state the compatibility of Fourier transform with tensor products, we use the convolution product.
Any (commutative) group scheme $G \in \Sch_S^\ft$ is a (commutative) group object in $\Corr(\Sch_S^\ft)$ in two distinct ways: 1) using the (cocommutative) coalgebra structure in $\Sch_S$ (with comultiplication given by the diagonal $X \r X \x_S X$) together with the obvious inclusion inclusion $(\Sch_S^\ft)^\opp \r \Corr(\Sch_S^\ft)$ and 2) using the group structure and the other inclusion $\Sch_S^\ft \r \Corr(\Sch_S^\ft)$.
As $\DA$ is a lax symmetric monoidal functor, $\DA(G)$ therefore carries two distinct structures of a (commutative) algebra object (in $\PrL$).
The underlying functors $\DA(G) \t \DA(G) \r \DA(G)$ are, respectively
\begin{enumerate}
  \item The usual tensor product
  $$F_1 \t F_2 := \Delta^*(F_1 \boxtimes F_2),$$
  with monoidal unit $\Z$.
  \item The \emph{!-convolution product}, denoted here $\star$, given by 
  \begin{equation}
    F_1 \star F_2 := \add_! (F_1 \boxtimes F_2).\label{convolution}
  \end{equation}
  Here, $\add$ is the sum map on $G$.
  Its monoidal unit is $i_! \Z$, where $i: S \r G$ is the neutral element.
\end{enumerate}
The following is well-known, e.g., mentioned in \cite[Remark~3.9]{BoyarchenkoDrinfeld:Character}.

\coro
\thlabel{T monoidal}
The Fourier transform is symmetric monoidal with respect to the convolution product on $E$ (cf.~\eqref{convolution}) and the usual tensor product on $E^\dual$:
$$T := T_\chi : (\DA(E), \star) \r (\DA(E^\dual), \t).$$
\xcoro

\pf
This follows from the above discussion, since the dual of $\add : E \x E \r E$ and the zero section $i : S \r E$ are, respectively, $\Delta_{E^\vee}$ and the structural map $q^\vee$ of the dual bundle $E^\vee$.
The lax symmetric monoidality of $T$ now implies that $T_E$ is a map of commutative algebra objects in $\PrL$, i.e., a functor of presentably symmetric monoidal categories.
\xpf

\lemm
\thlabel{Fourier transform basic objects}
Let $E = \A^1_S$ be a trivial bundle of rank 1. We have 
\begin{align}
T(1_! \Lambda) & = \chi, \label{T1!}\\
T(i_! \Lambda) & = q^{\vee*} i^* \chi \stackrel{\eqref{i^* chi}} = \Lambda, \label{T0!} \\
T(\Lambda) & \stackrel \cong \r i_! i^* T(\Lambda) = i_! \Lambda(-1)[-2] \text{ for } \chi \ne 1. \label{T Lambda}
\end{align}
\xlemm

\pf
Indeed, more generally than \eqref{T1!} and \eqref{T0!}, we have
$T(a_! \Lambda) = \mu_a^* \chi$ for any point $a \in \A^1_S(S)$, where $\mu_a : E^\vee = \A^1 \r \A^1$ is the multiplication with $a$. This follows from base change and the projection formula.

The computation \eqref{T Lambda} is similar to \cite[Lemma~I.5.9]{KiehlWeissauer}; an analogous statement in the context of the homogeneous Fourier transform appears in \cite[Proposition~8.2.12]{FengKhan:Modularity}.
We briefly review the argument:
the first isomorphism holds by the following computation, where $j : \Gm \r \A^1 = E^\dual$:
$$j^* T_\chi(\Lambda) = j^* \pr^\vee_! \mu^* \chi = \pr^\vee_! \act^* \chi = \pr^\vee_! \pr^* \chi =  p^* q_! \chi \stackrel{\eqref{p! chi}} = 0.$$
We have used the following diagram, where the square is cartesian and $\act$ denotes the $\Gm$-action on $\A^1$:
$$\xymatrix{
\A^1 \x \Gm \ar[r]^{(\act, \id)}_\cong \ar[dr]_{\act}
& \A^1 \x \Gm \ar[r]^(.6){\pr^\vee} \ar[d]^{\pr} & \Gm \ar[d]^p \\
& \A^1 \ar[r]^q & S.
}$$
The second isomorphism in \eqref{T Lambda} (and also the formula \eqref{T0!}) holds since the zero section $i: S \r E$ is dual to the structural map $q^\vee : E^\vee \r S$, so we conclude by \thref{Fourier functoriality}, homotopy invariance and \thref{T rk 0}.
\xpf

\subsubsection{The Fourier transform is an equivalence}

\theo
\thlabel{Fourier transform}
Suppose $\chi \ne 1$.
Then we have
$$T_{\chi} \circ T_{\chi^{-1}} = \id(-1)[-2].$$
Hence we have an \emph{equivalence} of symmetric monoidal \ii-categories
$$T_{\chi^{-1}} : (\DA(E), \star) \stackrel \cong \rightleftarrows (\DA(E^\dual), \t) : T_{\chi(1)[2]}.$$
\xtheo

\pf
The proof that of the Fourier equivalence for $\ell$-adic sheaves in \cite[Theorem~I.5.8]{KiehlWeissauer} uses only homotopy invariance and six functor formalism.
We briefly revisit the argument: $T_\chi T_{\chi^{-1}(1)[2]}$ is the endofunctor on $\DA(E^\vee)$ given by the kernel $w_! (\chi \boxtimes_E \chi^{-1}(1)[2])$, where $w : E^\vee \x E \x_E E \x E^\vee \r E^\vee \x E^\vee$ is the natural projection and $\boxtimes_E$ denotes the exterior product over $E$.
We claim that this kernel is isomorphic to $\Delta^\vee_! \Lambda$, whence the proposition.
This follows from the next diagram, whose left square is cartesian:
$$\xymatrix{
E^\vee \ar[d]^{q^\vee} \ar[r]^{\Delta^\vee} & E^\vee \x E^\vee \ar[d]^{\sub} & \ar[l]_w E^\vee \x E \x E^\vee \ar[d]^{\alpha = \pr_E \x \sub_{E^\vee}} \ar[rr]^{\mu_{12} \x \mu_{23}}  & & \A^1 \x \A^1 \ar[d]^\sub \\
S \ar[r]^{i=0} & E^\vee & E \x E^\vee \ar[rr]^\mu \ar[l]^{\pr^\vee} & & \A^1
}$$ 
We then have $\sub^* \chi = \chi \boxtimes \chi^{-1}$ by \eqref{add chi}, so that
$$
w_! (\mu^* \chi \boxtimes_E \mu^* \chi^{-1}(1)[2]) = w_! \alpha^* \mu^* \chi(1)[2] = \sub^* \pr^\vee_! \mu^* \chi(1)[2] =: \sub^* T_\chi(\Lambda(1)[2]) \stackrel{\eqref{T Lambda}} = \sub^* i_! \Lambda = \Delta^\vee_! \Lambda.$$
\xpf

\subsection{The Whittaker category}
\label{Whittaker category}
We continue working with a coefficient ring $\Lambda$ and a base scheme $S / \Fp$ as in \thref{notation Artin--Schreier}.
We fix a \emph{non-trivial character} $\chi : G = \Fp \r \Lambda^\x$.
We are going to use the Fourier transform equivalence (\thref{Fourier transform}) in order to construct the category of Whittaker motives.
If $\Ga$ acts on a prestack $X$, the category $\DA(X)$ becomes a $\DA(\Ga)$-module. The action is by convolution along the action map $\act : \Ga \x X \r X$, similarly to \eqref{convolution}.

\defilemm
\thlabel{DM Whit definition}
Let $X$ be a prestack with a $\Ga$-action.
\begin{enumerate}
\item
The object $\chi(1)[2] \in \DA(\Ga)$ is an idempotent monoid with respect to the convolution product.
We can therefore define the \emph{Whittaker category} as 
$$\DAW(X) := \DA_{\mathrm{W}, \chi} (X) := \Mod_{\chi(1)[2]} \DA(X, \Lambda).$$
We have an adjunction with a fully faithful forgetful functor $U$:
$$\xymatrix{\DA(X) \ar@<.5ex>[r]^{\chi \star -} & \DAW(X). \ar@{^{(}->}@<.5ex>[l]^U}$$
\item
As a mere full subcategory of $\DA(X)$, the Whittaker category admits also the following descriptions:
\begin{align}
\DAW(X) &= \{\chi(1)[2] \star A \ | \ A \in \DA(X) \}\label{bla1} \\
& = \{B \in \DA(X) \ | \ \exists \text{ an isomorphism } B \cong \chi(1)[2] \star B \} \label{bla2}
\end{align}
 \end{enumerate}
\xdefilemm

\pf
The object $1_! \Lambda$ is an idempotent algebra object in $(\DA(\Ga), \t)$.
Under the (symmetric monoidal!) Fourier equivalence $T_\chi(1)[2]$ it is mapped, by \eqref{T1!}, to $\chi(1)[2]$, which is therefore an idempotent monoid in $\DA(\Ga)$ with respect to $\star$.
The full faithfulness of $U$ and the equalities \eqref{bla1} and \eqref{bla2} are then a generality about idempotent monads \cite[Proposition~4.8.2.10]{Lurie:HA}.
\xpf

\rema
\thlabel{independence base scheme}
The category $\DAW(X)$ is independent of the choice of base scheme $S$, i.e., if $X / S$ and we momentarily write for clarity $\chi_S \in \DA(\GaX S)$ as constructed above, i.e., $\chi_S = f^* \chi_{\Spec \Fp}$, where $f : S \r \Spec \Fp$ is the structural map.
We also write $\star_S$ for the convolution of sheaves on $S$-schemes, i.e.,
$$F \star_S G = \act_{S!} (F \boxtimes_S G),$$
where for $F \in \DA(\GaX S)$ and $G \in \DA(X)$, $F \boxtimes_S G$ is the exterior product over $S$ (it can be computed as the pullback along $\GaX S \x_S X \r \GaX S \x X$ of the exterior product $F \boxtimes G$), and $\act_S : \GaX S \x_S X \r X$ is the action map, which is isomorphic to $\act_{\Fp} : \GaX {\Fp} \x_{\Fp} X \r X$.
Then
$$\chi_S(1)[2] \star_S A = \act_{S!} (\chi_S(1)[2] \boxtimes_S A) = \act_{S!}(\chi_{\Fp}(1)[2] \boxtimes \Lambda_S \boxtimes_S A) = \act_{\Fp}(\chi_{\Fp}(1)[2] \boxtimes A) = \chi_{\Fp}(1)[2] \star A.$$
\xrema

\rema
\thlabel{Whit functoriality}
If $f : X \r X'$ is a $\Ga$-equivariant map, the compatibility of $\bx$ with !-pushforward and *-pullback implies that $f_! \DAW(X) \subset \DAW(X')$ and $f^* \DAW(X') \subset \DAW(X)$. 
\xrema

\lemm
\thlabel{Whit other descriptions}
As a mere full subcategory of $\DA(X)$, the Whittaker category admits the following further descriptions:
\begin{align}
\DAW(X) & = \{C \in \DA(X) \ | \ \exists \text{ an isomorphism } \act^* C \cong \chi \boxtimes C\} \label{bla3} \\
& = \{D \in \DA(X) \ | \ \exists \text{ an isomorphism } \act^! D \cong \chi(1)[2] \boxtimes D\}. \label{bla4}
\end{align}
\xlemm

\pf
Let us write $\DAW'(X)$ for one of the (obviously identical) other full subcategories.
For $D \in \DAW'(X)$, homotopy invariance of $\DA$ gives an $D \stackrel \cong \r \act_! \act^! D \cong \act_! (\chi(1)[2] \boxtimes D) = \chi(1)[2] \star C$. Hence, $\DAW'(X) \subset \DAW(X)$.

To show the converse inclusion, we use the following preliminary observation: given a $\Ga$-equivariant map $\pi : X' \r X$, the functors $\pi^*$ and $\pi_!$ preserve the $\DAW'$-subcategories.
In particular, for $\pi = \act : X' := \Ga \x X \r X$, the full faithfulness of $\pi^*$ implies that $F \in \DAW'(X)$ iff $\pi^* F \in \DAW'(X')$.
The same permanence property also holds for the $\DAW$-subcategory, by \thref{Whit functoriality}.
Thus, for $F' \in \DAW(X')$, it suffices to see that $\act^* F' \in \DAW'(X')$.

We use that $\DAW(X)$ and $\DAW'(X)$ are independent of the base scheme (\thref{independence base scheme}).
We can then apply \thref{Fourier transform} (to $X'$, with the base scheme being $X$):
any object in $\DAW(X')$ is of the form $T_\chi(1_! M) = \chi \t M$, with $M \in \DA(X)$.
Now $\act^* (\chi \t M) = \act^* \chi \t M \stackrel{\eqref{add chi}} = (\chi \boxtimes \chi) \t M$ shows that $\chi \t M$ lies in $\DAW'(X)$.
\xpf

\lemm
\thlabel{Whit right adjoint}
The inclusion $\DAW(X) \subset \DA(X)$ admits a \emph{right} adjoint given by *-convolution with $\chi$, i.e., 
\begin{equation}
		F \mapsto \chi \star^* F := \act_* (\chi \boxtimes F)
		\label{convolution star}
\end{equation}
\xlemm

\pf
In order to see that $\act_* (\chi \boxtimes F)$ lies in $\DAW(X)$, we reduce, as in the proof of \thref{Whit other descriptions}, to the case where $X = \Ga \x S$.
We now use the characterization in \eqref{bla4} and construct an isomorphism $\add^! \add_*(\chi \boxtimes F) = \chi(1)[2] \boxtimes \add_* (\chi \boxtimes F)$ using the cartesian diagram
$$\xymatrix{
\Ga \x \Ga \x \Ga \ar[r]^(.6){\add \x \id} \ar[d]^{\id \x \add} &
\Ga \x \Ga \ar[d]^{\add} \\
\Ga \x \Ga \ar[r]^{\add} &
\Ga
}$$
which yields, using the Künneth formula for *-pushforward and !-pullback \cite[Theorem~2.4.6]{JinYang:Kuenneth}, as well as the additivity of $\chi$, cf.~\eqref{add chi}:
\begin{align*}
\add^! \add_* (\chi \boxtimes F) &  = (\id \x \add)_* (\add \x \id)^! (\chi \boxtimes F) \\
& = (\id \x \add)_* (\add^! \chi \boxtimes F) \\
& = (\id \x \add)_* (\chi \boxtimes \chi (1)[2] \boxtimes F) \\
& = \chi (1)[2] \boxtimes \add_*(\chi \boxtimes F).
\end{align*}
This proves that $\chi \star^* F \in \DAW(X)$.
To see that $\chi \star^* -$ is right adjoint to $U$, we compute for $A \in \DAW(X)$:
$$\Map_X(A, \chi \star^* C) = \Map_{\Ga \x X}(\act^*A, \chi \boxtimes C) = \Map_{\Ga \x X} (\chi \boxtimes A, \chi \boxtimes C).$$
The object $\chi$ is invertible with respect to the usual tensor product in $\DA(\Ga)$. 
Hence we can compute this mapping space as $\Map_{\Ga \x X}(\Lambda \boxtimes A, \Lambda \boxtimes B)$, which is isomorphic to $\Map_X(A, B)$ by homotopy invariance.
\xpf

\syno
\thlabel{synopsis Fourier}
The Fourier transform provides the following dictionary (where left adjoints are depicted at the left of their right adjoints)
$$\xymatrix{
\DAW(\A^1_S) \ar@{^{(}->}[d]|U
& & &
\DA(S)  \ar@{^{(}->}[d]|{1_* = 1_!}
\\
\DA(\A^1_S) \ar@/^1pc/[u]^{\chi \star -} \ar@/_1pc/[u]_{\chi \star^* -}="U" 
& & &
\DA(\A^1_S) \ar@/_1pc/[u]_{1^!} \ar@/^1pc/[u]^{1^*}="V"
\ar@{=>}@/^1pc/ "U";"V"^{T_{\chi(1)[2]}}
\ar@{} "U";"V"|\cong
\ar@{=>}@/^1pc/ "V";"U"^{T_{\chi^{-1}}}
}$$
\xsyno

\coro
\thlabel{Whit four functors}
Let $f$ be a $\Ga$-equivariant map.
Then the four functors $f^*$, $f_*$, $f^!$ and $f_!$ preserve the subcategory $\DAW \subset \DA$ (and form adjoints as usual).
\xcoro

\pf
For $f_!$ and $f^*$, this was noted in \thref{Whit functoriality}.
For $f_*$ and $f^!$ it follows the same way by \thref{Whit right adjoint}.
\xpf

\prop
\thlabel{descent Whit}
Let $f : X \r Y$ be a smooth surjective $\Ga$-equivariant map of finite type (between two prestacks).
Then there is an equivalence
\begin{equation}
\DAW(X) \stackrel[f^*]{\cong} \r \lim \left (\DAW(Y) \stackrel[p_2^*]{p_1^*} \rightrightarrows \DAW(Y \x_X Y) \dots \right ).
\end{equation}
\xprop

\pf
We immediately reduce to $Y$ being a scheme.
Given the conservativity of $f^* : \DA(Y) \r \DA(X)$ \cite[Proposition~3.3.33]{CisinskiDeglise:Triangulated}, this is proven the same way as \thref{descent DMexp}: the descent for $\DA(-)$ implies it also for $\DAW(-)$ since all functors in sight preserve that subcategory (\thref{Whit four functors}).
\xpf

\subsection{Exponential motives vs.~Whittaker motives}

Using the setup from \thref{notation Artin--Schreier}, we will now prove a motivic variant of \cite[Proposition~A.3.2]{GaitsgoryLysenko:Metaplectic}: in the presence of a $\Ga \rtimes \Gm$-action, Whittaker motives are equivalent to exponential motives.
It is worth noting that the latter category works over any base scheme and with integral coefficients.
We write $\DAe$ etc. for exponential étale motives, i.e., the theory obtained by replacing $\DM$ by $\DA$ in \refsect{exponential motives}.

\theo
\thlabel{Whit vs Kir}
Suppose $\Ga \rtimes \Gm$ acts on a prestack $X$.
Let $\chi$ be a \emph{non-trivial} character.
Then there is a commutative diagram as follows, with the top horizontal functor being an equivalence:
$$\xymatrix{
\DAW(X) \ar@{^{(}->}[d] \ar[r]^(.4)\cong & 
\DAesh(X / \Gm) \ar@{^{(}->}[d] \ar[dr]^\cong \\
\DA(X) \ar[r]_{\av^{\Gm}_!} & 
\DA(X / \Gm) \ar[r]^\can &
\DAe(X / \Gm).
}$$
In particular, the category $\DAW(X)$ is independent (as an abstract category, albeit not as a full subcategory of $\DA(X)$) of the choice of the non-trivial character $\chi$.
\xtheo

\pf
We first check that $\av^{\Gm}_!$ maps $\DAW(X)$ to the full subcategory $\DAesh(X /\Gm)$, i.e., that $\av_!^{\Ga} \av_!^{\Gm} (\chi \star F) = 0$.
It suffices to show that $\chi \star F$ is colocal w.r.t.~the full subcategory $\DA(X / \Ga) \subset \DA(X)$.
For any $\Ga$-equivariant sheaf $G$, we compute (where $\pr = q \x \id : \Ga \x X \r X$ is the projection map)
\begin{align}
\Map_{\DA(X)}(\chi \star F, u^! G) & = \Map_{\DA(\Ga \x X)}(\chi \boxtimes F, \act^! u^! G) \\
& = \Map_{\DA(\Ga \x X)}(\chi \boxtimes F, \pr^! u^! G)  \nonumber \\
& = \Map_{\DA(\Ga \x X)}((q \x id)_! \chi \boxtimes F, u^! G) \nonumber \\
& = \Map_{\DA(\Ga \x X)}(q_! \chi \boxtimes F, u^! G) \label{noteworthy} \\
& = 0 & \text{by \thref{q chi}}. \nonumber 
\end{align}
The isomorphism \eqref{noteworthy} holds by compatibility of $\bx$ with !-pushforward.

To show $\av_!|_{\DAW(X)}$ is an equivalence, we write $G := \Ga \rtimes \Gm$ and consider the natural smooth action map 
$$f = \act: X' := G \x X \r X.$$
It is $G$-equivariant, where $G$ acts on $X'$ only from the left (and trivially on $X$).
By construction, the natural transformation $\DAW(-) \r \DAesh(-)$ is functorial with respect to *- and !-pullback along a smooth map and !-pushforwards.
Using descent for both categories (\thref{descent DMexp}, \thref{descent Whit}, using the unit section of $f$), we reduce the claim for $X$ to the one for $Y_n := X' \x_X \dots \x_X X'$, with $n \ge 1$ factors of $X'$.
There is a $G$-equivariant isomorphism
$$Z := G \x \underbrace{G^{n-1} \x X}_{=:S} \stackrel \cong \r X' \x_X \dots \x_X X',$$
where at the right we have $n \ge 1$ factors $X'$, and on the left hand term $G$ acts on trivially on $G^{n-1} \x X$.

Note that on $Z = \Ga \x \Gm \x S$, $\Ga$ acts from the left, and $\Gm$ acts by multiplication on $\Ga \x \Gm$.
The actions on $S$ are trivial, so in the sequel we often omit $S$ from the notation.
Thus
$$\DAesh(Z/\Gm) = \ker (q_! : \DA(\A^1) \r \DA(S)).$$
Moreover, the functor $\av^{\Gm}_!$ is simply $\act_! : \DA(\A^1 \x \Gm) \r \DA(\A^1)$, where $\act$ is the (usual multiplication) action of $\Gm$ on $\A^1$.

We apply the Fourier transform $T := T_{\chi(1)[2]}$ (\thref{Fourier transform}, \thref{synopsis Fourier}) to $Z$ (with $\GmX S$ being the base scheme; note that $\Ga$ acts trivially on those factors) and also to $Z / \Gm = \A^1$ (with $S$ as the base scheme).
The bottom right square commutes by \thref{Fourier functoriality}, so the top right category is $\ker 0^* = \im j_!$, where $j : \GmX S\subset \A_S^1$.

$$\xymatrix{
\DA(\Gm) \ar@{^{(}->}[d]_{1_!} & 
\DAW(\A^1 \x \Gm) \ar@{^{(}->}[d] \ar[r] \ar[l]^\cong_T & 
\DAesh(\A^1) = \ker (q_!) \ar@{^{(}->}[d]  \ar[r]_\cong &
\ker (0^*) = \DA(\Gm) \ar@{^{(}->}[d]^{j_!} 
\\
\DA(\A^1 \x \Gm) \ar@<-1.5ex>[u]_{(\id \x 1)^!}  & 
\DA(\A^1 \x \Gm) \ar@<.5ex>[r]^{\act_!} \ar[l]^\cong_{T_{\A^1 \x \Gm}} & 
\DA(\A^1) \ar[r]_\cong^{T_{\A^1}} \ar[d]^{q_!} \ar@<.5ex>[l]^{\act^!} &
\DA(\A^1) \ar[d]^{0^*} \\
& & 
\DA(S) \ar[r]_\cong^{T_S = \id} &
\DA(S).
}$$
In order to check that $\act_!|_{\DAW(\A^1 \x \Gm)}$ is an equivalence, we consider the right adjoints and prove that $(\id \x 1)^! \circ T_{\A^1 \x \Gm} \act^! T_{\A^1}^{-1} \circ j_!$ is the identity functor $\ker 0^* = \DA(\Gm) \r \DA(\Gm)$.

We have a cartesian diagram, where $\act$ denotes the (diagonal) multiplication:
$$\xymatrix{
\A^1 \x \Gm \ar[d]^\cong &
\A^2 \x \Gm \ar[d]^\cong \ar[r] \ar[l] &
\A^1 \x \Gm \ar[d]^\cong & 
\Gm \ar[l]_(.4)\Delta \ar@{=}[d] \\
\A^1 \x \Gm \ar[d]^\act &
\A^2 \x \Gm \ar[d]^\act \ar[r] \ar[l] &
\A^1 \x \Gm \ar[d]^\act &
\Gm \ar[l]_(.4){\id \x 1} \\
\A^1 &
\A^2 \ar[l]_{p_1} \ar[r]^{p_2} &
\A^1. 
}$$
The inverse of the Fourier transform at the bottom row $T_{\chi(1)[2]}^{-1} = T_{\chi^{-1}}$ is given by convolution with $\mu^* \chi^{-1}$, where $\mu : \A^2 \r \A^1$ is the multiplication map.
The isomorphisms from the top row to the middle row are of the form $(\lambda, x) \mapsto (\lambda, \lambda^{-1} x)$, so that the composites to the bottom row are the projection maps.
Using this, we compute the above composite:
$$(\id \x 1)^! T_{\A^1 \x \Gm, \act^* \mu^* \chi} = \Delta^! T_{\A^1 \x \Gm, \pr^* \mu^* \chi} \pr^! = \Delta^! \pr^! T_{\A^1, \mu^* \chi} = j^! T_{\A^1, \mu^* \chi}.$$
Thus
$$(\id \x 1)^! \circ T_{\A^1 \x \Gm} \act^! T_{\A^1}^{-1} \circ j_! = j^! T_{\A^1, \mu^* \chi} T_{\A^1, \mu^* \chi}^{-1} j_! = j^! j_! = \id.$$ 
\xpf

\rema
\thlabel{Whit av*}
\begin{enumerate}
	\item Using the *-convolution (cf.~\eqref{convolution star}), the inverse of the equivalence is the following composite:
\begin{equation}
  \DAesh(X / \Gm) \subset \DA(X / \Gm) \stackrel {u^!} \r \DA(X) \stackrel{\chi \star^* -} \r \DAW(X).\label{adjoint}
\end{equation}
Indeed, by \thref{Whit right adjoint}, this is the right adjoint of $\av^{\Gm}_!|_{\DAW(X)}$.
\item
There is also an equivalence $\DAW(X) \stackrel \cong \r \DAest(X/\Gm)$, given by $\av^{\Gm}_*$.
This is proven along similar lines, replacing the argument in \eqref{noteworthy} by the compatibility of *-pushforward and $\bx$ for $\Fp$-schemes \cite[Proposition~2.1.20]{JinYang:Kuenneth} (which is applicable since $p$ is by assumption invertible in the coefficient ring $\Lambda$).
The following composite is \emph{another} equivalence of categories
$$\DAW(X) \stackrel{\av^{\Gm}_*} \r \DA(X/\Gm) \stackrel \can \r \DAe(X/\Gm).$$
It is related to the one established above by 
\begin{equation}
  \can \circ \av_!^{\Gm} = \can \circ \av_*^{\Gm}  (-1)[-2]. \label{can av! av*}
\end{equation}
This follows from \thref{av chi}, which gives $\av^{\Gm}_* T_{\chi(1)[2]}(1_! \Lambda) = \av^{\Gm}_* \chi(1)[2] \stackrel{\text{\eqref{av! chi}}} = \ol j_! \Lambda(-1)[-1]$ while $\av^{\Gm}_! \chi(1)[2] \stackrel{\text{\eqref{av* chi}}} = \ol j_* \Lambda(1)[1]$.
In $\DAe(X/\Gm)$, as in \thref{j! j* DMexp}, the boundary maps of the localization sequences are isomorphisms as shown below, confirming the assertion:
$$\xymatrix{
& \ol i_* \Lambda(-1)[-2] \ar[r]^\cong \ar@{=}[d] & \ol j_!\Lambda(-1)[-1] \\
\ol j_* \Lambda[-1] \ar[r]^(.4)\cong & \ol i_* \ol i^! \Lambda = \ol i_* \Lambda (-1)[-2] .
}$$
\end{enumerate}
\xrema

Recall the functoriality of $\DAe$ from \thref{functoriality DMexp}(\ref{functor Kir}), and the functoriality of $\DAW$ from \thref{Whit four functors}.

\coro
\thlabel{Whit vs DMexp functoriality}
For a $\Ga \rtimes \Gm$-equivariant map $f : X \r Y$ of prestacks, the equivalence
$$\can \circ \av^{\Gm}_! : \DAW(-) \stackrel \cong \r \DAe(-/\Gm)$$
established in \thref{Whit vs Kir} is compatible with the four functors $f^*$, $f_*$ (provided that $f$ is schematic), and $f^!$ and $f_!$ (any $f$).
\xcoro

\pf
The canonical functor commutes with all four functors, and $\av^{\Gm}_!$ commutes with $f_!$ and $f^*$ (since their right adjoints do, cf.~\thref{functoriality equivariant}). Likewise $\av^{\Gm}_*$ commutes with $f_*$ and $f^!$.
We conclude using \eqref{can av! av*}.
\xpf

\rema
\thlabel{box product DM acts on DMexp}
The action of $\DA$ on $\DAe$ mentioned in \thref{DM acts on DMexp} admits a variant for the Whittaker category. I.e., $\boxtimes$ restricts to a functor $\DA(X) \t \DAW(Y) \r \DAW(X \x Y)$, as follows from \eqref{bla3}.
By \thref{av} (see also the computation of $\av_! T 1_! (M)$ in the proof of \thref{Whit vs Kir}), the equivalence in \thref{Whit vs DMexp functoriality} is compatible with taking exterior products, i.e.,
$$\av^{\Gm}_! (- \boxtimes \chi \star F) = - \boxtimes \av^{\Gm}_! (\chi \star F).$$
\xrema

\rema
All the arguments in this section carry over to any homotopy invariant six functor formalism (with coefficients in a ring $\Lambda$ as in \thref{notation Artin--Schreier}) in place of $\DA$. In \thref{Equivalences hold}, we will apply it to the formalism of \(\ell\)-adic sheaves.
\xrema

\section{Exponential geometry of affine flag varieties}
\label{sect--Kirillov}
We now construct the stratifications of affine flag varieties to which we will apply the formalism of exponential motives, culminating in \thref{affinestrat}. Additionally, we prove that fibers of certain convolution morphisms admit cellular filtrable decompositions in \thref{fibers of exponential convolution}. All statements in this section are purely geometric; we will add motives to the picture in \refsect{exponential motives Fl}. 

\subsection{Reductive groups}
\label{sect--setup}
Let $G$ be a split reductive group over $\Z$. We assume that $G$ is not a torus. Fix a pinning of $G$,  e.g.~as in \cite[\S 2]{MarshRietsh:Flag}.
Then we have a maximal torus and Borel $T \subset B = U \rtimes T \subset G$.
Let $\Phi \supset \Phi^+ \supset \Delta$,  be the roots, positive roots, and simple roots determined by $B$. Let $X^*(T), X_*(T),$ and $X_*(T)^+$ be the characters, cocharacters, and dominant cocharacters determined by $T$, respectively. Let $W_0 = N_G(T)/T$ be the Weyl group.

For each root $\alpha$ let $U_\alpha \subset G$ be the corresponding root subgroup, and choose a root isomorphism $x_\alpha \colon \Ga \r U_\alpha$. If $\pm \alpha \in \Delta$ then $x_\alpha$ is already determined by the pinning, and since we work over $\Z$ then all the $x_\alpha$ are unique up to signs.
Consider the surjective group homomorphisms
\begin{equation} \label{eq--simple} U \to \prod_{\alpha \in \Delta} U_\alpha \cong \prod_{\alpha \in \Delta} \Ga \xrightarrow{\text{add}} \Ga.
\end{equation} Let $\phi \colon U \rightarrow \Ga$ be the composition, and let $U_0 = \ker \phi$.
Note that $U_0$ depends on the sign choice in the $x_\alpha$ for $\alpha \in \Delta$.
The group $U_0$ contains the non-simple positive root subgroups and is isomorphic to an affine space of dimension $\dim U - 1$.

Let $\hat{\rho}$ be half the sum of the positive coroots; this gives a homomorphism $\hat{\rho} \colon \Gm \r T_{\adj}$, where $T_{\adj}$ is the adjoint torus.
Since $T_{\adj}$ acts on $U$ by conjugation, we can form the semi-direct product $U \rtimes \Gm$ using $\hat{\rho}$.
The action of $\Gm$ on $U$ preserves each root subgroup.
Since $\hat{\rho}$ is also the sum of the fundamental coweights, the resulting action of $\Gm$ on the quotient $\prod_{\Delta} \Ga$ is by the standard scaling action on each factor. 
Hence we can also form the semi-direct product $U_0 \rtimes \Gm$. We define
$$U^{\exp} := U_0 \rtimes \Gm.$$

\subsection{Affine Weyl groups}
We follow similar conventions for affine Weyl groups as in \cite{dHL:Frobenius}, except we change the sign of translation elements to align with the Kottwitz homomorphism.
Let $\mathbf{a}_0$ be the alcove in the apartment $X_*(T) \otimes \mathbb{R}$ which is contained in the chamber associated to $B$ and whose closure contains the origin. The affine roots are $\mathcal{R} = \Phi + \Z$. The positive affine roots $\mathcal{R}^+$ are those of the form $\tilde{\alpha} = \alpha + n$ where $n \geq 1$ and $\alpha \in \Phi$, or $n = 0$ and $\alpha \in \Phi^+$. Equivalently, $\tilde{\alpha}$ is positive if $\restr{\tilde{\alpha}}{\mathbf{a}_0} > 0$. 

Let $W$ be the Iwahori--Weyl group, and for $\lambda \in X_*(T)$ let $t_{\lambda} \in W$ denote the translation element which acts on $X_*(T) \otimes \mathbb{R}$ by $x \mapsto x - \lambda$. We can write an arbitrary element of $W$ as $w = t_{\lambda} v$, where $\lambda \in X_*(T)$ and $v \in W_0$, and we have $W = X_*(T) \rtimes W_0$.
The action of $W$ on $\mathcal{R}$ given by $(w \tilde{\alpha}) (p) = \tilde{\alpha}(w^{-1} p)$ for $p \in X_*(T) \otimes \mathbb{R}$ can then be written as
$$(t_\lambda v) (\alpha+ n) = v \alpha + \langle v\alpha, \lambda \rangle + n.$$

Let $\mathcal{S} \subset W$ be the set of simple reflections in the walls of $\mathbf{a}_0$; this consists of reflections across the the zero loci of the elements of $\Delta$ and those $\alpha + 1 \in \mathcal{R}^+$ such that $-\alpha$ is a highest positive root. We denote the length function on $W$ associated to $\mathcal{S}$ by $l$.
Let $\Omega \subset W$ be the subgroup of elements of length zero; equivalently, $\Omega$ is the normalizer in $W$ of $\mathbf{a}_0$. Then $W= W_{\text{af}} \rtimes \Omega$ where $W_{\text{af}}$ is the affine Weyl group. 

For any facet $\mathbf{f}$ in the closure of $\mathbf{a}_0$ (henceforth called a standard facet), let $W_{\mathbf{f}} \subset W$ be the subgroup which fixes $\mathbf{f}$ pointwise. From now on we identify $W/W_{\mathbf{f}} \subset W$ (resp.~$W_0  \backslash W \subset W)$  with the subset of elements having \emph{minimal} length in their right $W_{\mathbf{f}}$ coset (resp.~\emph{maximal} length in their left $W_0$ coset).

\defi
Let $${_0}W_{\mathbf{f}} : = W_0 \backslash W \cap W/W_{\mathbf{f}}  \subset W$$ be the subset of elements which are simultaneously minimal in their right $W_{\mathbf{f}}$-orbit and maximal in their left $W_0$-orbit. Also, let
$$W^{\exp}_{\mathbf{f}} :=  W/W_{\mathbf{f}} \sqcup {_0}W_{\mathbf{f}}.$$
\xdefi

We emphasize that $W^{\exp}_{\mathbf{f}}$ is a (formal) disjoint union; while an element of ${_0}W_{\mathbf{f}}$ can also be viewed as an element of $W/W_{\mathbf{f}}$, we wish to have two copies of each element of ${_0}W_{\mathbf{f}}$, and we will carefully distinguish which copy we are working with at all times.
We will show in \thref{orbit-structure} below that exponential strata in the partial affine flag variety associated to $\mathbf{f}$ are enumerated by ${_0}W_{\mathbf{f}}$.

\lemm \thlabel{lemm-exampleW}
We explicate the structure of ${_0}W_{\mathbf{f}}$ in the following important cases.
\begin{enumerate}
    \item If $\mathbf{f} = \mathbf{a}_0$ then $W_{\mathbf{a}_0} = \{1\}$ and ${_0}W_{\mathbf{a}_0}$ consists of those elements $w = v t_\lambda \in W$ for $v \in W_0$ and $\lambda \in X_*(T)$ satisfying any of the following equivalent conditions.
    \begin{enumerate}
        \item $w$ is maximal in its left $W_0$-orbit.
        \item $v(-\lambda + \mathbf{a}_0) \subset X_*(T) \otimes \mathbb{R}$ belongs to the Weyl chamber opposite to the one containing $\mathbf{a}_0$.
        \item $\restr{w^{-1} \alpha}{\mathbf{a}_0} < 0$ for all $\alpha \in \Phi^+$ (equivalently, for all $\alpha \in \Delta$).
    \end{enumerate}

    \item If $\mathbf{f} = \mathbf{f}_0$ is the unique facet containing the origin then $W_{\mathbf{f}_0} = W_0$ and ${_0}W_{\mathbf{f}_0}$ consists of those elements $t_{\lambda}$ for $\lambda \in X_*(T)^{++} \subset X_*(T)$ strictly dominant, i.e.~$\langle \alpha, \lambda \rangle > 0$ for all $\alpha \in \Phi^+$ (equivalently, for all $\alpha \in \Delta$).
\end{enumerate}
\xlemm

\pf 
Recall that for $w \in W$, we have $$l(w) = \# \: \{\tilde{\alpha} \in \mathcal{R}^+ \: : \: \restr{w^{-1} (\tilde{\alpha})}{\mathbf{a}_0} < 0\}.$$ For (1), we note that the condition in (a) is the definition of ${_0}W_{\mathbf{a}_0}$. A left $W_0$-orbit in $W$ consists of all elements of the form $w = vt_{\lambda}$, where $\lambda \in X_*(T)$ is fixed and $v$ varies over $W_0$. Write an arbitrary positive affine root as $\tilde{\alpha} = \alpha + n$, where $n \geq 1$ or $\alpha \in \Phi^+$ and $n = 0$. Then for $p \in \mathbf{a}_0$, we have
$$(w^{-1}\tilde{\alpha})(p) = (\alpha+n)(wp) = \langle \alpha, v(-\lambda + p) \rangle + n.$$
For fixed $\lambda$ and $v$ varying over $W_0$, we claim this is negative for the greatest number of $\tilde{\alpha} \in \mathcal{R}^+$ precisely when the condition in (b) is satisfied. Indeed, by $W_0$-equivariance of the pairing, for fixed $n > 1$ and $\lambda$, the number of $\alpha \in \Phi$ for which the above quantity is negative is constant along the $W_0$-orbit. For $n=0$ the number of $\alpha \in \Phi^+$ for which the above quantity is negative does vary along the $W_0$-orbit, it this number is maximal (and equal to the order of $\Phi^+)$ precisely when the condition in (b) is satisfied. This also implies that (b) is equivalent to (c). The fact that the condition in (c) for all $\alpha \in \Delta$ implies the same for all $\alpha \in \Phi^+$ follows since every element of $\Phi^+$ is a positive integer combination of elements of $\Delta$. 

For (2), we consider a right $W_0$-orbit, which consists of elements $w = t_\lambda v$ for $\lambda \in X_*(T)$ fixed and $v$ varying over $W_0$. Then $t_{\lambda} v = v t_{v^{-1}\lambda}$, and it follows from the first part of the proof that this element is maximal in its left $W_0$-orbit if and only if $-\lambda + v \mathbf{a}_0$ belongs to the opposite Weyl chamber. This condition implies in particular that $\lambda$ is dominant. On the other hand, by \cite[Lemma 9.3]{Zhu:Coherence}, if $\lambda$ is dominant then $t_{\lambda}$ is the unique minimal element in its right $W_0$-orbit. Hence, if $w \in {_0}W_{\mathbf{f}_0}$ then $v = 1$ and $w = t_{\lambda}$ for some dominant $\lambda$. The condition that this element is maximal in its left $W_0$-orbit then amounts to strict dominance. 
\xpf

\subsection{Loop groups} Let $LG$ and $L^+G$ be the loop group and positive loop group, defined on $\Z$-algebras by $$LG(R) = G(R(\!(t)\!)), \quad L^+G(R) = G(R[\![t]\!]).$$ For any standard facet $\mathbf{f}$, let $L^+G_{\mathbf{f}} \subset LG$ be parahoric loop group associated to the facet $\mathbf{f}$ (see \cite[\S 4.1.2]{CassvdHScholbach:Central} or \cite[\S 4.2]{RicharzScholbach:Intersection} for more details). In the case $\mathbf{f} = \mathbf{f}_0$ is the facet containing the origin, we have $L^+G_{\mathbf{f}_0} = L^+G$. We denote the Iwahori subgroup by $\calI = L^+G_{\mathbf{a}_0} \subset L^+G$, which is also the preimage of $B$ under $L^+G \r G$, $t \mapsto 0$. For any $\tilde{\alpha} = \alpha + n \in \mathcal{R}$, the associated affine root group $U_{\tilde{\alpha}} \subset LG$ is given by the image of $U_{\alpha} \to LG$, $u \mapsto ut^n$.

Let $\calU \subset \calI$ be the pro-unipotent radical; this is generated by $\mathcal{T}^{>0} := \ker (L^+T \xrightarrow{t \mapsto 0} T)$ and the affine root groups $U_{\tilde{\alpha}}$ for $\tilde{\alpha} \in \mathcal{R}^+$.  We have $\calI = \calU \rtimes T$ and $\calU = \calU^{>0} \rtimes U$, where $\calU^{>0} = \ker (\calU \xrightarrow{t \mapsto 0} U)$. Let $L^iG(R) = G(R[t]/t^{i+1})$ for $i \geq 0$. This leads to a presentation $\calU = \lim_i \calU_i$ for $i \geq 0$, where $\calU_i \subset L^iG$ is the preimage of $U$ under $L^iG \rightarrow G$. Each $\calU_i$ is isomorphic as a scheme to an affine space over $\Z$, the transition maps in $\lim_i \calU_i$ are relative affine spaces, and $\ker(\calU_i \rightarrow \calU_{i-1})$ is isomorphic to a vector group, i.e.,~a product of copies of $\Ga$, for each $i \geq 0$ (see \cite[Lemma 4.2.7]{RicharzScholbach:Intersection} for justification of these claims). In particular, $\ker(\calU \rightarrow U)$ is a split pro-unipotent group as in \cite[Definition A.4.5]{RicharzScholbach:Intersection}.

\defi
By composing the mod $t$ surjection $\calU \r U$ with the composition $\phi$ in \eqref{eq--simple}, we get a surjection 
\begin{equation}\label{map from U to A1}\calU \r \Ga\end{equation}
 with kernel $\calU^{>0} \rtimes U_0$.
Using the action of $T_{\adj}$ on $\calU$ by conjugation and the cocharacter $\hat{\rho} \colon \Gm \r T_{\adj}$, we may form the semi-direct products $\calU \rtimes \Gm$ and $(\calU^{>0} \rtimes U_0) \rtimes \Gm$. 
We define
$$\calU^{\exp} : = (\calU^{>0} \rtimes U_0) \rtimes \Gm.$$
\xdefi

The following lemma describes the structure of $\calU^{\exp}$.

\lemm \thlabel{UexpStructrure}
There is a presentation $\calU^{\exp} = \lim_i \calU_i^{\exp}$ where $\calU_0^{\exp} = U^{\exp}$, each $\calU_i^{\exp}$ is isomorphic as a scheme to the product of an affine space over $\Z$ with $\Gm$, the transition maps are relative affine spaces, and $\ker( \calU_i^{\exp} \rightarrow \calU_{i-1}^{\exp})$ is a vector group for each $i \geq 1$.
\xlemm

\pf
Let
$$\calU_i^{\exp} = (\calU_i^{>0} \rtimes U_0) \rtimes \Gm,$$ where $\calU_i^{>0} = \ker (\calU_i \xrightarrow{t \mapsto 0} U)$ and $i \geq 0$.
Since $U_0$ is an affine space, the claims follow from the analogous facts for $\calU = \lim_i \calU_i$.
\xpf

\subsection{Affine flag varieties}
\label{sect--flagnotation}
For any standard facet $\mathbf{f}$, the partial affine flag variety is the ind-scheme $\Fl_{\mathbf{f}} = (LG/L^+G_{\mathbf{f}})_{\et}$. This agrees with $(LG/L^+G_{\mathbf{f}})_{\Zar}$, which is important when working with Nisnevich motives (see \cite[\S 4.1.4]{CassvdHScholbach:Central} for more details). We denote $\Fl = \Fl_{\mathbf{a}_0}$ and $\Gr = \Fl_{\mathbf{f}_0}$, where $\mathbf{f}_0$ is the unique facet containing the origin. For any other standard facet $\mathbf{f}'$, the left $L^+G_{\mathbf{f}'}$-orbits in $\Fl_{\mathbf{f}}$ are enumerated by $W_{\mathbf{f}'} \backslash W / W_{\mathbf{f}}$, which can be viewed as a subset of $W$ as in \cite[Lemma 1.6]{Richarz:Schubert} with length function $l$.

For each $v \in W_0$, fix a lift $\dot{v} \in N_G(T)(\Z)$. For $\lambda \in X_*(T)$, we have the element $\lambda(t) \in T(\Z(\!(t)\!))$.  Then for any $w \in W$ (or in the subset $W_{\mathbf{f}'} \backslash W / W_{\mathbf{f}}$), write $w = t_{\lambda} \dot{v}$, so that we have the lift $\dot{w} = \lambda(t) \dot{v} \in N_{G(\Z(\!(t)\!))} T(\Z[\![t]\!])$. This choice of lift ensures that
\begin{equation} \label{eq--conj} \dot{w} U_{\tilde{\alpha}} \dot{w}^{-1} = U_{w \tilde{\alpha}} \end{equation} for any $\tilde{\alpha} \in \mathcal{R}$ (see \cite[Eqn.~(3.5)]{dHL:Frobenius}, and note that we use $\lambda(t)$ instead of $\lambda(t^{-1})$ according to our sign convention for $t_{\lambda}$).
Let $\overline{e} \in \Fl_{\mathbf{f}}(\Z)$ be the basepoint. The point $\dot{w} \cdot \overline{e} \in \Fl_{\mathbf{f}}(\Z)$ is independent of the lift $\dot{w}$. 

For $w \in W_{\mathbf{f}'} \backslash W / W_{\mathbf{f}}$, the affine Schubert scheme $\Fl_w(\mathbf{f}', \mathbf{f})$ is the scheme-theoretic image of the map
$L^+G_{\mathbf{f}'} \rightarrow \Fl_{\mathbf{f}}$, $g \mapsto g \cdot \dot{w} \cdot e$. The orbit $\Fl_w^\circ(\mathbf{f}', \mathbf{f})$ is the \'etale sheaf-theoretic image of the same map. The schemes $\Fl_w^\circ(\mathbf{f}', \mathbf{f})$  are cellular by \cite[Lemma 4.3.7]{RicharzScholbach:Intersection}, and their closure relations are determined by the Bruhat order; see \cite[Proposition 4.4]{CassvdHScholbach:Central} for a more precise statement summarizing some of the results in \cite[\S 4.3, \S 4.4]{RicharzScholbach:Intersection}.

If $\mathbf{f}' = \mathbf{f}$, we often shorten the notation and write $\Fl_w^\circ(\mathbf{f}, \mathbf{f}) = \Fl_w^\circ(\mathbf{f})$ and $\Fl_w(\mathbf{f}, \mathbf{f}) = \Fl_w(\mathbf{f})$. If $\mathbf{f} = \mathbf{f}_0$, there is a canonical isomorphism $W_{\mathbf{f}_0} \backslash W / W_{\mathbf{f}_0} \cong X_*(T)^+$ with the monoid of dominant cocharacters. For $\mu \in X_*(T)^+$, we denote
$$\Gr^{\mu} = \Fl_{t_{\mu}}^\circ(\mathbf{f}_0), \quad \Gr^{\leq \mu} = \Fl_{t_{\mu}}(\mathbf{f}_0).$$

\subsection{Exponential orbits}  Our next aim is to show that ${_0}W_{\mathbf{f}}$ parametrizes the $\calU^{\exp}$-orbits in $\Fl_{\mathbf{f}}$. If $w \in {_0}W_{\mathbf{f}}$, we have the lift $\dot{w}$ chosen above by viewing $w$ as an element in $W$. Also, if $\alpha \in \Delta$, recall that we fixed a trivialization $x_{\alpha} \colon \Ga \r U_\alpha$. Then for such $w$ and $\alpha$, we may form the point $x_{\alpha}(1) \cdot \dot{w} \cdot \overline{e} \in \Fl_{\mathbf{f}}(\Z)$. This point depends on $\alpha$ and $x_\alpha$, but as we show below, its orbit under $\calU^{\exp}$ is independent of these choices.

The lemma below implies there is indeed a well-defined action of $\calU^{\exp}$ on $\Fl_{\mathbf{f}}$, which entitles us to define the $\calU^{\exp}$-orbits afterward.

\lemm \thlabel{lemm--UexpAction}
For any $w \in W/W_{\mathbf{f}}$, consider the inclusions $$\Fl^\circ_w(\mathbf{a}_0, \mathbf{f}) \rightarrow \Fl_w(\mathbf{a}_0, \mathbf{f}) \rightarrow \Fl_{\mathbf{f}}.$$
The left action of $\calI = \calU \rtimes T$ on each of these (ind-)schemes factors through an action of $\calU \rtimes T_{\adj}$, and the inclusions are  $\calU \rtimes T_{\adj}$-equivariant. 
\xlemm

\pf
The action of $\calI$ on $\Fl_{\mathbf{f}}$ comes from the \'{e}tale-sheafification of the left action of $\calI \subset LG$ on the presheaf quotient $LG/L^+G_{\mathbf{f}}$. 
Since the center $Z(G)$ acts trivially on the presheaf quotient, and the \'{e}tale and fppf-sheafifications of $LG/L^+G_{\mathbf{f}}$ agree, then we get an action of $\calU \rtimes T_{\adj} = \calU \rtimes (T/Z(G))_{\text{fppf}}$ on $\Fl_{\mathbf{f}}$. 
Next, note that $T_{\adj}$ fixes $\dot{w} \cdot \overline{e}$, 
so that $\Fl^\circ_w(\mathbf{a}_0, \mathbf{f})$ agrees with the \'{e}tale and fppf-sheafifications of the orbit maps for each of $\calI$, $\calU$, and $\calU \rtimes T_{\adj}$ acting on $\dot{w} \cdot \overline{e}$. Similarly, $\Fl_w(\mathbf{a}_0, \mathbf{f})$ agrees with the scheme-theoretic images of each of these orbit maps. Hence we also get an action of $\calU \rtimes T_{\adj}$ on $\Fl_w(\mathbf{a}_0, \mathbf{f})$, and the above inclusions are equivariant.
\xpf

\defi \thlabel{def-KirillovOrbits} Fix a simple root $\alpha_0 \in \Delta$. For $w \in W^{\exp}_{\mathbf{f}}$, we define $\Fl^{\exp}_w(\mathbf{f}) \subset \Fl_{\mathbf{f}}$ to be the \'etale sheaf-theoretic image of the orbit map

\begin{equation}\label{orbit map}\calU^{\exp} \rightarrow \Fl_{\mathbf{f}}, \quad  \begin{cases}
g \mapsto g \cdot \dot{w} \cdot \overline{e}, & \text{if } w \in W/W_{\mathbf{f}},\\
g \mapsto g \cdot  x_{\alpha_0}(1) \cdot \dot{w}  \cdot \overline{e},  & \text{if } w \in {_0}W_{\mathbf{f}}.
\end{cases}
\end{equation}
\xdefi

Note that since $\calU^{\exp} \subset \calU \rtimes T_{\adj}$ and $x_{\alpha}(1) \subset \calU(\Z)$, there is a natural inclusion $$\Fl_w^{\exp}(\mathbf{f}) \rightarrow \Fl_w^{\circ}(\mathbf{a}_0, \mathbf{f}),$$ where we view $w$ as an element of $W/W_{\mathbf{f}}$ to make sense of the $\calI$-orbit $\Fl_w^{\circ}(\mathbf{a}_0, \mathbf{f})$.

We now recall the structure of the the $\calI$-orbits, which agree with the $\calU$-orbits since $T$ fixes $\dot{w} \cdot \overline{e}$ for all $w \in W/W_{\mathbf{f}}$. Since $\dot{w} U_{\tilde{\alpha}} \dot{w}^{-1} = U_{w \tilde{\alpha}}$ for any $\tilde{\alpha} \in \mathcal{R}$ \cite[Eqn.~(3.5)]{dHL:Frobenius}, the stabilizer $\calU_w := (\calU \cap \dot{w}\calI \dot{w}^{-1})$ of $\dot{w} \cdot \overline{e}$ in $\calU$ is
generated by $\mathcal{T}^{>0}$ and the $U_{\tilde{\alpha}}$ such that $\restr{\tilde{\alpha}}{\mathbf{a}_0} > 0$ and  $\restr{w^{-1}\tilde{\alpha}}{\mathbf{a}_0} > 0$.
By omitting the root subgroups in each $\calU_i$ not appearing in $\calU_w$, we get a presentation $\lim_i \calU_{w,i} = \calU_w$, where each $\calU_{w,i} \subset \calU_i$ is isomorphic to an affine space over $\Z$, each transition map is a relative affine space, and $\ker(\calU_{w,i} \rightarrow \calU_{w,i-1})$ is a vector group for each $i \geq 1$ (cf.~\cite[Lemma 4.3.7]{RicharzScholbach:Intersection}).

Note that by definition $\Fl^\circ_{w}(\mathbf{a}_0, \mathbf{f}) = (\calU/ \calU_w)_{\et}$, but in fact  $\Fl^\circ_{w}(\mathbf{a}_0, \mathbf{f})$ agrees with the presheaf quotient $\calU/ \calU_w$. Indeed, the product of affine root groups $\prod U_{\tilde{\alpha}}$ in any order, for $\restr{\tilde{\alpha}}{\mathbf{a}_0} > 0$ and $\restr{w^{-1}\tilde{\alpha}}{\mathbf{a}_0} < 0$, is an affine space of dimension $l(w)$ which is a retraction of $\calU$ and which maps isomorphically onto $\Fl^\circ_{w}(\mathbf{a}_0, \mathbf{f})$ in the quotient.

\lemm \thlabel{orbit-structure} Let $w \in W^{\exp}_{\mathbf{f}}$.
 \begin{enumerate}
\item Let $\calU^{\exp}_w \subset \calU^{\exp}$ be the stabilizer of 
$$\begin{cases}
\dot{w} \cdot \overline{e}, & \text{if } w \in W/W_{\mathbf{f}},\\
x_{\alpha_0}(1) \cdot \dot{w}  \cdot \overline{e},  & \text{if } w \in {_0}W_{\mathbf{f}}.
\end{cases}$$
Then $\calU^{\exp}_w$ is a closed $\Z$ subgroup of $\calU^{\exp}$.
\item There is a presentation $\calU^{\exp}_w  =\lim_i \calU^{\exp}_{w,i}$ as a pro-algebraic group, where the transition maps are relative affine spaces and $\ker(\calU_{w,i}^{\exp} \rightarrow \calU_{w,i-1}^{\exp})$ is a vector group for each $i \geq 1$. Furthermore, if 
$w \in W/W_{\mathbf{f}}$ (resp. $w \in {_0}W_{\mathbf{f}}$) then each $\calU^{\exp}_{w,i}$ is isomorphic as a scheme to the product of an affine space with $\Gm$ (resp.~ just an affine space).
\item\label{it--global section} We have $$\Fl^{\exp}_w(\mathbf{f}) = (\calU^{\exp} / \calU^{\exp}_w)_{\et} = (\calU^{\exp}_i / \calU^{\exp}_{w,i})_{\et}$$
for $i \gg 0$. Furthermore, these \'{e}tale quotients agree with the Zariski and presheaf quotients.
In other words, there exists a subscheme \(H_w\subseteq \calU^{\exp}\) such that the orbit map \eqref{orbit map} induces an isomorphism \(H_w\cong \Fl_w^{\exp}(\bbf)\).
\item Consider the natural inclusion $\Fl_w^{\exp}(\mathbf{f}) \rightarrow \Fl_w^{\circ}(\mathbf{a}_0, \mathbf{f})$. 
\begin{enumerate}
    \item If $w \in W/W_{\mathbf{f}}$ and $w$ has maximal length in its left $W_0$-orbit in $W$, this inclusion is, after trivializing affine root groups, isomorphic to the inclusion of a hyperplane $\A^{l(w)-1} \subset \A^{l(w)}$.

    \item  If $w \in W/W_{\mathbf{f}}$ and $w$ does not have maximal length in its left $W_0$-orbit in $W$, this inclusion is an isomorphism of affine spaces.

    \item If $w \in {_0}W_{\mathbf{f}}$, this inclusion is the open complement of the orbit of $\dot{w} \cdot \overline{e}$ as in part (a), and $\Fl_w^{\exp}(\mathbf{f})$ is isomorphic to $\A^{l(w)-1} \times \Gm$.
    \end{enumerate}
\end{enumerate}
\xlemm

\pf We consider two cases.

{\em Case 1:} $w \in W/W_{\mathbf{f}}$. Recall that $\calU = \calU^{>0} \rtimes U$. Since $\Gm$ fixes $\dot{w} \cdot \overline{e}$, we have $$\calU^{\exp}_w = (\calU_w \cap (\calU^{>0} \rtimes U_0)) \rtimes \Gm,$$ which proves (1).

For (2), we consider two subcases. If $w$ is maximal in its left $W_0$-orbit in $W$, then none of the $U_\alpha$ for $\alpha \in \Phi^+$ belong to $\calU_w$ by \thref{lemm-exampleW}. Hence $\calU_w \subset \calU^{>0}$, and we may take $\calU^{\exp}_{w,i} = \calU_{w,i} \rtimes \Gm$. 
If $w$ is not maximal, there is some $\alpha \in \Delta$ with $U_\alpha \subset \calU_w$. In this case, write $\calU_{w,i} = \calU_{w,i}^{>0} \rtimes U_w$, where $U_w \subset U$ and $\calU_{w,i}^{>0} = \ker (\calU_{w,i} \rightarrow U)$.
Note that $U_0 \cap U_w \subset U_w$ is a hyperplane cut out by a nonzero linear homogeneous equation with coefficients all equal to $1$ (after trivializing the root groups), so it isomorphic to an affine space. Thus, we may take $\calU^{\exp}_{w,i} = (\calU_{w,i}^{>0} \rtimes (U_0 \cap U_w)) \rtimes \Gm$. In both subcases, the claims in (2) now follow from the analogous facts for the $\calU_{w,i}$.

To prove (3), we first note that $\ker (\calU^{\exp} \rightarrow \calU_i^{\exp}) \subset \calU_w$ for $i \gg 0$, so that by definition $\Fl^{\exp}_w(\mathbf{f})$ agrees with the \'{e}tale quotients. To prove that it agrees with the presheaf quotients, we exhibit a global section $\Fl^{\exp}_w(\mathbf{f}) \rightarrow \calU^{\exp}$ (and also $\Fl^{\exp}_w(\mathbf{f}) \rightarrow \calU_i^{\exp}$ for $i \gg 0$), as follows.

Recall that $\Fl_w^{\circ}(\mathbf{a}_0, \mathbf{f})$ is given by the product of positive affine root groups $\prod_{\tilde{\alpha}} U_{\tilde{\alpha}}$ for $\restr{w^{-1}\tilde{\alpha}}{\mathbf{a}_0} < 0$. 
If $w$ is maximal, this product contains $U$, and by separating these $\tilde{\alpha}$ into $\Phi^+$ and the remaining 
 affine roots $\tilde{\alpha} = \alpha' \notin \Phi^+$, we can write $\Fl_w^{\circ}(\mathbf{a}_0, \mathbf{f}) = \prod_{\alpha'} U_{\alpha'} \times U$. Then by considering the stabilizer $\calU_w$ computed in (2), the inclusion $\Fl_w^{\exp}(\mathbf{f}) \rightarrow \Fl_w^{\circ}(\mathbf{a}_0, \mathbf{f})$ is identified with the hyperplane \begin{equation} \label{eq--hyerplane} \prod_{\alpha'} U_{\alpha'} \times U_0 \rightarrow \prod_{\alpha'} U_{\alpha'} \times U. \end{equation} The product $\prod_{\alpha'} U_{\alpha'} \times U_0$ provides the necessary section in (3), and this also completes the proof of (4) when $w$ is maximal.

If $w$ is not maximal, write the product $\Fl_w^{\circ}(\mathbf{a}_0, \mathbf{f}) = \prod_{\tilde{\alpha}} U_{\tilde{\alpha}}$ for $\restr{w^{-1}\tilde{\alpha}}{\mathbf{a}_0} < 0$ as $\prod_{\alpha'} U_{\alpha'} \times \prod_{\alpha''} U_{\alpha''}$ where $\alpha' \notin \Phi^+$ and $\alpha'' \in \Phi^+$. While the product $\prod_{\alpha'} U_{\alpha'} \times \prod_{\alpha''} U_{\alpha''}$ does not lie in $\calU^{\exp}$, we can build a section that does lie in $\calU^{\exp}$ as as follows. Fix $\alpha \in \Delta$ such that $\restr{w^{-1}\alpha}{\mathbf{a}_0} > 0$, i.e.,~$U_{\alpha}$ stabilizes $\dot{w} \cdot \overline{e}$. Let $f \colon \prod_{\alpha''} U_{\alpha''} \rightarrow \Ga$ be the inclusion into $U$ followed by \eqref{eq--simple}. Consider the map
$$\prod_{\alpha'} U_{\alpha'} \times \prod_{\alpha''} U_{\alpha''} \rightarrow \calU^{>0} \times U_0, \quad (u_1, u_2) \mapsto (u_1, u_2 \cdot x_{\alpha}(-f(u_2)).$$
This provides the necessary section in (3) and shows that $\Fl_w^{\exp}(\mathbf{f}) = \Fl_w^{\circ}(\mathbf{a}_0, \mathbf{f})$, finishing the proof of (4) when $w \in W/W_{\mathbf{f}}$.

{\em Case 2:} $w \in {_0}W_{\mathbf{f}}$. Note that
$$\calU_w^{\exp} = \calU^{\exp} \cap x_{\alpha_0}(1) \dot{w} \mathcal{I}_{\adj} \dot{w}^{-1} x_{\alpha_0}(1)^{-1},$$ where $\calI_{\adj} := \calU \rtimes T_{\adj}$. We will compute $x_{\alpha_0}(1)^{-1} \calU_w^{\exp} x_{\alpha_0}(1)$. Starting with $\calU^{\exp}$, note that conjugation by $x_{\alpha_0}(1)^{-1} \in \calU(\Z)$ induces an automorphism of $\calU^{>0} \rtimes U_0$. Furthermore, for any $\Z$-algebra $R$ and $\gamma \in \Gm(R) \subset T_{\adj}(R)$, viewed as a subgroup via $\hat{\rho}$, we have
$x_{\alpha_0}(1)^{-1}  \gamma x_{\alpha_0}(1) = x_{\alpha_0}(\gamma-1) \gamma$. It follows that $(x_{\alpha_0}(1)^{-1} \calU^{\exp} x_{\alpha_0}(1))(R)$ is given by
\begin{equation} \label{eq--stab_int} \left\{(u_1, u_2, \gamma) \in ((\calU^{>0} \rtimes U) \rtimes \Gm)(R) \: : \: \phi(u_2) =  \gamma - 1 \right\},\end{equation}
where $\phi$ is the composition \eqref{eq--simple}. To compute $\dot{w} \mathcal{I}_{\adj} \dot{w}^{-1}$, note that conjugation by $\dot{w}$ is an automorphism of $T_{\adj}$, but it may move the subgroup $\Gm \subset T_{\adj}$. Regardless, we have $\dot{w} \calI_{\adj} \dot{w}^{-1} \cap (\calU \rtimes \Gm) \subset \calU_w \rtimes \Gm$. Since $\calU_w \subset \calU^{>0}$, then taking the intersection with \eqref{eq--stab_int} gives
$x_{\alpha_0}(1)^{-1} \calU_w^{\exp} x_{\alpha_0}(1) = \calU_w$, and consequently
$$\calU_w^{\exp} = x_{\alpha_0}(1) \calU_w x_{\alpha_0}(1)^{-1}.$$
Thus, we may take $\calU_{w,i}^{\exp} = x_{\alpha_0}(1) \calU_{w,i} x_{\alpha_0}(1)^{-1}$, proving (1) and (2).

We now finish the proof of (3) and (4). Viewing $w$ as an element of $W/W_{\mathbf{f}}$, we can write $\Fl_w^\circ(\mathbf{a}_0, \mathbf{f}) = \prod_{\alpha'} U_{\alpha'} \times U $ as above, where $\alpha'$ ranges over the positive affine roots with $\restr{w^{-1}\alpha'}{\mathbf{a}_0} < 0$. Consider the product $\prod_{\alpha'} U_{\alpha'} \times U_0 \times \Gm \subset \calU_i^{\exp}$ for $i \gg 0$. It suffices to show that by taking the orbit of $x_{\alpha_0}(1) \cdot \dot{w}  \cdot \overline{e}$, this product maps isomorphically onto the complement of $\Fl_w^{\exp}(\mathbf{f}) \subset \Fl_w^{\circ}(\mathbf{a}_0, \mathbf{f})$, which is given in terms of affine root groups by \eqref{eq--hyerplane}. For an $R$-point $(u_1, u_2, \gamma)$ of  $\prod_{\alpha'}U_{\alpha'} \times U_0 \times \Gm$, we have
$$u_1 u_2 \gamma x_{\alpha_0}(1) \dot{w} \cdot \overline{e} = u_1 u_2 x_{\alpha_0}(\gamma) \dot{w} \cdot \overline{e}.$$
Hence, we can conclude by observing that $$\{u_2 x_{\alpha_0}(\gamma) \: : \: u_2 \in U_0(R), \gamma \in \Gm(R)\} = (U \setminus U_0)(R).$$
\xpf

\subsection{Exponential stratification of $G/B$}
By \cite[Lemma 4.3.7]{RicharzScholbach:Intersection}, the $\Fl^\circ_{w}(\mathbf{a}_0, \mathbf{f})$ provide a stratification of $\Fl_{\mathbf{f}}$ with closure relations determined by the Bruhat order on $W/W_{\mathbf{f}}$. When combined with \thref{orbit-structure}(4), we have a set-theoretic equality 
$$\bigsqcup_{w \in W^{\exp}_{\mathbf{f}}} \Fl_w^{\exp}(\mathbf{f}) = \Fl_{\mathbf{f}}.$$
We also note that \thref{orbit-structure}(4) implies the $\Fl_w^{\exp}(\mathbf{f})$ for $w \in {_0}W_{\mathbf{f}}$ are independent of the choice of $\alpha_0 \in \Delta$ we made in \thref{def-KirillovOrbits}.

\defi \thlabel{def-KirillovClosures} For $w \in W^{\exp}_{\mathbf{f}}$, we define $\overline{\Fl}^{\exp}_w(\mathbf{f}) \subset \Fl_{\mathbf{f}}$ to be the scheme-theoretic image of the orbit map

$$\calU^{\exp} \rightarrow \Fl_{\mathbf{f}}, \quad  \begin{cases}
g \mapsto g \cdot \dot{w} \cdot \overline{e}, & \text{if } w \in W/W_{\mathbf{f}},\\
g \mapsto g \cdot  x_{\alpha_0}(1) \cdot \dot{w}  \cdot \overline{e},  & \text{if } w \in {_0}W_{\mathbf{f}}.
\end{cases}$$
\xdefi

Our next goal is to show that the $\Fl_w^{\exp}(\mathbf{f})$ refine the Bruhat stratification, and in particular, each $\overline{\Fl}^{\exp}_w(\mathbf{f})$ is a union of $\calU^{\exp}$-orbits. We would like the flexibility to work over an arbitrary base $S$, so we first define these objects in this generality.

\defi Let $S$ be a noetherian scheme. When working over $S$ instead of $\Z$, we make the following definitions.

\begin{enumerate} \thlabel{generalS}
\item The groups $G$, $LG$, and their various subgroups and affine flag varieties in \refsect{setup} - \refsect{flagnotation} are the corresponding object obtained by base change from $\Z$.
\item The orbits (resp.~orbit closures), e.g., $\Fl_w^\circ( \mathbf{f}, \mathbf{f}')$, $\Fl_w^{\exp}(\mathbf{f})$, $\overline{\Fl}_w^{\exp}(\mathbf{f})$, etc., are the \'etale- (resp.~scheme-) theoretic images of the corresponding orbit maps over $S$.
\end{enumerate}
When we wish to compare the above objects over $\Z$ and $S$ we sometimes use a subscript to denote the latter, e.g. $\Fl_w^{\exp}(\mathbf{f})_S$.
\xdefi

\rema We refer to \cite[\S 4.4]{RicharzScholbach:Intersection} and \cite[\S 4.1.2]{CassvdHScholbach:Central} for further discussion on affine flag varieties over arbitrary bases. Since \'{e}tale-sheafification commutes with arbitrary base change, there is little difference with the theory over $\Z$, with the exception of the scheme-theoretic images, e.g.~$\Fl_w(\mathbf{f}, \mathbf{f}')$. It is shown in \cite[Proposition 4.4.3]{RicharzScholbach:Intersection} their formation commutes with base change in $S$ up to a nilpotent thickening, and we will show the same is true for the $\overline{\Fl}^{\exp}_w(\mathbf{f})$. Since categories of motives are insensitive to passing to reduced subschemes, this will be good enough for our applications.
\xrema

We start with the ordinary flag variety $G/B$, from which we will bootstrap the affine case. Let $w_0 \in W_0$ be the longest element. The action of $\calI$ on $\Fl_{w_0}(\mathbf{a}_0)$ factors through $B$, and there is a canonical $B$-equivariant isomorphism $\Fl_{w_0}(\mathbf{a}_0) \cong G/B$ identifying $U$-orbits, both of which are indexed by $W_0$. By \thref{lemm-exampleW}, $W_0 \cap ({_0}W_{\mathbf{a}_0}) = \{w_0\}$, so that only the largest 
$U$-orbit $\Fl^\circ_{w_0}(\mathbf{a}_0)$ in $G/B$ splits into two $U^{\exp}$-orbits.

We now recall several tools which we will use to analyze the closures, starting with the Demazure schemes.
Let $w_0 = s_{\alpha_1} \cdots s_{\alpha_n}$ be a reduced decomposition by simple reflections in $\mathcal{S} \cap W_0$.
Let $P_{\alpha_i} \supset B$ be the minimal parabolic associated to $\alpha_i$.
The Demazure scheme is
$$D(\alpha_1, \ldots, \alpha_n) = P_{{\alpha_1}} \overset{B}{\times} \cdots \overset{B}{\times} P_{{\alpha_n}} / B \xrightarrow{\pi} G/B,$$ where $\pi$ is induced by multiplication in the group $G$. The map $\pi$ is a proper surjection, and it is an isomorphism over the Bruhat cell $\Fl^\circ_{w_0}(\mathbf{a}_0)$. We note that the Demazure schemes and $\pi$ can also be defined over a general base $S$. By smoothness,  $D(\alpha_1, \ldots, \alpha_n)_S = D(\alpha_1, \ldots, \alpha_n) \times S.$

We have
\begin{equation} \label{eq--Dem} \pi^{-1}(\Fl^\circ_{w_0}(\mathbf{a}_0)) = U_{\alpha_1} \dot{s}_{\alpha_1} \cdots U_{\alpha_n} \dot{s}_{\alpha_n} \cdot B
\end{equation}
where each $\dot{s}_{\alpha_i} \in P_{{\alpha_i}}(\Z)$ is our chosen lift.
Our chosen trivializations $x_{\alpha_i}$ of the $U_{\alpha_i}$ then give an isomorphism
\begin{equation} \label{eqn--DemazureFiber}
    \pi^{-1}(\Fl^\circ_{w_0}(\mathbf{a}_0)) \cong \A^n.
\end{equation}
Let us denote by $t_i$ the coordinate function on $\A^n$ corresponding to $U_{\alpha_i}$ under \eqref{eqn--DemazureFiber}.

 \lemm \thlabel{DemEqn} View $w_0$ as an element of $W_0$ (and not ${_0}W_{\mathbf{a}_0}$). Suppose $l(w_0) > 1$. Then under \eqref{eqn--DemazureFiber}, over any base $S$ the fiber $\pi^{-1}(\Fl_{w_0}^{\exp}(\mathbf{a}_0))$ is a hyperplane in $\A^n$ given by an equation of the form 
    $$\sum_{i=1}^n c_i t_i = 0, \quad c_i \in \{0, 1, -1\}.$$
Furthermore, the coefficients $c_1$ and $c_n$ are nonzero.
\xlemm

\pf
By \eqref{eq--conj} we have
$$U_{\alpha_1} s_{\alpha_1} \cdots U_{\alpha_n} s_{\alpha_n} \cdot B = U_{\alpha_1} \cdot U_{\alpha_1(s_{\alpha_2})} \cdots U_{s_{\alpha_1} \cdots s_{\alpha_{n-1}} (\alpha_n)} \dot{w}_0 \cdot B.$$
The collection of roots $s_{\alpha_1} \cdots s_{\alpha_{i-1}}(\alpha_i)$ for $1 \leq i \leq n$ is precisely $\Phi^+$, and the first and last of these are simple (both claims are standard and can be deduced, e.g., by the arguments in \cite[\S 5.6 Proposition]{Humphreys:Representations}).  The result now follows from the definition of $U_0$.
\xpf

\lemm  \thlabel{coordinate-change}
Let $\alpha \in \Delta$. For a suitable choice of lift $\dot{s}_{\alpha}$, we have 
$$x_{-\alpha}(t) = x_{\alpha}(t^{-1}) \dot{s}_{\alpha} x_{\alpha}(t) \alpha^\vee(t)$$ for all $t \in \Gm$. 
 \xlemm

 \pf
This follows from the identity $$\begin{pmatrix} 1 & 0 \\ t & 1 \end{pmatrix} = \begin{pmatrix} 1 & t^{-1} \\ 0 & 1 \end{pmatrix} \begin{pmatrix} 0 & -1 \\ 1 & 0 \end{pmatrix} \begin{pmatrix} 1 & t \\ 0 & 1 \end{pmatrix} \begin{pmatrix} t & 0 \\ 0 & t^{-1} \end{pmatrix}$$ in $\text{SL}_2$.
 \xpf

\prop \thlabel{Dem-fiber}
Let $w_0 = s_{\alpha_1} \cdots s_{\alpha_n}$ be a reduced decomposition and let $v = w_0 s_{\alpha_n}$. Suppose that $l(w_0) > 1$. Then over any base $S$, the scheme-theoretic image of $\pi^{-1}(\Fl_{w_0}^{\exp}(\mathbf{a}_0))$ in $ D(\alpha_1, \ldots, \alpha_n)$ has empty intersection with 
$\pi^{-1}(\Fl^\circ_{v}(\mathbf{a}_0))$.
\xprop

 \rema
Before we prove \thref{Dem-fiber}, we note that it will be used to show that $\overline{\Fl}^{\exp}_{w_0}(\mathbf{a}_0)$ does not meet any Schubert cell of codimension one in $G/B$. Over an algebraically closed field, the latter statement can be deduced more easily as follows. The complement $\overline{\Fl}^{\exp}_{w_0}(\mathbf{a}_0) \setminus \Fl^{\exp}_{w_0}(\mathbf{a}_0)$ is a union of $U^{\exp}$-orbits, all of dimension strictly less than $\dim \Fl^{\exp}_{w_0}(\mathbf{a}_0) = l(w_0)-1$, which are thus $U$-orbits.
\xrema

\pf
First, we claim that the only way to obtain an expression for $v$ by deleting factors from $s_{\alpha_1} \cdots s_{\alpha_n}$ is by deleting $s_{\alpha_n}$. Indeed, if $v' = s_{\alpha_1} \cdots \hat{s}_{\alpha_i} \cdots s_{\alpha_n}$ for some $1 \leq i \leq n-1$, then $v' \neq v$ since $v' s_{\alpha_n} \neq w_0$. It follows that 
\begin{equation} \label{eq--Dem3} \pi^{-1}(\Fl^\circ_{v}(\mathbf{a}_0)) = U_{\alpha_1} \dot{s}_{\alpha_1} \cdots U_{\alpha_{n-1}} \dot{s}_{\alpha_{n-1}} \cdot B, \end{equation} which is the open cell in $D(\alpha_1, \ldots, \alpha_{n-1}) \subset D(\alpha_1, \ldots, \alpha_n)$. An open affine neighborhood of $\pi^{-1}(\Fl^\circ_{v}(\mathbf{a}_0))$ is therefore given by
\begin{equation} \label{eq--Dem2} U_{\alpha_1} \dot{s}_{\alpha_1} \cdots U_{\alpha_{n-1}}  \dot{s}_{\alpha_{n-1}}  \cdot U_{-\alpha_{n}} \cdot B. \end{equation}

Our chosen trivializations of the root subgroups give an isomorphism
between $\eqref{eq--Dem2}$ and $\A^n$. 
Let $u_i$ denote the coordinate function on $\A^n$ corresponding to $U_{\alpha_i}$ for $i \in \{1, \ldots, n-1\}$, and $U_{-\alpha_n}$ for $i = n$.
Note that the intersection of $\eqref{eq--Dem2}$ with $\pi^{-1}(\Fl^\circ_{w_0}(\mathbf{a}_0))$ is isomorphic to $\A^{n-1} \times \Gm$. 
By \thref{coordinate-change}, the change of variables isomorphism over this intersection is given by
$$u_i \mapsto t_i, \: i \in \{1, \ldots, n-1\}, \quad u_n \mapsto t_{n}^{-1}.$$
Here the coordinates $t_i$ on $\pi^{-1}(\Fl^\circ_{w_0}(\mathbf{a}_0))$ are defined in \eqref{eqn--DemazureFiber}.

Let $X$ denote the intersection of $\eqref{eq--Dem2}$ with the scheme-theoretic image of $\pi^{-1}(\Fl_{w_0}^{\exp}(\mathbf{f}))$ in $ D(\alpha_1, \ldots, \alpha_n)$. To compute $X$, we use \cite[Tag 01R8]{StacksProject}, and we work in the coordinates $u_i$. We have the natural inclusion $\A^{n-1} \times \Gm \rightarrow \A^n$, and by \thref{DemEqn}, $X$ is given by the scheme-theoretic image under this map of the closed subscheme in $\A^{n-1} \times \Gm$ defined by an equation of the form
$$f = c_1u_1 + \cdots c_{n-1} u_{n-1} + c_n u_n^{-1} = 0$$ where  $c_i \in \{0, 1, -1\}$ and $ c_1, c_n \neq 0$. 
Affine-locally on $S$, again using \cite[Tag 01R8]{StacksProject}, it follows that $X$ is the closed subscheme defined by the kernel of the unique $\mathcal{O}_S$-algebra homomorphism
$$\mathcal{O}_S[u_1, \ldots, u_n] \rightarrow \mathcal{O}_S[u_1, \ldots, u_{n-1}, u_n^{\pm 1}]/(f)$$ which preserves the $u_i$. Using the conditions on the $c_i$, it is straightforward to show that this kernel is the principal ideal generated by $$u_n f = c_1 u_1 u_n + \cdots + c_{n-1} u_{n-1}u_n + c_n,$$ independently of the base $S$. Hence, affine-locally on $S$, $X = \Spec \mathcal{O}_S[u_1, \ldots, u_n]/(u_n f)$. Thus, $X$ has trivial intersection with the locus in $\A^n$ where $u_n = 0$. By \eqref{eq--Dem3}, the locus in $\A^n$ where $u_n = 0$ is precisely $\pi^{-1}(\Fl^\circ_{v}(\mathbf{a}_0))$, so this completes the proof.
\xpf

For a regular dominant cocharacter of $T_{\adj}$, consider the fixed points $(G/B)^0$ and repellers $(G/B)^-$ \cite[Eqn.~(0.1)]{RIcharz:Spaces} for the induced $\Gm$-action on $G/B$ by conjugation. We have $G^0 = T$ and $G^- = B^{\opp}$, the Borel opposite to $B$, from which it follows that, at least on reduced loci, \begin{equation} \label{repeller} (G/B)^0_{\text{red}} = \bigsqcup_{w \in W_0} \dot{w} \cdot \overline{e}, \quad (G/B)^{-}_{\text{red}} = \bigsqcup_{w \in W_0} U^{\opp} \cdot \dot{w} \cdot \overline{e}.\end{equation}
Here $U^{\opp} \cdot \dot{w} \cdot \overline{e}$ is an opposite Schubert cell, defined as an \'{e}tale sheaf-theoretic image as in \thref{def-KirillovOrbits}. This description of the fixed points and repellers holds over a general base $S$ instead of $\Z$, cf.~\cite[Theorem 3.17]{HainesRicharz:TestFunctionsWeil}, \cite[Lemma 3.5]{CassvdHScholbach:MotivicSatake} for similar statements for the affine Grassmannian.

For $v, w \in W_0$, let $R_{v,w} \subset G/B$ be the reduced locus of
$$ (U^{\opp} \cdot \dot{v} \cdot \overline{e}) \cap (U \cdot \dot{w} \cdot \overline{e}).$$ Over a field, this intersection is called a Richardson variety. Deodhar \cite{Deodhar:GeometricDecomposition} has given a decomposition of $R_{v,w}$ in terms of products of root groups and punctured root groups. Marsh and Rietsch \cite{MarshRietsh:Flag} made Deodhar's decomposition more explicit. While these results were stated over a field, the proofs are purely combinatorial and uniform across all fields, and are thus also valid over a general base $S$. We will apply \cite{MarshRietsh:Flag} to $R_{v,w_0}$ for $v$ chosen as in the lemma below.

\lemm \thlabel{lemm--Weyl}
Let $v \in W_0$ be such that $l(v) =n - 2$, where $n = l(w_0)$.
\begin{enumerate}
    \item  There exist $\alpha_1, \alpha_n \in \Delta$ such that $w_0 = s_{\alpha_1} v s_{\alpha_n}$.
    \item For  $\alpha_1, \alpha_n \in \Delta$ chosen as above, the root $v(\alpha_n)$ is positive, and the root subgroup $U_{\alpha_1}$ commutes with $U_{v(\alpha_n)}$.
\end{enumerate}
\xlemm

\pf
Recall that $l(w) = | \{\Phi^+ \cap w^{-1}(\Phi^-)\}|$.
Since $v \neq w_0$ there exists $\alpha_1 \in \Delta$ such that $v^{-1}(\alpha_1) \in \Phi^+$. Since $s_{\alpha_1}(\alpha_1) = - \alpha_1$ and $s_{\alpha_1}$ permutes $\Phi^+ \setminus \{\alpha_1\}$, it follows that $l(s_{\alpha_1} v) = l(v) + 1$. Now apply the same argument to $(s_{\alpha_1} v)^{-1}$ to get $\alpha_n$. This proves (1).

For (2), choose a reduced expression $v = s_{\alpha_2} \cdots s_{\alpha_{n-1}}$. The collection of roots $s_{\alpha_1} \cdots s_{\alpha_{i-1}}(\alpha_i)$ for $1 \leq i \leq n$ is precisely $\Phi^+$, with no repetitions, and the last of these, $s_{\alpha_1} \cdots s_{\alpha_{n-1}}(\alpha_n)$, is simple. Thus, $s_{\alpha_1}v(\alpha_n)$ is a simple root distinct from $\alpha_1$, so that $v(\alpha_n)$ is positive.  To show that $U_{\alpha_1}$ commutes with $U_{v(\alpha_n)}$, it suffices to check that there are no positive roots of the form $i \alpha_1 + j v(\alpha_n)$, for $i, j > 0$. If $i \alpha_1 + j v(\alpha_n)$ is a root, then applying $s_{\alpha_1}$, it follows that $-i \alpha_1 + j (s_{\alpha_1}v(\alpha_n))$ is also a root. But this is impossible since $\alpha_1$ and $s_{\alpha_1}v(\alpha_n)$ are simple and distinct, while $-i$ and $j$ have different signs.
\xpf

Recall that by \thref{orbit-structure}, if we view $w_0$ as an element of $W_0$ (and not ${_0}W_{\mathbf{a}_0}$), then $\Fl^{\exp}_{w_0}(\mathbf{a}_0) = U^{\exp} \cdot \dot{w}_0 \cdot \overline{e} \subset \Fl_{w_0}^\circ(\mathbf{a}_0)$ is the smaller of the two $U^{\exp}$-orbits in  $\Fl_{w_0}^\circ(\mathbf{a}_0)$. The following
\thref{Deodhar} will be used to show that $\overline{\Fl}^{\exp}_{w_0}(\mathbf{a}_0)$ meets every Schubert cell of codimension at least two in $G/B$.

\prop \thlabel{Deodhar}
View $w_0$ as an element of $W_0$ (and not ${_0}W_{\mathbf{a}_0}$). Let $v \in W_0$ be such that $l(v) =l(w_0) - 2$. Choose as in \thref{lemm--Weyl} a reduced expression $w_0 = s_{\alpha_1} s_{\alpha_2} \cdots s_{\alpha_n}$ such that $v = s_{\alpha_2} \cdots s_{\alpha_{n-1}}$.
Then over any base $S$, the intersection
$$(U^{\opp} \cdot \dot{v} \cdot \overline{e}) \cap \Fl^{\exp}_{w_0}(\mathbf{a}_0)$$ contains 
$$\{x_{\alpha_1}(t_1) \cdot x_{ s_{\alpha_1} \cdots s_{\alpha_{n-1}}(\alpha_n)}(t_n) \: : \: t_1, t_n \in \Gm, \: t_1 + t_n = 0 \} \cdot \dot{w}_0B$$
as a locally closed subscheme isomorphic to $\Gm$.
\xprop

\pf
In the sense of \cite[Definition 3.1]{MarshRietsh:Flag} we have the \emph{expression} $$\textbf{w}_0 = (w_{(0)}, \ldots,  w_{(n)}) = (1, s_{\alpha_1}, s_{\alpha_1} s_{\alpha_2}, \ldots, s_{\alpha_1} s_{\alpha_2} \cdots s_{\alpha_n})$$ for $w_0$, with \emph{sequence of factors} $(s_{\alpha_1}, s_{\alpha_2}, \ldots, s_{\alpha_n})$. The expression $$\mathbf{v} = (1, 1, s_{\alpha_2}, s_{\alpha_2} s_{\alpha_3}, \ldots, s_{\alpha_1} \cdots s_{\alpha_{n-1}}, s_{\alpha_1} \cdots s_{\alpha_{n-1}})$$ is easily verified to be a \emph{subexpression} for $v$ in $\textbf{w}$ in the sense of \cite[Definition 3.3]{MarshRietsh:Flag}, meaning that
$$v_{(j)} \in \left\{ v_{(j-1)},v_{(j-1)} s_{\alpha_j}\right\} \text{ for } j \in \{1, \ldots, n\},$$ and $v_{(n)} = v$. It is moreover \emph{distinguished}, meaning that
$$v_{(j)} \leq v_{(j-1)} s_{\alpha_j} \text{ for } j \in \{1, \ldots, n\}.$$
(In fact, $\mathbf{v}$ is the unique \emph{positive subexpression} \cite[Lemma 3.5]{MarshRietsh:Flag} for $v$ in $\textbf{w}_0$.) 

It is straightforward to compute the sets in \cite[Definition 3.2]{MarshRietsh:Flag} as
\begin{align*}J_{\mathbf{v}}^+ := \{j \in \{1, \ldots, n\} \: : \: v_{(j-1)} < v_{(j)} \} & = \{2, \ldots, n-1\} \\
J_{\mathbf{v}}^\circ := \{j \in \{1, \ldots, n\} \: : \: v_{(j-1)} = v_{(j)} \} &= \{1, n\} \\ J_\mathbf{v}^- := \{j \in \{1, \ldots, n\} \: : \: v_{(j-1)} > v_{(j)} \} &= \emptyset.\end{align*}
Then by \cite[Proposition 5.2]{MarshRietsh:Flag}, we have a locally closed subscheme
$$\mathcal{R}_{\mathbf{v}, \mathbf{w}_0} := \{x_{-\alpha_1}(t_1)\dot{s}_{\alpha_2} \cdots \dot{s}_{\alpha_{n-1}} x_{-\alpha_{n}}(t_n) \: : \: t_1, t_n \in \Gm \} \cdot B \subset R_{v, w_0}.$$ This is in fact a dense subscheme, and it is isomorphic to $\Gm^2$.

We now compute $\mathcal{R}_{\mathbf{v}, \mathbf{w}_0} \cap \Fl^{\exp}_{w_0}(\mathbf{a}_0)$. For this purpose we are free to pick the lifts $\dot{s}_{\alpha_i}$. Then using \thref{coordinate-change}
to eliminate the  $x_{-\alpha_i}(t_i)$ for $i \in \{1, n\}$, it follows that an arbitrary element of $\mathcal{R}_{\mathbf{v}, \mathbf{w}_0}$ can be written as
$$x_{\alpha_1}(t_1^{-1}) \dot{s}_{\alpha_1} x_{\alpha_1}(t_1) x_{v(\alpha_n)} (\pm t_1^{\langle v(\alpha_n), \alpha_1^\vee \rangle} t_n^{-1}) \dot{v} \dot{s}_{\alpha_n} \cdot B,$$ for $t_1^{-1}, t_n^{-1} \in \Gm$. Here it is not necessary to resolve the ambiguity in the sign, as we are always free to replace $t_n$ with $-t_n$.
By \thref{lemm--Weyl}, $U_{\alpha_1}$ commutes with $U_{v(\alpha_n)}$. Note that $x_{\alpha_1}(t_1) \dot{v} \dot{s}_{\alpha_n} = \dot{v} \dot{s}_{\alpha_n} x_{s_{\alpha_n} \cdots s_{\alpha_{2}}(\alpha_1)}(\pm t_1)$ and $s_{\alpha_n} \cdots s_{\alpha_{2}}(\alpha_1)$ is positive, so that the above element of $\mathcal{R}_{\mathbf{v}, \mathbf{w}_0}$ is equal to
$$x_{\alpha_1}(t_1^{-1}) \dot{s}_{\alpha_1} x_{v(\alpha_n)} (\pm t_1^{\langle v(\alpha_n), \alpha_1^\vee \rangle} t_n^{-1}) \dot{v} \dot{s}_{\alpha_n} \cdot B.$$ By moving $\dot{s}_{\alpha_1}$ to the right, this is equal to
$$x_{\alpha_1}(t_1^{-1}) x_{s_{\alpha_1} \cdots s_{\alpha_{n-1}}(\alpha_n)} (\pm t_1^{\langle v(\alpha_n), \alpha_1^\vee \rangle} t_n^{-1}) \dot{w}_0 \cdot B.$$ It follows that
$$\mathcal{R}_{\mathbf{v}, \mathbf{w}_0} = \{x_{\alpha_1}(t_1) x_{s_{\alpha_1} \cdots s_{\alpha_{n-1}}(\alpha_n)}(t_n) \: : \: t_1, t_n \in \Gm \} \cdot \dot{w}_0B.$$
Hence, $\mathcal{R}_{\mathbf{v}, \mathbf{w}_0} \cap \Fl^{\exp}_{w_0}(\mathbf{a}_0)$ is the closed subscheme of the above defined by $t_1 + t_n = 0$, which is isomorphic to $\Gm$ since $\alpha_1 \neq s_{\alpha_1} \cdots s_{\alpha_{n-1}}(\alpha_n)$.
\xpf

We will use the following lemma in the proof of \thref{theo--strat} below, and after that, often without comment.

\lemm \thlabel{imageofclosure}
Let $f \colon X \r Y$ and $g \colon Y \r Z$ be morphisms of noetherian schemes. Let $\overline{X}$ be the scheme-theoretic image of $X$ in $Y$. Then the scheme-theoretic image of $\overline{X}$ in $Z$ agrees with the scheme-theoretic image of $X$ in $Z$.
\xlemm

\pf
This follows from the description of scheme-theoretic images in terms of ideal sheaves in \cite[Tag  01R8]{StacksProject}.
\xpf

\theo \thlabel{theo--strat}
View $w_0$ as an element of $W_0$ (and not ${_0}W_{\mathbf{a}_0}$).
\begin{enumerate}
    \item The closure $\overline{\Fl}^{\exp}_{w_0}(\mathbf{a}_0)$ is set-theoretically the disjoint union of the following locally closed subschemes:
    $$\overline{\Fl}^{\exp}_{w_0}(\mathbf{a}_0) = \Fl^{\exp}_{w_0}(\mathbf{a}_0) \sqcup \bigsqcup_{v \in W_0, \: l(v) \leq l(w_0) - 2} \Fl_{v}^
    {\circ}(\mathbf{a}_0).$$
        \item The natural map $\Fl^{\exp}_{w_0}(\mathbf{a}_0) \rightarrow \overline{\Fl}^{\exp}_{w_0}(\mathbf{a}_0)$ is an open immersion.
    \item For any base $S$, the natural map
    $$\overline{\Fl}^{\exp}_{w_0}(\mathbf{a}_0)_S \rightarrow \overline{\Fl}^{\exp}_{w_0}(\mathbf{a}_0) \times S$$ is a nilpotent thickening.
\end{enumerate}
\xtheo

\rema
Since the other $U^{\exp}$-orbit $\Fl^{\exp}_{w_0}(\mathbf{a}_0) \setminus U$ is dense in $\Fl_{w_0}^{\circ}(\mathbf{a}_0)$ by \thref{orbit-structure}, then \thref{theo--strat} and the Bruhat decomposition imply that the $U^{\exp}$-orbits give a stratification of $G/B$, independent of the base $S$. 
\xrema

\pf If $l(w_0) = 1$ then $\Fl^{\exp}_{w_0}(\mathbf{a}_0) \subset G/B \cong \P$ is a point and all of the statements are clear. Otherwise, we necessarily have $l(w_0) \geq  2$. Then let $v \in W_0$ be an arbitrary element such that $l(v) =l(w_0) - 2$.
Taking the orbit under $\Gm \subset U^{\exp}$ of the point  $x_{\alpha_1}(1) \cdot x_{ s_{\alpha_1} \cdots s_{\alpha_{n-1}}(\alpha_n)}(-1)$ as in \thref{Deodhar} gives a map 
$$\Gm \r (U^{\opp} \cdot \dot{v} \cdot \overline{e}) \cap \Fl^{\exp}_{w_0}(\mathbf{a}_0).$$
As $U^{\opp} \cdot \dot{v} \cdot \overline{e}$ is a component of the repeller \eqref{repeller}, this extends to a map $\P \setminus \{0\} \r G/B$ whose value at $\infty$ is $\dot{v} \cdot \overline{e}$. The map $\P \setminus \{0\} \r G/B$ necessarily factors through $\overline{\Fl}^{\exp}_{w_0}(\mathbf{a}_0)$ (as this is the closure of a larger group action), so that $\overline{\Fl}^{\exp}_{w_0}(\mathbf{a}_0)$ contains the $\Z$-point $\dot{v} \cdot \overline{e}$. Since $\overline{\Fl}^{\exp}_{w_0}(\mathbf{a}_0)$ is necessarily $U^{\exp}$-stable, it follows from \thref{orbit-structure} and the Bruhat decomposition that $\overline{\Fl}^{\exp}_{w_0}(\mathbf{a}_0)$ contains the locally closed subschemes in (1).

Since $\Fl^{\exp}_{w_0}(\mathbf{a}_0) \subset \Fl_{w_0}^\circ(\mathbf{a}_0) $ is closed and $\Fl_{w_0}^\circ(\mathbf{a}_0) \subset G/B$ is open, then  $\overline{\Fl}^{\exp}_{w_0}(\mathbf{a}_0) \cap \Fl_{w_0}^\circ(\mathbf{a}_0) = \Fl^{\exp}_{w_0}(\mathbf{a}_0)$ and (2) holds. Now to finish the proof of (1) it suffices to show that  if $v \in W_0$ is such that $l(v) = l(w_0)-1$, then $\overline{\Fl}^{\exp}_{w_0}(\mathbf{a}_0) \cap \Fl_{v}^\circ(\mathbf{a}_0) = \emptyset$. For this, we may argue as in the proof of \thref{lemm--Weyl} to find $\alpha_i \in \Delta$ such that $w_0 = s_{\alpha_1} \cdots s_{\alpha_{n}}$ is a reduced expression and $v =  w_0 s_{\alpha_n}$. Consider the Demazure resolution $\pi \colon D(\alpha_1, \ldots, \alpha_n) \r G/B$. Since our schemes are noetherian and $\pi$ is proper, then by \thref{imageofclosure} the set-theoretic locus of $\overline{\Fl}^{\exp}_{w_0}(\mathbf{a}_0)$ is the image under $\pi$ of the scheme-theoretic image of $\pi^{-1}(\Fl_{w_0}^{\exp}(\mathbf{f}))$ in $ D(\alpha_1, \ldots, \alpha_n)$. Thus, $\overline{\Fl}^{\exp}_{w_0}(\mathbf{a}_0) \cap \Fl_{v}^\circ(\mathbf{a}_0) = \emptyset$ by \thref{Dem-fiber}. Finally, (3) holds since all of the arguments we have given apply over any base $S$, which in particular implies the map in (3) is an equality on topological spaces.
\xpf

\subsection{Exponential stratification of the affine flag variety}
In this section we show that the $\Fl_w^{\exp}(\mathbf{f})$ give a stratification of $\Fl_{\mathbf{f}}$. We will frequently apply \thref{imageofclosure} to certain convolution maps, which we now recall. For any standard facets $\mathbf{f}$, $\mathbf{f}'$ we let
$$\Fl_{\mathbf{f}'}  \widetilde\times  \Fl_{\mathbf{f}} := LG \overset{L^+G_{\mathbf{f}'} }{\times} \Fl_{\mathbf{f}} \xrightarrow{m} \Fl_{\mathbf{f}},$$ where $L^+G_{\mathbf{f}'} $ acts diagonally by $g \cdot (h_1, h_2) = (h_1 g^{-1}, g h_2)$ and $m$ is induced by multiplication in $LG$. The ind-scheme $\Fl_{\mathbf{f}'}  \widetilde\times  \Fl_{\mathbf{f}}$ and the map $m$ are both ind-proper. The above process can be iterated as well, so that we have $n$-fold convolution Grassmannians $\Fl_{\mathbf{f}_1}  \widetilde\times  \Fl_{\mathbf{f}_2} \widetilde\times \cdots \widetilde\times \Fl_{\mathbf{f}_n}$.

If $X \subset \Fl_{\mathbf{f}'}$ is an ind-subscheme and $Y \subset \Fl_{\mathbf{f}}$ is $L^+G_{\mathbf{f}'}$-stable, we may form the twisted product $X \widetilde\times Y \subset \Fl_{\mathbf{f}'}  \widetilde\times  \Fl_{\mathbf{f}}$. In this section we will only need to consider the case $\mathbf{f}'=\mathbf{a}_0$, so that $L^+G_{\mathbf{f}'} = \calI$.

\lemm \thlabel{conv-GB}
View $w_0$ as an element of $W_0$ (and not ${_0}W_{\mathbf{a}_0}$). Let $\alpha \in \Delta$, with associated simple reflection $s_{\alpha} \in W_0$. Then over any base $S$, the scheme-theoretic image of the convolution map
$$\Fl^{\exp}_{w_0}(\mathbf{a}_0) \widetilde\times \Fl_{s_\alpha}^\circ(\mathbf{a}_0) \r \Fl$$
is $\Fl_{w_0}(\mathbf{a}_0) = G/B$.
\xlemm

\pf
The convolution map clearly factors through $G/B$. Choose a reduced expression $w_0 = s_{\alpha_1} \cdots s_{\alpha_n}$ such that $s_{\alpha_n} = s_{\alpha}$. Similarly to \eqref{eq--Dem}, let 
$$X: = \Fl^\circ_{w_0}(\mathbf{a}_0) \widetilde\times \Fl_{s_\alpha}^\circ(\mathbf{a}_0) = U_{\alpha_1}  \dot{s}_{\alpha_1} \cdots U_{\alpha_n} \dot{s}_{\alpha_n}  \cdot U_{\alpha_n} \dot{s}_{\alpha_n}.$$ 
Here an element of the right side is viewed as a formal product, i.e., we do not (yet) multiply the factors $U_{\alpha_i} \dot{s}_{\alpha_i}$ in $G$.
Then $X \cong \A^{n+1}$ with coordinate functions $t_i$ corresponding to the trivializations of the leftmost $U_{\alpha_i}$ for $i \in \{1, \ldots, n\}$, and $t_{n+1}$ corresponding to the last copy of $U_{\alpha_n}$.

We now compute a lift to $G$ of the image of a point $(p_1, \ldots, p_{n+1}) \in \A^{n+1} \cong X$ under convolution, such that $p_{n+1} \neq 0$. Using $\dot{s}_{\alpha_n} \cdot U_{\alpha_n} \cdot \dot{s}_{\alpha_n}^{-1} = U_{-\alpha_{n}}$ and \thref{coordinate-change}, such a lift is given by
\begin{equation} \label{lift} U_{\alpha_1}(p_1) \dot{s}_{\alpha_1} \cdots U_{\alpha_{n-1}}(p_{n-1})  \dot{s}_{\alpha_{n-1}} \cdot U_{\alpha_n}(p_n - p_{n+1}^{-1})  \dot{s}_{\alpha_n} \in G. \end{equation}

Let $V \subset X \cong \A^{n+1}$ be the locally closed subscheme defined by  $c_1t_1 + \cdots + c_n t_n = 0$ as in \thref{DemEqn} and the condition $t_{n+1} \neq 0$. Since $\Fl^{\exp}_{w_0}(\mathbf{a}_0) \widetilde\times \Fl^\circ_{s_\alpha}(\mathbf{a}_0) \subset X$ is cut out by the latter equation, then $V$ is a dense open subscheme of $\Fl^{\exp}_{w_0}(\mathbf{a}_0) \widetilde\times \Fl_{s_\alpha}^\circ (\mathbf{a}_0)$. By \eqref{lift}, the convolution map from $V$ to $G/B$ factors through $\Fl^\circ_{w_0}(\mathbf{a}_0)$. Since $c_n \in \{1, -1\}$, then $V \cong \A^{n-1} \times \Gm$ with coordinate functions $t_1, \ldots, t_{n-1}, t_{n+1}^{\pm 1}$ inherited from $V$. Using \eqref{lift}, in terms of the isomorphism $\Fl^\circ_{w_0}(\mathbf{a}_0) \cong \A^n$ in \eqref{eqn--DemazureFiber}, the convolution map $V \r \Fl^\circ_{w_0}(\mathbf{a}_0)$ is given in coordinates by the map $\A^{n-1} \times \Gm \r \A^n$ such that
$$\quad (p_1, \ldots, p_{n-1}, p_{n+1}) \mapsto (p_1, \ldots, p_{n-1}, -(c_nc_1 p_1 + \cdots + c_n c_{n-1}p_{n-1}) - p_{n+1}^{-1}).$$
It is not difficult to check that the scheme-theoretic image of the above map is all of $\A^n \cong \Fl^\circ_{w_0}(\mathbf{a}_0)$, for any base $S$. Since $G/B = \Fl_{w_0}(\mathbf{a}_0)$  for any base $S$ (as both are covered by $W_0$-translates of $\Fl_{w_0}^\circ(\mathbf{a}_0) \cong \A^n$), the lemma follows.
\xpf

\lemm \thlabel{convlem}
Let $w \in W$ have maximal length in its left $W_0$-orbit, and write $w = w_0 v$ for $v \in W$. Then the convolution map
$$\Fl_{w_0}^{\exp}(\mathbf{a}_0) \widetilde\times \Fl_{v}^\circ(\mathbf{a}_0) \rightarrow \Fl_{w}^{\exp}(\mathbf{a}_0) $$ is an isomorphism over any base $S$.
\xlemm

\pf
Since $l(w_0) + l(v) = l(w)$ then convolution gives an isomorphism from $\Fl_{w_0}^\circ(\mathbf{a}_0) \widetilde\times \Fl_{v}^\circ(\mathbf{a}_0)$ to $ \Fl_{w}^\circ(\mathbf{a}_0)$. We have $\Fl_{w_0}^\circ(\mathbf{a}_0) = U \cdot \dot{w}_0 \cdot \overline{e}$ and $\Fl_{v}^\circ(\mathbf{a}_0) = \prod_{\tilde{\alpha}} U_{\tilde{\alpha}} \cdot \dot{v} \cdot \overline{e}$, where $\tilde{\alpha}$ runs over the positive affine roots such that  $\restr{v^{-1}\tilde{\alpha}}{\mathbf{a}_0} < 0$. Note that none of these $\tilde{\alpha}$ lie in $\Phi^+$ since $v$ is minimal in its left $W_0$-orbit.
Then $\Fl_w^\circ(\mathbf{a}_0) = U \cdot \prod_{\tilde{\alpha}} U_{w_0(\tilde{\alpha})} \cdot \dot{w} \cdot \overline{e}$, and the lemma now follows from the definition of the exponential orbits.
\xpf

For any standard facet $\mathbf{f}$, we construct a canonical map
\begin{equation} \label{orbit-projection} W_{\mathbf{a}_0}^{\exp} = W \sqcup {_0}W_{\mathbf{a}_0} \r W/W_{\mathbf{f}} \sqcup {_0}W_{\mathbf{f}} =  W_{\mathbf{f}}^{\exp} \end{equation} as follows. If $w \in W \subset W_{\mathbf{a}_0}^{\exp}$, we map $w$ to the unique element $w_{\min} \in W/W_{\mathbf{f}}$ of minimal length in the orbit $w W_{\mathbf{f}}$. If $w \in {_0}W_{\mathbf{a}_0}$, there are two cases to consider. If $w_{\min}$ also has maximal length in its left $W_0$-orbit then we map $w$ to $w_{\min} \in {_0}W_{\mathbf{f}}$, and if $w_{\min}$ does not have maximal length then we map $w$ to $w_{\min} \in W/W_{\mathbf{f}}$.

\prop \thlabel{Uexp-projection}
Fix a standard facet $\mathbf{f}$ and consider the projection $\pi \colon \Fl \rightarrow \Fl_{\mathbf{f}}$. Let $w \in W_{\mathbf{a}_0}^{\exp}$ and let $w_{\min} \in  W_{\mathbf{f}}^{\exp}$ be its image under \eqref{orbit-projection}. Write $w = w_{\min} v$ for a unique $v \in W_{\mathbf{f}}$. 
\begin{enumerate}
    \item At the level of $\calU$-orbits, if $w \in W$ then $\pi$ restricts to a relative affine space $\Fl_w^\circ(\mathbf{a}_0) \r 
\Fl^\circ_{w_{\min}}(\mathbf{a}_0, \mathbf{f})$ of relative dimension $l(v)$.

\item At the level of $\calU^{\exp}$-orbits, one of the following holds.
    \begin{enumerate}
        \item If $w \in W$ and $w$ is not maximal in its left $W_0$-orbit, then $\pi$ restricts to a relative affine space $\Fl^{\exp}_{w}(\mathbf{a}_0) \r \Fl^{\exp}_{w_{\min}}(\mathbf{f})$, which is identified with the map in (1) under \thref{orbit-structure}(4b).
        \item If $w \in W$ and $w$ is maximal in its left $W_0$-orbit, then $\pi$ restricts to a relative affine space $\Fl_w^{\exp}(\mathbf{a}_0) \r \Fl^{\exp}_{w_{\min}}(\mathbf{f})$ of relative dimension
            $$\begin{cases}
            l(v)-1, & \text{ if } w_{\min} \text{ is not maximal in its left } W_0\text{-orbit} ,\\
            l(v),  & \text{ if } w_{\min} \text{ is maximal in its left } W_0\text{-orbit}.
            \end{cases}$$

        \item If $w \in {}_0W_{\mathbf{a}_0}$, then $\pi$ restricts to a map $\Fl_w^{\exp}(\mathbf{a}_0) \r \Fl^{\exp}_{w_{\min}}(\mathbf{f})$, which, in suitable coordinates, is isomorphic to a projection of the following form:
        $$\begin{cases}
            \A^{l(w_{\min})} \times \A^{l(v)-1} \times \Gm\r \A^{l(w_{\min})},  & \text{ if } w_{\min} \text{ is not maximal},\\
            \A^{l(w_{\min})-1} \times \Gm \times \A^{l(v)} \r \A^{l(w_{\min})-1} \times \Gm,  & \text{ if } w_{\min} \text{ is maximal}.
        \end{cases}$$
        
    \end{enumerate}
\end{enumerate}
\xprop

\pf
Part (1) is shown in the proof of \cite[Proposition 4.3.13]{RicharzScholbach:Intersection}. More precisely, separate the set of positive affine roots $\tilde{\alpha}$ such that $\restr{w^{-1}\tilde{\alpha}}{\mathbf{a}_0} < 0$ into those $\alpha$ such that $\restr{w_{\min}^{-1}\alpha}{\mathbf{a}_0} < 0$ and those $\alpha'$ such that the latter inequality does not hold. Then $\Fl_w^{\circ}(\mathbf{a}_0) = \prod_{\alpha} U_{\alpha} \times \prod_{\alpha'} U_{\alpha'} \cdot \dot{w} \cdot \overline{e}$ and $\Fl_{w_{\min}}^{\circ}(\mathbf{a}_0, \mathbf{f}) = \prod_{\alpha} U_{\alpha} \times \dot{w}_{\min} \cdot \overline{e}$. The map $\pi$ sends $\dot{w} \cdot \overline{e}$ to $\dot{w}_{\min} \cdot \overline{e}$ and forgets the $U_{\alpha'}$. (Here we are free to fix once and for all an order of the affine root groups under each of the two $\prod$ symbols, but we must position the $U_{\alpha'}$ directly adjacent to $\dot{w} \cdot \overline{e}$ for this description of $\pi$ to hold.)

For part (2), all that remains to prove in case (a) is that $w_{\min}$ is also not maximal in its left $W_0$-orbit, which is straightforward to check. In case (b), write $\Fl_w^{\circ}(\mathbf{a}_0) = \prod_{\alpha} U_{\alpha} \times \prod_{\alpha'} U_{\alpha'} \cdot \dot{w} \cdot \overline{e}$ as in part (1). Since $w$ is maximal, then $\Delta$ is contained in the union of the $\{\alpha\}$ and $\{\alpha'\}$ by  \thref{lemm-exampleW}. Moreover, the collection of just the $\{\alpha\}$ contains $\Delta$ if and only if $w_{\min}$ is maximal. Now case (b) follows from part (1) and the fact that the $\calU^{\exp}$-orbits are defined by the vanishing of the sum over $\Delta$ of root subgroups. Finally, case (c) follows from the previous cases and \thref{orbit-structure}(c).
\xpf

\coro \thlabel{Wexp-lift}
For $w \in W_{\mathbf{f}}^{\exp}$, there is a canonical lift $\tilde{w} \in W_{\mathbf{a}_0}^{\exp}$ with respect to \eqref{orbit-projection} such that the projection $\pi \colon \Fl \rightarrow \Fl_{\mathbf{f}}$ restricts to an isomorphism $\Fl_{\tilde{w}}^{\exp}(\mathbf{a}_0) \cong \Fl_{w}^{\exp}(\mathbf{f})$.
\xcoro

\pf
If $w \in W/W_{\mathbf{f}} \subset W_{\mathbf{f}}^{\exp}$, we let $\tilde{w} = w \in W$ viewed as an element of minimal length in its right $W_{\mathbf{f}}$-orbit. If $w \in {}_0W_{\mathbf{f}}$, we let $\tilde{w} = w \in {}_0W_{\mathbf{a}_0}$. This choice of $\tilde{w}$ satisfies the required property by \thref{Uexp-projection}.
\xpf

\lemm \thlabel{Iconvolution}
Let $w_1 \in W$ and $w_2 \in W/W_{\mathbf{f}}$. Then the scheme-theoretic image of the convolution map $$\Fl^{\circ}_{w_1}(\mathbf{a}_0) \widetilde\times \Fl^{\circ}_{w_2}(\mathbf{a}_0, \mathbf{f}) \r \Fl_{\mathbf{f}}$$
is set-theoretically a disjoint union of locally closed subschemes of the form $\Fl^{\circ}_{v}(\mathbf{a}_0, \mathbf{f})$, for certain $v \in W/W_{\mathbf{f}}$ which are independent of the base $S$.
\xlemm

\pf
Factor the convolution map as $\Fl_{w_1}^\circ(\mathbf{a}_0) \widetilde\times \Fl_{w_2}^\circ(\mathbf{a}_0) \r \Fl \r \Fl_{\mathbf{f}}$. The projection $\Fl \r \Fl_{\mathbf{f}}$ maps each $\calI$-orbit onto another combinatorially defined $\calI$-orbit \cite[Proposition 4.3.13]{RicharzScholbach:Intersection}, so we reduce to the case $\mathbf{f} = \mathbf{a}_0$. If $l(w_1) + l(w_2) = l(w_1w_2)$ then the convolution map is an isomorphism onto $\Fl^{\circ}_{w_1w_2}(\mathbf{a}_0, \mathbf{f})$, so the result follows from the Bruhat decomposition in this case. In general, we proceed using induction on the lengths of $w_1$ and $w_2$. By factoring the convolution map according to reduced decompositions of $w_1$ and $w_2$, it is not difficult to further reduce to the case where $w_1 = w_2 =s$ are the same simple reflection, in which case the scheme-theoretic image is the closure $\Fl_{s}(\mathbf{a}_0)$ by the first part of the proof of \cite[Proposition 3.19]{RicharzScholbach:Intersection}.
\xpf

\theo \thlabel{conv-image}
Let $w_1 \in W_{\mathbf{a}_0}^{\exp}$ and $w_2 \in W/W_{\mathbf{f}}$. 
Then the scheme-theoretic image of the convolution map $$\Fl^{\exp}_{w_1}(\mathbf{a}_0) \widetilde\times \Fl^{\circ}_{w_2}(\mathbf{a}_0, \mathbf{f}) \r \Fl_{\mathbf{f}}$$
is set-theoretically a disjoint union of locally closed subschemes of the form $\Fl^{\exp}_{v}(\mathbf{f})$, for certain $v \in W_{\mathbf{f}}^{\exp}$ which are independent of the base $S$.
\xtheo

\pf
Identifying $w_2$ with its canonical lift to $W$ under \thref{Wexp-lift}, we may factor the convolution map as $\Fl^{\exp}_{w_1}(\mathbf{a}_0) \widetilde\times \Fl^{\circ}_{w_2}(\mathbf{a}_0) \r \Fl \r \Fl_{\mathbf{f}}$. Hence by \thref{Uexp-projection} we may assume $\mathbf{f} = \mathbf{a}_0$. If $w_1 \in W$ and $w_1$ is not maximal in its left $W_0$-orbit, then $\Fl^{\exp}_{w_1}(\mathbf{a}_0)$ is an $\calI$-orbit by \thref{orbit-structure}, so \thref{Iconvolution} applies. If $w_1 \in {}_0W_{\mathbf{a}_0}$ then $\Fl^{\exp}_{w_1}(\mathbf{a}_0)$ is dense in the $\calI$-orbit $\Fl^{\circ}_{w_1}(\mathbf{a}_0)$, so that \thref{Iconvolution} applies again. 

Finally, if $w_1 \in W$ and $w_1$ is maximal in its left $W_0$-orbit, write $w_1 = w_0 v$ for some $v \in W$, where $w_0 \in W_0$ is the longest element. Then using \thref{convlem}, we may rewrite the convolution map as an iterated convolution $\Fl_{w_0}^{\exp}(\mathbf{a}_0) \widetilde\times \Fl_v^\circ(\mathbf{a}_0) \widetilde\times \Fl_{w_2}^\circ(\mathbf{a}_0) \r \Fl$. By convolving the middle and right factors first and using \thref{Iconvolution}, we may further reduce to the case where $w_1 = w_0 \in W$ and $w_2 = v$ for some $v \in W$. In this case we may first compute the scheme-theoretic image of $\Fl_{w_0}^{\exp}(\mathbf{a}_0) \widetilde\times \Fl_v^\circ(\mathbf{a}_0)$ in $\Fl \widetilde\times \Fl$, which is $\overline{\Fl}^{\exp}_{w_0}(\mathbf{a}_0) \widetilde\times \Fl_v(\mathbf{a}_0)$, and then by properness of the convolution map $m$, it suffices to show the latter scheme maps onto a union of $\calU^{\exp}$-orbits under $m$. 

By \thref{theo--strat},  $\overline{\Fl}^{\exp}_{w_0}(\mathbf{a}_0) \setminus \Fl^{\exp}_{w_0}(\mathbf{a}_0)$ is a union of $\calI$-orbits. Hence by \thref{convlem}, it suffices to show that $\Fl_{w_0}^{\exp}(\mathbf{a}_0) \widetilde\times \Fl_v^\circ(\mathbf{a}_0)$ maps onto a union of $\calU^{\exp}$-orbits under $m$, for any $v \in W$. If $v$ is minimal in its left $W_0$-orbit, then
$l(w_0) + l(v) = l(w_0v)$ and the latter scheme maps isomorphically onto a $\calU^{\exp}$-orbit by \thref{convlem}. If $v$ is not minimal, write $v = v' v_{\min}$ where $v' \in W_0$ is nontrivial and $v_{\min}$ is minimal. At this point, it suffices to show that the scheme-theoretic image of $\Fl_{w_0}^{\exp}(\mathbf{a}_0) \widetilde\times \Fl_v^\circ(\mathbf{a}_0)$ under $m$ is a union of $\calU^{\exp}$-orbits. Rewrite the convolution map as  $\Fl_{w_0}^{\exp}(\mathbf{a}_0) \widetilde\times \Fl_{v'}^\circ(\mathbf{a}_0) \widetilde\times \Fl_{v_{\min}}^\circ(\mathbf{a}_0) \r \Fl$, and convolve the left and middle factors first. 
By \thref{conv-GB}, the scheme-theoretic agrees with that of $\Fl_{w_0}(\mathbf{a}_0) \widetilde\times \Fl_{v_{\min}}(\mathbf{a}_0) \r \Fl$, which is a union of $\calI$-orbits by \thref{Iconvolution}.
\xpf

\coro \thlabel{affinestrat}
Let $w \in W_{\mathbf{f}}^{\exp}$.
\begin{enumerate}
    \item The scheme $\overline{\Fl}^{\exp}_{w}(\mathbf{f})$ is set-theoretically a disjoint union of locally closed subschemes of the form $\Fl^{\exp}_{v}(\mathbf{f})$, for certain $v \in W_{\mathbf{f}}^{\exp}$ which are independent of the base $S$.
        \item The natural map $\Fl^{\exp}_{w}(\mathbf{f}) \rightarrow \overline{\Fl}^{\exp}_{w}(\mathbf{f})$ is an open immersion.
    \item For any base $S$, the natural map
    $$\overline{\Fl}^{\exp}_{w}(\mathbf{f})_S \rightarrow \overline{\Fl}^{\exp}_{w}(\mathbf{f}) \times S$$ is a nilpotent thickening.
\end{enumerate}
\xcoro

\pf
If $w \in W/W_{\mathbf{f}}$ and $w$ is not maximal in its left $W_0$-orbit, then $\Fl^{\exp}_{w}(\mathbf{f}) = \Fl^{\circ}_{w}(\mathbf{a}_0, \mathbf{f})$ by \thref{orbit-structure}. In this case, all of the claims follow from the Bruhat decomposition \cite[Lemma 4.3.7, Proposition 4.4.3]{RicharzScholbach:Intersection}. If $w \in {}_0W_{\mathbf{f}}$ then $\Fl^{\exp}_{w}(\mathbf{f}) \subset \Fl^{\circ}_{w}(\mathbf{a}_0, \mathbf{f})$ is a dense open immersion so the claims again follow from loc.~cit. Finally, suppose $w \in W/W_{\mathbf{f}}$ and $w$ is maximal in its left $W_0$-orbit. Then (1) follows from \thref{conv-image} by taking $w_1 \in W$ the canonical lift under \thref{Wexp-lift}, and $w_2$ trivial. As in the proof of \thref{theo--strat}, this immediately implies (3), and (2) holds since $\Fl^{\exp}_{w}(\mathbf{f}) \subset \Fl^{\circ}_{w}(\mathbf{a}_0, \mathbf{f})$ is closed and $\Fl^{\circ}_{w}(\mathbf{a}_0, \mathbf{f}) \subset \Fl_{w}(\mathbf{a}_0, \mathbf{f})$ is open.
\xpf

\subsection{Fibers of exponential convolution morphisms}

Haines \cite{Haines:Pavings} showed that the fibers of convolution morphisms of 
the form $\Fl_{w_1}^{\circ}(\bbf) \widetilde{\times} \Fl_{w_2}^{\circ}(\bbf) \widetilde{\times} \cdots \widetilde{\times} \Fl_{w_n}^{\circ}(\bbf) \to \Fl_{\bbf}$, for $w_i \in W_{\mathbf{f}}\backslash W/W_{\mathbf{f}}$, can be decomposed into products of \(\A^1\)'s and \(\Gm\)'s. In this section, we consider analogous exponential convolution morphisms where $\Fl_{w_1}^{\circ}(\bbf)$ is replaced with $\Fl_{w_1}^{\exp}(\bbf)$ for some $w_1 \in W_{\mathbf{f}}^{\exp}$. A product of \(\A^1\)'s and \(\Gm\)'s is too much to hope for, as the following example shows.

\exam\thlabel{example punctured Gm}
Let us consider convolution on the finite flag variety for \(G=\mathrm{SL}_2\).
Let \(s\in W_0\) be the unique simple reflection, and let \(z\in {}_0W_{\mathbf{a}_0}^{\exp}\) be the representative of the open exponential orbit in $G/B$. Identify the flag variety as $$\P = \Gm\sqcup \{0\} \sqcup \{\infty\} \cong \Fl_z^{\exp}(\mathbf{a}_0) \sqcup \Fl_s^{\exp}(\mathbf{a}_0) \sqcup \Fl_e^{\exp}(\mathbf{a}_0) = \Fl_{s}(\mathbf{a}_0) \subset \Fl.$$
In particular, we identify \(\Fl_s^{\circ}(\mathbf{a}_0) \subseteq \Fl_{s}(\mathbf{a}_0)\) with \(\A^1\subset \P\).

The proof of \cite[Proposition 3.19]{RicharzScholbach:Motivic} shows that the convolution map \(\Fl_s^{\circ}(\mathbf{a}_0)\widetilde{\times}\Fl_s^{\circ}(\mathbf{a}_0)\to \Fl_s(\mathbf{a}_0)\) can be identified with
\((\A^1\times \P)\setminus \Delta(\A^1) \to \P\), via the projection onto the second factor.
Similarly, since convolving with \(\Fl_s^{\exp}(\mathbf{a}_0) = \{\dot{s} \cdot \overline{e}\}\subset \Fl\) on the left is simply a translation, we can identify the convolution \(\Fl_s^{\exp}(\mathbf{a}_0) \widetilde{\times}  \Fl_s^{\circ}(\mathbf{a}_0) \to \Fl_{s}(\mathbf{a}_0)\) with a map \(\{\infty\} \times \A^1\to \P\), whose image is \(\Gm \sqcup \{\infty\}\subset \P\).
Using this, one can compute the fibers of \(m\colon \Fl_z^{\exp}(\mathbf{a}_0) \widetilde{\times} \Fl_s^{\circ}(\mathbf{a}_0)\to \Fl_{s}(\mathbf{a}_0)\cong \P\) as follows.
\[m^{-1}(x)\cong\begin{cases}
	\Gm & \text{ if } x=0\\
	\Gm & \text{ if } x=\infty\\
	\Gm\setminus\{1\} & \text{ else. }
\end{cases}\]
The union of these fibers yields \(\Gm\sqcup \Gm\sqcup ((\Gm\setminus \{1\})\times \Gm) \cong \A^1\times \Gm\).
Thus, the fibers of exponential convolution morphisms do not necessarily admit a filtrable decomposition by products of \(\A^1\)'s and \(\Gm\)'s.
\xexam

In the example, we denoted \(m^{-1}(x):=m^{-1}(\dot{x} \cdot \overline{e})\); we will use this shorthand notation from now on.

Of course, having cells isomorphic to the complement of finitely many rational points of \(\A^1\) does not yield problems with respect to Tateness.
Thus, throughout this paper we will use the following definition.

\defi \thlabel{defi--filtrable}
A scheme \(X/S\) admits a \emph{cellular filtrable decomposition} if it admits a sequence of closed subschemes \(\varnothing = X_0 \subsetneq X_1\subsetneq \ldots \subsetneq X_n=X\), such that the reduced subscheme of each successive complement \(X_i\setminus X_{i-1}\) is isomorphic to a product of schemes in \(\{\A^1,\Gm,\Gm\setminus \{1\} \}\).
\xdefi

For a scheme \(X/S\) as above, any cell \(X_i\setminus X_{i-1}\) is \emph{admissible} in the sense of \cite[Definition 2.9]{CassvdHScholbach:MotivicSatake}.
The main result concerning fibers of exponential convolution morphisms is the following:

\theo\thlabel{fibers of exponential convolution}
Let \(\bbf, \bbf_1,\ldots,\bbf_n\) be standard facets, with corresponding partial affine flag varieties \(\Fl_{\bbf_i}\).
Let \(v_i\in W_{\bbf_{i-1}}\backslash W/W_{\bbf_i}\) for \(i=1,\ldots,n\), as well as \(v_0\in W^{\exp}_{\bbf}\), and consider the convolution morphism
\[m\colon \Fl_{v_0}^{\exp}(\bbf) \widetilde{\times} \Fl_{v_1}^\circ(\bbf,\bbf_1) \widetilde{\times} \ldots \widetilde{\times} \Fl_{v_n}^\circ(\bbf_{n-1},\bbf_n) \to \Fl_{\bbf_n}.\]
Then for any \(w\in W_{\bbf_n}^{\exp}\subseteq \Fl_{\bbf_n}\), the fiber \(m^{-1}(w)\) admits a cellular filtrable decomposition.
\xtheo

Here, we implicitly write \(\bbf_{i-1}=\bbf\) when \(i=1\); we have refrained from using \(\bbf_0\), since this already denotes the facet containing the origin.

The above theorem refines the main result of \cite{Haines:Pavings} by allowing exponential orbits, but we also allow different parahorics, rather than setting \(\bbf=\bbf_1=\ldots=\bbf_n\) as in loc.~cit.
Before we prove the theorem, we need some preparations.

\lemm\thlabel{triviality of convolution}
Fix notation as in \thref{fibers of exponential convolution}.
\begin{enumerate}
	\item The \(\calI = \calU\rtimes T\)-action on the convolution product \(\Fl_{\bbf} \widetilde{\times} \ldots \widetilde{\times} \Fl_{\bbf_n}\) factors through \(\calU\rtimes T_{\adj}\).
	In particular, \(\Fl_{\bbf} \widetilde{\times} \ldots \widetilde{\times} \Fl_{\bbf_n}\) admits a natural \(\calU^{\exp}\)-action, for which the convolution map is equivariant.
	\item The convolution morphism
	\[m\colon\Fl_{v_0}^{\exp}(\bbf) \widetilde{\times} \Fl_{v_1}^\circ(\bbf,\bbf_1) \widetilde{\times} \ldots \widetilde{\times} \Fl_{v_n}^\circ(\bbf_{n-1},\bbf_n) \to \Fl_{\bbf_n}\]
	is trivial over every exponential orbit in its image, in the sense of \cite[Definition 5.1]{Haines:Pavings}.
\end{enumerate}
\xlemm

\pf
(1): the \(\calU \rtimes T\)-action on \(\Fl_{\bbf}\) factors through \(\calU \rtimes T_{\adj}\) by \thref{lemm--UexpAction}.
Consider the isomorphism
\[\Fl_{\bbf} \widetilde{\times} \ldots \widetilde{\times} \Fl_{\bbf_n}\cong \Fl_{\bbf}\times \ldots\times \Fl_{\bbf_n}\colon (g_0,\ldots,g_n) \mapsto (g_0, \ldots, g_0\cdot \ldots \cdot g_n).\]
This is \(\calU\rtimes T\)-equivariant for the diagonal action on the target.
Since the action on the target factors through \(\calU \rtimes T_{\adj}\), we get the factorization on the source as well.
It is clear that the resulting convolution map \(\Fl_{\bbf} \widetilde{\times} \ldots \widetilde{\times} \Fl_{\bbf_n} \to \Fl_{\bbf_n}\) is equivariant for the \(\calU \rtimes T_{\adj}\)-action, and hence for the \(\calU^{\exp}\)-action.

(2): let \(\Fl_w^{\exp}(\mathbf{f}_n)\) be an exponential orbit in the image of \(m\).
Recall from \thref{orbit-structure} \eqref{it--global section} that we can find a subscheme \(H_w\subseteq \calU^{\exp}\) such that the action map induces an isomorphism \(\phi\colon H_w\cong \Fl_w^{\exp}(\mathbf{f}_n)\colon h\mapsto h\cdot w\).
This yields an isomorphism
\[m^{-1}(w)\times \Fl_w^{\exp}(\mathbf{f}_n) \cong m^{-1}(w)\times H_w \xrightarrow{\cong} m^{-1}(\Fl_w^{\exp}(\mathbf{f}_n)) \colon (x,h)\mapsto h\cdot x.\]
Indeed, the inverse is given by sending \(y\to (h^{-1}y,h)\), with \(h=\phi^{-1}(m(y))\).
\xpf

The following lemma is the key to upgrade the results from \cite{Haines:Pavings} to our setting.
Following \cite[§6.1]{Haines:Pavings}, for \(g_1,g_2\in LG\) and \(v\in W\), we write \(g_1\mathcal{I} \overset{v}{\text{---}} g_2\mathcal{I}\) if \(g_1^{-1}g_2\in \mathcal{I}v\mathcal{I}\).

\lemm\thlabel{key lemma cellularity}
Let \(v,w\in W_{\bba_0}^{\exp}\), and \(s\) an affine simple reflection.
Then the intersection \(\{\mathcal{I}'\in \Fl\mid w\mathcal{I}\overset{s}{\text{---}} \mathcal{I}'\}\cap \Fl_v^{\exp}(\bba_0)\) inside \(\Fl=\Fl_{\bba_0}\) is isomorphic to a scheme in \(\{\varnothing,\A^0,\A^1,\Gm,\Gm\setminus \{1\}\}\).
\xlemm
\pf 
Note that \(\A^1\cong \{\mathcal{I}'\in \Fl\mid w\mathcal{I}\overset{s}{\text{---}} \mathcal{I}'\} \subseteq \Fl_w^{\circ}(\bba_0) \cup \Fl_{ws}^{\circ}(\bba_0)\); it thus suffices to consider exponential orbits contained in this union.
We distinguish several cases.

\textbf{(I) First, we assume that \(w\in W\subseteq W_{\bba_0}^{\exp}\), and that \(w<ws\).}
The subscheme \(\{\mathcal{I}'\in \Fl\mid w\mathcal{I}\overset{s}{\text{---}} \mathcal{I}'\}\) is isomorphic to \(\A^1\), and contained in \(\Fl_{ws}^{\circ}(\bba_0)\) by \cite[§6.1]{Haines:Pavings}; we may thus assume \(\Fl_v^{\exp}(\bba_0)\subseteq \Fl_{ws}^{\circ}(\bba_0)\).
More precisely, the intersection is the subset of \(\Fl_{ws}^{\circ}(\bba_0)\) given by \(U_{w(\alpha)}\cdot ws\), where \(\alpha\) is the simple affine root corresponding to \(s\).
Indeed, we have 
\[\calI s \calI = \calI U_{\alpha} s \calI = \calI w^{-1} w U_{\alpha} w^{-1}ws \calI = \calI w^{-1} U_{w(\alpha)} ws \calI.\]
Now, the map \(U_{w(\alpha)}\to \Ga\) obtained by restricting \eqref{map from U to A1} is either an isomorphism or trivial, depending on whether \(w(\alpha)\in \Delta\).
Thus, when \(\Fl_{ws}^{\circ}(\bba_0)\) splits up into two exponential orbits, \(U_{w(\alpha)}\cdot ws\subseteq \Fl_{ws}^{\circ}(\bba_0)\) intersects these exponential orbits as \(\A^0\sqcup \Gm\) in the first case, and as \(\A^1\sqcup \varnothing\) in the second case (if \(\Fl_{ws}^{\circ}(\bba_0)\) is itself an exponential orbit, the desired intersection does not split up, and is hence isomorphic to \(\A^1\)).

\textbf{(II) Next, we still assume \(w\in W\subseteq W_{\bba_0}^{\exp}\), but that \(ws<w\).}
Then \cite[Lemma 6.2]{Haines:Pavings} implies that \(\{\mathcal{I}'\in \Fl\mid w\mathcal{I}\overset{s}{\text{---}}\mathcal{I}'\} \cap \Fl_w^{\circ}(\bba_0)\cong \Gm\), so that we may assume \(\Fl_v^{\exp}(\bba_0)\subseteq \Fl_{w}^{\circ}(\bba_0)\).
More specifically, this intersection agrees with \(U_{-w(\alpha)}^\times \cdot w\). Indeed, we have
\(\calI w^{-1} U_{-w(\alpha)} w\calI = \calI U_{-\alpha}\calI\), for which \(\calI \cdot 1\cdot \calI = \calI\) and \(\calI U_{-\alpha}^\times \calI \subseteq \calI s \calI\).
Moreover, \(U_{-w(\alpha)}^\times\) either lies in the kernel of \eqref{map from U to A1}, or does not intersect it, so that the intersections with the exponential orbits are isomorphic to \(\Gm\) or empty.

\textbf{(III) Finally, let \(w\in {}_0W_{\bba_0} \subseteq  W_{\bba_0}^{\exp}\),} and let \(z\in W\) be such that \(\Fl_w^{\exp}(\mathbf{a}_0)\sqcup \Fl_z^{\exp}(\mathbf{a}_0) = \Fl_z^{\circ}(\bba_0)\).
Let \(a\in \mathcal{U}\) be an element such that \(a\cdot z = w\in \Fl\); then \(a\) maps to \(1\in \Ga\) under \eqref{map from U to A1}.
Moreover, the set \(\{\mathcal{I}'\in \Fl\mid w\mathcal{I}\overset{s}{\text{---}} \mathcal{I'}\}\) can be obtained by translating the set \(\{\mathcal{I}'\in \Fl\mid z\mathcal{I}\overset{s}{\text{---}} \mathcal{I'}\}\).
Although translating by \(a\) preserves the Iwahori-orbits, it does not preserve the exponential orbits when they do not agree with Iwahori-orbits; we may thus assume that \(\Fl_v^{\exp}(\bba_0)\) is not an Iwahori-orbit by appealing to \cite[Lemma 6.2]{Haines:Pavings}.
The intersection of \(\{\mathcal{I}'\in \Fl\mid w\mathcal{I}\overset{s}{\text{---}} \mathcal{I'}\}\) with exponential orbits is then determined by the preimage of \(-1\in \Ga\) under \(\calU\to \Ga\), and more specifically the restrictions of this map to the group schemes appearing above.
In order to make this explicit, we distinguish further cases.

\textbf{(IIIa)} If \(z<zs\), we may assume \(\Fl_v^{\exp}(\bba_0)\subseteq \Fl_{zs}^{\circ}(\bba_0)\) as above. 
Let us write \(\Fl_{zs}^{\circ}(\bba_0) = \Fl_{zs}^{\exp}(\bba_0) \sqcup \Fl_y^{\exp}(\bba_0)\), where \(\Fl_y^{\exp}(\bba_0)\) is the open exponential orbit.
Then, the intersection of \(\{\mathcal{I}'\in \Fl\mid w\mathcal{I}\overset{s}{\text{---}} \mathcal{I'}\}\cong \A^1\subseteq \Fl_{zs}^{\circ}(\bba_0)\) with \(\Fl_{zs}^{\exp}(\bba_0) \sqcup \Fl_y^{\exp}(\bba_0)\) yields
\[\begin{cases}
	\A^1\sqcup \varnothing &\text{ if } U_{z(\alpha)} \to \Ga \text{ is trivial,}\\
	\Gm\sqcup \A^1 & \text{ if } U_{z(\alpha)} \to \Ga \text{ is an isomorphism}.
\end{cases}\]

\textbf{(IIIb)} If \(zs<z\), we may assume \(\Fl_v^{\exp}(\bba_0) \subseteq \Fl_z^{\circ}(\bba_0)\), and we now write \(\Fl_{z}^{\circ}(\bba_0) = \Fl_{z}^{\exp}(\bba_0) \sqcup \Fl_y^{\exp}(\bba_0)\), with \(\Fl_y^{\exp}(\bba_0)\) the open exponential orbit.
Then the intersection of
\(\{\mathcal{I}'\in \Fl\mid z\mathcal{I}\overset{s}{\text{---}} \mathcal{I'}\}\) with \(\Fl_{z}^{\exp}(\bba_0) \sqcup \Fl_y^{\exp}(\bba_0)\) yields
\[\begin{cases}
	\Gm\sqcup \varnothing &\text{ if } U_{-z(\alpha)}^{\times} \to \Ga \text{ is trivial,}\\
	\Gm\setminus\{-1\}\sqcup \A^1 & \text{ if } U_{-z(\alpha)}^{\times} \to \Ga \text{ is an open immersion}.
\end{cases}\]
This concludes all the necessary case distinctions, and provides the desired descriptions of the intersections.
\xpf

\exam
We generalize \thref{example punctured Gm} to general \(G\) as follows. 
We consider the convolution in the finite flag variety \(G/B\), which we view as a subscheme \(\Fl_{w_0}(\bba_0)\subseteq \Fl\) of the full affine flag variety.
Then there is a unique Iwahori-orbit splitting up into two exponential orbits, corresponding to the longest element \(w_0\in W_0\subseteq W\) in the finite Weyl group.
Let \(s\in W_0\) be a simple reflection of the finite Weyl group, and \(z\) correspond to the unique open exponential orbit in \(\Fl_{w_0}(\bba_0)\).
Then the proof of the lemma above yields the following description of the desired intersections in the interesting cases:
\begin{center}
	\begin{tabular}{ c || c | c | c|} 
		\(\cap\) & \(\Fl_{w_0}^{\exp}(\bba_0)\) & \(\Fl_z^{\exp}(\bba_0)\) & \(\Fl_{w_0s}^{\exp}(\bba_0)\)\\
		\hline
		\hline
		\(\{\mathcal{I}'\in \Fl\mid w_0\mathcal{I} \overset{s}{\text{---}} \mathcal{I'}\}\) & \(\varnothing\) & \(\Gm\) & \(\A^0\) \\
		\hline
		\(\{\mathcal{I}'\in \Fl\mid z\mathcal{I} \overset{s}{\text{---}} \mathcal{I'}\}\) & \(\A^0\) & \(\Gm \setminus \{1\}\) & \(\A^0\) \\
		\hline
		\(\{\mathcal{I}'\in \Fl\mid w_0s\mathcal{I} \overset{s}{\text{---}} \mathcal{I'}\}\) & \(\A^0\) & \(\Gm\) & \(\varnothing\) \\
		\hline
	\end{tabular}
\end{center}
\xexam

We can now prove \thref{fibers of exponential convolution} in a special case.

\prop\thlabel{fibers for full flag variety}
\thref{fibers of exponential convolution} holds for the full affine flag variety, i.e., when \(\bbf=\ldots =\bbf_n = \bba_0\).
\xprop
\pf
We keep the notation of \thref{fibers of exponential convolution}.
For \(i=1,\ldots,n\), let us fix reduced expressions \(v_i = \tau_i s_{i1}\ldots s_{ir}\), with \(\tau_i\) of length 0.
This yields isomorphisms \(\Fl_{v_i}^{\circ}(\bba_0) \cong \Fl_{\tau_i}^{\circ}(\bba_0)\widetilde{\times} \Fl_{s_{i1}}^{\circ}(\bba_0) \widetilde{\times} \ldots \widetilde{\times} \Fl_{s_{ir}}^{\circ}(\bba_0)\). 
Since conjugation by \(\tau_i\) normalizes \(\calU\) and permutes the simple affine reflections, we obtain a further isomorphism
\[\Fl_{v_i}^{\circ}(\bba_0) \cong \Fl_{\tau_i s_{i1} \tau_i^{-1}}^{\circ}(\bba_0) \widetilde{\times} \ldots \widetilde{\times} \Fl_{\tau_i s_{ir} \tau_i^{-1}}^{\circ}(\bba_0) \widetilde{\times} \Fl_{\tau_i}^{\circ}(\bba_0),\]
compatibly with the natural maps to \(\Fl\).
Thus, we can reduce the proposition to the case where each \(v_i=s_i\) is a simple reflection.
On the other hand, by \thref{convlem}, any Iwahori-orbit that splits up into two exponential orbits can be obtained by convolving \(\Fl_{w_0}^\circ(\bba_0)\) with a (smaller) Iwahori-orbit.
Thus, we may assume that \(v_0=w_0\) is the longest element in the finite Weyl group \(W_0\) of \(G\), or that \(v_0\) is such that \(\Fl_{v_0}^{\exp}(\bba_0)\sqcup \Fl_{w_0}^{\exp}(\bba_0) = \Fl_{w_0}^{\circ}(\bba_0)\).

Next, we proceed by induction, the case \(n=0\) being trivial.
Note that the composition 
\[m^{-1}(w) \subseteq \Fl_{v_0}^{\exp}(\bba_0) \widetilde{\times} \Fl_{s_1}^\circ(\bba_0) \widetilde{\times} \ldots \widetilde{\times} \Fl_{s_{n}}^\circ(\bba_0) \to \Fl_{v_0}^{\exp}(\bba_0) \widetilde{\times} \Fl_{s_1}^\circ(\bba_0) \widetilde{\times} \ldots \widetilde{\times} \Fl_{s_{n-1}}^\circ(\bba_0)\] 
induces an isomorphism onto its image, where the second morphism is the projection away from the last factor.
Using this, consider the convolution morphism
\[m'\colon \Fl_{v_0}^{\exp}(\bba_0) \widetilde{\times} \Fl_{s_1}^\circ(\bba_0) \widetilde{\times} \ldots \widetilde{\times} \Fl_{s_{n-1}}^\circ(\bba_0) \to \Fl ,\]
as well as the map 
\[\xi\colon m^{-1}(w)\to \Fl_{ws_n}^{\circ}(\bba_0) \sqcup \Fl_w^{\circ}(\bba_0),\]
given by the restriction of \(m'\).
We claim that if the intersection of \(\im(\xi)\) with an exponential orbit in \(\Fl_{ws_n}^{\circ}(\bba_0)\sqcup \Fl_w^{\circ}(\bba_0)\) is nonempty, then it is isomorphic to \(\A^0\), \(\A^1\), \(\Gm\), or \(\Gm\setminus\{1\}\).
Indeed, convolution morphisms are \(\calU^{\exp}\)-equivariant by \thref{triviality of convolution}, so an exponential orbit in \(\Fl_{ws_n}^{\circ}(\bba_0)\sqcup \Fl_w^{\circ}(\bba_0)\) is either contained in \(\im(m')\), or does not meet \(\im(m')\).
In the former case, the intersection is described by \thref{key lemma cellularity}, and in the latter case the intersection is empty.
This implies the proposition, as the map \(m'\) has cellular fibers by induction, and \(m'\) is trivial over the exponential orbits in its image by \thref{triviality of convolution}.
\xpf

To deduce \thref{fibers of exponential convolution} in general, we will use the following lemma.

\lemm\thlabel{intersections and convolution fibers}
Let \(\bbf,\bbf'\) be standard facets, and choose elements \(v\in W_{\bbf}^{\exp}\), \(v'\in W_{\bbf}\backslash W/W_{\bbf'}\) and \(w\in W_{\bbf'}^{\exp}\). 
Consider the convolution morphism \(m_{v,v'}\colon \Fl_{v}^{\exp}(\bbf) \widetilde{\times} \Fl_{v'}^\circ(\bbf,\bbf') \to \Fl_{\bbf'}\).
Then there is a canonical isomorphism of schemes
\[\Fl_{v}^{\exp}(\bbf)\cap w\cdot \Fl_{(v')^{-1}}^\circ(\bbf',\bbf)\cong m_{v,v'}^{-1}(w),\]
where we have implicitly chosen a lift of \(w\) to \(LG\).
\xlemm
\pf
Similarly to \cite[Lemma 7.3]{Haines:Pavings}, both sides represent the étale sheafifications of the presheaf sending an \(S\)-algebra \(R\) to those \(gL^+G_{\bbf}(R)\in \Fl_{v}^{\exp}(\bbf)(R)\) such that \(g^{-1}w\in L^+G_{\bbf}(R)v_1L^+G_{\bbf'}(R)\).
\xpf

\pf[Proof of \thref{fibers of exponential convolution}]
Let \(v_0,\ldots,v_n,w\) be as in the statement of \thref{fibers of exponential convolution}, and 
\[m\colon \Fl_{v_0}^{\exp}(\bbf) \widetilde{\times} \Fl_{v_1}^\circ(\bbf,\bbf_1) \widetilde{\times} \ldots \widetilde{\times} \Fl_{v_n}^\circ(\bbf_{n-1},\bbf_n)\to \Fl_{\bbf_n}\]
the corresponding convolution morphism.
We again consider the restricted convolution
\[m'\colon \Fl_{v_0}^{\exp}(\bbf) \widetilde{\times} \Fl_{v_1}^\circ(\bbf,\bbf_1) \widetilde{\times} \ldots \widetilde{\times} \Fl_{v_{n-1}}^\circ(\bbf_{n-2},\bbf_{n-1})\to \Fl_{\bbf_{n-1}},\]
and let \(\xi\colon m^{-1}(w)\to \Fl_{\bbf_{n-1}}\) be its restriction to \(m^{-1}(w)\).
By \thref{intersections and convolution fibers}, the image of \(\xi\) agrees with \(\im(m')\cap w\Fl_{v_n^{-1}}(\bbf_n,\bbf_{n-1})\).
As in the proof of \thref{fibers for full flag variety}, it will suffice to show that for any \(y\in W^{\exp}_{\bbf_{n-1}}\), the intersection \(\im(\xi)\cap \Fl_y^{\exp}(\bbf_{n-1})\) admits a cellular filtrable decomposition.
Since convolution morphisms are \(\calU^{\exp}\)-equivariant by \thref{triviality of convolution}, if such an intersection is nonempty we have \(\Fl_y^{\exp}(\bbf_{n-1})\subseteq \im(m')\).
Thus, it suffices to consider the intersections \(\Fl_y^{\exp}(\bbf_{n-1})\cap w\Fl_{v_n^{-1}}(\bbf_n,\bbf_{n-1})\).
Clearly, it even suffices to consider the intersections \(\Fl_y^{\exp}(\bbf_{n-1})\cap w\Fl_{z}(\bba_0,\bbf_{n-1})\) for any lift \(z\in W/W_{\bbf_{n-1}}\) of \(v_{n-1}\). 

Below, we will write \(y\in W_{\bba_0}^{\exp}\) for the minimal length representative of \(y\), and similarly for \(z\in W\).
By \cite[Lemma 7.2]{Haines:Pavings} (which works over any base \(S\)), the projection \(\Fl\to \Fl_{\bbf_{n-1}}\) induces the following isomorphism, after passing to reduced subschemes:
\[\Fl_y^\circ(\bba_0)\cap \left(\bigsqcup_{u\in W_{\bbf_{n-1}}} w\Fl_{zu}^\circ(\bba_0)\right) \cong \Fl_y^\circ (\bba_0,\bbf_{n-1})\cap w\Fl_z^\circ(\bba_0,\bbf_{n-1}).\]
Since \(\Fl\to \Fl_{\bbf_{n-1}}\) is \(\calU^{\exp}\)-equivariant, this restricts to a similar isomorphism
\[\Fl_y^{\exp}(\bba_0)\cap \left(\bigsqcup_{u\in W_{\bbf_{n-1}}} w\Fl_{zu}^\circ(\bba_0)\right) \cong \Fl_y^{\exp} (\bbf_{n-1})\cap w\Fl_z^\circ(\bba_0,\bbf_{n-1})\]
after passing to reduced subschemes.
Thus, we are reduced to showing that \(\Fl_y^{\exp}(\bba_0)\cap w\Fl_{zu}^\circ(\bba_0)\) admits a cellular filtrable decomposition.
By \thref{intersections and convolution fibers}, this is isomorphic to the fiber of an exponential convolution morphism for the facet \(\bba_0\), so we can conclude by \thref{fibers for full flag variety}.
\xpf

Let us now record the following consequence of \thref{fibers of exponential convolution}.
Although we will only need the result below for Iwahori-orbits (cf.~the proof of \thref{Equivalent conditions}), which can be deduced from \cite{Haines:Pavings}, the following generalization is immediate.

\coro\thlabel{cellularity of intersections of translates} 
Let \(\bbf,\bbf',\bbf''\) be standard facets, and choose elements \(v\in W^{\exp}_{\bbf}\), \(w\in W_{\bbf'}\backslash W/W_{\bbf}\), and \(u\in W_{\bba_0}^{\exp}\), for which we also fix a lift \(u\in LG(\Z)\).
Then the intersection \(\Fl_v^{\exp}(\bbf)\cap u\cdot \Fl_w^\circ(\bbf',\bbf)\subseteq \Fl_{\bbf}\) admits a cellular filtrable decomposition.

Consequently, if \(x\in W_{\bbf''}\backslash W/W_{\bbf}\), the intersection \(\Fl_x^{\circ}(\bbf'',\bbf) \cap u\cdot\Fl_w^{\circ}(\bbf',\bbf)\subseteq \Fl_{\bbf}\) also admits a cellular filtrable decomposition.
\xcoro
\pf
The intersection \(\Fl_v^{\exp}(\bbf)\cap u\cdot \Fl_w^\circ(\bbf',\bbf)\subseteq \Fl_{\bbf}\) is isomorphic to a fiber of a convolution morphism by \thref{intersections and convolution fibers}, so the first statement follows from \thref{fibers of exponential convolution}.

The second statement follows, since any parahoric orbit in \(\Fl_{\bbf}\) is stratified by exponential orbits.
\xpf

\section{Exponential motives on affine flag varieties}
\label{sect--exponential motives Fl}
We continue using the notation from \refsect{Kirillov}. We now assume the scheme $S$ is connected, smooth, and of finite-type over Dedekind ring or a field, so that we have the category $\DM(X)$ of integral motivic sheaves on finite-type $S$-schemes $X$.

In this section we show that the formalism of exponential stratified Tate motives applies to the $\calU^{\exp}$-orbit stratification. Highlights include Whitney--Tateness (\thref{exponential orbits WT}) and the existence of a t-structure (\thref{nunu} and \thref{defi exp on Fl}). We also investigate the t-exactness properties of exponential (co)standard and convolution functors. We conclude in §\ref{subsec--twisted exponential motives} with a slight variant of the $\calU^{\exp}$-orbit stratification. 

\subsection{Whitney--Tateness} \label{sect--WTExp}

For \(w\in W_{\mathbf{f}}^{\exp}\), let \(\iota_w\colon \mathcal{U}^{\exp}\backslash\Fl_w^{\exp}(\mathbf{f}) \to 
\mathcal{U}^{\exp}\backslash\Fl_{\mathbf{f}}\) be the inclusion of prestacks, and let \(h_w\colon \mathcal{U}^{\exp}\backslash\Fl_w^{\exp}(\mathbf{f}) \to \mathcal{U}^{\exp}\backslash S\) be the prestack quotient of the structure map, over a fixed base $S$ as in \thref{generalS}.

\defi\thlabel{Definition CoStandard}
For any \(w\in W_{\mathbf{f}}^{\exp}\), the standard functor is defined as
\[\Delta_w^{\exp} \colon \DM(\mathcal{U}^{\exp}\backslash S) \to \DM(\mathcal{U}^{\exp}\backslash \Fl_{\mathbf{f}})\colon \mathcal{F} \mapsto \iota_{w!}h_w^*\mathcal{F} [\dim \Fl_w^{\exp}(\mathbf{f})],\]
Similarly, the costandard functor for \(w\) is defined as 
\[\nabla_w^{\exp} \colon \DM(\mathcal{U}^{\exp}\backslash S) \to \DM(\mathcal{U}^{\exp} \backslash \Fl_{\mathbf{f}})\colon \mathcal{F} \mapsto \iota_{w*}h_w^*\mathcal{F} [\dim \Fl_w^{\exp}(\mathbf{f})].\]
\xdefi

We refer to \cite[\S 2.5]{CassvdHScholbach:MotivicSatake} for a discussion on why the above equivariant $*$ and $!$ functors exist, and that after forgetting the equivariance, they agree with the usual functors.
We also have the analogously defined (co)standard functors for $\calI$-orbits,
$$\Delta_w, \nabla_w \colon \DTM(\mathcal{I}\backslash S) \to \DTM(\mathcal{I\backslash \Fl_{\mathbf{f}}}), \quad w \in W/W_{\mathbf{f}}.$$

We will discuss properties of the (co)standard functors later. In this section we will only apply them to the motive $\Z$, and we will also not need to keep track of the equivariance data. This is because our goal is to show that the stratification in \thref{affinestrat} is Whitney--Tate, i.e., that after forgetting equivariance, for all $v, w \in W_{\mathbf{f}}^{\exp}$ the motive $\iota_v^* \nabla_w^{\exp}(\Z)$ lies in the presentable stable subcategory $\DTM(\Fl^{\exp}_v(\mathbf{f})) \subset \DM(\Fl^{\exp}_v(\mathbf{f}))$ generated by the Tate twists $\Z(n)$ for $n \in \mathbb{Z}$. In fact, we will show that $\iota_v^* \nabla_w^{\exp}(\Z)$ lies in the smaller subcategory  $\DTM^{\anti}(\Fl^{\exp}_v(\mathbf{f}))$ of anti-effective motives, generated by $\Z(n)$ for $n \leq 0$, and that the formation of this motive is universal in the sense that it commutes with base change in $S$ (see \cite[Eqn.~(2.4)]{CassvdHScholbach:MotivicSatake} for the precise definition). We start with the finite flag variety.

\theo \thlabel{WT-finite}
The stratification of $G/B = \Fl_{w_0}(\mathbf{a}_0)$ by $U^{\exp}$-orbits is anti-effective universally Whitney--Tate in the sense of \cite[Definition 2.6]{CassvdHScholbach:MotivicSatake}.
\xtheo

\pf
The same result for $U$-orbits (which are the same as $B$-orbits) follows from the more general affine case in \cite[Proposition 3.7]{CassvdHScholbach:MotivicSatake} (in the non-affine case this is essentially \cite[Proposition 4.10]{SoergelWendt:Perverse}). It remains to deal with the two $U^{\exp}$-orbits in the big open cell $\Fl_{w_0}^\circ(\mathbf{a}_0)$. Since $\Fl_{w_0}^\circ(\mathbf{a}_0)$ is dealt with in loc.~cit., then by relative purity and excision for $*$-pullback it suffices to consider only the smaller $U^{\exp}$-orbit $\Fl_{w_0}^{\exp}(\mathbf{a}_0)$, where we view $w_0$ as an element of $W \subset W_{\mathbf{a}_0}^{\exp}$. 

We need to show that if $v \in W_0$ is such that $l(v) < l(w_0)$, then $\iota_v^* \nabla^{\exp}_{w_0}( \Z) \in \DTM^{\anti}(\Fl_{v}^{\circ}(\mathbf{a}_0))$, and that the formation of this motive commutes with base change in $S$ as in \cite[Eqn. (2.4)]{CassvdHScholbach:MotivicSatake}. We will prove this by induction on $l(v)$, where the base case $l(v) = l(w_0) - 1$ holds since $\iota_v^* \nabla^{\exp}_{w_0}( \Z) = 0$ by \thref{theo--strat}. Now suppose $l(v) \leq l(w_0) - 2$. Let $\alpha \in \Delta$ be such that $l(vs_{\alpha}) = l(v) + 1$. Let $\mathbf{f}$ be the wall of $\mathbf{a}_0$ passing through the origin and associated to the reflection $s_{\alpha}$, so that $W_{\mathbf{f}} = \{1, s_{\alpha}\}$. Then $G/P_{\alpha} = \Fl_{w_0s_{\alpha}}(\mathbf{a}_0, \mathbf{f})$, and the projection $\Fl \r \Fl_{\mathbf{f}}$ restricts to a $\P$-fibration $\pi \colon G/B \r G/P{_{\alpha}}$.  

In the following commutative diagram the square is cartesian and the horizontal maps are the natural inclusions.
$$\xymatrix{
\Fl_{vs_{\alpha}}^\circ(\mathbf{a}_0) \sqcup \Fl_{v}^\circ(\mathbf{a}_0) \ar[r] \ar[rd] & \pi^{-1}(\Fl_{v}^\circ(\mathbf{a}_0, \mathbf{f})) \ar[r]^-{\iota} \ar[d]^{\pi'} &  G/B \ar[d]^\pi \\
& \Fl_{v}^\circ(\mathbf{a}_0, \mathbf{f}) \ar[r]^-{\iota'} & G/P_{s_\alpha}
}$$
Since $*$-pullback along any affine space $
\A^n \r S$ induces an equivalence $\DTM^{\anti}(\A^n) \cong \DTM^{\anti}(S)$, then by induction we may identify $\iota_{v s_{\alpha}}^*\nabla^{\exp}_{w_0}( \Z)$ with an element of $\DTM^{\anti}(S)$, and moreover we may assume its formation commutes with base change in $S$. The projection $\Fl_{vs_{\alpha}}^\circ(\mathbf{a}_0) \r \Fl_{v}^\circ(\mathbf{a}_0, \mathbf{f})$ is a relative one-dimensional affine space, so we may similarly identify the $!$-pushforward of $\iota_{v s_{\alpha}}^*\nabla^{\exp}_{w_0}( \Z)$ along this map with $\iota_{v s_{\alpha}}^*\nabla^{\exp}_{w_0}( \Z)(-1)[-2] \in \DTM^{\anti}(S)$. Additionally, we note that the projection $\Fl_{v}^\circ(\mathbf{a}_0) \r \Fl_{v}^\circ(\mathbf{a}_0, \mathbf{f})$ is an isomorphism.

Consider the excision sequence for $\iota^* \nabla^{\exp}_{w_0}( \Z)$ with respect to the top left map (for $*$-pullback). Applying $\pi_!'$ and the above considerations gives a fiber sequence
\begin{equation} \label{WTfiberG/B} \iota_{v s_{\alpha}}^*\nabla^{\exp}_{w_0}( \Z)(-1)[-2] \r \pi'_!\iota^* \nabla^{\exp}_{w_0}( \Z) \r \iota_{v}^* \nabla^{\exp}_{w_0}( \Z).\end{equation} By base change, $\pi'_!\iota^* \cong \iota'^*\pi_!$. The map $\pi$ sends $\Fl_{w_0}^{\exp}(\mathbf{a}_0)$ isomorphically onto $\Fl^\circ_{w_0s_{\alpha}}(\mathbf{a}_0, \mathbf{f})$ by \thref{Uexp-projection}. Since the stratification of $G/P_{s_{\alpha}}$ by $U$-orbits is anti-effective universally Whitney--Tate \cite[Proposition 3.7]{CassvdHScholbach:MotivicSatake} and $\pi_! = \pi_*$, then we may identify $\iota'^*\pi_! \nabla^{\exp}_{w_0}( \Z)$ with an element $\DTM^{\anti}(S)$, and the formation of this motive commutes with base change in $S$. Thus, the same is true for $\iota_{v}^* \nabla^{\exp}_{w_0}( \Z)$.
\xpf

\lemm \thlabel{Omega-translate}
Let $\tau \in \Omega \subset W$. Then right multiplication by $\dot{\tau}$ in $LG$ induces an $LG$-equivariant automorphism of $\Fl$ which maps $\calU^{\exp}$-orbits isomorphically onto $\calU^{\exp}$-orbits.
\xlemm

\pf
This is immediate from the fact that $\Omega$ normalizes $\calI$, and right multiplication by $\tau$ in $W$ preserves elements which are maximal in their left $W_0$-orbits. 
\xpf

\rema
 Left multiplication by $\dot{\tau}$ also induces an automorphism of $\Fl$ which maps $\calI$-orbits isomorphically onto $\calI$-orbits. However, in general $\tau$ does not normalize $\calU^{\exp}$, so that  this map does not preserve $\calU^{\exp}$-orbits. For example, let $G = \PGL_2$ with the standard upper triangular $B$ and $\Delta = \{\alpha\}$.
 Then  $W = \langle s_0, s_1, \tau \: | \: s_0^2=s_1^2=\tau^2=1, \tau s_0 = s_1 \tau \rangle$, where $s_0 = t_{-\alpha^\vee} s_{\alpha}$, $s_1 = s_{\alpha}$, and $\tau = t_{-\frac{\alpha^\vee}{2}} s_{\alpha}$. Then $s_1$ is maximal in its left $W_0 = \{1, s_1\}$-orbit, but $\tau s_1 = s_0\tau$ is not maximal. 
\xrema

\lemm \thlabel{piAnti}
Let $M \in \DM(\Fl)$ be a motive such that for each $w \in W_{\mathbf{a}_0}^{\exp}$, the $*$-restriction of $M$ to $\Fl^{\exp}_{w}(\mathbf{a}_0)$ lies in $\DTM^{\anti}(\Fl^{\exp}_{w}(\mathbf{a}_0))$. Let $\pi \colon \Fl \r \Fl_{\mathbf{f}}$ be the projection. Then for each $v \in W_{\mathbf{f}}^{\exp}$, the $*$-restriction of $\pi_!M$ to $\Fl^{\exp}_{v}(\mathbf{f})$ lies in $\DTM^{\anti}(\Fl^{\exp}_{v}(\mathbf{f}))$.
\xlemm

\pf
By \thref{Uexp-projection}, the fiber of any $\calU^{\exp}$-orbit in $\Fl_{\mathbf{f}}$ is set-theoretically a union of locally closed $\calU^{\exp}$-orbits in $\Fl$. The same result, together with base change and excision for $*$-pullback, reduces the lemma to the following straightforward statement: the $!$-pushforward of $\Z$ along a projection of the form $\A^n \times \Gm \r S$ or $\A^n \r S$ is anti-effective Tate.
\xpf

\theo\thlabel{exponential orbits WT}
The stratification of $\Fl_{\mathbf{f}}$ by $\calU^{\exp}$-orbits is anti-effective universally Whitney--Tate.
\xtheo

\pf
We first suppose that $\mathbf{f} = \mathbf{a}_0$. The same result is known for $\calI$-orbits by \cite[Proposition 3.7]{CassvdHScholbach:MotivicSatake} (which is essentially \cite[Theorem 5.1.1]{RicharzScholbach:Intersection}). As in the proof of \thref{WT-finite}, by \thref{orbit-structure} it suffices to consider only those $\calU^{\exp}$-orbits of the form $\Fl_{w}^{\exp}(\mathbf{a}_0)$, where $w \in W \subset W_{\mathbf{a}_0}^{\exp}$ is maximal in its left $W_0$-orbit. For such $w$, write $w = w_0 w_{\min}$ where $w_{\min}$ is minimal in its left $W_0$-orbit. We will prove the necessary facts about $\nabla_{w}^{\exp}(\Z)$ by induction on $l(w_{\min})$.
For the base case, if $l(w_{\min}) = 0$ then $w_{\min} \in \Omega$ and hence $\Fl_{w}^{\exp}(\mathbf{a}_0)$ is the right $\dot{w}_{\min}$-translate of $\Fl_{w_0}^{\exp}(\mathbf{a}_0) \subset G/B$ by \thref{Omega-translate}. Thus, the case $l(w_{\min}) = 0$ reduces to \thref{WT-finite}.

If $l(w_{\min}) > 0$, we may let $s \in W$ be a simple affine reflection such that $l(w_{\min}s) < l(w_{\min})$. The maximality of $w$ is equivalent to the statement that any product of reduced words for $w_0$ and $w_{\min}$ gives a reduced word for $w = w_0 w_{\min}$. Hence, by choosing a reduced word for $w_{\min}$ which ends in $s$, it follows that $w_0 w_{\min}s = ws$ is also maximal in its left $W_0$-orbit.  Let $\mathbf{f}'$ be the wall of $\mathbf{a}_0$ passing through the origin and associated to $s$, so that $W_{\mathbf{f}'} = \{1, s\}$ and we have the $\P$-fibration $\pi \colon \Fl \r \Fl_{\mathbf{f'}}$. By \thref{Uexp-projection} we have the following commutative diagram, where the square is cartesian and the horizontal maps are the natural inclusions.
$$\xymatrix{
\Fl_{w}^{\exp}(\mathbf{a}_0) \sqcup \Fl_{ws}^{\exp}(\mathbf{a}_0) \ar[r] \ar[rd] & \pi^{-1}(\Fl_{ws}^{\exp}(\mathbf{f}')) \ar[r]^-{\iota} \ar[d]^{\pi'} &  \Fl \ar[d]^\pi \\
& \Fl_{ws}^{\exp}(\mathbf{f}') \ar[r]^-{\iota'} & \Fl_{\mathbf{f}'}
}$$
Additionally, by \thref{Uexp-projection}, the map $\Fl_{w}^{\exp}(\mathbf{a}_0) \r \Fl_{ws}^{\exp}(\mathbf{f}')$ is a relative one-dimensional affine space, and the map $\Fl_{ws}^{\exp}(\mathbf{a}_0) \r \Fl_{ws}^{\exp}(\mathbf{f}')$ is an isomorphism.

We are now in a geometric setup analogous to the proof of \cite[Proposition 3.7, Eqn.~(3.2)]{CassvdHScholbach:MotivicSatake}, so that in particular, excision with respect to the top left horizontal map leads to a fiber sequence
\begin{equation} \label{FlWTseq} \nabla^{\exp}_{ws} (\Z)(-1)[-1] \r \pi^! \pi_! \nabla^{\exp}_{ws} (\Z)(-1)[-1] \r \nabla^{\exp}_{w} (\Z).\end{equation}
By induction, the $*$-pullback of $\nabla^{\exp}_{ws} (\Z)(-1)[-1]$ to each $\calU^{\exp}$-orbit is anti-effective Whitney--Tate, and the formation of this motive commutes with base change in $S$. Since $\pi$ is smooth and proper, the formation of $\pi^! \pi_! \nabla^{\exp}_{ws} (\Z)(-1)[-1]$ also commutes with base change in $S$, so the same is true of $\nabla^{\exp}_{w} (\Z)$. Next, note that $\pi^! \pi_! = \pi^* \pi_!(1)[2]$.
Thus, to complete the induction, it remains to see that $\pi^* \pi_!$ preserves the class of motives having the property that $*$-pullback to each $\calU^{\exp}$-orbit is anti-effective Whitney--Tate. Since $\pi$ maps $\calU^{\exp}$-orbits onto $\calU^{\exp}$-orbits, it is obvious that $\pi^*$ has this property, and $\pi_!$ is handled in \thref{piAnti}.

Now suppose that $\mathbf{f}$ is arbitrary, and let $w \in W_{\mathbf{f}}^{\exp}$. Let $\tilde{w} \in W_{\mathbf{a}_0}^{\exp}$ be the lift of $w$ in \thref{Wexp-lift}, so that the inclusion of $\Fl_w^{\exp}(\mathbf{f})$ factors as $\Fl_{\tilde{w}}^{\exp}(\mathbf{a}_0) \r \Fl \rightarrow \Fl_{\mathbf{f}}$. By base change for the proper map $\pi$, the universality isomorphism \cite[Eqn. (2.4)]{CassvdHScholbach:MotivicSatake} for $\nabla^{\exp}_{\tilde{w}} (\Z)$ implies the analogous map for $\nabla^{\exp}_{w} (\Z)$ is an isomorphism. The rest of the necessary statements follow from \thref{piAnti} applied to $\nabla^{\exp}_{\tilde{w}} (\Z)$.
\xpf

\subsection{Tate motives on affine flag varieties} 
In this section we define various categories of Tate motives on affine flag varieties, and convolution functors between them, beginning with properties which do not rely on the existence of a t-structure.

\subsubsection{Tate motives without a t-structure}
Let $\Fl_{\mathbf{f}}^{\exp}$ be the disjoint union of the $\calU^{\exp}$-orbits in $\Fl_{\mathbf{f}}$. By \thref{exponential orbits WT}, this is a Whitney--Tate stratification, and we have the full subcategory $\DTM(\Fl_{\mathbf{f}}, \Fl_{\mathbf{f}}^{\exp}) \subset \DM(\Fl_{\mathbf{f}})$ of stratified Tate motives. We also have the anti-effective and reduced variants, denoted $\DTM^{\anti}(\Fl_{\mathbf{f}}, \Fl_{\mathbf{f}}^{\exp})$ and  $\DTM_{\red}(\Fl_{\mathbf{f}}, \Fl_{\mathbf{f}}^{\exp})$, respectively, or $\DTM_{\redx}^{\xanti}(\Fl_{\mathbf{f}}, \Fl_{\mathbf{f}}^{\exp})$ to allow any combination of these (see \cite[\S 2.1.3, \S 2.2]{CassvdHScholbach:MotivicSatake} for more details). Unless otherwise stated, everything we prove from now on also holds verbatim for reduced motives, but we will usually not mention this.

\defi \thlabel{non-exp-Tate}
We define $\DTM(\mathcal{U}^{\exp}\backslash\Fl_{\mathbf{f}}) \subset \DM(\mathcal{U}^{\exp}\backslash\Fl_{\mathbf{f}})$ as the full subcategory consisting of equivariant motives whose underlying motive lies in $\DTM(\Fl_{\mathbf{f}}, \Fl_{\mathbf{f}}^{\exp})$.
\xdefi

We refer to \cite[\S 2.5]{CassvdHScholbach:MotivicSatake} for further details on equivariant motives. The category $\DTM(\mathcal{U}^{\exp}\backslash\Fl_{\mathbf{f}})$ admits anti-effective and reduced variants. We also have the subcategory $\DTM(\mathcal{U}^{\exp}\backslash\Fl_{\mathbf{f}})^{\text{lc}} \subset \DTM(\mathcal{U}^{\exp}\backslash\Fl_{\mathbf{f}})$ of locally compact motives \cite[Definition 5.52]{CassvdHScholbach:MotivicSatake}, consisting of motives whose image in $\DTM(\Fl_{\mathbf{f}}, \Fl_{\mathbf{f}}^{\exp})$ is compact. Equivalently, an object in $\DTM(\mathcal{U}^{\exp}\backslash\Fl_{\mathbf{f}})$ is locally compact if and only if it is supported on finitely many strata, and its $*$- (or $!$-) pullback to each stratum is compact after forgetting the equivariance.

For any standard facets $\mathbf{f}, \mathbf{f}'$ we have the category $\DTM(L^+G_{\mathbf{f}'} \backslash \Fl_{\mathbf{f}})$ thanks to \cite[Proposition 3.7]{CassvdHScholbach:MotivicSatake}. The (co)standard functors (for both $\calI$- and $\calU^{\exp}$-orbits) from \refsect{WTExp} preserve Tateness, and from now on we will only apply them to Tate motives.

We note that we always have a K\"unneth formula for $*$-pullback and $!$-pushforward \cite[Lemma 2.2.3]{JinYang:Kuenneth}. For $*$-pushforward, the following lemma is the exponential analogue of \cite[Corollary 3.10]{CassvdHScholbach:MotivicSatake}.

\lemm
\thlabel{WT.prod}
Let $\iota_X\colon  \coprod X_v \r X$ be any stratified Whitney--Tate (ind-)scheme $X$, and let $\iota_{\mathbf{f}} \colon \coprod \Fl^{\exp}_w(\mathbf{f}) \r \Fl_{\mathbf{f}}$ be the exponential stratification.
Then the product stratification on $X \x \Fl_{\mathbf{f}}$ is again Whitney--Tate. In particular, the inclusions $\iota_X$ and $\iota_{\mathbf{f}}$ satisfy the K\"unneth formula
$$(\iota_{X} \times \iota_{\mathbf{f}})_*\Z = \iota_{X *}\Z \boxtimes \iota_{\mathbf{f} *}\Z.$$
\xlemm

\pf
This follows exactly as in \cite[Corollary 3.10]{CassvdHScholbach:MotivicSatake}, using the universality of the stratification $\iota_{\mathbf{f}}$. 
\xpf

\rema \thlabel{rema--Compact} 
By revisiting the fiber sequences  in the proof of \thref{exponential orbits WT}, one can check that the functors $$\Delta_w^{\exp}, \nabla_w^{\exp} \colon \DTM(\mathcal{U}^{\exp}\backslash S) \to \DTM(\mathcal{U}^{\exp}\backslash \Fl_{\mathbf{f}})$$ preserve locally compact objects. Alternatively, this property is automatic for $\Delta_w^{\exp}$ as $!$-pushforward is left adjoint to the colimit-preserving functor of $!$-pullback. The case of $\nabla_w^{\exp}$ then follows by Verdier duality, which is applicable since the former induces an anti-equivalence of $\DTM(S)^{\comp}$. Similar remarks apply to the Iwahori (co)standard functors, reduced motives, and $*$- and $!$-pullback along equivariant inclusions of strata in $\Fl_{\mathbf{f}}$.
\xrema

Let \(\bbf,\bbf'\) be two standard facets, and consider the convolution diagram
\[LG /L^+G_{\bbf'} \times L^+G_{\bbf'} \backslash LG/L^+G_{\bbf} \xleftarrow{p} LG \overset{L^+G_{\bbf'}}{\times} LG/L^+G_{\bbf} \xrightarrow{m} LG/L^+G_{\bbf},\]
where $p$ and $m$ are the projection and multiplication maps, respectively. Since $
\DM$ is invariant under Nisnevich sheafification, we may define the convolution
 product
\[m_!p^!(-\boxtimes -) \colon \DM(\Fl_{\bbf'}) \times \DM (L^+G_{\bbf'} \backslash \Fl_{\bbf}) \to \DM(\Fl_{\bbf}).\]
Here $p^!$ is the functor which forgets equivariance. The map $p$ is equivariant for the natural left action of $\calU^{\exp}$ (the action on the target is only on the left factor $\Fl_{\bbf'}$), and $m$ is $\calU^{\exp}$-equivariant by \thref{triviality of convolution}, so that we may view convolution as a functor 
\[(-)\star(-) \colon \DM(\calU^{\exp} \backslash \Fl_{\bbf'}) \times \DM (L^+G_{\bbf'} \backslash \Fl_{\bbf}) \to \DM(\calU^{\exp} \backslash \Fl_{\bbf}).\]

By \cite[Propositions 4.17, 4.18]{CassvdHScholbach:Central}, $\DM(L^+G_{\mathbf{f}} \backslash \Fl_{\mathbf{f}}) = \DM(L^+G_{\mathbf{f}} \backslash LG/ L^+G_{\mathbf{f}})$ is a monoidal $\infty$-category under convolution. 

\lemm
The convolution product $\star$ endows $ \DM(\calU^{\exp} \backslash \Fl_{\bbf})$ with a module structure for the monoidal $\infty$-category $\DM(L^+G_{\mathbf{f}} \backslash LG/ L^+G_{\mathbf{f}})$.
\xlemm

\pf
By the above remarks, $\calU^{\exp} \setminus LG / L^+G_{\bbf}$ is a right module over the Hecke prestack $L^+G_\bbf \setminus LG / L^+G_\bbf$, which is an algebra object in the category denoted by $\Corr(\IndSch^{\mathrm{pl}}_S(\mathrm{sift}))$ in \cite[§2.2.8]{CassvdHScholbach:Central}. The claim results by applying the lax symmetric monoidal functor $\DM^!$ on this category of correspondences.
\xpf

\subsubsection{Tate motives with a t-structure}
From now on, unless we are working with reduced motives, we assume that the base scheme $S$ satisfies the Beilinson--Soul\'e vanishing condition, as discussed in \cite[§2.3]{CassvdHScholbach:MotivicSatake}. Then there is a motivic t-structure on $\DTM(S)$ whose heart $\MTM(S)$ is compactly generated by the motives $\Z(k)$ for $k \in \Z$. 

Recall that by \cite[Lemma 2.26, Proposition 3.4]{CassvdHScholbach:MotivicSatake}, the t-structures on the strata $\DTM(L^+G_{\mathbf{f}'} \backslash \Fl^{\circ}_w(\mathbf{f}',\mathbf{f}))$ for $w \in W_{\mathbf{f}'} \backslash W / W_{\mathbf{f}}$, normalized so that $\Z[\dim \Fl^{\circ}_w({\mathbf{f}}', \mathbf{f})] \in \MTM(L^+G_{\mathbf{f}'} \backslash \Fl^{\circ}_w(\mathbf{f}',\mathbf{f}))$, glue to a t-structure on 
$\DTM(L^+G_{\mathbf{f}'} \backslash \Fl_{\mathbf{f}})$, and the same is also true for non-equivariant Tate motives. The following lemma says that the same statement holds for the stratification by $\calU^{\exp}$-orbits.

\lemm
The ind-scheme $\Fl_{\mathbf{f}}$, equipped with the action of $\calU^{\exp}$ and the stratification by $\calU^{\exp}$-orbits in \thref{exponential orbits WT}, satisfies all the assumptions of \cite[Lemma 2.26]{CassvdHScholbach:MotivicSatake}. In particular, we have the abelian subcategory
$$\MTM(\calU^{\exp} \backslash \Fl_{\mathbf{f}}) \subset \DTM(\calU^{\exp} \backslash \Fl_{\mathbf{f}})$$
of mixed Tate motives. The same is true for non-equivariant motives, reduced motives, and anti-effective motives.
\xlemm

\pf
We apply \cite[Lemma 2.26]{CassvdHScholbach:MotivicSatake} in the case where $X = \Fl_{\mathbf{f}}$ and $Y = S$. The admissibility condition in loc.~cit. was introduced to deal with certain non-cellular stratifications, and it is satisfied in the present situation since all orbits are products of \(\A^1\)'s and \(\Gm\)'s by \thref{orbit-structure} (see \cite[Remark 2.10]{CassvdHScholbach:MotivicSatake}). The remaining hypotheses all deal with the cellular nature of $\calU^{\exp}$ and the stabilizers of points in $\Fl_{\mathbf{f}}$, which are also dealt with in \thref{orbit-structure}. The non-equivariant case is \cite[Lemma 2.15]{CassvdHScholbach:MotivicSatake}, and the anti-effective case follows from the anti-effectivity of the stratification (see \cite[\S 2.3]{CassvdHScholbach:MotivicSatake}). 
\xpf

At the level of $\calU^{\exp}$-orbits, we have the following result.

\lemm \thlabel{lemm--DMorbit}
For $w \in W_{\mathbf{f}}^{\exp}$, let $e \colon S / \calU_{w}^{\exp} \rightarrow  \calU^{\exp} \backslash \Fl_{w}^{\exp}(\mathbf{f})$ be the induced map of pre-stacks (this is the inclusion of the basepoint modulo stabilizer into the orbit). Then there is a canonical equivalence
$$\DM( \calU^{\exp} \backslash \Fl_{w}^{\exp}(\mathbf{f}) ) \xrightarrow[\sim]{e^![\dim \Fl_{w}^{\exp}(\mathbf{f})]} \DM(S/\calU_{w}^{\exp}).$$ Furthermore, the following hold.
\begin{enumerate}
\item The above equivalence restricts to $\DTM^{\xanti}$, and also holds for the finite-type quotients $\calU^{\exp}_{i}$, $\calU^{\exp}_{w,i}$, for $i \gg 0$. 
\item Any $\calU^{\exp}$-equivariant Tate motive in $\DTM(\Fl_{w}^\circ(\mathbf{f}))$ is anti-effective if and only if its *-pullback to the basepoint $S \subset \Fl_w^\circ(\mathbf f)$ is anti-effective.
\item The above equivalence is t-exact on $\DTM$, and additionally $$\MTM(S/\calU^{\exp}) = \MTM(S/\calU_{w}^{\exp}) = \MTM(S).$$
\item We also have equivalences
$$
\DTM(S/\calU_{w}^{\exp}) = \begin{cases}
    \DTM(S/ \Gm), & w \in W/W_{\mathbf{f}} \\
    \DTM(S), & w \in {}_0W_{\mathbf{f}}.
    \end{cases}
$$
\end{enumerate}
\xlemm

\pf
Parts (1)-(3) are the analogue of \cite[Lemma 4.28]{CassvdHScholbach:Central} for $\calU^{\exp}$ rather than parahorics (see also \cite[Proposition 3.4]{CassvdHScholbach:MotivicSatake}). The proof is the same, making use of \thref{orbit-structure} for all the necessary conditions. Part (4) follows from the triviality of all equivariance data for trivial actions of split pro-unipotent groups \cite[Proposition 2.2.11]{RicharzScholbach:Intersection}, since $\Gm \subset \calU_w^{\exp}$ if and only if $w \in W/W_{\mathbf{f}}$.
\xpf

\theo\thlabel{convDTM}
Let \(\bbf,\bbf'\) be two standard facets.
Then the convolution product preserves stratified Tate motives, i.e., it restricts to a functor 
\[(-)\star(-) \colon \DTM(\calU^{\exp} \backslash \Fl_{\bbf'}) \times \DTM (L^+G_{\bbf'} \backslash \Fl_{\bbf}) \to \DTM(\calU^{\exp} \backslash \Fl_{\bbf}).\]
Furthermore, convolution preserves anti-effective motives and locally compact motives.
\xtheo

\pf We start with the preservation of Tateness.
As in the parahoric case \cite[Theorem 3.17]{RicharzScholbach:Motivic} (see also \cite[Definition and Lemma 4.11]{CassvdHScholbach:MotivicSatake}), by  \thref{MTM equivariant} and \thref{lemm--DMorbit}(3) it is enough to consider a convolution product of the form $M:=\Delta^{\exp}_w(\Z) \star \Delta_v(\Z)$, for $w \in W_{\mathbf{f}'}^{\exp}$ and $v \in W_{\mathbf{f}'} \backslash W / W_{\mathbf{f}}$. 
For $w' \in W_{\mathbf{f}}^{\exp}$, let $M_{w'} \in \DM(S)$ be the $*$-pullback to the basepoint of the stratum $\Fl_{w'}^{\exp}(\mathbf{f})$ (as in \thref{def-KirillovOrbits}). 

Since Tateness may be detected by $*$-pullback to strata \cite[Definition and Lemma 3.1.11]{RicharzScholbach:Intersection}, then by \thref{triviality of convolution} (or \thref{lemm--DMorbit}) it suffices to show that $M_{w'}$ is Tate. By base change and ind-properness of \(m\), the motive $M_{w'}$ is the $!$-pushforward of $\Z$ along a fiber of the convolution morphism $\Fl_{w}^{\exp}(\mathbf{f}') \widetilde{\times} \Fl_{v}(\mathbf{f}', \mathbf{f}) \r \Fl_{\mathbf{f}}$. Any such fiber admits a cellular filtrable decomposition by \thref{fibers of exponential convolution}, which has Tate cohomology  by inductively applying the localization sequence: the argument is identical to the first part of \cite[Lemma 2.20]{CassvdHScholbach:MotivicSatake}, allowing for $\Gm \setminus \{1\}$ in addition to $\Gm$ and $\A^1$. In fact, this shows that $M_{w'} \in \DTM^{\anti}(S)$, so that $\star$ preserves anti-effective Tate motives by \thref{lemm--DMorbit}(2). Finally, the preservation of locally compact Tate motives follows by the same argument as in \cite[Lemma 4.18]{CassvdHScholbach:Central}. 
\xpf

\subsection{t-exactness properties of (co)standards and convolution}

By \cite[Proposition 5.26]{CassvdHScholbach:Central}, the (co)standard functors 
$$\Delta_w, \nabla_w \colon \DTM(\mathcal{I}\backslash S) \to \DTM(\mathcal{I\backslash \Fl_{\mathbf{f}}}), \quad w \in W/W_{\mathbf{f}}$$
are t-exact. The proof of this fact intersects in several places with the left (resp.~right) t-exactness of convolution by $\Delta_w(\Z)$ (resp.~$\nabla_w(\Z)$) on $\DTM(\calI \backslash \Fl_{\mathbf{f}})$, for $w \in W$ \cite[Lemma 5.27]{CassvdHScholbach:Central}. The difficulty in proving all of these results is the lack of a motivic Artin vanishing theorem. In this section we consider the analogous problems for exponential stratified motives. 

Since t-exactness properties can be checked after forgetting the equivariance (see \cite[Lemma 2.15, Lemma 2.26]{CassvdHScholbach:MotivicSatake}), we will often ignore equivariance data on co(standard) objects. Additionally, the remark below explains why we can often ignore equivariance data in the domain as well.

\rema \thlabel{rema--Exact}
Similarly to \thref{lemm--DMorbit}(4), we have equivalences
$$\DTM(\calI \backslash \Fl_w^{\circ}(\mathbf{a}_0, \mathbf{f})) = \DTM(\calI \backslash S) = \DTM(T \backslash S), \quad w \in W/W_{\mathbf{f}}.$$ Thus, any of these can be viewed as the domain for the Iwahori co(standard) functors. 
On the other hand, $\DTM(\calU^{\exp} \backslash \Fl_w^{\exp}(\mathbf{f}))$ is isomorphic to either $\DTM(S/\calU^{\exp}) = \DTM(S / \Gm)$ or $\DTM(S)$. This will not be a problem, but we record here that the co(standard) functors are insensitive to the equivariance datum if $w \in {}_0W_{\mathbf{f}}$. 

We shall also make extensive use of \thref{lemm--DMorbit}(3) which is true regardless of whether $w \in {}_0W_{\mathbf{f}}$, and its parahoric analogue \cite[Lemma 4.28]{CassvdHScholbach:Central}. In particular, to prove that a functor with domain $\DTM(\calU^{\exp} \backslash S)$ or $\DTM(L^+G_{\mathbf{f}} \backslash S)$ is left t-exact, it suffices to show this on elements of $\MTM(S)$, since the t-structures are accessible and right-complete (see the proof of \cite[Proposition 5.10]{CassvdHScholbach:MotivicSatake} for more details). Furthermore, to show that a functor with domain $\DTM(\calU^{\exp} \backslash \Fl_{\mathbf{f}})$ (resp.~$\DTM(L^+G_{\mathbf{f}'} \backslash \Fl_{\mathbf{f}}$)) is left t-exact, by excision and the previous comments, it suffices to consider elements of the form $\nabla_w^{\exp}(M)$ (resp.~the $*$-extension of $M[\dim \Fl_w^\circ(\mathbf{f}', \mathbf{f})]$ along $\Fl_w^\circ(\mathbf{f}', \mathbf{f}) \subset \Fl_{\mathbf{f}}$) for $M \in \MTM(S)$.

For right t-exactness the situation is even simpler by \thref{MTM equivariant}, as we only have to consider Tate twists of the $!$-extension of $\Z$ along orbits (shifted by the dimension). Also by \thref{MTM equivariant}, $\DTM(\calU^{\exp} \backslash \Fl_{\mathbf{f}})$ and $\DTM(L^+G_{\mathbf{f}'} \backslash \Fl_{\mathbf{f}})$ are generated by shifts and twists of either $!$- or $*$-extensions of $\Z$ along orbits, so that when proving a functor preserves Tate motives we only have to consider such objects (as was already noted in the proof of \thref{convDTM}).
\xrema

\subsubsection{The full affine flag variety.} We start with the case $\mathbf{f}=\mathbf{a}_0$.

\prop \thlabel{expFLstandards}
If $\mathbf{f}= \mathbf{a}_0$ and $w \in W^{\exp}_{\mathbf{a}_0}$, the functors $$\Delta_w^{\exp}, \nabla_w^{\exp} \colon \DTM(\mathcal{U}^{\exp}\backslash S) \to \DTM(\mathcal{U}^{\exp}\backslash \Fl)$$ are t-exact.
\xprop

\pf We consider the functors separately. 

{\em The functor $\nabla_w^{\exp}$}: 
This is clearly left t-exact, so it suffices to prove that $\nabla_w^{\exp}(\Z) \in \DTM(\Fl)^{\leq 0}$.
By \thref{orbit-structure} and the Iwahori case \cite[Proposition 5.2]{CassvdHScholbach:Central}, we can assume that either $w \in W$ and $w$ has maximal length in its left $W_0$-orbit, or $w \in {}_0W_{\mathbf{a}_0}$. If $w \in W$ has maximal length and $w' \in {}_0W_{\mathbf{a}_0}$ is the same element regarded as as a member of ${}_0W_{\mathbf{a}_0}$, then by excision and relative purity we have a fiber sequence
$$\nabla_w^{\exp} (\Z)(-1)[-1] \r \nabla_w(\Z) \r \nabla_{w'}^{\exp}(\Z).$$
The middle term lies in degree $0$, so it suffices to consider the case $w \in W$ has maximal length. Write $w=w_0 w_{\min}$. We proceed by induction on $l(w_{\min})$.

If $l(w_{\min}) = 0$, then by \thref{Omega-translate} we can assume that $w_{\min}$ is trivial. Then we are dealing with $\nabla_{w_0}^{\exp}(\Z)$, which is supported on $G/B$. In this case, we prove by descending induction on $l(v)$, for $v \in W_0$ and $v \neq w_0$, that $\iota_v^* \nabla_{w_0}^{\exp}(\Z) \in \DTM(\Fl_{v}^\circ(\mathbf{a}_0))^{\leq 0}$. If $l(v) = l(w_0) -1$ then $\iota_v^* \nabla_{w_0}^{\exp}(\Z) = 0$ by \thref{theo--strat}. Now if $l(v) \leq l(w_0) -2$, we apply the fiber sequence \eqref{WTfiberG/B}. The result is that if $\alpha \in \Delta$ is such that $l(v s_{\alpha}) = l(v) + 1$, there is a fiber sequence
\begin{equation} \label{seq1} \iota_{v s_{\alpha}}^*(\nabla_{w_0}^{\exp}(\Z))(-1)[-2] \r \iota_{v}'^*\nabla^{\mathbf{f}}_{w_0 s_{\alpha}}(\Z) \r \iota_{v}^* \nabla_{w_0}^{\exp}(\Z) \end{equation} in $\DTM(\Fl_v^\circ(\mathbf{a}_0))$. Here $\mathbf{f}$ is the wall of $\mathbf{a}_0$ associated to $s_{\alpha}$, $\nabla^{\mathbf{f}}_{w_0 s_{\alpha}}(\Z)$ is the Iwahori costandard object associated to $\Z$ and $w_0 s_{\alpha} \in W/W_{\mathbf{f}}$ on $G/P_{s_{\alpha}}$, and $\iota_v' \colon \Fl_v^\circ(\mathbf{a}_0) \cong \Fl_v^\circ(\mathbf{a}_0,\mathbf{f}) \r G/P_{s_{\alpha}}$ is the inclusion. The term on the left is associated with an element of  $\DTM(\Fl_v^\circ(\mathbf{a}_0))$ by the canonical equivalences $\DTM(S) \cong \DTM(\Fl_v^\circ(\mathbf{a}_0)) \cong \DTM(\Fl_{v s_{\alpha}}^{\circ}(\mathbf{a}_0))$ induced by $*$-pullback and homotopy invariance. By descending induction on $l(v)$, the left term lies in degrees $\leq 1$. The middle term lies in degrees $\leq 0$ by \cite[Proposition 5.26]{CassvdHScholbach:Central}, so that the right term lies in degrees $\leq 0$, as needed.

If $l(w_{\min}) > 0$, we adapt the proof of \cite[Proposition 5.2]{CassvdHScholbach:Central}. Let $s \in \mathcal{S}$ be such that $l(w_{\min}s) < l (w_{\min})$. Let $\mathbf{f}$ be the facet such that $W_{\mathbf{f}} = \{1, s\}$ and let $\pi \colon \Fl \r \Fl_{\mathbf{f}}$ be the $\P$-fibration. The fiber sequence \eqref{FlWTseq} gives 
$$\nabla_w^{\exp}(\Z) \r \nabla_{ws}^{\exp}(\Z) (-1) \r \pi^! \pi_! \nabla_{ws}^{\exp} (\Z)(-1).$$
By induction $\nabla_{ws}^{\exp}(\Z) (-1) \in \MTM(\mathcal{U}^{\exp}\backslash \Fl)$, so it 
suffices to prove that the endofunctor $\fib (\id \r \pi^! \pi_!)$ on $\DTM(\mathcal{U}^{\exp}\backslash \Fl)$ is right t-exact. For this, let $v \in W_{\mathbf{f}}^{\exp}$ be arbitrary and consider the pullback diagram
\[\begin{tikzcd}
	\pi^{-1}(\Fl_{\mathbf{f}}^{\exp}(v)) \arrow[d, "p"] \arrow[r, "k"] & \Fl \arrow[d, "\pi"] \\
	\Fl_{\mathbf{f}}^{\exp}(v)  \arrow[r, "i"] & \Fl_{\mathbf{f}}.
\end{tikzcd}\]
By base change, $k^* \fib (\id \r \pi^! \pi_!) = \fib (k^* \r p^! p_! k^*)$. All of these maps are $\calU^{\exp}$-equivariant, so we are reduced to proving that $\fib (\id \r p^! p_!)$ is right t-exact as an endofunctor on $\DTM(\mathcal{U}^{\exp}\backslash \pi^{-1}(\Fl_{\mathbf{f}}^{\exp}(v)))$.

Over each $\calU$-orbit in $\Fl_{\mathbf{f}}$ the map $\pi$ is a relative $\P$, so the same is true of the map $p$. Let $\Fl_v^{\circ}(\mathbf{a}_0, \mathbf{f})$ be the $\calU$-orbit containing $\Fl_{\mathbf{f}}^{\exp}(v)$. Then $\pi^{-1}(\Fl_v^{\circ}(\mathbf{a}_0, \mathbf{f}))$ decomposes as the closed subscheme $\Fl_v^{\circ}(\mathbf{a}_0)$ mapping isomorphically onto $\Fl_v^{\circ}(\mathbf{a}_0, \mathbf{f})$ and the open complement  $\Fl_{vs}^{\circ}(\mathbf{a}_0)$ which is a relative $\A^1$ over $\pi^{-1}(\Fl_v^{\circ}(\mathbf{a}_0, \mathbf{f}))$. The next step is to understand the possible stratifications of $\pi^{-1}(\Fl_{\mathbf{f}}^{\exp}(v))$. The intersection $\pi^{-1}(\Fl_{\mathbf{f}}^{\exp}(v)) \cap \Fl_v^{\circ}(\mathbf{a}_0, \mathbf{f})$ is always identified with the section at $\{\infty\}$ in the relative $\P$ map $p$, and it remains to understand $\pi^{-1}(\Fl_{\mathbf{f}}^{\exp}(v)) \cap \Fl_{vs}^{\circ}(\mathbf{a}_0, \mathbf{f})$. There are several cases to consider.

First, suppose that $v \in W/W_{\mathbf{f}}$ and that $v$ is not left maximal. Then $\Fl^{\exp}_{\mathbf{f}}(v) =  \Fl_v^{\circ}(\mathbf{a}_0, \mathbf{f})$ is a $\calU$-orbit. There are two further subscases.
If $vs \in W$ is also not left maximal then $\Fl^{\exp}_{vs}(\mathbf{a}_0) = \Fl_{vs}^{\circ}(\mathbf{a}_0, \mathbf{f})$, so that at the level of strata $p$ is identified with the map $\P = \A^1 \sqcup \{\infty\} \r S$, where we have replaced our original base with $S = \Fl_{\mathbf{f}}^{\exp}(v)$. If $vs \in W$ is left maximal then $\Fl_{vs}^{\circ}(\mathbf{a}_0)$ decomposes as a union of two $\calU^{\exp}$-orbits, and by \thref{Uexp-projection}, $p$ is identified with the map $\P = \Gm \sqcup \{0\} \sqcup \{\infty\} \r S$.

Next, suppose that $v \in W/W_{\mathbf{f}}$ and that $v$ is left maximal. Then $vs$ is also left maximal. Indeed, if $v = w_0 v_{\min}$ with $v_{\min} \in W$ and $l(v) = l(w_0) + l(v_{\min})$, then $l(vs) = l(v) + 1 = l(w_0) + l(v_{\min}) +1$, so that $l(w_0 v_{\min}s) = l(w_0) + l(v_{\min} s)$. Now by \thref{Uexp-projection}, $\pi^{-1}(\Fl_{\mathbf{f}}^{\exp}(v))$ meets only the smaller of the two $\calU^{\exp}$-orbits in $\Fl_{vs}^{\circ}(\mathbf{a}_0, \mathbf{f})$, and $p$ is identified with the map $\P = \A^1 \sqcup \{\infty\} \r S$.

Finally, suppose $v \in {}_0 W_{\mathbf{f}}$. Then arguing as in the previous paragraph, $\pi^{-1}(\Fl_{\mathbf{f}}^{\exp}(v))$ meets only the larger of the two $\calU^{\exp}$-orbits in $\Fl_{vs}^{\circ}(\mathbf{a}_0, \mathbf{f})$, and $p$ is again identified with the map $\P = \A^1 \sqcup \{\infty\} \r S$.

Thus, the most general case we need to consider is the endofunctor $\fib (\id \r p^! p_!)$ for the stratified map $\P = \Gm \sqcup \{0\} \sqcup \{\infty\} \r S$, which is handled by the first part of \thref{fiblem} below. This completes the proof that $\nabla_w^{\exp}$ is t-exact.

{\em The functor $\Delta_w^{\exp}$}: The proof is similar, so we will be brief. By \thref{rema--Exact}, we need to prove that if $M \in \MTM(S)$, then $\Delta_w^{\exp}(M) \in \DTM(\Fl_{\mathbf{f}})^{\geq 0}$. By the Iwahori case \cite[Proposition 5.2]{CassvdHScholbach:Central}, we can assume that either $w \in W$ and $w$ has maximal length in its left $W_0$-orbit, or $w \in {}_0W_{\mathbf{a}_0}$. In the former case, there is a fiber sequence $\Delta_w^{\exp}(M) \r \Delta_{w'}^{\exp}(M) \r \Delta_w(M)$ where $w' = w$ regarded as an element of ${}_0 W_{\mathbf{a}_0}$, so we reduce the case $w \in W$ is left maximal. Writing $w = w_0 w_{\min}$, we again induct on $l(w_{\min})$. The base case reduces to $\Delta_{w_0}^{\exp}(M)$ on $G/B$, where we argue as follows. If $v \in W_0$ and $l(v) \leq l(w_0) -2$, a dual argument to the construction of \eqref{WTfiberG/B} gives a fiber sequence
$$\iota_v^! \Delta_{w_0}^{\exp}(M) \r \iota'^! \Delta_{w_0 s_{\alpha}}^{\mathbf{f}}(M) \r \iota_{v s_{\alpha}}^! \Delta_{w_0}^{\exp}(M)$$ where $l(vs_{\alpha}) = l(v) + 1$. 
By descending induction on $l(v)$ and \cite[Proposition 5.26]{CassvdHScholbach:Central} applied to the middle term, we find that $\iota_v^! \Delta_{w_0}^{\exp}(M)$ lies in non-negative degrees, completing the base case. 

If $l(w_{\min}) > 0$, an argument dual to the construction of \eqref{FlWTseq} gives a fiber sequence
$$\pi^* \pi_* \Delta_{ws}^{\exp}(M) \r \Delta_{ws}^{\exp}(M) \r \Delta_w^{\exp}(M)$$ where $l(ws) = l(w) -1$. Using base change, we are reduced to proving the left t-exactness of $\text{cof}(p^* p_* \r \id)$ for $p \colon \P \r S$, where $\P$ is stratified as $\P = \Gm \sqcup \{0\} \sqcup \{\infty\} \r S$. We do not claim left t-exactness in this full generality, but since the fiber sequence above is $\calU^{\exp}$-equivariant, then by \thref{rema--Exact} we are reduced to the setup in the second part of \thref{fiblem} below.
\xpf

\lemm \thlabel{fiblem}
Let $p \colon \P \r S$ be the relative projective line, and give $\P$ the stratification $\P = \Gm \sqcup \{0\} \sqcup \{\infty\}$. Then the endofunctor $\fib (\id \r p^! p_!)$ of $\DTM(\P)$ is right t-exact. Additionally, let $i_0, i_\infty$, and $j$ be the inclusions of $\{0\}, \{\infty\},$ and $\Gm$ into $\P$. If $L \in \MTM(S)$ and $M$ is one of the motives $i_{0*} L$, $i_{\infty *} L$, or $j_* L[1]$, then $\text{cof}(p^*p_*M \r M) \in \DTM(\P)^{\geq 0}$.
\xlemm

\pf
We start with the claim about $\fib (\id \r p^! p_!)$. It suffices to let $M$ be one of the motives $i_{0!} \Z$, $i_{\infty !} \Z$, or $j_! \Z[1]$, and we must show that $\fib (M \r p^! p_! M) \in \DTM(\P)^{\leq 0}$. If $M = i_{0!} \Z$ then $p^! p_! M = \Z(1)[2] \in \DTM(\P)$ lies in degree $-1$. Since $M$ lies in degree $0$ then $\fib (M \r p^! p_! M) \in \DTM(\P)^{\leq 0}$. The case $M = i_{\infty!} \Z$ is identical. Finally, if $M = j_! \Z[1]$ then $M \r p^! p_! M$ is a map $j_! \Z[1] \r \Z(1)[2] \oplus \Z[1]$. Since $\Z(1)[2]$ lies in degree $-1$, it suffices to check that $\fib(j_!\Z[1] \r \Z[1])$ lies in non-positive degrees. However, the map $j_!\Z[1] \r \Z[1]$ restricts to the identity map over $\Gm$ so that the latter fiber is just $i_{0!} \Z \oplus i_{\infty!} \Z$, as needed.

For the claim about  $\text{cof}(p^*p_*(M) \r M)$, first suppose that $M = i_{0*}L$. Then $p^*p_*M \r M$ is a map $p^*L \r i_{0*} L$. By relative purity, $p^*L$ lies in degree $1$, so that the cofiber lies in non-negative degrees. The case $M = i_{\infty*} L$ is identical. Finally, if $M = j_* L[1]$ then $p^*p_*M = L[1] \oplus L(-1)$ (where we have used the shorthand $L=p^*L$), and the map $L[1] \r j_* L[1]$ is characterized by the property that its restriction to $\Gm$ is the identity map. This is just an excision map for $L[1]$ (and this is where we use that $j^*M = L[1]$ is pulled back from $S$, as opposed to being an arbitrary object of $\MTM(\Gm)$), so that by relative purity the cofiber $i_{0*} L(-1) \oplus i_{\infty*} L(-1)$ of the latter map lies in degree $0$.
\xpf

\prop \thlabel{Flconvtexact}
For any $w \in W^{\exp}_{\mathbf{a}_0}$, convolution  by $\Delta_w^{\exp}(
\Z)$ (resp.~$\nabla_w^{\exp}(\Z)$) defines a left (resp.~right) $t$-exact functor $\DTM(
\calI \backslash \Fl) \r \DTM(\calU^{\exp} \backslash \Fl)$.
\xprop

\pf {\em Convolution by $\Delta_w^{\exp}(\Z)$}:
 By \thref{rema--Exact}, we need to show that if $v \in W$ and $M \in \MTM(S)$, then $\Delta_w^{\exp}(\Z) \star \nabla_v(M) \in \DTM(\calU^{\exp} \backslash \Fl)^{\geq 0}$. If the $
\calU^{\exp}$-orbit of $w$ agrees with the $\calU$-orbit, then $\Delta_w^{\exp}(\Z) = \Delta_w(\Z)$ and the claim follows from the Iwahori case \cite[Lemma 5.4]{CassvdHScholbach:Central}. It remains to consider those $\calU^{\exp}$-orbits which are a proper subscheme of a $\calU$-orbit. Suppose that $w \in W$ is maximal on the left, and let $w' \in {}_0W_{\mathbf{a}_0}$ be the same element. Then from the fiber sequence $\Delta_w^{\exp}(\Z) \r \Delta_{w'}^{\exp}(\Z) \r \Delta_w(\Z)$ and the Iwahori case, this is the only type of $w$ we need to consider. If we write $w = w_0 w_{\min}$ then $\Delta_{w}^{\exp}(\Z) \cong \Delta_{w_0}^{\exp}(\Z) \star \Delta_{w_{\min}}(\Z)$ by \thref{convlem} (and the K\"unneth formula in \thref{WT.prod}), so that we can further assume $w = w_0$.

We induct on the difference $l(v) - l(v_{\min})$. The convolution $\Delta_{w_0}^{\exp}(\Z) \star \nabla_v(M)$ is isomorphic to the $*$-pushforward of $\Delta_{w_0}^{\exp}(\Z) \widetilde{\boxtimes} M[l(v)]$ along the restricted convolution map $\overline{\Fl}_{w_0}^{\exp}(\mathbf{a}_0) \widetilde{\times} \Fl^{\circ}_v(\mathbf{a}_0) \r \Fl$. Now assuming $v = v_{\min}$, this restricted convolution map is an isomorphism onto its image since $l(uv_{\min}) = l(u) + l(v_{\min})$ for all $u \in W_0$. By base change, our task is to show that the $!$-pullback of $\Delta_{w_0}^{\exp}(\Z) \widetilde{\boxtimes} M[l(v_{\min})]$ to each $
\Fl_u^{\circ}(\mathbf{a}_0) \widetilde{\times} \Fl^{\circ}_{v_{\min}}(\mathbf{a}_0)$ lies in nonnegative degrees, where $u \in W_0$, $u \neq w_0$.

Consider the following cartesian diagram, where the horizontal maps are the projections and the vertical maps are the inclusions.
\[\begin{tikzcd}
	\Fl_{w_0}^{\exp} \widetilde{\times} \Fl^{\circ}_{v_{\min}}(\mathbf{a}_0)  \arrow[r] \arrow[d, "j"] & \Fl_{w_0}^{\exp} \arrow[d] \\
	\overline{\Fl}_{w_0}^{\exp} \widetilde{\times} \Fl^{\circ}_{v_{\min}}(\mathbf{a}_0) \arrow[r, "p"] & \overline{\Fl}_{w_0}^{\exp}
\end{tikzcd}\]
By the K\"unneth formula for $!$-pushforward, $\Delta_{w_0}^{\exp}(\Z) \widetilde{\boxtimes} M[l({v_{\min}})] \cong j_!(M[l(w_0 v_{\min})]$. By base change, since $L$ is pulled back from $\MTM(S)$, $j_!(M[l(w_0 v_{\min})] \cong p^* \Delta_{w_0}^{\exp}(M) [l(v_{\min})]$. Since $p$ is a relative affine space, the necessary non-negativity of the $!$-pullbacks follows from \thref{expFLstandards} applied to $\Delta_{w_0}^{\exp}(M)$.

For the inductive step, suppose that $v$ is not left minimal. Let $\alpha \in \Delta$ be such that $l(s_{\alpha} v) < l(v)$, and set $s=s_{\alpha}$. Let $\mathbf{f}$ be the facet such that $W_{\mathbf{f}} = \{1, s\}$ and let $\pi \colon \Fl \r \Fl_{\mathbf{f}}$ be the projection. The $L^+G_{\mathbf{f}}$-orbit $\Fl_{v}^{\circ}(\mathbf{f}, \mathbf{a}_0)$ of $\dot{v} \cdot \overline{e}$ in $\Fl$ consists of $\Fl_v^{\circ}(\mathbf{a}_0)$ and $\Fl_{sv}^{\circ}(\mathbf{a}_0)$. Let $\nabla_{v}^{\mathbf{f}}(M) \in \DTM(L^+G_{\mathbf{f}} \backslash \Fl)$ be the $*$-extension of $M[l(v)]$ from $\Fl_v^{\circ}(\mathbf{f}, \mathbf{a}_0)$. This lies in nonnegative degrees, and by relative purity, there is a fiber sequence $\nabla_{v}^{\mathbf{f}}(M) \r \nabla_{v}(M) \r \nabla_{sv}(M)(-1)$. By induction, it suffices to show that $\Delta_{w_0}^{\exp}(\Z) \star \nabla_{v}^{\mathbf{f}}(M) \in \DTM(\calU^{\exp} \backslash \Fl)^{\geq 0}$. For this, note that the convolution map factors as
\begin{equation} \label{convFactor}
\begin{tikzcd}
	LG \times^{\calI} \Fl \arrow[rr] \arrow[rd, "\pi \widetilde{\times} \id"] & & \Fl. \\
    & LG \times^{L^+G_{\mathbf{f}}} \Fl \arrow[ru, "m"] &
\end{tikzcd}
\end{equation}
Since $\nabla_{v}^{\mathbf{f}}(M)$ is $L^+G_{\mathbf{f}}$-equivariant, it follows easily (cf.~\cite[Lemma 2.5]{BGMRR:IwahoriWhittaker}) that
$$\Delta_{w_0}^{\exp}(\Z) \star \nabla_{v}^{\mathbf{f}}(M) \cong m_! (\pi_!\Delta_{w_0}^{\exp}(\Z) \widetilde{\boxtimes} \nabla_{v}^{\mathbf{f}}(M)).$$ In other words, the left side only depends on $\pi_!\Delta_{w_0}^{\exp}(\Z)$. Since $\pi_!\Delta_{w_0}^{\exp}(\Z) \cong \pi_! \Delta_{w_0s}(\Z)$ by \thref{Uexp-projection} (both pushforwards are the $!$-extension of $\Z[l(w_0s)]$ supported on $\Fl_{w_0s}^{\circ}(\mathbf{a}_0, \mathbf{f})$), then by the same argument, 
$$\Delta_{w_0}^{\exp}(\Z) \star \nabla_{v}^{\mathbf{f}}(M) \cong \Delta_{w_0s}(\Z) \star \nabla_{v}^{\mathbf{f}}(M).$$ This lies in nonnegative degrees by \cite[Lemma 5.4]{CassvdHScholbach:Central}.

{\em Convolution by $\nabla_w^{\exp}(\Z)$}:
The proof is analogous, and it is slightly simpler since we need only show that $\nabla_w^{\exp}(\Z) \star \Delta_v(\Z) \in \DTM(\calU^{\exp} \backslash \Fl)^{\leq 0}$ for $v \in W$. In short, the fiber sequence $\nabla_{w}(\Z) \r \nabla_{w'}^{\exp}(\Z) \r \nabla_w^{\exp}(\Z)(-1)$ dual to the one in the first paragraph reduces us (with \cite[Lemma 5.4]{CassvdHScholbach:Central}) to the case $w \in W$ is left minimal. For such $w$, $\nabla_{w}^{\exp}(\Z) \cong \nabla_{w_0}^{\exp}(\Z) \star \nabla_{w_{\min}}^{\exp}(\Z)$, so we can assume $w = w_0$. If $v = v_{\min}$, then by writing $\nabla_{w_0}^{\exp}(\Z) \star \Delta_{v_{\min}}^{\exp}(\Z)$ as the $!$-pushforward of $\nabla_{w_0}^{\exp}(\Z) \widetilde{\boxtimes} \Z[l(v_{\min})]$ along $\overline{\Fl}_{w_0}^{\exp}(\mathbf{a}_0) \widetilde{\times} \Fl^{\circ}_{v_{\min}}(\mathbf{a}_0) \r \Fl$, a similar argument as above (and using the K\"unneth formula in \thref{WT.prod}) reduces the problem to checking that $\nabla_{w_0}^{\exp}(\Z)$ lies in nonpositive degrees, which follows from \thref{expFLstandards}. For the inductive step, if $s_{\alpha}v < v$ then similarly as above, we reduce to the case of convolution with the $L^+G_{\mathbf{f}}$-equivariant standard object $\Delta_{v}^{\mathbf{f}}(\Z)$ and ultimately appeal to \cite[Lemma 5.4]{CassvdHScholbach:Central}.
\xpf

\subsubsection{Partial affine flag varieties} We now let $\mathbf{f}$ be a standard facet.

\lemm \thlabel{lemma-cons+texact}
Let $\pi \colon \Fl \r \Fl_{\mathbf{f}}$ be the projection, and let $d = \dim L^+G_{\mathbf{f}}/\calI$. Then the functor $$\pi^*[d] = \pi^![-d](-d) \colon \DTM(\calU^{\exp} \backslash \Fl_{\mathbf{f}}) \r \DTM(\calU^{\exp} \backslash \Fl)$$ is conservative and t-exact.
\xlemm

\pf
The map $\pi$ is Zariski locally isomorphic to a projection with fibers $L^+G_{\mathbf{f}}/\calI$, so $\pi^*$ is conservative. It also follows that $\pi$ is smooth and proper, and $\pi^*[d] = \pi^![-d](-d)$. Moreover, $\pi$ maps strata onto strata by \thref{Uexp-projection}, so that $\pi^*$ preserves Tate motives. Right t-exactness then follows from base change, since this can be checked on the generators $\Delta_{w}^{\exp}(\Z)$ for $w \in W_{\mathbf{f}}^{\exp}$. For left t-exactness, by \thref{rema--Exact} we need to show that $\pi^*[d]\nabla_{w}^{\exp}(M) \in \DTM(\calU^{\exp} \backslash \Fl)^{\geq 0}$ if $M \in \MTM(S)$, which also follows from base change.
\xpf

\theo \thlabel{expFLfstandards}
For any standard facet $\mathbf{f}$ and $w \in W_{\mathbf{f}}^{\exp}$, the functors $$\Delta_w^{\exp}, \nabla_w^{\exp} \colon \DTM(\mathcal{U}^{\exp}\backslash S) \to \DTM(\mathcal{U}^{\exp}\backslash \Fl_{\mathbf{f}})$$ are t-exact.
\xtheo

\pf
Granting our work above, the proof is very similar to \cite[Proposition 4.7]{Achar:ModularPerverseSheavesII} and \cite[Proposition 5.26]{CassvdHScholbach:Central}. Briefly, if $w_{\mathbf{f}} \in W_{\mathbf{f}}$ is the longest element, then we have a cartesian diagram as follows.
\[\begin{tikzcd}
	\Fl\overset{\mathcal{I}}{\times}\Fl_{w_{\mathbf{f}}}(\mathbf{a}_0) \arrow[d] \arrow[r, "m"] & \Fl \arrow[d, "\pi"] \\
	\Fl \arrow[r, "\pi"] & \Fl_{\mathbf{f}}
\end{tikzcd}\]
Here $m$ is the multiplication map, and the left vertical map is the projection onto the first factor.
Let $\tilde{w} \in W_{\mathbf{a}_0}^{\exp}$ be the lift of $w$ as in \thref{Wexp-lift}, and
let $d = \dim L^+G_{\mathbf{f}}/\calI = \dim \Fl_{w_{\mathbf{f}}}(\mathbf{a}_0)$. By base change and the K\"unneth formula for $!$-pushforward, 
$$\pi^*[d]\Delta_{w}^{\exp}(M) \cong \pi^*[d]\pi_! \Delta_{\tilde{w}}^{\exp}(M) \cong \Delta_{\tilde{w}}^{\exp}(\Z) \star \pi^*[d]\Delta_{e}(M),$$ where $\Delta_{e}(M)$ is supported on the basepoint of $\Fl_{\mathbf{f}}$.

By \thref{lemma-cons+texact}, if $M \in \DTM(\mathcal{U}^{\exp}\backslash S)^{\geq 0}$, then $\pi^*[d]\Delta_{e}(M)\in \DTM(\mathcal{U}^{\exp}\backslash \Fl)^{\geq 0}$, and thus the above convolution product lies in $\DTM(\mathcal{U}^{\exp}\backslash \Fl)^{\geq 0}$ by \thref{Flconvtexact}. Hence $\pi^*[d]\Delta_{w}^{\exp}(M) \in \DTM(\mathcal{U}^{\exp}\backslash \Fl)^{\geq 0}$, which is enough to conclude that $\Delta_w^{\exp}$ is left t-exact (and hence t-exact) by \thref{lemma-cons+texact}. The proof that $\nabla_w^{\exp}$ is t-exact is analogous.
\xpf

\theo \thlabel{leftright}
For any standard facets $\mathbf{f}, \mathbf{f}'$, consider the convolution functor $$\DTM(\calU^{\exp} \backslash \Fl_{\mathbf{f}'}) \times \DTM(L^+G_{\mathbf{f'}} \backslash \Fl_{\mathbf{f}}) \xrightarrow{\star}  \DTM(\calU^{\exp} \backslash \Fl_{\mathbf{f}}).$$ 
Then for any $w \in W_{\mathbf{f}'}^{\exp}$, convolution by $\Delta_w^{\exp}(
\Z)$ (resp.~$\nabla_w^{\exp}(\Z)$) defines a left (resp.~right) $t$-exact functor $\DTM(
L^+G_{\mathbf{f}'} \backslash \Fl_{\mathbf{f}}) \r \DTM(\calU^{\exp} \backslash \Fl_{\mathbf{f}})$.
\xtheo

\pf
{\em The case $\mathbf{f}'=\mathbf{a}_0$}: 
Consider the following cartesian diagram, where the horizontal maps are convolution.
\[\begin{tikzcd}
	\Fl \times^{\calI} \Fl \arrow[r] \arrow[d, "\id \widetilde{\times} \pi"] & \Fl\arrow[d, "\pi"]\\
	\Fl \times^{\calI} \Fl_{\mathbf{f}} \arrow[r] & \Fl_{\mathbf{f}}
\end{tikzcd}\] Let $M \in \DTM(\calI \backslash S)$ and $v \in W/W_{\mathbf{f}}$.
By base change, $\pi^*[d](\Delta_w^{\exp}(\Z) \star \nabla_{v}(M)) \cong \Delta_w^{\exp}(\Z) \star \pi^*[d]\nabla_{v}(M)$, where the $\star$ on the right side refers to the convolution product along the top horizontal map. The same is true if ones switches $\Delta$ and $\nabla$. Now the result follows from \thref{lemma-cons+texact}, and the case $\mathbf{f}=\mathbf{a}_0$ in \thref{Flconvtexact}.

{\em The case of arbitrary $\mathbf{f}'$}:
Let $\tilde{w} \in W_{\mathbf{a}_0}^{\exp}$ be the lift of $w$ as in \thref{Wexp-lift}. Since $\pi_! \Delta_{\tilde{w}}^{\exp}(\Z) = \Delta_{w}^{\exp}(\Z)$ and likewise for $\nabla$, this follows from the case $\mathbf{f}'=\mathbf{a}_0$ and a factorization of the convolution map similar to \eqref{convFactor}.
\xpf

\subsection{$\Ga$-averageability} \label{sect--Ga}
Recall that the pro-unipotent radical of $\calI$ is $\calU = \calU^{>0} \rtimes U$, and that $\calU^{\exp} = (\calU^{>0} \rtimes U_0) \rtimes \Gm$.
In \refsect{setup}, we fixed a short exact sequence $0 \r U_0 \r U \r \Ga \r 0$. This leads to a short exact sequence
$$0 \r \calU^{>0} \rtimes U_0 \r \calU \rtimes \Gm  \r \Ga \rtimes \Gm \r 0.$$
For any standard facet $\mathbf{f}$, the prestack $(\calU^{>0} \rtimes U_0) \backslash \Fl_{\mathbf{f}}$ is equipped with an action of $\Ga \rtimes \Gm$, induced by the action of $\calU \rtimes \Gm$ on $\Fl_{\mathbf{f}}$. Moreover, ignoring the $\Ga$-action, we have (cf.~\eqref{quotient semidirect}) $$\Gm \backslash ((\calU^{>0} \rtimes U_0) \backslash \Fl_{\mathbf{f}}) = \calU^{\exp} \backslash \Fl_{\mathbf{f}}.$$

Consider the forgetful functor $\DM((\calU \rtimes \Gm) \backslash \Fl_{\mathbf{f}}) \r \DM(\calU^{\exp} \backslash \Fl_{\mathbf{f}})$ with its left and right adjoints $\av^{\Ga}_!$ and $\av^{\Ga}_*$. These may be computed as follows.

\lemm \thlabel{avg-comp}
Let $\beta \in \Delta$ be any simple root, and let $m_{\beta} \colon U_{\beta} \times \Fl_{\mathbf{f}} \r \Fl_{\mathbf{f}}$ be the restriction of the multiplication map $LG \times \Fl_{\mathbf{f}} \r \Fl_{\mathbf{f}}$ to $ U_{\beta}$. Let $u^! \colon \DM((\calU \rtimes \Gm) \backslash \Fl_{\mathbf{f}}) \r \DM(\Fl_{\mathbf{f}})$ be the forgetful functor. 
Then for any $M \in \DM(\calU^{\exp} \backslash \Fl_{\mathbf{f}})$, there are canonical isomorphisms
$$u^! \av^{\Ga}_! M = m_{\beta !} (\Z(1)[2] \boxtimes M), \quad u^! \av^{\Ga}_* M = m_{\beta *} (\Z \boxtimes M).$$
\xlemm

\pf By arguments similar to explicit computation of the (co)averaging functors in \cite[Lemma 2.22]{CassvdHScholbach:MotivicSatake} (see also \thref{av} or \cite[\S 2.3]{BGMRR:IwahoriWhittaker}), the functor $u^! \av^{\Ga}_!$ may be computed as follows. Using $\calU^{\exp}$-equivariance of $M$, form the twisted product $\Z[1](2)  \widetilde \boxtimes M$ on $(\calU \rtimes \Gm) \times^{\calU^{\exp}} \Fl_{\mathbf{f}}$ (cf.~\cite[\S 4.2.4]{CassvdHScholbach:MotivicSatake}), and take the $!$-pushforward along the multiplication map to $\Fl_{\mathbf{f}}$. (The shift and twist in $\Z[1](2)$ appear since $\av^{\Ga}_!$ involves $!$-pullback along a smooth projection map, and $\calU^{\exp} \subset \calU \rtimes \Gm$ has codimension $1$.) For any simple root $\beta$, the multiplication map $U_\beta \times \calU^{\exp} \r \calU \rtimes \Gm$ is an isomorphism, so that $(\calU \rtimes \Gm) \times^{\calU^{\exp}} \Fl_{\mathbf{f}} = U_{\beta} \times \Fl_{\mathbf{f}}$, and the result for  $u^! \av^{\Ga}_!$ follows. The argument for $u^! \av^{\Ga}_*$ is similar.
\xpf

\rema
Recall that $\DTM(\calU^{\exp} \backslash \Fl_{\mathbf{f}}) \subset \DM(\calU^{\exp} \backslash \Fl_{\mathbf{f}})$ consists of motives which are Tate with respect to the stratification by $\calU^{\exp}$-orbits (\thref{non-exp-Tate}). We likewise have the category $\DTM((\calU \rtimes \Gm) \backslash \Fl_{\mathbf{f}}) \subset \DM((\calU \rtimes \Gm) \backslash \Fl_{\mathbf{f}})$. While the latter category could also be defined using the stratification by $\calU^{\exp}$-orbits, by $\calU$-equivariance it is in fact equivalent to define Tateness with respect to the $\calU$-orbits $\Fl_{w}^{\circ}(\mathbf{a}_0, \mathbf{f})$. This follows since Tateness can be detected by pullback to the basepoints of orbits, cf.~\thref{lemm--DMorbit} and \cite[Lemma 4.28]{CassvdHScholbach:Central}. The same remark applies to anti-effective Tate motives.
\xrema

\prop \thlabel{GaAverageable}
The stratification of $\Fl_{\mathbf{f}}$ by $\calU^{\exp}$ is anti-effectively $\Ga$-averageable, i.e., the functors $$\av^{\Ga}_!, \av^{\Ga}_* \colon \DM(\calU^{\exp} \backslash \Fl_{\mathbf{f}}) \r \DM((\calU \rtimes \Gm) \backslash \Fl_{\mathbf{f}}),$$ preserve Tate motives, and $\av^{\Ga}_*$ furthermore preserves anti-effectivity.
\xprop

{\em{The functor $\av^{\Ga}_!$}}: Let $w \in W_{\mathbf{f}}^{\exp}$, and let $\Fl_{w}^{\circ}(\mathbf{a}_0, \mathbf{f})$ be the $\calU$-orbit containing $\Fl_w^{\exp}(\mathbf{f})$. By \thref{avg-comp}, $u^!\av^{\Ga}_! \Delta_w^{\exp}(\Z) = m_{\beta !}(\Z(1)[2] \boxtimes  \Delta_w^{\exp}(\Z))$. This is the same as the $!$-pushforward of $\Z(1)[2 + \dim \Fl_w^{\exp}(\mathbf{f})]$ along the composition $U_{\beta} \times \Fl_w^{\exp}(\mathbf{f}) \r U_{\beta} \times \Fl_{\mathbf{f}} \r \Fl_{\mathbf{f}}$, consisting of the inclusion followed by $m_{\beta}$. Since $U_\beta$ acts on $\Fl_{w}^{\circ}(\mathbf{a}_0, \mathbf{f})$, this composition can be refactored as
$$U_{\beta} \times \Fl_w^{\exp}(\mathbf{f}) \xrightarrow{a} \Fl_{w}^{\circ}(\mathbf{a}_0,\mathbf{f}) \xrightarrow{\iota} \Fl_{\mathbf{f}},$$ where $a$ is the action map and $\iota$ is the inclusion. Thus,
\begin{equation}
\label{avGa!}
u^!\av^{\Ga}_! \Delta_w^{\exp}(\Z)  = \iota_! a_! \Z(1)[2 + \dim \Fl_w^{\exp}(\mathbf{f})].
\end{equation} 
We now consider three cases. 

If $w \in W/W_{\mathbf{f}}$ and $w$ does not have maximal length in its left $W_0$-orbit, then $\Fl_w^{\exp}(\mathbf{f}) = \Fl^\circ_w(\mathbf{a}_0, \mathbf{f})$ by \thref{orbit-structure}. Thus $\Delta_w^{\exp}(\Z) = \Delta_w(\Z)$, and hence $u^!\av^{\Ga}_! \Delta_w^{\exp}(\Z) = \Delta_w^{\exp}(\Z)$ by \thref{example-avg}. This lies in $\DTM(\calU^{\exp} \backslash \Fl_{\mathbf{f}})$ by \thref{exponential orbits WT}.

If $w \in W/W_{\mathbf{f}}$ and $w$ has maximal length in its left $W_0$-orbit, then the inclusion $\Fl_w^{\exp}(\mathbf{f}) \subset \Fl^\circ_w(\mathbf{a}_0, \mathbf{f})$ is described in \eqref{eq--hyerplane}. We are free to switch the order of the factors, and it follows that for any $\beta \in \Delta$, the multiplication map $a \colon U_\beta \times \Fl_w^{\exp}(\mathbf{f})\r \Fl^\circ_w(\mathbf{a}_0, \mathbf{f})$ is an isomorphism. Then by \eqref{avGa!}, $u^!\av^{\Ga}_! \Delta_w^{\exp}(\Z)  = \Delta_w(\Z)(1)[1] \in \DTM(\calU^{\exp} \backslash \Fl_{\mathbf{f}})$. 

Finally, if $w \in {}_0W_{\mathbf{f}}$, let $v \in W/W_{\mathbf{f}}$ be the same element. Applying $u^!\av^{\Ga}_!$ to the fiber sequence $\Delta_v^{\exp}(\Z) \r \Delta_{w}^{\exp}(\Z) \r \Delta_v(\Z)$ and using the previous cases, we get a fiber sequence $\Delta_v(\Z)(1)[1] \r u^!\av^{\Ga}_! \Delta_w^{\exp}(\Z) \r \Delta_v(\Z)$.

{\em{The functor $\av^{\Ga}_*$}}: Similarly to \eqref{avGa!} (and using \thref{WT.prod}), if $w \in W_{\mathbf{f}}^{\exp}$ we have
$$u^!\av^{\Ga}_* \nabla_w^{\exp}(\Z)  = \iota_* a_* \Z[\dim \Fl_w^{\exp}(\mathbf{f})].$$ By arguments analogous to those above, $u^! \av^{\Ga}_* \nabla_w^{\exp}(\Z) = \nabla_w^{\exp}(\Z)$ if $w \in W/W_{\mathbf{f}}$ is not maximal, $u^! \av^{\Ga}_* \nabla_w^{\exp}(\Z) = \nabla_w(\Z)[-1]$ if $w \in W/W_{\mathbf{f}}$ is maximal, and we have a fiber sequence $\nabla_v(\Z) \r u^!\av^{\Ga}_* \nabla_w(\Z)^{\exp} \r \nabla_v(\Z)(-1)[-1]$ if $w \in {}_0W_{\mathbf{f}}$ and $v \in W/W_{\mathbf{f}}$ is the same element. Thus, $u^!\av^{\Ga}_* \nabla_w^{\exp}(\Z) \in \DTM^{\anti}(\calU^{\exp} \backslash \Fl_{\mathbf{f}})$ in each case by \thref{exponential orbits WT}.

\defi
We define $\fib_!^{\Ga} := \fib (\id \r v^! \av^{\Ga}_!)$ to be the endofunctor of 
$\DTM(\calU^{\exp} \backslash \Fl_{\mathbf{f}})$ as in \thref{t-structure DTMexp}, where $v^! \colon \DM((\calU \rtimes \Gm) \backslash \Fl_{\mathbf{f}}) \r \DM(\calU^{\exp} \backslash \Fl_{\mathbf{f}})$ is the forgetful functor.
\xdefi

\lemm \thlabel{nunuR} Let $u^! \colon \DM(\calU^{\exp}  \backslash \Fl_{\mathbf{f}}) \r \DM(\Fl_{\mathbf{f}})$ be the forgetful functor.  Then for any simple root $\beta \in \Delta$ and $M \in \DTM(\calU^{\exp} \backslash \Fl_{\mathbf{f}})$, we have
$$u^! \fib_!^{\Ga} M = m_{\beta !} j_* j^* (\Z(1)[1] \boxtimes M)$$
where $j$ is the inclusion of the complement of the zero section $\{u_\beta(1)\} \times \Fl_{\mathbf{f}} \r U_\beta \times  \Fl_{\mathbf{f}}$, and $m_\beta \colon U_\beta \times \Fl_{\mathbf{f}}\r \Fl_{\mathbf{f}}$ is the multiplication map.
\xlemm

Note that we have used the notation $u^!$ in \thref{avg-comp} and \thref{nunuR} to denote forgetful functors with slightly different domains, but there should be no risk of confusion.

\pf 
The map $0_! M \r \pr^!M$ in the notation of \thref{t-structure DTMexp} is part of a localization sequence for $\pr^!M$ on $U_\beta \times \Fl_{\mathbf{f}}$. The cofiber of this localization sequence is (after forgetting equivariance) given by $j_* j^* \pr^!M$. Taking the fiber then introduces a shift by $[-1]$. By our work in \thref{avg-comp}, $u^! \fib_!^{\Ga}$ is then obtained by further applying $m_{\beta !} = \act_!$.
\xpf

\prop \thlabel{nunu}
The endofunctor $\fib_!^{\Ga}$ of
$\DTM(\calU^{\exp} \backslash \Fl_{\mathbf{f}})$ is t-exact.
\xprop

\pf
{\em{Right t-exactness}}: It suffices to show that if $w \in W_{\mathbf{f}}^{\exp}$, then $u^! \fib_!^{\Ga} \Delta_w^{\exp}(\Z) \in \DTM(\Fl_{\mathbf{f}})^{\leq 0}$. Since $\fib_!^{\Ga}$ vanishes on $(\calU \rtimes \Gm)$-equivariant motives by \thref{example-avg}, this completes the case where $w \in W/W_{\mathbf{f}}$ is not left $W_0$-maximal. As in the proof of \thref{GaAverageable}, by considering an excision sequence on any $\calU$-orbit which splits into two $\calU^{\exp}$-orbits, we reduce to the case where $w \in W/W_{\mathbf{f}}$ is left $W_0$-maximal.

If  $w \in W/W_{\mathbf{f}}$ is left $W_0$-maximal, let $\Fl_w^\circ(\mathbf{a}_0, \mathbf{f})$ be the $\calU$-orbit containing $\Fl^{\exp}_w(\mathbf{f})$. By \thref{nunuR} (and the K\"unneth formula \thref{WT.prod} applied to $j_*$), we have
\begin{equation}
\label{fibGa!}
u^! \fib_!^{\Ga} \Delta_w^{\exp}(\Z) = \iota_! a_! (L \boxtimes \Z[\dim \Fl^{\exp}_w(\mathbf{f})]),
\end{equation}
where $a \colon U_\beta \times \Fl^{\exp}_w(\mathbf{f}) \r \Fl_w^{\circ}(\mathbf{a}_0, \mathbf{f})$ is the action map, $\iota \colon \Fl_w^{\circ}(\mathbf{a}_0, \mathbf{f}) \r \Fl_{\mathbf{f}}$ is the inclusion, and $L \in \DM(U_\beta)$ is the $*$-extension of $\Z(1)[1]$ from $U_\beta \setminus \{u_{\beta}(1)\} \cong \Gm$. By excision on $U_\beta$, $L$ sits in the middle of a fiber sequence $L_1 \r L \r L_2$, where $L_1$ is the $!$-extension of $\Z(1)[1]$ from $U_\beta \setminus \{u_{\beta}(1)\} \cong \Gm$ to $U_\beta$, and $L_2$ is supported on $u_{\beta}(1) \cong S$. By relative purity, $L_2 \in \DTM(S)^{\leq 0}$. Note that $a$ is a stratum-preserving isomorphism, where the domain is stratified by $U_\beta \setminus \{u_{\beta}(1)\} \times \Fl_w^{\exp}(\mathbf{a}_0, \mathbf{f})$ and $\{u_{\beta}(1)\} \times \Fl_w^{\exp}(\mathbf{a}_0, \mathbf{f})$, and the codomain is stratified by the two $\calU^{\exp}$-orbits. Since $\iota_!$ is right t-exact, we conclude that $u^! \fib_!^{\Ga} \Delta_w^{\exp}(\Z) \in \DTM(\Fl_{\mathbf{f}})^{\leq 0}$. 

{\em{Left t-exactness}}: By \thref{rema--Exact}, it suffices to show that if $M \in \MTM(S)$ and $w \in W_{\mathbf{f}}^{\exp}$, then $u^! \fib_!^{\Ga} \nabla_w^{\exp}(M) \in \DTM(\Fl_{\mathbf{f}})^{\geq 0}$. Similarly to the proof of right t-exactness, by considering \thref{example-avg} and excision, we reduce to the case where $w \in W/W_{\mathbf{f}}$ is left $W_0$-maximal. Write $w = w_0 w_{\min}$ where $w_{\min} \in W$ is left $W_0$-minimal. By \thref{convlem} and \thref{WT.prod}, $\nabla_w^{\exp}(M) = \nabla_{w_0}^{\exp}(\Z) \star \nabla_{w_{\min}}(M)$. Then by \thref{avg-comp}, $u^! \av_!^{\Ga} \nabla_w^{\exp}(M) = u^! \av_!^{\Ga} \nabla_{w_0}^{\exp}(\Z) \star \nabla_{w_{\min}}(M)$. It follows that
$$u^! \fib_!^{\Ga} (\nabla_w^{\exp}(M)) = u^! \fib_!^{\Ga} (\nabla_{w_0}^{\exp}(\Z)) \star \nabla_{w_{\min}}(M).$$

The fiber of the natural map $\Delta_{w_0}^{\exp}(\Z) \r  \nabla_{w_0}^{\exp}(\Z)$ lies in the subcategory generated by the forgetful image of $\DTM (\calU \rtimes \Gm \backslash \Fl)$, so that by \thref{example-avg}, $u^! \fib_!^{\Ga} (\nabla_{w_0}^{\exp}(\Z)) = u^! \fib_!^{\Ga} (\Delta_{w_0}^{\exp}(\Z))$. Let $w_0' \in {}_0W_{\mathbf{a}}$ be the element $w_0$. By \eqref{fibGa!} and the remarks following it, $u^! \fib_!^{\Ga} (\Delta_{w_0}^{\exp}(\Z))$ is obtained taking the $*$-extension of $\Z(1)[\dim \Fl^{\exp}_{w_0'}(\mathbf{a}_0)]$ along the open immersion $\Fl^{\exp}_{w_0'}(\mathbf{a}_0) \r \Fl^{\circ}_{w_0}(\mathbf{a}_0)$, followed by the $!$-extension along $\Fl^{\circ}_{w_0}(\mathbf{a}_0) \r \Fl$. Hence by excision for $\Z$ on $\Fl^{\circ}_{w_0}(\mathbf{a}_0)$, there is a fiber sequence
$$\Delta_{w_0}(\Z)(1) \r u^! \fib_!^{\Ga} (\Delta_{w_0}^{\exp}(\Z)) \r \Delta^{\exp}_{w_0}(\Z).$$ Putting this all together, $u^! \fib_!^{\Ga} (\nabla_w^{\exp}(M))$ lies in the middle of a fiber sequence between $\Delta_{w_0}(\Z)(1) \star \nabla_{w_{\min}}(M)$ and $\Delta_{w_0}^{\exp}(\Z) \star \nabla_{w_{\min}}(M)$. These lie in nonnegative degrees by \cite[Lemma 5.27]{CassvdHScholbach:Central} and \thref{leftright}, respectively.
\xpf

\subsection{Exponential Tate motives} 
In this section we define exponential Tate motives on affine flag varieties and deduce some of their properties. The group actions recalled at the start of \refsect{Ga} apply to the setup of \thref{DMexp G-equivariant} (see also \thref{usage  exponential motives later on}), so that we have the category
$$\DM_{\exp}(\Fl_{\mathbf{f}}) : = \DM(\calU^{\exp} \backslash \Fl_{\mathbf{f}}) / \DM((\calU \rtimes \Gm) \backslash \Fl_{\mathbf{f}}),$$ along with the equivalent categories $\DM_{\exp!}(\Fl_{\mathbf{f}})$ and $\DM_{\exp*}(\Fl_{\mathbf{f}})$.

Recall that the stratification of $\Fl_{\mathbf{f}}$ by $\calU^{\exp}$-orbits is anti-effectively $\Ga$-averageable by \thref{GaAverageable}, and the functor $\fib^!_{\Ga}$ is t-exact by \thref{nunu}.  Thus, the setups of \thref{DTMexp}, \thref{anti-effective DTMexp}, and \thref{t-structure DTMexp} apply (or, more precisely, their extensions to the prestack quotients  as discussed before \thref{def--exp}). This entitles us to the following definition.

\defi \thlabel{defi exp on Fl}
The category of exponential stratified (anti-effective) Tate motives on $\Fl_{\mathbf{f}}$ is
$$\DTM_{\exp}^{\xanti}(\Fl_{\mathbf{f}}) : = \DTM^{\xanti}(\calU^{\exp} \backslash \Fl_{\mathbf{f}})  / \DTM^{\xanti}((\calU \rtimes \Gm) \backslash \Fl_{\mathbf{f}}).$$ Additionally, $\MTM_{\exp}^{\xanti}(\Fl_{\mathbf{f}}) \subset \DTM_{\exp}^{\xanti}(\Fl_{\mathbf{f}})$ is the heart of the t-structure afforded by \thref{t-structure DTMexp}.
\xdefi

We also have the equivalent categories $\DTM_{\exp!}(\Fl_{\mathbf{f}})$ and $\DTM_{\exp*}(\Fl_{\mathbf{f}})$ of colocal, resp.~local objects in $\DTM(\calU^{\exp} \backslash \Fl_{\mathbf{f}})$. The category $\DTM_{\exp}^{\xanti}(\Fl_{\mathbf{f}})$ is a full subcategory of $\DM_{\exp}(\Fl_{\mathbf{f}})$ by \thref{full exp subcat}. As a consequence, any functor between categories of the form
$\DTM^{\xanti}(\calU^{\exp} \backslash \Fl_{\mathbf{f}})$ also induces a functor on $\DTM^{\xanti}_{\exp}(\Fl_{\mathbf{f}})$ (cf.~\thref{functoriality DTMexp anti}).

The following result generalizes \thref{DMexp A^1} to the case of $\calI$-orbits in affine flag varieties.

\prop\thlabel{exponential motives on iwahori orbits}
Let $w \in W/W_{\mathbf{f}}$. 
\begin{enumerate}
    \item There are canonical $t$-exact equivalences
$$
\DTM_{\exp}^{\xanti}(\Fl_{w}^{\circ}(\mathbf{a}_0, \mathbf{f})) = \begin{cases}
    \DTM(S)^{\xanti}, & w \text{ is left } W_0\text{-maximal,} \\
    0, & \text{otherwise}.
    \end{cases}
$$
\item If $w$ is left $W_0$-maximal, regard $w$ as an element of ${}_0W_{\mathbf{f}}$, and let $j \colon \calU^{\exp} \backslash \Fl_{w}^{\exp}(\mathbf{f}) \r \calU^{\exp} \backslash \Fl_{w}^{\circ}(\mathbf{a}_0, \mathbf{f})$ be the prestack quotient of the open inclusion. Then the above equivalence is the top row of the following commutative diagram:
\[\begin{tikzcd}
	\DTM(S) \arrow[d, "\sim"] \arrow[r, "\sim"] & \DTM_{\exp}(\Fl_{w}^{\circ}(\mathbf{a}_0, \mathbf{f})) \\
	\DTM(\calU^{\exp} \backslash \Fl_{w}^{\exp}(\mathbf{f}) ) \arrow[r, "{j_!}"] &  \DTM(\calU^{\exp} \backslash \Fl_{w}^{\circ}(\mathbf{a}_0, \mathbf{f})) \arrow[u, "q"']
\end{tikzcd}\]
Here $q$ is the categorical quotient, and the left vertical equivalence equips an element $h^*M[\dim \Fl_{w}^{\exp}(\mathbf{f})]$ with the trivial $\calU^{\exp}$-equivariance datum, where $M \in \DTM(S)$ and $h \colon \Fl_w^{\exp}(\mathbf{f}) \r S$ is the structure map.

\item In the situation of part (2), we have an equivalence
$$q \circ j_*(-)= q \circ j_!(-)(-1) $$ of functors $\DTM(\calU^{\exp} \backslash \Fl_{w}^{\exp}(\mathbf{f}) ) \r \DTM_{\exp}(\Fl_{w}^{\circ}(\mathbf{a}_0, \mathbf{f}))$.
\end{enumerate}
\xprop

\pf
If $w$ is not left $W_0$-maximal, then $\Fl_{w}^{\circ}(\mathbf{a}_0, \mathbf{f}) = \Fl^{\exp}_{w}(\mathbf{f})$ by \thref{orbit-structure}. By \thref{lemm--DMorbit}, we have $\DTM(\calU^{\exp} \backslash \Fl^{\exp}_{w}(\mathbf{f})) = \DTM(S/ \Gm)$. By essentially the same argument, $\DTM((\calU \rtimes \Gm) \backslash \Fl^{\exp}_{w}(\mathbf{f})) = \DTM(S/ \Gm)$, with the forgetful functor from the latter category to the former category being identified with the identity functor on $\DTM(S/\Gm)$. This implies $\DTM_{\exp}(\Fl_{w}^{\circ}(\mathbf{a}_0, \mathbf{f})) = 0$ is the trivial category.

If $w$ is left $W_0$-maximal, let us relabel $w \in {}_0 W_{\mathbf{f}}$ to be the same element, and let $v$ be the element $w$ regarded as an element of $W/W_{\mathbf{f}}$, so that we have the closed immersion $i \colon \calU^{\exp} \backslash \Fl_{v}^{\exp}(\mathbf{f}) \r \calU^{\exp} \backslash \Fl_{w}^{\circ}(\mathbf{a}_0, \mathbf{f})$. By \thref{DTMexp}, the functor $q$ induces an isomorphism of $\ker (\av_!^{\Ga})$ with $\DTM_{\exp}(\Fl_{w}^{\circ}(\mathbf{a}_0, \mathbf{f}))$, where $\av_!^{\Ga} \colon \DTM(\calU^{\exp} \backslash \Fl_{w}^{\circ}(\mathbf{a}_0, \mathbf{f})) \r \DTM((\calU \rtimes \Gm) \backslash \Fl_{w}^{\circ}(\mathbf{a}_0, \mathbf{f}))$ is the $!$-averaging functor.  

We will determine $\ker (\av_!^{\Ga})$. For any $M \in \DTM(\calU^{\exp} \backslash \Fl_{w}^{\circ}(\mathbf{a}_0, \mathbf{f}))$, we have a fiber sequence
$$i_* i^! M \r M \r j_* j^*M$$
in $\DTM(\calU^{\exp} \backslash \Fl_{w}^{\circ}(\mathbf{a}_0, \mathbf{f}))$. We claim that $\av_!^{\Ga} (j_*j^*M) = 0$. Since $\DTM(\calU^{\exp} \backslash \Fl_{w}^{\exp}(\mathbf{f}) )  \cong \DTM(S)$ by \thref{lemm--DMorbit}, it suffices to prove the claim when $M= \Z$ supported on $\Fl_{w}^{\circ}(\mathbf{a}_0, \mathbf{f})$. By relative purity, $i_*i^!\Z = i_*i^* \Z(-1)[-2]$. Now using that the multiplication map $a \colon U_\beta \times \Fl_{v}^{\exp}(\mathbf{f}) \r \Fl_{w}^{\circ}(\mathbf{a}_0, \mathbf{f})$ is an isomorphism and arguing as in \eqref{avGa!}, it follows that $\av_!^{\Ga}(i_*i^!\Z) = \Z$ and the canonical map to $\av_!^{\Ga} \Z = \Z$ is an isomorphism. Thus, $\av_!^{\Ga} (j_*j^*M) = 0$. Hence $\av_!^{\Ga} M = 0$ if and only if $\av_!^{\Ga} (i_*i^!M) = 0$. Again, using that $a$ is an isomorphism, $\av_!^{\Ga} (i_*i^!M) = 0$ if and only if $i_*i^!M = 0$. Thus, $\ker (\av_!^{\Ga})$ is identified with the essential image of $j_* \colon\DTM(\calU^{\exp} \backslash \Fl_{w}^{\exp}(\mathbf{f})) \r \DTM((\calU \rtimes \Gm) \backslash \Fl_{w}^{\circ}(\mathbf{a}_0, \mathbf{f}))$. Since $j_*$ is fully faithful, $\DTM_{\exp}(\Fl_{w}^{\circ}(\mathbf{a}_0, \mathbf{f}))  \cong \DTM(S)$. This also proves that we have a diagram satisfying all the properties of (2), except where $j_!$ is replaced by $j_*$. However, the proof of (3) below then gives diagram in (2) with $j_!$.

To prove (3) we argue as in \thref{j! j* DMexp}. More precisely, view any element $M \in \DTM(S)$ as an element of $\DTM(\Fl_w^\circ(\mathbf{a}_0, \mathbf{f}))$ via $*$-pullback, and equip $M$ with the trivial $\calU^{\exp}$-equivariance datum. Note that $q(M) = 0$ since $M$ is also canonically $\calU \rtimes \Gm$-equivariant. Then by relative purity and excision, after applying $q$ we have isomorphisms
$$j_*j^*M \cong i_*i^! M[1] \cong i_*i^*M(-1)[-1] \cong j_! j^*M(-1)$$
in $\DTM_{\exp}(\Fl_{w}^{\circ}(\mathbf{a}_0, \mathbf{f}))$.

All that remains to prove is that the equivalence in (2) respects anti-effective motives. As in \thref{full exp subcat}, $\DTM_{\exp}^{\anti}(\Fl_{w}^{\circ}(\mathbf{a}_0, \mathbf{f}))$ is generated by the $*$- (equivalently, by the $!$-) extensions of $\Z(n)$ along the two $\calU^{\exp}$-orbits in $\Fl_{w}^{\circ}(\mathbf{a}_0, \mathbf{f})$, for $n \leq 0$. By the last isomorphism in the above chain, the $!$-extensions of $\Z(n)$ along the larger orbit $\Fl_{w}^{\exp}(\mathbf{f})$ suffice. (This also explains why we prefer to use $j_!$ rather than $j_*$ in (2), as $q \circ j_*\Z(1) = q \circ j_!\Z$, so that the analogous diagram with $j_*$ does not identify anti-effective subcategories.)
\xpf

We now define the exponential standard and costandard functors.
To distinguish the inclusion maps for exponential and parahoric orbits, we introduce the following notation.

\nota\thlabel{notation parahoric inclusions}
Let \(\bbf,\bbf'\) be standard facets, and \(w\in W_{\bbf'}\backslash W /W_{\bbf}\).
Then the inclusion of the \(L^+G_{\bbf'}\)-orbit in \(\Fl_{\bbf}\) through \(w\) is denoted
\[\iota_w^{\bbf'}\colon \Fl_w^{\circ}(\bbf',\bbf) \hookrightarrow \Fl_{\bbf}.\]
\xnota

The facet \(\bbf\) does not appear in the notation, since it will always be clear from the context.
In particular, any inclusion of the form \(\iota_w\), i.e., without superscript, will denote the inclusion of a \(\calU^{\exp}\)-orbit.

\defi\thlabel{defi-exponential co/standard}
Let \(w\in {}_0 W_{\bff}\).
Then the exponential standard functor is
\[\Delexp_w \colon \colon \DTM(S) \cong \DTM_{\exp}(\Fl_w^{\circ}(\bba_0,\bbf)) \xrightarrow{\iota_{w!}^{\bba_0}} \DTM_{\exp}(\Fl_{\bbf}).\]
Similarly, the exponential costandard functor is
\[\nabexp_w\colon \DTM(S) \cong \DTM_{\exp}(\Fl_w^{\circ}(\bba_0,\bbf)) \xrightarrow{\iota_{w*}^{\bba_0}} \DTM_{\exp}(\Fl_{\bbf}).\]
Here, the equivalence \(\DTM(S)\cong \DTM_{\exp}(\Fl_w^{\circ}(\bba_0,\bbf))\) is given by \thref{exponential motives on iwahori orbits} (2).
\xdefi

The fact that these functors indeed take values in \(\DTM_{\exp}(\Fl_\bbf) \subseteq \DM_{\exp}(\Fl_{\bbf})\) follows from \thref{exponential orbits WT}, and the following lemma.

\lemm
Let \(w\in {}_0 W_{\bbf}\), and denote by \(q\colon \DTM(\calU^{\exp}\backslash \Fl_{\bbf}) \to \DTM_{\exp}(\Fl_{\bbf})\) the quotient map.
Then there are natural equivalences
\[\Delexp_w \cong q\circ \Delta^{\exp}_w \cong \DTM(S) \to \DTM_{\exp}(\Fl_{\bbf})\]
and 
\[\nabexp_w \cong q \circ \nabla^{\exp}_w(-1) \colon \DTM(S) \to \DTM_{\exp}(\Fl_{\bbf}).\] Furthermore, $\DTM_{\exp}^{\xanti}(\Fl_\mathbf{f}) \subset \DM_{\exp}(\Fl_\mathbf{f})$ is the presentable stable subcategory generated by the $\Delexp_w \Z(n)$ for $n \in \Z$ (resp. $n \leq 0)$ and \(w\in {}_0 W_{\bff}\).
\xlemm
\pf
This follows from \thref{exponential motives on iwahori orbits} and the equivariant analogue of \thref{full exp subcat}.
\xpf

\rema \thlabel{rema--CompactExp}
The functors $\Delexp_w $ and $\nabexp_w$ preserve compact objects. Indeed, this follows from arguments similar to \thref{rema--Compact}, combined with \thref{exponential motives on iwahori orbits}. Additionally, an object of $\DTM_{\exp}(\Fl_{\bbf})$ is compact if and only if it is supported on finitely many Iwahori orbits, and its $*$- (or $!$-) pullback to each orbit is identified with a compact object of $\DTM(S)$ under \thref{exponential motives on iwahori orbits}. Similar remarks apply to reduced motives.
\xrema

We can now deduce the following t-exactness properties, since \(q\colon \DTM(\calU^{\exp} \backslash \Fl_{\bbf}) \to \DTM_{\exp}(\Fl_{\bbf})\) is t-exact by \thref{t-structure DTMexp}.

\prop\thlabel{exponential t-exactness}
For \(w\in {}_0 W_{\bbf}\), the functors \(\Delexp_w, \nabexp_w\colon \DTM(S) \to \DTM_{\exp}(\Fl_{\bbf})\) are t-exact.
\xprop
\pf
This follows from \thref{expFLfstandards} and the equivalence $\nabexp_w \cong q \circ \nabla^{\exp}_w(-1)$.
\xpf

\prop\thlabel{exponential convolution}
Let \(\bbf,\bbf'\) be standard facets. Then convolution induces a functor
\[(-)\star(-) \colon \DTM_{\exp}^{\xanti}(\Fl_{\bbf'}) \times \DTM^{\xanti}(L^+G_{\bbf'} \backslash \Fl_{\bbf}) \to \DTM_{\exp}^{\xanti}(\Fl_{\bbf}).\]
Moreover, it restricts to a functor
\[(-)\star(-) \colon \DTMexp(\Fl_{\bbf'})^{\comp} \times \DTM(L^+G_{\bbf'} \backslash \Fl_{\bbf})^{\locc} \to \DTMexp(\Fl_{\bbf})^{\comp},\]
i.e., right-convolution with a locally compact object preserves compact objects.

Finally, for any \(w\in {}_0 W_{\bbf'}\), convolution by $\Delexp_w(\Z)$ (resp.~$\nabexp_w(\Z)$) defines a left (resp.~right) $t$-exact functor $\DTM^{\xanti}(L^+G_{\mathbf{f}'} \backslash \Fl_{\mathbf{f}}) \r \DTM_{\exp}^{\xanti}(\Fl_{\bbf})$.
\xprop
\pf
The convolution product \[\star\colon \DTM^{\xanti}(\calU^{\exp}\backslash \Fl_{\bbf'}) \times \DTM^{\xanti}(L^+G_{\bbf'} \backslash \Fl_{\bbf}) \to \DTM^{\xanti}(\calU^{\exp} \backslash \Fl_{\bbf})\] from \thref{convDTM} respects $(\calU \rtimes \Gm)$-equivariant objects, so the same construction gives a functor
\[\DTM^{\xanti}((\calU \rtimes \Gm)\backslash \Fl_{\bbf'}) \times \DTM^{\xanti}(L^+G_{\bbf'} \backslash \Fl_{\bbf}) \to \DTM^{\xanti}((\calU\rtimes \Gm) \backslash \Fl_{\bbf}).\]
Passing to the quotient categories $\DTM_{\exp}^{\xanti}(\Fl_{\bbf'})$ and $\DTM_{\exp}^{\xanti}(\Fl_{\bbf})$ in \thref{convDTM} then gives the desired convolution functor. The t-exactness properties follow from \thref{leftright}.
Finally, the fact that right-convolution with a locally compact object preserves compact objects also follows \thref{convDTM}, since the quotient \(\DTM(\calU^{\exp} \backslash \Fl_{\bbf}) \to \DTMexp(\Fl_\bbf)\) sends locally compact objects to compact objects.
\xpf

Let us also record the following result, concerning certain Hom-groups in \(\DTMexp(\Fl_{\bbf})\).

\lemm\thlabel{exponential homs}
Let \(\bbf\) be a standard facet, \(v,w\in {}_0W_{\bbf}\), and \(\calF,\calG\in \DTM(S)\).
Then 
\[\Hom_{\DTM_{\exp}(\Fl_{\bbf})}(\Delexp_v(\calF),\nabexp_w(\calG)) \cong \begin{cases}
	\Hom_{\DTM(S)}(\calF,\calG) &\text{ if } v=w\\
	0 & \text{ if } v\neq w.
\end{cases}\]
Similar equivalences hold when replacing the Hom-groups \(\Hom\) by the mapping spectra \(\Maps\).
\xlemm
\pf
This follows immediately from adjunction and base change, and the description of exponential motives on Iwahori-orbits as in \thref{exponential motives on iwahori orbits}.
\xpf

\subsection{Twisted exponential motives}\label{subsec--twisted exponential motives}

As a final topic in this section, we introduce \(\hat{\rho}\)-twisted exponential motives.
These will be used to state \thref{Motivic CasselmanShalika} uniformly for all groups, thereby removing a restriction appearing in \cite{BGMRR:IwahoriWhittaker}.
This twist is motivated by \cite[Definition 2.2.1]{Raskin:Walgebras}.

Recall that \(\hat{\rho}\) defines a cocharacter of the adjoint torus \(T_{\adj}\).
Since the conjugation action of the center \(Z_G\) of \(G\) on \(LG\) is trivial, we get a conjugation action of \(LT_{\adj}\) on \(LG\).
Consider the subgroup \(\calI_{\hat{\rho}}\subseteq LG\), given by the image of the Iwahori \(\calI\subseteq LG\) under the automorphism \(LG\cong LG\colon x\mapsto \hat{\rho}(t^{-1}) x \hat{\rho}(t)\).
In particular, there is an isomorphism \(\calI_{\hat{\rho}} \cong \calI\colon x\mapsto \hat{\rho}(t)x\hat{\rho}(t^{-1})\).
Similarly, we define \((\calU^{>0} \rtimes U_0)_{\hat{\rho}} \subseteq \calU_{\hat{\rho}}\subseteq \calI_{\hat{\rho}}\) and \(\calU_{\hat{\rho}}^{\exp} := (\calU^{>0}\rtimes U_0)_{\hat{\rho}} \rtimes \Gm\).

\rema
Fix a standard facet \(\bbf\).
The \(\calI_{\hat{\rho}}\)- and \(\calU_{\hat{\rho}}^{\exp}\)-orbits on \(\Fl_{\bbf}\) satisfy the same properties as the \(\calI\)- and \(\calU^{\exp}\)-orbits.
Indeed, when \(G\) is adjoint (so that \(\hat{\rho}\) actually defines an element in \(X_*(T)^+\)), the isomorphism \(\Fl_{\bbf}\cong \Fl_{\bbf}\colon x\mapsto \hat{\rho}(t^{-1})\cdot x\) identifies the \(\calI\)-orbits with the \(\calI_{\hat{\rho}}\)-orbits, and the \(\calU^{\exp}\)-orbits with the \(\calU_{\rho}^{\exp}\)-orbits.
In general, this follows from the following lemma.
\xrema

Recall that the standard facets for the group \(G\) and \(G_{\adj}\) are naturally in bijection.
To state the lemma below, we denote by \(\Fl_{\adj,\bbf}\) the partial affine flag variety for the adjoint group \(G_{\adj}\), and the facet \(\bbf\) under this bijection.

\lemm\thlabel{reduction to adjoint case}
The map \(\Fl_{\bbf} \to \Fl_{\adj,\bbf}\) identifies \(\calI_{\hat{\rho}}\)- and \(\calU_{\hat{\rho}}^{\exp}\)-orbits (for the respective groups \(G\) and \(G_{\adj}\)).
Moreover, when restricted to each connected component, it induces an isomorphism on \(\DM\), which identifies the subcategory of (anti-effective) stratified Tate motives.
In particular, the stratification of \(\Fl_{\bbf}\) by \(\calI_{\hat{\rho}}\)- or \(\calU_{\hat{\rho}}^{\exp}\)-orbits is anti-effective universally Whitney--Tate.
\xlemm
\pf
The map \(\alpha\colon \Fl_{\bbf}\to \Fl_{\adj,\bbf}\) induces a universal homeomorphism on each connected component, which induces an isomorphism stratawise.
This is proven in \cite[Lemma 6.24]{CassvdHScholbach:MotivicSatake} for the affine Grassmannian, but the proof holds in general.
This shows that \(\alpha\) preserves the orbits under \(\calI_{\hat{\rho}}\) and \(\calU_{\hat{\rho}}^{\exp}\).
The second statement then follows immediately from the description of \(\DM(\Fl_{\bbf})\) and \(\DTM(\Fl_{\bbf})\) as a lax limit \cite[Remark 4.10]{EberhardtScholbach:Integral}, and the fact that proper pushforward along \(\alpha\) commutes with push-pull along the strata.
Finally, the fact that the stratification is anti-effective universally Whitney--Tate follows from the untwisted case, which is \thref{exponential orbits WT}.
\xpf

In particular, we can define (\(\hat{\rho}\)-)twisted exponential motives on \(\Fl_{\bbf}\), by mimicking \thref{defi exp on Fl}.

\defi\thlabel{defi twisted exponentials}
The category of (\(\hat{\rho}\)-)twisted exponential motives on \(\Fl_{\bbf}\) is
\[\DM_{\rhoexp}(\Fl_{\bbf}) := \DM(\calU_{\hat{\rho}}^{\exp} \backslash \Fl_{\bbf}) / \DM((\calU_{\hat{\rho}}\rtimes \Gm)\backslash \Fl_{\bbf}).\]
Similarly, the category of (anti-effective) stratified Tate (\(\hat{\rho}\)-)twisted exponential motives on \(\Fl_{\bbf}\) is
\[\DTM_{\rhoexp}^{\xanti}(\Fl_{\bbf}) := \DTM^{\xanti}(\calU_{\hat{\rho}}^{\exp}\backslash \Fl_{\bbf})/\DTM^{\xanti}((\calU_{\hat{\rho}}\rtimes \Gm)\backslash \Fl_{\bbf}).\]
By \thref{DTMexp full subcategory}, this is a full subcategory of \(\DM_{\rhoexp}(\Fl_{\bbf})\).
Moreover, by \thref{reduction to adjoint case}, \(\DTM_{\rhoexp}^{\xanti}(\Fl_{\bbf})\) admits a natural t-structure, compatible with the one on \(\DTM_{\exp}^{\xanti}(\Fl_{\bbf})\) in the adjoint case, whose heart is denoted \(\MTM_{\rhoexp}^{\xanti}(\Fl_{\bbf})\).
\xdefi

The fact that the convolution product involving twisted exponential motives preserves \(\DTM\) follows from the untwisted case. 

\coro\thlabel{twisted exponential convolution}
Let \(\bbf,\bbf'\) be two standard facets.
The natural convolution functor
\[\star \colon \DM_{\rhoexp}(\Fl_{\bbf'}) \times \DM(L^+G_{\bbf'} \backslash\Fl_{\bbf}) \to \DM_{\rhoexp}(\Fl_{\bbf})\]
restricts to a functor
\[\star \colon \DTM_{\rhoexp}(\Fl_{\bbf'}) \times \DTM(L^+G_{\bbf'}\backslash\Fl_{\bbf}) \to \DTM_{\rhoexp}(\Fl_{\bbf}).\]
\xcoro
\pf
By \thref{reduction to adjoint case} we may assume \(G\) is adjoint. 
Then the isomorphisms \(\Fl_{\bbf} \cong \Fl_{\bbf}\) and \(\Fl_{\bbf'}\cong \Fl_{\bbf'}\), both given by \(x\mapsto \hat{\rho}(t^{-1})x\), allow us to conclude by \thref{exponential convolution}.
\xpf

The main difference between the \(\calU_{\hat{\rho}}^{\exp}\)- and \(\calU^{\exp}\)-orbits on \(\Fl_{\bbf}\) is that they are indexed by different sets.
Let us focus on the case \(\bbf=\bbf_0\), so that \(\Fl_{\bbf_0}=\Gr\) is the affine Grassmannian.
Then the \(\calI\)-orbits on \(\Gr\) are naturally indexed by the cocharacters \(X_*(T)\), and those orbits splitting up into two \(\calU^{\exp}\)-orbits correspond to the strictly dominant cocharacters \(X_*(T)^{++}\) by \thref{lemm-exampleW}.

\prop\thlabel{twisted exponential orbits in Gr}
The \(\calI_{\hat{\rho}}\)-orbits on \(\Gr\) are naturally indexed by \(X_*(T)\).
Moreover, the orbits splitting up into two \(\calU_{\hat{\rho}}^{\exp}\)-orbits are naturally indexed by \(X_*(T)^+\).
\xprop
\pf
First, we assume \(G\) is adjoint.
In that case, the isomorphism \(\Gr \cong \Gr \colon x\mapsto \hat{\rho}(t^{-1})x\) identifies the \(\calI\)-orbits with the \(\calI_{\hat{\rho}}\)-orbits, so that the latter are also indexed by \(X_*(T)\).
Moreover, the \(\calI\)-orbits splitting up into two \(\calU^{\exp}\)-orbits are the orbits through \(t^{\lambda}\) for \(\lambda\in X_*(T)^{++}\), so that the \(\calI_{\hat{\rho}}\)-orbits in \(\Gr\) splitting up into two \(\calU_{\rho}^{\exp}\)-orbits are indexed by \(X_*(T)^{++} - \hat{\rho} = X_*(T)^+\); they are precisely the orbits through \(\lambda(t)\) for \(\lambda\in X_*(T)^+\).

Now, let \(G\) be general, with adjoint quotient \(G_{\adj}=G/Z_G\) and maximal torus \(T_{\adj} = T/Z_G\subseteq G_{\adj}\).
Note that there are natural isomorphisms
\[X_*(T) \cong X_*(T_{\adj}) \times_{\alpha,\pi_1(G_{\adj})} \pi_1(G),\]
identifying the subsets of (strictly) dominant cocharacters. 
Here \(\pi_1(-)\) denotes Borovoi's fundamental group, defined as the quotient of the cocharacter lattice by the coroot lattice; in particular there are natural maps \(X_*(T) \to \pi_1(G)\) and \(\alpha\colon X_*(T_{\adj}) \to \pi_1(G_{\adj})\).
Recall also that there is a canonical isomorphism \(\pi_1(G) \cong \pi_0(\Gr_G)\).
Now, the effect of conjugating by \(\hat{\rho}(t^{-1})\) is that the map \(\beta \colon X_*(T_{\adj})\to \pi_1(G_{\adj})\) sending an \(\calI_{\hat{\rho}}\)-orbit in \(\Gr_{G_{\adj}}\) to the connected component containing it, is given by shifting the natural map \(\alpha\colon X_*(T_{\adj}) \to \pi_1(G_{\adj})\) by \(-\hat{\rho}\).
Thus, the \(\calI_{\hat{\rho}}\)-orbits in \(\Gr_G\) are naturally indexed by
\[X_*(T_{\adj}) \times_{\beta,\pi_1(G_{\adj})} \pi_1(G) \cong X_*(T),\]
and the subset of orbits splitting up into two \(\calU_{\rho}^{\exp}\)-orbits is indexed by
\[X_*(T_{\adj})^{++} \times_{\beta,\pi_1(G_{\adj})} \pi_1(G) = X_*(T_{\adj})^+ \times_{\alpha,\pi_1(G_{\adj})} \pi_1(G) = X_*(T)^+,\]
as desired.
\xpf

\section{Exponential motives on the affine Grassmannian}\label{Testsection}

As before, we fix a split reductive group \(G\) over \(S\), which we assume is not a torus.

As an application of the machinery developed in the previous sections, we now provide a Langlands dual description of the category of exponential Tate motives on the affine Grassmannian, which can be viewed as a motivic Casselman--Shalika equivalence (compare \cite{FrenkelGaitsgoryVilonen:Whittaker,ABBGM:Modules,BGMRR:IwahoriWhittaker}).
We will prove this by first considering reduced motives, for which we will show that \(\MTM_{\red,\rhoexp}(\Gr)\) and \(\MTM_{\red,L^+G}(\Gr)\) are equivalent, following the methods of \cite{BGMRR:IwahoriWhittaker}.
In particular, this is an equivalence of abelian categories.
Motivated by the categorical Langlands philosophy, we will then pass to stable \(\infty\)-categories.
There, we can \emph{unreduce} the reduced category \(\DTM_{\red,\rhoexp}(\Gr)\), to obtain \(\DTM_{\rhoexp}(\Gr)\), cf.~\thref{Motivic CasselmanShalika} for a precise statement.
We note that this procedure of unreducing reduced motives is not formal, rather it crucially uses the fact that the Langlands dual group controlling the motivic Satake category from \cite{CassvdHScholbach:MotivicSatake} is reduced.
In particular, this gives a partial answer to a question asked in \cite[§1.6.1]{EberhardtScholbach:Integral}.

Throughout this section, we will need to distinguish regular motives \(\DM\) from reduced motives \(\DMr\) \cite{EberhardtScholbach:Integral}.
When we can treat both motivic theories uniformly, we will write \(\DMrx\), and similar for the subcategories of stratified (mixed) Tate motives.
In order to get an equivalence for all groups and independent of any choices (in contrast to \cite{BGMRR:IwahoriWhittaker}), we will work with \(\hat{\rho}\)-twisted exponential motives, as in \thref{defi twisted exponentials}.
For this purpose, let us introduce some notation for the (co)standard functors in this context.

Recall from \thref{twisted exponential orbits in Gr} that the \(\calI_{\hat{\rho}}\)-orbits in \(\Gr\) that split into two \(\calU_{\hat{\rho}}^{\exp}\)-orbits are indexed by \(X_*(T)^+\).
More precisely, they are the \(\calI_{\hat{\rho}}\)-orbits through \(\mu(t)\) for \(\mu\in X_*(T)^+\), and we will denote them by \(\Fl_\mu^{\rhocirc}(\bba_0,\bbf_0)\).
Similarly, we will denote the \(\calU_{\hat{\rho}}^{\exp}\)-orbit in \(\Gr\) through \(\mu(t)\) by \(\Fl_\mu^{\rhoexp}(\bbf_0)\).
The following lemma follows from \thref{twisted exponential orbits in Gr} and the usual structure of orbits in affine Grassmannians.

\lemm\thlabel{lemm-dimensions of twisted exponential orbits}
Let \(\mu\in X_*(T)^+\).
Then there are isomorphisms \(\Fl_\mu^{\rhocirc}(\bba_0,\bbf_0) \cong \A^{\langle 2\rho,\mu+\hat{\rho}\rangle}\) and \(\Fl_\mu^{\rhoexp}(\bbf_0) \cong \A^{\langle 2\rho,\mu+\hat{\rho}\rangle-1}\).
\xlemm

For \(\mu\in X_*(T)^+\), we denote the standard and costandard functors associated with the corresponding stratum, defined similarly as in \thref{defi-exponential co/standard}, by 
\[\Delrhoexp_{\mu}\colon \DTMrx(S) \to \DTM_{\redx,\rhoexp}(\Gr)\]
and
\[\nabrhoexp_{\mu}\colon \DTMrx(S) \to \DTM_{\redx,\rhoexp}(\Gr).\]
For example, if \(G\) is adjoint, the isomorphism \(\Gr\cong \Gr\colon x\mapsto \hat{\rho}(t^{-1})x\) induces an equivalence \(\DTM_{\redx,\exp}(\Gr) \cong \DTM_{\redx,\rhoexp}(\Gr)\), under which \(\Delexp_{\mu+\hat{\rho}}(-)\) corresponds to \(\Delrhoexp_{\mu}(-)\).

By \thref{exponential t-exactness} and \thref{reduction to adjoint case}, \(\Delrhoexp_{\mu}(-)\) and \(\nabrhoexp_{\mu}(-)\) are t-exact functors.
Similarly, by \thref{exponential convolution} and \thref{reduction to adjoint case}, (left) convolution with \(\Delrhoexp_{\mu}(\Z)\) is left t-exact, whereas (left) convolution with \(\nabrhoexp_{\mu}(\Z)\) is right t-exact.

Since \(\Fl_0^{\rhoexp}(\bba_0,\bbf_0)\subseteq \Gr\) is a minimal \(\calI_{\hat{\rho}}\)-orbit splitting up into two \(\calU_{\hat{\rho}}^{\exp}\)-orbits by \thref{twisted exponential orbits in Gr}, the following lemma is immediate.
\lemm
The natural transformation \(\Delrhoexp_{0}(-) \to \nabrhoexp_{0}(-)\) is an isomorphism of functors.
In particular, left convolution with \(\Delrhoexp_{0}(\Z)\cong \nabrhoexp_{0}(\Z)\) is a t-exact functor \(\DTM_{\redx,L^+G}(\Gr) \to \DTM_{\redx,\rhoexp}(\Gr)\).
\xlemm

We denote the resulting t-exact functor by
\[\Phi:= \Delrhoexp_{0}(\Z) \star (-) \colon \DTM_{\redx,L^+G}(\Gr) \to \DTM_{\redx,\rhoexp}(\Gr),\]
and its restriction to the hearts by
\[\Phi^{\heartsuit} \colon \MTM_{\redx,L^+G}(\Gr) \to \MTM_{\redx,\rhoexp}(\Gr).\] 

One of the goals of this section is to prove that, for reduced motives, \(\Phi^{\heartsuit}\) is an equivalence of categories.
Before we prove this, we need some preliminary results on standard and costandard objects in the Satake category \(\MTM_{\redx,L^+G}(\Gr)\).

\subsection{Complements on (co)standard objects in the Satake category}
We now extend the study of the motivic Satake category from \cite{CassvdHScholbach:MotivicSatake}.
For a dominant cocharacter \(\mu\in X_*(T)^+\), recall from \refsect{flagnotation} that \(\Gr^\mu\) denotes the corresponding \(L^+G\)-orbit in \(\Gr\), and let \(p_\mu\colon L^+G\backslash \Gr^\mu\to S\) and \(\iota_\mu^{\bbf_0}\colon L^+G\backslash \Gr^\mu\to L^+G\backslash \Gr\) be the natural maps.
Recall the (truncated) standard and costandard functors from \cite[Definition 5.40]{CassvdHScholbach:MotivicSatake}:
\[\calJ^\mu_!:=\pe \iota_{\mu,!}^{\bbf_0}(p_\mu^*(-)[\langle 2\rho,\mu\rangle])\colon \MTMrx(S)\to \MTM_{\redx,L^+G}(\Gr),\]
\[\calJ^\mu_*:=\pe \iota_{\mu,*}^{\bbf_0}(p_\mu^*(-)[\langle 2\rho,\mu\rangle])\colon \MTMrx(S)\to \MTM_{\redx,L^+G}(\Gr).\]

Recall also that for a prestack \(f\colon X\to S\), we denote the Verdier duality functor by \(\Du\colon \DMrx(X)^\opp \to \DMrx(X) := \IHom(-,f^!\Z)\) (but we emphasize that it is not necessarily an anti-equivalence, even on compact objects).

\lemm\thlabel{duality co/standards}
Let \(\mu\in X_*(T)^+\).
Then there are canonical isomorphisms 
\[\Du(\calJ_!^\mu(\Z)) \cong \calJ_*^\mu(\Z(\langle 2\rho,\mu\rangle)) \quad \text{ and } \quad \Du(\calJ_*^\mu(\Z)) \cong \calJ_!^\mu(\Z(\langle 2\rho,\mu\rangle))\]
in \(\MTM_{\redx,L^+G}(\Gr)\).
\xlemm
\pf
It suffices to show the first isomorphism, the second will follow from \cite[Proposition 5.54]{CassvdHScholbach:MotivicSatake}.
First, we claim that \(\Du(\calJ_!^\mu(\Z))\) lies in \(\MTM_{\redx,L^+G}(\Gr)\).
Indeed, by a similar argument as \cite[Proposition IV.6.13]{FarguesScholze:Geometrization}, \(\Du\) is compatible with constant term functors, so we can reduce to the case of a torus by \cite[Proposition 5.10]{CassvdHScholbach:MotivicSatake}.
In that case, the Schubert cells are just copies of \(S\) (after taking reduced subschemes), so the claim is clear.

Next, since \(\Du(\calJ_!^\mu(\Z))\) lies in degree 0, the canonical map
\[\Du(\calJ_!^\mu(\Z)) \to \Du(\iota_{\mu,!}^{\bbf_0}(\Z[\langle 2\rho,\mu\rangle])) \cong \iota_{\mu,*}^{\bbf_0}(\Z(\langle 2\rho,\mu\rangle)[\langle 2\rho,\mu\rangle])\]
(where the isomorphism follows from \cite[Theorem 2.4.50 (5)]{CisinskiDeglise:Triangulated}) factors through a map \(\Du(\calJ_!^\mu(\Z)) \to \calJ_*^\mu(\Z)\).
To check it is an isomorphism, we may apply the conservative fiber functor \(F\) \cite[Remark 5.35]{CassvdHScholbach:MotivicSatake}. 
After applying \(F\), both terms are direct sums of Tate twists, indexed by the top-dimensional irreducible components of Mirkovic--Vilonen cycles \cite[Lemma 5.42 and Remark 5.44]{CassvdHScholbach:MotivicSatake}, so that the natural map is indeed an isomorphism.
\xpf

Now, let \(\hat{G}\) be the Langlands dual group of \(G\), with maximal torus and Borels \(\hat{T} \subseteq \hat{B}\subseteq \hat{G}\), defined over \(\Z\).
These groups admit a natural \(\Gm\)-action via $\rho$ \cite[(6.11)]{CassvdHScholbach:MotivicSatake}, and the motivic Satake equivalence from \cite[Theorem 1.3]{CassvdHScholbach:MotivicSatake} yields an equivalence
\[\MTM_{\redx,L^+G}(\Gr) \cong \Rep_{\hat{G}}(\MTMrx(S));\]
we refer to \cite{CassvdHScholbach:MotivicSatake} for the precise meaning of the right hand side.
In particular, for reduced motives, we get
\begin{equation}\label{reduced Satake}
	\MTM_{\red,L^+G}(\Gr)\cong \Rep_{\hat{G}\rtimes \Gm}(\Ab).
\end{equation}

Let us give an explicit description for the \(\hat{G}\rtimes \Gm\)-representations corresponding to \(\calJ_!^{\mu}(\Z)\) and \(\calJ_*^{\mu}(\Z)\) in \(\MTM_{\red,L^+G}(\Gr)\).
The reductive group \(\hat{G}\rtimes \Gm\) contains \(\hat{T}\times \Gm\) as a maximal torus, with character group \(X^*(\hat{T})\times \Z =X_*(T)\times \Z\).
For any \((\mu,n)\in X_*(T)\times \Z\) with \(\mu\in X_*(T)^+\) dominant, we can view it as a character of \(\hat{B}\rtimes \Gm\) via
\[(\mu,n)\colon \hat{B}\rtimes \Gm \twoheadrightarrow \hat{T}\times \Gm \xrightarrow{(\mu,n)}\Gm.\]
Inducing this to a representation of \(\hat{G}\rtimes \Gm\) then yields the \emph{Schur module} 
\[S(\mu,n):=\operatorname{Ind}_{\hat{B}\rtimes \Gm}^{\hat{G}\rtimes \Gm}(\mu,n).\]
Dually, we have the \emph{Weyl module}
\[W(\mu,n):=S(-w_0(\mu),-n)^{\vee},\]
where \(w_0\) is the longest element in the finite Weyl group of \(G\).
Then the same proof as \cite[Proposition 13.1]{MirkovicVilonen:Geometric} shows the following:

\prop\thlabel{Weyl and Schur}
Let \(\mu\in X_*(T)^+\) and \(n\in \Z\).
Then under the equivalence \eqref{reduced Satake}, the following objects correspond:
\[\calJ_!^{\mu}(\Z(n)) \leftrightarrow W(\mu,n)\]
and
\[\calJ_*^{\mu}(\Z(n)) \leftrightarrow S(\mu,n).\]
\xprop

\rema\thlabel{reducedness of co/standards}
A priori, the above proposition only describes the images of \(\calJ_*^{\mu}(\Z(n))\) and \(\calJ_!^{\mu}(\Z(n))\) under the Satake equivalence for reduced motives.
However, by \cite[Lemma 5.42, Remark 5.44]{CassvdHScholbach:MotivicSatake}, the corresponding objects in the nonreduced category \(\MTM_{L^+G}(\Gr)\) are already reduced, in the sense that under the functor
\begin{equation}\label{reducedness of Satake category}
	\MTM_{\red,L^+G}(\Gr) \to \MTM_{L^+G}(\Gr)
\end{equation}
from \cite[Corollary 6.16]{CassvdHScholbach:MotivicSatake}, the reduced (co)standard objects are mapped to the corresponding (co)standard objects in \(\MTM_{L^+G}(\Gr)\).
We note that the existence of \eqref{reducedness of Satake category} crucially uses that the dual group \(\hat{G}\) is reduced \cite[Theorem 6.15]{CassvdHScholbach:MotivicSatake}.
\xrema

\subsection{Dimensions of fibers of exponential convolution morphisms}

A key ingredient in \cite{BGMRR:IwahoriWhittaker} is the geometric Casselman--Shalika formula \cite{FrenkelGaitsgoryVilonen:Whittaker,NgoPolo:Resolutions}, which is used to show that \(\Phi^{\heartsuit}\) maps (co)standard objects in the Satake category to (co)standard objects in the Iwahori--Whittaker category.
Concretely, this formula describes certain compactly supported \(\ell\)-adic cohomology groups (over the base scheme \(S=\Fqq\)) as
\[\H^i_c(\mathcal{S}_\mu^+,\sigma_\mu^*(\IC_{\lambda,\Qlq}) \otimes_{\Qlq} \chi_\mu^*(\mathcal{L}_\psi)) \cong \begin{cases}
	\Qlq & \text{ if } \lambda=\mu \text{ and } i=\langle2\rho,\mu\rangle\\
	0 & \text{ else.}
\end{cases}\]
Here, \(\lambda,\mu\in X_*(T)^+\) are dominant cocharacters, \(\mathcal{S}_\mu^+=LU\cdot \mu(t) L^+G/L^+G \xrightarrow{\sigma_\mu} \Gr\) is a semi-infinite orbit, \(\IC_{\lambda,\Qlq}\) is the intersection complex associated with the Schubert cell \(\Gr^\lambda\), the sheaf \(\mathcal{L}_\psi\) is the Artin--Schreier local system on \(\Ga\) associated with a character \(\psi\colon \Fq\to \Qlq^\times\), and \(\chi_\mu\colon \mathcal{S}_\mu^+\to \Ga\) is a natural map \cite[§3.4]{BGMRR:IwahoriWhittaker}.
In particular, this formula requires the choice of an Artin--Schreier local system, and since \(\mathcal{S}_\mu^+\) does not admit a natural \(\mathcal{I}\)-action, it is not clear how to formulate an analogue for exponential motives over a general base scheme \(S\).
The goal of this subsection is to find an alternative geometric property, which is equivalent to the geometric Casselman--Shalika formula when the latter makes sense, and which would similarly allow us to show \(\Phi^{\heartsuit}\) is compatible with (co)standard functors.
We will then deduce this property indeed holds from the results of \cite{BGMRR:IwahoriWhittaker}, and hence indirectly from the geometric Casselman--Shalika formula.

Throughout this subsection only, we will assume that \(G=G_{\adj}\) is adjoint, and work with usual (i.e., untwisted) exponential motives.
In particular, \(\hat{\rho}\in X_*(T)^+\) is an actual cocharacter.
By \thref{reduction to adjoint case}, the results below will also hold for \(\hat{\rho}\)-twisted exponential motives and for general \(G\), but this will significantly lighten the notation.
To further lighten the notation, we will denote the \(\mathcal{I}\)-orbit through \(\lambda(t)\), for \(\lambda\in X_*(T)^+\) a dominant cocharacter, by \(\Fl_\lambda^{\circ}(\bba_0,\bbf_0):= \Fl_{t_\lambda}^\circ(\bba_0,\bbf_0)\), and similarly for exponential orbits.
We will also still write \(\Phi(-) := \Delexp_{\hat{\rho}}(\Z)\star(-)\colon \DTM_{\redx,L^+G}(\Gr) \to \DTM_{\redx,\exp}(\Gr)\), and denote its restriction to the hearts by \(\Phi^{\heartsuit}\).
Finally, recall from \thref{notation parahoric inclusions} that we will denote the inclusion of \(L^+G\)-orbits by \(\iota_\mu^{\bbf_0}\colon \Gr^\mu\to \Gr\), whereas the \(\iota_{\lambda+\hat{\rho}}\), i.e., without superscript, denote the inclusion of exponential orbits.

\prop\thlabel{Equivalent conditions}
Let \(\lambda\neq \mu\in X_*(T)^+\), and consider the convolution morphism 
\[m_{\hat{\rho},\mu}\colon \Fl_{\hat{\rho}}^{\exp}(\bbf_0)\widetilde{\times} \Gr^\mu \to \Gr.\] 
Then the following statements are equivalent:
\begin{enumerate}
	\item\label{fibral dimension} The fiber \(m^{-1}_{\hat{\rho},\mu}(t^{\lambda+\hat{\rho}})\) has relative dimension \(<\langle\rho,\mu-\lambda \rangle\) (where the empty scheme has dimension \(-\infty\)).
	\item\label{degrees of convolution} For any \(M\in \MTMrx(S)\), the motive \[\iota_{\lambda+\hat{\rho}}^*(\Delexp_{\hat{\rho}}(\Z)\star \iota_{\mu,!}^{\bbf_0}(p_\mu^*(M)[\langle 2\rho,\mu\rangle]))\in \DTM_{\redx,\exp}(\Fl_{\lambda+\hat{\rho}}^{\circ}(\bba_0,\bbf_0))\]
	lies in degrees \(\leq -1\).
	\item\label{vanishing of all homs} For any \(M,N\in \MTMrx(S)\), we have \[\Hom_{\MTM_{\redx,\exp}(\Gr)} (\Phi^{\heartsuit}(\calJ_!^{\mu}(M)),\nabexp_{\lambda+\hat{\rho}}(N))=0.\]
\end{enumerate}
\xprop
The semismallness of the convolution morphism \cite[Lemma 4.4]{MirkovicVilonen:Geometric} implies that \(\dim_S m^{-1}_{\hat{\rho},\mu}(t^{\lambda+\hat{\rho}})\leq \langle\rho,\mu-\lambda \rangle\).
The content of the proposition will be to ensure that this equality cannot be achieved.
\pf
Let \(M,N\in \MTMrx(S)\).
When \(\lambda +\hat{\rho} \nleq \mu+\hat{\rho}\), the image of \(m_{\hat{\rho},\mu}\) intersects \(\Fl_{\lambda+\hat{\rho}}^{\circ}(\bba_0,\bbf_0)\) trivially, so all statements vacuously hold.
We may thus assume that \(\lambda<\mu\).

Since \(\Phi\) is t-exact and \(\iota_{\mu,!}^{\bbf_0}\colon \DTM_{\redx,L^+G}(\Gr^\mu)\to \DTM_{\redx,L^+G}(\Gr)\) is right t-exact, there are natural isomorphisms
\[\Hom_{\MTM_{\redx,\exp}(\Gr)}(\Phi^{\heartsuit}(\calJ_!^\mu(M)),\nabexp_{\lambda+\hat{\rho}}(N)) \cong \Hom_{\DTM_{\redx,\exp}(\Gr)}(\Phi(\iota_{\mu,!}^{\bbf_0}(p_\mu^*(M[\langle 2\rho,\mu\rangle]))),\nabexp_{\lambda+\hat{\rho}}(N))\]
\[\cong \Hom_{\DTM_{\redx,\exp}(\Fl_{\lambda+\hat{\rho}}^{\circ}(\bba_0,\bbf_0))}(\iota_{\lambda+\hat{\rho}}^*(\Delexp_{\hat{\rho}}(\Z)\star \iota_{\mu,!}^{\bbf_0}(p_\mu^*(M)[\langle 2\rho,\mu\rangle])),p^*_{\lambda+\hat{\rho}}(N[\langle 2\rho,\lambda+\hat{\rho}\rangle])),\]
using \thref{lemm-dimensions of twisted exponential orbits}.
Since \(p^*_{\lambda}(N[\langle 2\rho,\lambda+\hat{\rho}\rangle])\in \DTM_{\redx,\exp}(\Fl_{\lambda+\hat{\rho}}^{\circ}(\bba_0,\bbf_0))\) is concentrated in degree 0, we see that \(\eqref{degrees of convolution} \Rightarrow \eqref{vanishing of all homs}\).

Recall that we have a t-exact equivalence \(\DTM_{\redx,\exp}(\Fl_{\lambda+\hat{\rho}}^{\circ}(\bba_0,\bbf_0))\cong \DTMrx(S)\), \thref{exponential motives on iwahori orbits}.
Thus, if we knew that \(\iota_{\lambda+\hat{\rho}}^*(\Delexp_{\hat{\rho}}(\Z)\star \iota_{\mu,!}^{\bbf_0}(p_\mu^*(M)[\langle 2\rho,\mu\rangle]))\) was concentrated in degrees \(\leq 0\), we would also have the converse implication \(\eqref{vanishing of all homs}\Rightarrow \eqref{degrees of convolution}\).
This will be shown below.

Now, consider the diagram
\begin{equation}\label{diagram for dimension of fibers}\begin{tikzcd}[column sep=small]
		\Fl_{\hat{\rho}}^{\exp}(\bbf_0) \times (\Fl_{\lambda+\hat{\rho}}^{\circ}(\bba_0,\bbf_0)\cap t^{\hat{\rho}}\Gr^\mu) \arrow[d] \arrow[r] & \Fl_{\lambda+\hat{\rho}}^{\circ}(\bba_0,\bbf_0) \arrow[d]\\
		\Fl_{\hat{\rho}}^{\exp}(\bbf_0) \widetilde{\times} \Gr^\mu \arrow[r] & \Gr.
\end{tikzcd}\end{equation}
Cartesianness follows by considering the intersection \(\Fl_{\lambda+\hat{\rho}}^{\circ}(\bba_0,\bbf_0)\cap t^{\hat{\rho}}\Gr^\mu\), and taking the product with the subscheme \(\calU^{\exp} \supseteq H_{\hat{\rho}} \cong \Fl_{\hat{\rho}}^{\exp}(\bbf_0)\) from \thref{orbit-structure} \eqref{it--global section}. 

By \cite[Lemma 7.3]{Haines:Pavings}, $\Fl_{\lambda+\hat{\rho}}^{\circ}(\bba_0,\bbf_0)\cap t^{\hat{\rho}}\Gr^\mu$ is a subscheme of the fiber of the convolution morphism $\Gr^{\lambda + \hat{\rho}} \widetilde{\times} \Gr^{-w_0(\mu)} \r \Gr$ over $t^{\hat{\rho}}$. By semismallness of this convolution morphism \cite[Lemma 4.4]{MirkovicVilonen:Geometric}, the top left scheme in \eqref{diagram for dimension of fibers} has dimension at most 
\begin{equation}\label{dimension estimate}
	 \langle 2\rho,\hat{\rho}\rangle -1 + \langle \rho,\mu+\lambda\rangle = \langle\rho,\mu+\lambda+2\hat{\rho} \rangle-1.
\end{equation}

Now, let \(z\in {}_0W_{\mathbf{f}}\subseteq W_{\bbf_0}^{\exp}\) be such that \(\Fl_{\lambda+\hat{\rho}}^\circ(\bba_0,\bbf_0) = \Fl_{\lambda+\hat{\rho}}^{\exp}(\bbf_0)\sqcup \Fl_z^{\exp}(\bbf_0)\), with \(\Fl_z^{\exp}(\bbf_0)\subseteq \Fl_{\lambda+\hat{\rho}}^{\circ}(\bba_0,\bbf_0)\) open.
We will denote its basepoint by \(t^z:=x_{\alpha_0}(1)\cdot t^{\lambda+\hat{\rho}}\in \Fl_z^{\exp}(\bbf_0)\).
Then we can refine the above diagram to get two cartesian diagrams 
\[\begin{tikzcd}
	\Fl_{\lambda+\hat{\rho}}^{\exp}(\bbf_0) \times m_{\hat{\rho},\mu}^{-1}(t^{\lambda+\hat{\rho}}) \arrow[d] \arrow[r, "\widetilde{m}_{\lambda+\hat{\rho}}"] & \Fl_{\lambda+\hat{\rho}}^{\exp}(\bbf_0) \arrow[d]\\
	\Fl_{\hat{\rho}}^{\exp}(\bbf_0) \widetilde{\times} \Gr^\mu\arrow[r, "m_{\hat{\rho},\mu}"'] & \Gr,
\end{tikzcd} 
\quad \text{and} \quad
\begin{tikzcd}
	\Fl_z^{\exp}(\bbf_0) \times m_{\hat{\rho},\mu}^{-1}(t^z) \arrow[d] \arrow[r, "\widetilde{m}_z"] & \Fl_z^{\exp}(\bbf_0) \arrow[d]\\
	\Fl_{\hat{\rho}}^{\exp}(\bbf_0) \widetilde{\times} \Gr^\mu\arrow[r, "m_{\hat{\rho},\mu}"'] & \Gr,
\end{tikzcd}\]
using \thref{triviality of convolution}.
These fiber products admit a cellular filtrable decomposition by \thref{fibers of exponential convolution}. 

Below we will use \cite[Lemma 2.20]{CassvdHScholbach:MotivicSatake}, which is stated there for schemes with a filtrable decomposition by products of $\A^1$ and $\Gm$, but same proof works if we also allow $\Gm \setminus \{1\}$. In particular, for a scheme $X \r S$ as in \thref{defi--filtrable} with $\dim_S X = d$, the functor $f_!f^*$ maps $\DTMrx(S)^{\leq 0}$ to $\DTMrx(S)^{\leq 2d}$, and $\pe H^{2d}(f_!f^* \Z)$ is a direct sum of copies of $\Z(-d)$ indexed by the top-dimensional cells of $X$.

We will now do some computations in the category \(\DTM_{\redx,\exp}(\Gr)\).
More precisely, we will write down objects in the category \(\DTMrx(\calU^{\exp}\backslash \Gr)\), and implicitly consider their images under the t-exact quotient functor to \(\DTM_{\redx,\exp}(\Gr)\).
In particular, we can view \(\Delexp_{\hat{\rho}}(\Z)\) as the !-pushforward of \(\Z[\langle 2\rho,\hat{\rho}\rangle-1]\in \DTMrx(\Fl_{\hat{\rho}}^{\exp}(\bbf_0))\) to \(\Gr\); cf.~the proof of \thref{exponential motives on iwahori orbits}.
To further simplify the notation, we will denote by \(\underline{M}\) the *-pullback of \(M\) to any scheme via the structure map to \(S\). 

Then we want to compute
\[\iota_{\lambda+\hat{\rho}}^*(\Delexp_{\hat{\rho}}(\Z)\star \iota_{\mu,!}^{\bbf_0}(p_\mu^*(M)[\langle 2\rho,\mu\rangle])) \cong \iota_{\lambda+\hat{\rho}}^* m_{\hat{\rho},\mu,!}(\underline{M}[\langle2\rho,\mu+\hat{\rho} \rangle-1])\in \DTM_{\redx,\exp}(\Fl_{\lambda+\hat{\rho}}^\circ(\bba_0,\bbf_0)).\]
This is an extension of the !-pushforwards of 
\[\widetilde{m}_{\lambda+\hat{\rho},!}\underline{M}[\langle2\rho,\mu+\hat{\rho} \rangle-1] \text{ and } \widetilde{m}_{z,!}\underline{M}[\langle2\rho,\mu+\hat{\rho} \rangle-1]\]
to \(\DTM_{\redx,\exp}(\Fl_{\lambda+\hat{\rho}}^{\circ}(\bba_0,\bbf_0))\) and we will consider both terms separately.

First, by \eqref{dimension estimate}, \(\widetilde{m}_z\) is of relative dimension
$$d_z\leq \langle\rho,\mu+\lambda+2\hat{\rho} \rangle-1 - \langle 2\rho,\lambda+\hat{\rho}\rangle= \langle \rho,\mu-\lambda\rangle-1.$$
Thus, by \cite[Lemma 2.20]{CassvdHScholbach:MotivicSatake} we see that \(\widetilde{m}_{z,!}\underline{M}[\langle2\rho,\mu+\hat{\rho} \rangle-1]\) lies in degree at most
\[1-\langle2\rho,\mu+\hat{\rho}\rangle + 2(\langle \rho,\mu-\lambda\rangle-1) + \langle 2\rho,\lambda+\hat{\rho} \rangle = -1\] in \(\DTM_{\redx,\exp}(\Fl_{\lambda+\hat{\rho}}^\circ(\bba_0,\bbf_0))\).

On the other hand, \(\widetilde{m}_{\lambda+\hat{\rho}}\) is of relative dimension $d_{\lambda+\hat{\rho}}\leq \langle \rho,\mu-\lambda\rangle$, again by \eqref{dimension estimate}.
Thus, \(\widetilde{m}_{\lambda+\hat{\rho},!}\underline{M}[\langle2\rho,\mu+\hat{\rho} \rangle-1]\) lies in degree at most $0$,
with equality achieved for $M=\Z$ exactly when \(d_{\lambda+\hat{\rho}} = \langle \rho,\mu-\lambda\rangle\), i.e., when \(\dim_S m_{\hat{\rho},\mu}^{-1}(t^{\lambda+\hat{\rho}}) = \langle \rho,\mu-\lambda\rangle.\)
Since the degree 0 part arising from \(\widetilde{m}_{\lambda+\hat{\rho},!}\underline{M}[\langle2\rho,\mu+\hat{\rho} \rangle-1] \) cannot be canceled out by \(\widetilde{m}_{z,!}\underline{M}[\langle2\rho,\mu+\hat{\rho} \rangle-1]\), we get the desired equivalence \(\eqref{fibral dimension}\Leftrightarrow \eqref{degrees of convolution}\).
Note that the above paragraph, together with the argument in the beginning of the proof, also shows \(\eqref{vanishing of all homs} \Rightarrow \eqref{degrees of convolution}\).
\xpf

\prop\thlabel{Equivalences hold}
The equivalent conditions from \thref{Equivalent conditions} hold true.
\xprop
\pf
Let \(\lambda\) and \(\mu\) be as in \thref{Equivalent conditions}.

The proof of \thref{Equivalent conditions} also holds when working with étale sheaves over \(\Fqq\) with coefficients in \(\Qlq\).
In that case, we can use the equivalence of exponential sheaves with the Whittaker model, \thref{Whit vs Kir}, to see that, in the étale setting, condition \eqref{vanishing of all homs} of \thref{Equivalent conditions} was shown in \cite[Proposition 3.11]{BGMRR:IwahoriWhittaker}, by using the geometric Casselman--Shalika formula. 

Thus condition \eqref{fibral dimension} holds for \(S=\Fqq\).
But by \thref{fibers of exponential convolution}, the fiber \(m_{\hat{\rho},\mu}^{-1}(t^{\lambda+\hat{\rho}})\) is cellular over \(\Spec \Z\), so that its relative dimension is independent of the base.
This shows that condition \eqref{fibral dimension} holds in general, and then the other conditions follow from \thref{Equivalent conditions}.
\xpf

\subsection{The motivic Casselman--Shalika equivalence}

Finally, we can describe the category \(\DTM_{\redx,\rhoexp}(\Gr)\) in terms of Langlands dual information.
We now switch back to using \(\hat{\rho}\)-twisted exponential motives.

First, we show that rationally, and up to motivic cohomology coming from the base scheme, the category \(\MTM_{\redx,\rhoexp}(\Gr,\Q)\) is semisimple.

\prop\thlabel{rationalexp}
For \(\mu\in X_*(T)^+\), the natural map
\[\Delrhoexp_{\mu}(\Q) \to \nabrhoexp_{\mu}(\Q)\]
is an isomorphism.
In particular, these are simple objects of \(\MTM_{\redx,\rhoexp}(\Gr,\Q)\).
\xprop
\pf
Recall that the reduction functor \(\MTM_{\rhoexp}(\Gr,\Q) \to \MTM_{\red,\rhoexp}(\Gr,\Q)\) is exact and conservative, \thref{conservativity}.
Since the resulting category \(\MTM_{\red,\rhoexp}(\Gr,\Q)\) is independent of \(S\) \cite[Proposition 4.25]{EberhardtScholbach:Integral}, we may assume \(S=\Fq\) is a finite field.
In that case, regular rational motives agree with reduced rational motives by \cite[Proposition 5.3]{EberhardtScholbach:Integral}.

Now, fix a prime \(\ell\neq p\). 
The construction of the exponential category works for \(\ell\)-adic sheaves as well, and repeating \thref{defi twisted exponentials} yields a category of (twisted) exponential \(\ell\)-adic sheaves \(D_{\rhoexp}(\Gr,\Ql)\), together with a t-structure with heart \(\Perv_{\rhoexp}(\Gr,\Ql)\), and a realization functor \(\DTM_{\rhoexp}(\Gr,\Q)\to D_{\rhoexp}(\Gr,\Ql)\).
By \cite[Lemma 3.2.8]{RicharzScholbach:Intersection}, this realization functor is t-exact, and it is conservative when restricted to \(\MTM_{\rhoexp}(\Gr,\Q)\) by \thref{conservativity} (since the property of an object being trivial can be checked after pullback to the strata).
We may also enlarge our coefficients from \(\Ql\) to \(\Qlq\), and fix a non-trivial character \(\chi \colon \Fp\to \Qlq^\times\).
In that situation, we may work with the (Iwahori--)Whittaker category \(\Perv_{IW}(\Gr,\Qlq)\) attached to \(\chi\) by the \(\ell\)-adic analogue of \thref{Whit vs Kir}.
In that case, the proposition follows from \cite[(3.3)]{BGMRR:IwahoriWhittaker}.
\xpf

Next, we consider the t-exact functor \(\Phi^{\heartsuit}:=\Delrhoexp_{0}(\Z)\star (-)\).

\prop\thlabel{averaging of standard}
Let \(\mu\in X_*(T)^+\) and \(\calF\in \MTMrx(S)\).
Then there is a natural isomorphism 
\[\Delrhoexp_{0}(\Z) \star \calJ_!^\mu(\calF) \cong \Delrhoexp_{\mu}(\calF)\]
in \(\MTM_{\redx,\rhoexp}(\Gr)\).
\xprop
\pf
By \thref{reduction to adjoint case}, we may assume that \(G\) is adjoint, so that \(\hat{\rho}\in X_*(T)^+\) defines a cocharacter, and we may work with usual (untwisted) exponential motives.
Since the functors \(\calJ_!^\mu\) and \(\Delexp_{\mu+\hat{\rho}}\) are \(\MTMrx(S)\)-linear (by \cite[Proposition 5.41]{CassvdHScholbach:MotivicSatake} and the projection formula respectively), it suffices to consider \(\calF=\Z\).

Note that the convolution morphism \(\Fl_{\hat{\rho}}^\circ(\bba_0,\bbf_0) \widetilde{\times} \Gr^\mu \to \Gr\) restricts to an isomorphism over \(\Fl_{\mu+\hat{\rho}}^\circ(\bba_0,\bbf_0)\).
This yields a map
\[\alpha \colon \Delexp_{\mu+\hat{\rho}}(\Z) \to \Delexp_{\hat{\rho}}(\Z) \star \calJ_!^\mu(\Z) = \Phi^\heartsuit(\calJ_!^\mu(\Z))\]
in \(\MTM_{\redx,\exp}(\Fl_{\mu+\hat{\rho}}(\bba_0,\bbf_0))\), which restricts to an isomorphism over \(\Fl_{\mu+\hat{\rho}}^{\circ}(\bba_0,\bbf_0)\).
We claim that \(\coker(\alpha)=0\). 
Assume this is not the case, and let \(\nu\in {}_0W_{\bbf_0}\) be maximal such that \(\Fl_\nu^{\circ}(\bba_0,\bbf_0)\) lies in the support of \(\coker(\alpha)\); in particular \(\nu<\mu\).
Then \(\Phi^{\heartsuit}(\calJ_!^{\mu}(\Z))\) would admit a nonzero map to \(\iota_{\nu,*}^{\bba_0}\iota_{\nu}^{\bba_0,*}\coker(\alpha))\) (via the unit of the \((\iota_{\nu}^{\bba_0,*},\iota_{\nu,*}^{\bba_0})\)-adjunction).
This contradicts \thref{Equivalent conditions} \eqref{vanishing of all homs}, which holds by \thref{Equivalences hold}.
Thus, \(\alpha\) is surjective.

Now, by \thref{rationalexp}, \(\Delexp_{\mu+\hat{\rho}}(\Q)\) is a simple object, so that \(\alpha\) is an isomorphism after rationalizing.
In particular, the kernel of \(\alpha\) is a torsion object.
But since \(\Delexp_{\mu+\hat{\rho}}(\Z)\) is torsionfree (in the sense that the multiplication-by-\(n\) map has trivial kernel for any \(n\geq 1\)), as the image of the torsionfree object \(\Z\) (\thref{kernel n-multiplication in MTM}) under the exact functor \(\Delexp_{\mu+\hat{\rho}}(-)\), it does not have any torsion subobjects, so that \(\alpha\) has trivial kernel.
We conclude that \(\alpha\) is an isomorphism.
\xpf

\coro\thlabel{averaging of costandard}
Let \(\mu\in X_*(T)^+\), and \(\calF\in \MTMrx(S)\).
When working with regular motives, assume that \(\calF\) is flat.
Then there is a natural isomorphism
\[\nabrhoexp_{\mu}(\calF) \cong \Delrhoexp_{0}(\Z) \star \calJ_*^{\mu}(\calF).\]
\xcoro
\pf
Note that \(\Phi\) is compatible with Verdier duality on both \(\DTM_{\redx,L^+G}(\Gr)\) and \(\DTM_{\redx,\rhoexp}(\Gr)\): this follows from \cite[Theorem 2.4.50 (5)]{CisinskiDeglise:Triangulated}, and the compatibility of Verdier duality with exterior products of Tate motives (which can be proven as in \cite[Proposition 2.4.3]{JinYang:Kuenneth}, using \thref{WT.prod} as a replacement for \cite[Proposition 2.3.5]{JinYang:Kuenneth}). 
Combined with \thref{averaging of standard}, this yields isomorphisms \(\nabrhoexp_{\mu}(\Z) \cong \Delrhoexp_{0}(\Z) \star \calJ_*^{\mu}(\Z)\), so that the corollary holds for \(\calF=\Z\).

For general \(\calF\), we always have \(\nabrhoexp_{\mu}(\Z) \otimes \calF \cong \nabrhoexp_{\mu}(\calF)\) by \thref{WT.prod}.
In the Satake category, we have
\[\calJ_*^{\mu}(\Z)\otimes \calF \cong \calJ_*^{\mu}(\calF)\]
for general \(\calF\) when using reduced motives, and for flat \(\calF\) when using regular motives \cite[Remark 5.44]{CassvdHScholbach:MotivicSatake}.
Thus, the corollary follows from \(\MTMrx(S)\)-linearity of \(\Phi^{\heartsuit}\).
\xpf

We now prove a motivic refinement of the main result of \cite{BGMRR:IwahoriWhittaker} (for reduced motives).

\prop\thlabel{Equivalence in the reduced case}
For reduced motives, \(\Phi^{\heartsuit}\) induces an equivalence
\[\Phi^{\heartsuit}\colon \MTM_{\red,L^+G}(\Gr)\cong \MTM_{\red,\rhoexp}(\Gr).\]
Moreover, this equivalence identifies the anti-effective motives on both sides, and hence induces
\[\Phi^{\heartsuit}\colon \MTM_{\red,L^+G}^{\anti}(\Gr) \cong \MTM_{\red,\rhoexp}^{\anti}(\Gr).\]
\xprop
\pf

Both \(\MTM_{\red,L^+G}(\Gr)\) and \(\MTM_{\red,\rhoexp}(\Gr)\) are compactly generated \cite[Lemma 2.26]{CassvdHScholbach:MotivicSatake}.
Moreover, it follows from  the argument of \thref{convDTM} that \(\Phi^{\heartsuit}\) preserves compact objects, i.e., restricts to a functor
\[\Phi^{\heartsuit}\colon \MTM_{\red,L^+G}(\Gr)^{\comp} \to \MTM_{\red,\rhoexp}(\Gr)^{\comp}.\]
Let \(Z\subseteq X_*(T)^+\) be any finite subset, closed under the dominance order.
We identify this finite subset with a quasi-compact \(L^+G\)-stable closed subscheme \(Z\subseteq \Gr\), and let \(Z'\subseteq \Gr\) denote the closure of the image of \(\Fl_0^{\rhocirc}(\bba_0,\bbf_0) \widetilde{\times} Z \to \Gr\). 
We will show that \(\Phi^{\heartsuit}\) induces equivalences
\begin{equation}\label{bounded Casselman-Shalika}
	\MTM_{\red,L^+G}(Z)^{\comp} \cong \MTM_{\red,\rhoexp}(Z')^{\comp}
\end{equation}
for all \(Z\). 
Passing to the colimit and taking ind-completions will then imply the equivalence
\[\MTM_{\red,L^+G}(\Gr)\cong \MTM_{\red,\rhoexp}(\Gr).\]
Let us now fix such a \(Z\subseteq X_*(T)^+\) as above.

Consider the exact derived functor
\[\D(\Phi^{\heartsuit})\colon \D(\MTM_{\red,L^+G}(Z)^{\comp}) \to \D(\MTM_{\red,\rhoexp}(Z')^{\comp}),\]
where \(\D(-)\) denotes the derived category of an abelian category.
By \thref{averaging of standard} and \thref{averaging of costandard}, \(\D(\Phi^{\heartsuit})\) sends \(\calJ_!^{\mu}(\calF)[n]\) to \(\Delrhoexp_{\mu}(\calF)[n]\) and \(\calJ_*^{\mu}(\calF)[n]\) to \(\nabrhoexp_{\mu}(\calF)[n]\), for any \(\calF\in \MTMr(S)^{\comp}\), \(\mu\in Z\), and \(n\in \Z\).
It suffices to show that \(\D(\Phi^{\heartsuit})\) is an equivalence of triangulated categories.

The sets \(\left\{\calJ_!^{\mu}(\Z(m))[n]\mid \mu\in Z \text{ and }m,n\in \Z\right\}\) and \(\left\{\calJ_*^{\mu}(\Z(m))[n]\mid \mu\in Z \text{ and }m,n\in \Z\right\}\) each generate \(\D(\MTM_{\red,L^+G}(Z)^{\comp})\).
Similarly, the sets \(\left\{\Delrhoexp_{\mu}(\Z(m))[n],\mu\in Z\text{ and }m,n\in\Z\right\}\) and \(\left\{\nabrhoexp_{\mu}(\Z(m))[n],\mu\in Z\text{ and }m,n\in\Z\right\}\) each generate \(\D(\MTM_{\red,\rhoexp}(Z')^{\comp})\).
Moreover, we clearly have
\[\Hom_{\MTM_{\red,L^+G}(Z)}(\calJ_!^{\mu}(\Z),\calJ_*^{\mu}(\Z(m))) \cong \Hom_{\MTM_{\red,\rhoexp}(Z')}(\Delrhoexp_{\mu}(\Z),\nabrhoexp_{\mu}(\Z(m))),\]
as both are \(\Z\) when \(m=0\), and vanish otherwise.
Since \(\Fl_0^{\rhocirc}(\bba_0,\bbf_0) \widetilde{\times} \Gr^\mu \to \Gr\) restricts to an isomorphism over \(\Fl_{\mu}^{\rhocirc}(\bba_0,\bbf_0)\subseteq \Gr\), we see that \(\Phi^{\heartsuit}\) induces the above isomorphism on Hom-groups for \(m=0\).
Thus, it will suffice to show that for any \(\lambda,\mu\in Z\), we have
\begin{equation}\label{ext vanishing}
	\Ext_{\MTM_{\red,L^+G}(Z)^{\comp}}^n(\calJ_!^{\mu}(\Z),\calJ_*^{\lambda}(\Z(m))) = 0 = \Ext_{\MTM_{\red,\rhoexp}(Z')^{\comp}}^n(\Delrhoexp_{\mu}(\Z),\nabrhoexp_{\lambda}(\Z(m)))
\end{equation}
if \(\lambda\neq \mu\) or \(n\neq 0\).
In \(\MTM_{\red,L^+G}(Z)^{\comp}\), this follows from \thref{Weyl and Schur} and \cite[Proposition B.4]{Jantzen:Representations}, in view of the reduced motivic Satake equivalence \eqref{reduced Satake}.
So it remains to consider the category \(\MTM_{\red,\rhoexp}(Z')^{\comp}\).
This can be handled as in \cite[Remark 3.7 (2)]{BGMRR:IwahoriWhittaker}; we recall the argument.

When working with coefficients in a field \(\mathbf{F}\), \thref{exponential homs} (for reduced motives) and the arguments of \cite[§3.2]{Beilinson:Koszulduality} show that \(\MTM_{\red,\rhoexp}(Z',\mathbf{F})^{\comp}\) is a graded highest weight category in the sense of \cite[Definition A.1]{Achar:ModularPerverseSheavesII} (cf.~\cite[Remark 7.2 (2)]{Riche:Geometric} for a remark on the terminology).
Thus, the analogue of \eqref{ext vanishing} for \(\mathbf{F}\)-coefficients follows from \cite[(A.9)]{Achar:ModularPerverseSheavesII}.
Now, for integral coefficients, the arguments of \cite[§2]{RicheSoergelWilliamson:Modular} show that \(\MTM_{\red,\rhoexp}(Z')\) has enough projectives, each of which admits a filtration by objects of the form \(\Delrhoexp_{\nu}(\Z(l))\), for \(\nu\in Z\) and \(l\in \Z\).
(We note that since our (co)standard functors preserve compact objects by \thref{rema--CompactExp}, all required extension groups are indeed finitely generated.)
Thus, we may consider the right derived functor of \(\Hom_{\MTM_{\red,\rhoexp}(Z')^{\comp}}(-,\nabrhoexp_{\lambda}(\Z(m)))\).
Let \(M(m) := \RHom_{\MTM_{\red,\rhoexp}(\Z')^{\comp}}(\Delrhoexp_{\mu}(\Z),\nabrhoexp_{\lambda}(\Z(m)))\), which is a complex of abelian groups.
Since \(\Delrhoexp_{\mu}\) and \(\nabrhoexp_{\lambda}\) are \(\DTMr(S)\)-linear, we see that for any field \(\mathbf{F}\) we have
\[\mathbf{F}\otimes M(m) \cong \RHom_{\MTM_{\red,\rhoexp}(\Z',F)^{\comp}}(\Delrhoexp_{\mu}(\mathbf{F}),\nabrhoexp_{\lambda}(\mathbf{F}(m))).\]
Indeed, the natural transformation
\[\mathbf{F}\otimes \RHom_{\MTM_{\red,\rhoexp}(\Z')^{\comp}}(-,\nabrhoexp_{\lambda}(\Z(m))) \to \RHom_{\MTM_{\red,\rhoexp}(\Z',F)^{\comp}}(-,\nabrhoexp_{\lambda}(\mathbf{F}(m)))\]
is an isomorphism, since any projective object admits a filtration by \(\Delrhoexp_{\nu}(\Z(l))\), and both \[\Hom_{\MTM_{\red,\rhoexp}(\Z')^{\comp}}(\Delrhoexp_{\nu}(\mathbf{\Z}(l)),\nabrhoexp_{\lambda}(\mathbf{\Z}(m)))\] and \[\Hom_{\MTM_{\red,\rhoexp}(\Z',F)^{\comp}}(\Delrhoexp_{\nu}(\mathbf{F}(l)),\nabrhoexp_{\lambda}(\mathbf{F}(m)))\] are free of the same rank as a \(\Z\)-, resp.~\(\mathbf{F}\)-module, and all $\Ext^1$ groups vanish.
This implies that
\[M(m)\otimes_{\Z} \mathbf{F} \cong \begin{cases}
	\mathbf{F} &\text{ if } m=0 \text{ and } \lambda=\mu,\\
	0 & \text{ else.}
\end{cases}\]
Since we already know that
\[\Hom_{\MTM_{\red,\rhoexp}(Z')^{\comp}}(\Delrhoexp_{\mu}(\Z),\nabrhoexp_{\lambda}(\Z(m))) \cong \begin{cases}
	\Z & \text{ if } m=0 \text{ and } \lambda=\mu,\\
	0 & \text{ else,}
\end{cases}\]
we can deduce \eqref{ext vanishing} from the case of field coefficients.
(We emphasize that the above argument does not work for regular motives, as \(\MTM(S)\) does not necessarily have nontrivial projective objects, even with rational coefficients.)
This gives the equivalence \eqref{bounded Casselman-Shalika}.

Finally, note that \(\D(\MTM_{\red,L^+G}^{\anti}(\Gr))\) is generated by \(\calJ_!^{\mu}(\Z(m))[n]\) for \(\mu\in X_*(T)^+\) and \(m\leq 0\).
Similarly, \(\D(\MTM_{\red,\rhoexp}^{\anti}(\Gr))\) is generated by \(\Delrhoexp_{\mu}(\Z(m))\) for the same \(\mu\) and \(m\).
Since \(\D(\Phi^{\heartsuit})\) identifies these objects, we see that \(\Phi^{\heartsuit}\) restricts to an equivalence between the anti-effective categories.
\xpf

To describe \(\DTM_{\rhoexp}(\Gr)\) for regular motives, let us recall the stable presentable category \(\IndCoh(X)\) of ind-coherent sheaves on an algebraic stack \(X\), defined as the ind-completion of the category of coherent sheaves on \(X\) \cite[§9]{Zhu:Tame}.
For an algebraic group \(H/\Spec \Z\), let us denote by \(\B(H)\) the classifying stack of \(H\), for the fpqc-topology.
The following lemma is well known; we include a proof for convenience of the reader.
\lemm\thlabel{Indcoh on classifying stacks}
Let \(H/\Spec \Z\) be a smooth affine algebraic group.
Then there is a natural equivalence
\begin{equation}\label{IndCoh as Ind of Coh}
	\IndCoh(\B(H)) \cong \Ind(\D^b(\Rep_H(\Ab^{\comp}))),
\end{equation}
where the target is the ind-completion of the bounded derived category of \(H\)-representations on finitely generated abelian groups.
\xlemm
\pf
Recall that the category of quasi-coherent sheaves \(\QCoh(-)\) on a stack is defined by descent \cite[(9.3)]{Zhu:Tame}.
Since for (noetherian) schemes, the stable category \(\Coh(-)\subseteq \QCoh(-)\) consists of those complexes concentrated in finitely many degrees whose cohomology sheaves are coherent (i.e., finitely generated), we deduce that \(\Coh(-)\) also satisfies descent.
On the other hand, perfect complexes \(\Perf(-)\) also satisfy descent \cite[Tag 09UG]{StacksProject}, so that \(\Perf(-)\) can be extended to any prestack, cf.~\cite[(9.4)]{Zhu:Tame}.
Equivalently, \(\Perf(-)\subseteq \QCoh(-)\) consists of the dualizable objects.
Now, since \(H\) is smooth, we have \(\Coh(H^{\times n}) = \Perf(H^{\times n})\) as subcategories of \(\QCoh(H^{\times n})\), for any \(n\geq 0\).
Summarizing the above, we see that 
\[\Coh(\B(H)) = \Perf(\B(H)) \subseteq \QCoh(\B(H)),\]
which consists of representations of \(H\) on perfect complexes of abelian groups.
In particular, \(\Perf(\B(H)) \cong \D^b(\Rep_H(\Ab^{\comp}))\),  
and we conclude by taking ind-completions.
\xpf

The main result of this paper, the motivic Casselman--Shalika equivalence, is the following.
Note that \(\IndCoh\left(\B(\hat{G}\rtimes \Gm)\right)\) is naturally a module under \(\IndCoh(\B(\Gm)) \cong \DTMr(S)\).
\theo\thlabel{Motivic CasselmanShalika}
There is a natural commutative diagram 
\[\begin{tikzcd}
	\IndCoh\left(\B(\hat{G}\rtimes \Gm)\right) \otimes_{\DTMr(S)} \DTM(S) \arrow[d] \arrow[r, "\cong"', "\CS"] & \DTM_{\rhoexp}(\Gr)\arrow[d]\\
	\IndCoh\left(\B(\hat{G}\rtimes \Gm)\right) \arrow[r, "\cong", "\CSr"'] & \DTM_{\red,\rhoexp}(\Gr),
\end{tikzcd}\]
where the horizontal arrows are equivalences.
\xtheo
\pf
Recall that \(\DTM_{\red,\rhoexp}(\Gr)\) is compactly generated;  its compact objects are those supported on a bounded subscheme of \(\Gr\), and are compact there.
Then \eqref{ext vanishing} implies that \(\DTM_{\red,\rhoexp}(\Gr)^{\comp} \cong \D^b(\MTM_{\red,\rhoexp}(\Gr)^{\comp})\); compare \cite[§3.2]{BGMRR:IwahoriWhittaker}.
Thus, we see that
\begin{equation}\label{derived category of reduced exponential}
	\DTM_{\red,\rhoexp}(\Gr) \cong \Ind(\D^b(\MTM_{\red,\rhoexp}(\Gr)^{\comp})).
\end{equation}
Combining this with \eqref{IndCoh as Ind of Coh}, \eqref{reduced Satake}, and \thref{Equivalence in the reduced case} yields the equivalence \(\CSr\).
Below, we will sometimes abuse notation, and consider objects of \(\DTM_{\red,\rhoexp}(\Gr)\) as objects of \(\IndCoh\left(\B(\hat{G}\rtimes \Gm)\right)\) under \(\CSr\).

We now move to the case of regular (i.e., nonreduced) motives.
Consider the functors
\[\Rep_{\hat{G}}(\MTMr(S)) \to \Rep_{\hat{G}}(\MTM(S)) \xrightarrow{\Phi^{\heartsuit}} \MTM_{\rhoexp}(\Gr),\]
where the first one crucially use the reducedness of \(\hat{G}\) \cite[Theorem 6.15, Corollary 6.16]{CassvdHScholbach:MotivicSatake}.
By \eqref{IndCoh as Ind of Coh}, this extends to a \(\DTMr(S)\)-linear functor
\[\IndCoh\left(\B(\hat{G}\rtimes \Gm)\right) \to \DTM_{\rhoexp}(\Gr).\]
Let us abbreviate $I := \IndCoh\left(\B(\hat{G}\rtimes \Gm)\right)$ to simplify the notation.
Then the adjunction
\[\FunL_{\DTMr(S)}(I,\DTM_{\rhoexp}(\Gr)) \cong \FunL_{\DTMr(S)}(I,\FunL_{\DTM(S)}(\DTM(S),\DTM_{\rhoexp}(\Gr)))\] 
\[\cong \FunL_{\DTM(S)}(I\otimes_{\DTMr(S)} \DTM(S),\DTM_{\rhoexp}(\Gr))\]
gives a functor
\[\CS\colon \IndCoh\left(\B(\hat{G}\rtimes \Gm)\right) \otimes_{\DTMr(S)} \DTM(S) \to \DTM_{\rhoexp}(\Gr).\]
This yields the commutative diagram
\[\begin{tikzcd}
	\IndCoh\left(\B(\hat{G}\rtimes \Gm)\right) \otimes_{\DTMr(S)} \DTM(S) \arrow[d] \arrow[r, "\CS"] & \DTM_{\rhoexp}(\Gr)\arrow[d]\\
	\IndCoh\left(\B(\hat{G}\rtimes \Gm)\right) \arrow[r, "\cong", "\CSr"'] & \DTM_{\red,\rhoexp}(\Gr),
\end{tikzcd}\]
where the vertical arrows are induced by the reduction functors.
It remains to show the upper arrow \(\CS\) is an equivalence.
For \(\mu\in X_*(T)^+\) and \(\calF\in \DTM(S)\), it sends \(\Delrhoexp_{\mu}(\Z) \otimes \calF\in \IndCoh\left(\B(\hat{G}\rtimes \Gm)\right) \otimes_{\DTMr(S)} \DTM(S)\) to \(\Delrhoexp_{\mu}(\calF)\in \DTM_{\rhoexp}(\Gr)\), and \(\nabrhoexp_{\mu}(\Z)\otimes \calF\) to \(\nabrhoexp_{\mu}(\calF)\); this follows from \thref{reducedness of co/standards}, as well as \thref{averaging of standard} and \thref{averaging of costandard}.
In particular, it is essentially surjective as soon as we have full faithfulness.
To see the latter, it suffices to show that for \(\lambda,\mu\in X_*(T)^+\) and \(\calF,\calG\in \DTM(S)\), the natural maps induce equivalences between the mapping spectra
\begin{equation}\label{eq-mapping spectra}
	\Maps_{I \otimes_{\DTMr(S)} \DTM(S)}(\Delrhoexp_{\mu}(\Z)\otimes \calF,\nabrhoexp_{\lambda}(\Z) \otimes \calG) \cong 
\begin{cases}
	\Maps_{\DTM(S)}(\calF,\calG) & \text{ if } \lambda=\mu\\
	0 & \text{ else,}
\end{cases}\end{equation}
by \thref{exponential homs}.
Indeed, then \(\CS\) will induce the above equivalence between mapping spectra for \(\lambda=\mu\), since the convolution map \(\Fl_0^{\rhocirc}(\bba_0,\bbf_0)\widetilde{\times} \Gr^\mu\to \Gr\) restricts to an isomorphism over \(\Fl_{\mu}^{\rhocirc}(\bba_0,\bbf_0)\).
We may assume that $\mathcal F$ is a compact object.

To see \eqref{eq-mapping spectra}, we abbreviate $A:=\DTMr(S)$ and write $\Map_{C/A}(-,-) \in A$ for the relative mapping object, for any $A$-module $C$, such as $C=I$, $C=\DTM(S)$, or $C=I \otimes_A \DTM(S)$.
The relative mapping object is characterized by $\Map_{A}(a, \Map_{C/A}(x,y)) = \Map_C(a\otimes x, y)$ for $a\in A$ and $x,y\in C$, where $\Map_A$ and $\Map_C$ denote the usual mapping spectra.
We claim
$$\Map_{I / A}(\Delrhoexp_\mu(\Z), \nabrhoexp_\mu (\Z)) = \begin{cases} \Z  & \text{ if } \lambda = \mu,\\ 0 & \text{ else.} \end{cases}$$
Indeed, $\Map_{\DTMr(S)}(\Z(m), \Map_{I/A}(\Delrhoexp_\mu(\Z), \nabrhoexp_\lambda (\Z))) = \Map_I(\Delrhoexp_\mu(\Z(m)), \nabrhoexp_\lambda (\Z))$ is zero if $\lambda \ne \mu$ or if $\lambda = \mu$ and $m \ne 0$, by \thref{exponential homs}. For $\lambda = \mu$ and $m = 0$, it is $\Z$, by the same lemma. 

According to \cite[Lemma~2.7]{RicharzScholbach:Categorical}, the $A$-atomic objects in these categories $C$ are precisely the compact objects, and we may use Proposition~2.9 there to compute
$$\Map_{I \t_{\DTMr(S)} \DTM(S)}(\Delrhoexp_\mu(\Z)\otimes \calF, \nabrhoexp_\lambda(\Z) \otimes \calG) = \Map_{I/A}(\Delrhoexp_\mu(\Z), \nabrhoexp_\lambda (\Z)) \otimes \Map_{\DTM(S) / A}(\calF, \calG).$$
This is $\Map_{\DTM(S) / A}(\calF, \calG)$ if $\lambda = \mu$ and zero otherwise.
Applying $\Map_A(\Z,-)$ then confirms that \eqref{eq-mapping spectra} holds, concluding the proof.
\xpf

\subsection{Hecke modules}

To motivate this subsection, let us briefly change the setup, and assume \(G\) is a (split) reductive group over a nonarchimedean local field \(F\) of residue characteristic \(p\).
Then the geometric Casselman--Shalika formula \cite{FrenkelGaitsgoryVilonen:Whittaker,NgoPolo:Resolutions} is a geometrization of the classical Casselman--Shalika formula \cite{CasselmanShalika:Unramified}, which says that a certain Whittaker module is a free rank one module under the spherical Hecke algebra. 
This Whittaker module requires the choice of a nontrivial additive character of \(U(F)\),  where \(U\) is the unipotent radical of a Borel.
In particular, this requires sufficiently many roots of unity, and does not work with \(\Fp\)-coefficients.
Using the machinery developed in this paper, we can give a different construction of a module under the spherical Hecke algebra, which recovers the classical Whittaker module after fixing a choice of generic additive character.
This construction will be canonical (only depending on a pinning of \(G\)), and works for Hecke algebras and modules with \(\Z\)-coefficients (and hence also for \(\Fp\)-coefficients).
In fact, this construction will also work in the setting of generic Hecke algebras (cf.~e.g.~\cite{Vigneras:Algebres,PepinSchmidt:Generic}), in which the residue cardinality \(q\) is replaced by a formal parameter \(\mathbf{q}\), continuing the study from \cite{CassvdHScholbach:MotivicSatake,CassvdHScholbach:Central}.

Now, we go back to the usual setup, i.e., let \(G\) be defined over a suitable base \(S\), and we assume that \(G\) is not a torus.
Throughout this section, we will take Grothendieck groups of categories of (locally) compact stratified Tate motives on the affine Grassmannian.
Since these Grothendieck groups cannot differ between regular and reduced motives \cite[§2.2.6]{CassvdHScholbach:Central}, we may work with reduced motives, without loss of generality.
By \thref{exponential convolution} the convolution product
\[\DTM_{\red,\rhoexp}(\Gr) \times \DTM_{\red,L^+G}(\Gr) \to \DTM_{\red,\rhoexp}(\Gr)\]
from \thref{twisted exponential convolution} restricts to a functor
\begin{equation}\label{conv for module}
	\DTM_{\red,\rhoexp}^{\anti}(\Gr)^{\comp} \times \DTM_{\red,L^+G}^\anti(\Gr)^{\locc} \to \DTM_{\red,\rhoexp}^\anti(\Gr)^{\comp}
\end{equation}
(Here, \(\DTM_{\red,L^+G}(\Gr)^{\locc}\) denotes the category of locally compact objects, which are compact after forgetting the equivariance \cite[Definition 5.52]{CassvdHScholbach:MotivicSatake}.
This is relevant since the (co)standard functors \(\calJ_!^{\mu}\) and \(\calJ_*^{\lambda}\) do not preserve not compact objects, but they do after forgetting the equivariance.
When using exponential motives, the (co)standard functors \(\Delrhoexp_{\mu}\) and \(\nabrhoexp_{\lambda}\) do preserve compact objects by \thref{rema--CompactExp}.)
Recall also the generic spherical Hecke algebra \(\calH^{\sph}(\mathbf{q})\) from \cite[Example 6.2]{CassvdHScholbach:Central}, defined as the Grothendieck ring of \(\DTM_{\red,L^+G}^{\anti}(\Gr)^{\locc}\).

\defi\thlabel{genericExpMod}
The \emph{generic exponential module} is the Grothendieck group \[M_{\exp}(\mathbf{q}):=K_0(\DTM_{\red,\rhoexp}^{\anti}(\Gr)^{\comp})\cong K_0(\MTM_{\red,\rhoexp}^{\anti}(\Gr)^{\comp}),\]
where the equivalence follows from \eqref{derived category of reduced exponential}.
Using \eqref{conv for module}, it is naturally a (right) module under \(\calH^{\sph}(\mathbf{q})\).
\xdefi

This definition, combined with \thref{Equivalence in the reduced case}, immediately gives the following result. 

\prop
The generic exponential module \(M_{\exp}(\mathbf{q})\) is a free \(\calH^{\sph}(\mathbf{q})\)-module of rank 1, with a basis given by the class \([\Delrhoexp_{0}(\Z)] = [\nabrhoexp_{0}(\Z)]\).
\xprop

We can also describe \(M_{\exp}(\mathbf{q})\) more concretely, or at least its specializations.
For a prime power \(q\), we will write \(F:=\Fq\rpot{t}\) and \(\mathcal{O}:=\Fq\pot{t}\), to simplify the notations throughout this subsection.
Let \[C_c(\calU_{\hat{\rho}}^{\exp}(\Fq) \backslash G(F)/G(\mathcal{O}),\Z)\] be the space of compactly supported \(\calU_{\hat{\rho}}^{\exp}(\Fq) \times G(\mathcal{O})\)-bi-equivariant functions \(G(F) \to \Z\), which is naturally a module under the spherical Hecke algebra \(\calH^{\sph}:=\calH^{\sph}(q):=\calH^{\sph}(\mathbf{q}) \otimes_{\Z[\mathbf{q}], \mathbf{q}\mapsto q} \Z\).
Consider the submodule \(C_c\left((\calU_{\hat{\rho}}\rtimes \Gm)(\Fq) \backslash G(F)/G(\mathcal{O}),\Z\right)\) of functions \(G(F)/ G(\mathcal{O})\to \Z\) which are equivariant under \(\calU_{\hat{\rho}}(\Fq) \supseteq \calU_{\hat{\rho}}^{\exp}(\Fq)\).
Since we are just working with functions, this is isomorphic to \(C_c\left(\calU_{\hat{\rho}}(\Fq)\backslash G(F)/G(\mathcal{O}),\Z\right)\), or even to \(C_c\left(\calI_{\hat{\rho}}(\Fq)\backslash G(F)/G(\mathcal{O}),\Z\right)\).

\prop\thlabel{prop exponential module}
For any prime power \(q\), there is a natural isomorphism 
\[M_{\exp}(\mathbf{q}) \otimes_{\Z[\mathbf{q}], \mathbf{q}\mapsto q} \Z\cong \frac{C_c\left(\calU_{\hat{\rho}}^{\exp}(\Fq) \backslash G(F)/G(\mathcal{O}),\Z\right)}{C_c\left(\calU_{\hat{\rho}}(\Fq)\backslash G(F)/G(\mathcal{O}),\Z\right)}\]
of \(\calH^{\sph}\)-modules.
\xprop
\begin{proof}
	A similar argument as \cite[Proposition 6.4]{CassvdHScholbach:Central} shows that there are isomorphisms
	\[K_0(\DTMr^{\anti}(\calU_{\hat{\rho}}^{\exp}\backslash \Gr)^{\locc}) \otimes_{\Z[\mathbf{q}], \mathbf{q}\mapsto q} \Z \cong C_c(\calU_{\hat{\rho}}^{\exp}(\Fq) \backslash G(F)/G(\mathcal{O}),\Z)\]
	and
	\[K_0(\DTMr^\anti((\calU_{\hat{\rho}}\rtimes \Gm)\backslash \Gr)^{\locc}) \otimes_{\Z[\mathbf{q}], \mathbf{q}\mapsto q} \Z \cong C_c\left((\calU_{\hat{\rho}}\rtimes \Gm)(\Fq) \backslash G(F)/G(\mathcal{O}),\Z\right)\]
	of \(\calH^{\sph}\)-modules.
	Since \(\DTM_{\red,\rhoexp}^{\anti}(\Gr)^{\comp}\) is the quotient of the two categories appearing above, the proposition follows from the long exact sequence of K-groups, and the fact that the quotient splits by \thref{(co)localization nonsense}.
\end{proof}

Finally, we would like to compare these specializations to the usual Whittaker module, involving functions on \(G(F)/G(\mathcal{O})\) that are equivariant under the action of \(U(F)\).
The key is the following lemma, analogous to \cite[Lemma 2]{ArkhipovBezrukavnikov:Perverse}.

Let \(\Lambda\) be a \(\Z[\frac{1}{p}]\)-algebra admitting a nontrivial additive character \(\psi\colon \Fq\to \Lambda^\times\), and consider the induced maps
\begin{equation}\label{whittaker map}
	\Psi_U\colon U(F) \to \prod_{\alpha \in \Delta} U_\alpha(F) \xrightarrow{\sum_{i\gg -\infty} a_it^i \mapsto a_{-1}} \prod_{\alpha \in \Delta}U_{\alpha}(\Fq) \xrightarrow{\sum_\alpha u_\alpha} \Fq \xrightarrow{\psi} \Lambda^\times,
\end{equation}
and
\[\Psi_{\calU_{\hat{\rho}}} \colon \calU_{\hat{\rho}}(\Fq) \cong \calU(\Fq) \to U(\Fq) \to \prod_{\alpha \in \Delta} U_{\alpha}(\Fq) \xrightarrow{\sum_\alpha u_\alpha} \Fq \xrightarrow{\psi} \Lambda^\times,\]
where the isomorphism \(\calU_{\hat{\rho}}\cong \calU\) is given by \(x\mapsto \hat{\rho}(t)x\hat{\rho}(t^{-1})\) (as in §\ref{subsec--twisted exponential motives}).

\lemm\thlabel{combinatorics of Whittaker}
Let \(\lambda\in X_*(T)\).
Then the following conditions are equivalent:
\begin{enumerate}
	\item \(\lambda \in X_*(T)^+\) is dominant.
	\item \(\Stab_{U(F)}(\lambda(t) G(\mathcal{O})) \subseteq \Ker(\Psi_U)\).
	\item \(\Stab_{\calU_{\hat{\rho}}(\Fq)} (\lambda(t) G(\mathcal{O})) \subseteq \Ker(\Psi_{\calU_{\hat{\rho}}})\).
\end{enumerate}
Moreover, if these conditions are satisfied, the \(\calU_{\hat{\rho}}(\Fq)\)-orbit through \(\lambda(t) G(\mathcal{O})\) in \(G(F)/G(\mathcal{O})\) is contained in a single \(U(F)\)-orbit, i.e., 
\[\calU_{\hat{\rho}}(\Fq)\lambda(t) G(\mathcal{O}) \subseteq U(F) \lambda(t) G(\mathcal{O}).\]
\xlemm
\pf
Let \(\lambda\in X_*(T)\).
Then the stabilizer \[\Stab_{U(F)}(\lambda(t) G(\mathcal{O})) = U(F)\cap \lambda(t) G(\mathcal{O}) \lambda(t^{-1})\] is generated by the images of the affine root groups \(U_{\alpha+n}(\Fq)\), for \(\alpha\in \Phi^+\) a positive root, and \(n\geq \langle \alpha,\lambda\rangle\).
On the other hand, if \(\alpha\in \Delta\) is a simple root, the affine root group \(U_{\alpha+n}\) lies in the kernel of \(\Psi_U\) if and only if \(n\neq -1\).
Thus, we see that (1) \(\Leftrightarrow\) (2).

Next, the stabilizer \(\Stab_{\calU(\Fq)}(\lambda(t) G(\mathcal{O}))\) is generated by the affine root groups \(U_{\alpha+n}\) for \(n\geq \max(0,\langle \alpha,\lambda\rangle)\) if \(\alpha\in \Phi^+\) is a positive root, and for \(n\geq \max(1,\langle \alpha,\lambda\rangle)\) if \(\alpha\) is a negative root.
Moreover, if \(\alpha\in \Delta\) is a simple root, the kernel \(\Ker(\Psi_{\calU})\) of \(\calU(\Fq) \to U(\Fq) \to \Fq \to \Lambda\) contains \(U_{\alpha+n}\) exactly when \(n\neq 0\).
Thus, \(\Stab_{\calU(\Fq)}(\lambda(t) G(\mathcal{O})) \subseteq \Ker(\Psi_{\calU})\) exactly when \(\lambda\) is strictly dominant.
By an argument similar to \thref{twisted exponential orbits in Gr}, we can deduce that \(\Stab_{\calU_{\hat{\rho}}(\Fq)}(\lambda(t) G(\mathcal{O})) \subseteq \Ker(\Psi_{\calU_{\hat{\rho}}})\) if and only if \(\lambda\) is dominant, so that (1)\(\Leftrightarrow\)(3).

Finally, the fact that the \(\calU_{\hat{\rho}}\)-orbit through \(\lambda(t) G(\mathcal{O})\) is contained in a unique \(U(F)\)-orbit can be checked after passing to the adjoint group \(G_{\adj}\).
In that case, as explained in \cite[Lemma 2.2]{NgoPolo:Resolutions}, the action on \(t^{\lambda+\hat{\rho}}\) induces an isomorphism
\begin{equation}\label{inclusion of orbits}
	\prod_{\alpha\in \Phi^+} \prod_{j=0}^{\langle\alpha, \lambda+\hat{\rho}\rangle-1} U_{\alpha,j} \cong \Fl_{\lambda+\hat{\rho}}^{\circ}(\bba_0,\bbf_0).
\end{equation}
This implies that
\[\calU_{\hat{\rho}}(\Fq) \lambda(t) G(\mathcal{O}) = \hat{\rho}(t^{-1}) \calU(\Fq) (\lambda+\hat{\rho})(t) G(\mathcal{O}) \subseteq U(F) \lambda(t) G(\mathcal{O}),\]
as desired. 
\xpf

Let us now denote by \(C_c\left((U(F),\Psi_U) \backslash G(F)/G(\mathcal{O}),\Lambda\right)\) the set of functions \(f\colon G(F)/G(\mathcal{O}) \to \Lambda\) such that \(f(ug) = \Psi_U(u)f(g)\) for \(u\in U(F)\) and \(g\in G(F)\), for which \(f(\lambda(t))\neq 0\) only for finitely many \(\lambda\in X_*(T)\).
We also define \(C_c\left((\calU_{\hat{\rho}}(\Fq),\Psi_{\calU_{\hat{\rho}}}) \backslash G(F)/G(\mathcal{O}),\Lambda\right)\) as the set of compactly supported functions \(G(F)/G(\mathcal{O})\to \Lambda\), satisfying a similar equivariance condition.
Both of these are modules under the spherical Hecke algebra (with \(\Lambda\)-coefficients.)

\coro\thlabel{Whittaker vs baby}
There is a natural \(\mathcal{H}^{\sph}\)-equivariant isomorphism
\[C_c\left((U(F),\Psi_U) \backslash G(F)/G(\mathcal{O}),\Lambda\right) \cong C_c\left((\calU_{\hat{\rho}}(\Fq),\Psi_{\calU_{\hat{\rho}}}) \backslash G(F)/G(F),\Lambda\right).\]
\xcoro
\pf
Consider the map
\begin{equation}\label{averaging}
	C_c\left((U(F),\Psi_U) \backslash G(F)/G(\mathcal{O}),\Lambda\right) \to C_c\left((\calU_{\hat{\rho}}(\Fq),\Psi_{\calU_{\hat{\rho}}}) \backslash G(F)/G(F),\Lambda\right),
\end{equation}
given by sending a function \(f\) to the averaged function
\[\av(f) \colon G(F)\to \Lambda\colon g\mapsto \frac{1}{|\calU_{\hat{\rho}}(\Fq)|} \sum_{x\in \calU_{\hat{\rho}}(\Fq)} \Psi_{\calU_{\hat{\rho}}}^{-1}(x) f(xg).\]
(Strictly speaking, this does not make sense, since \(\calU_{\hat{\rho}}(\Fq)\) is an infinite group.
But on each bounded part \(\Gr^{\leq \mu}\subseteq \Gr\), the \(\calU_{\hat{\rho}}\)-action factors through a finite type quotient, with split pro-unipotent kernel.
On the corresponding parts of \(G(F)\), we may then replace the above sum by a finite one, and the result is independent of the choice of finite type quotient over which the \(\calI_{\hat{\rho}}\)-action factors.
Thus, this procedure extends to all of \(G(F)\), defining \eqref{averaging}.)

Then \eqref{averaging} is clearly \(\mathcal{H}^{\sph}\)-equivariant.
Moreover, for \(\lambda\in X_*(T)^+\), it sends a function supported on \(U(F)\lambda(t)G(\mathcal{O})\) to a function supported on \(\calU_{\hat{\rho}}(\Fq)\lambda(t)G(\mathcal{O})\) by \thref{combinatorics of Whittaker}.
To see that \eqref{averaging} is an isomorphism, it remains to see that for a function \(f\) supported on some \(U(F)\lambda(t)G(\mathcal{O})\), we have \(f(\lambda(t)) = \av(f)(\lambda(t))\).

For this, note that \(\Psi_U\) and \(\Psi_{\calU_{\hat{\rho}}}\) agree on the intersection \(U(F) \cap \calU_{\hat{\rho}}(\Fq)\).
Then for each \(x\in \calU_{\hat{\rho}}(\Fq)\), we may find \(x'\in U(F)\) such that \(x\lambda(t)G(\mathcal{O}) = x'\lambda(t)G(\mathcal{O})\) and for which the image of \(x\) under \(\calU_{\hat{\rho}}(\Fq) \cong \calU(\Fq) \to U(\Fq)\) agrees with the image of \(x'\) under \(U(F) \xrightarrow{\sum_{i\gg -\infty} a_it^i\mapsto a_{-1}} U(\Fq)\); this follows from \eqref{inclusion of orbits}.
Then we have \(\Psi_{\calU_{\hat{\rho}}}(x) = \Psi_U(x')\), and we compute
\[\av(f)(\lambda(t)) = \frac{1}{|\calU_{\hat{\rho}}(\Fq)|} \sum_{x\in \calU_{\hat{\rho}}(\Fq)}\Psi_{\calU_{\hat{\rho}}}^{-1}(x)f(x\lambda(t)) = \frac{1}{|\calU_{\hat{\rho}}(\Fq)|} \sum_{x\in \calU_{\hat{\rho}}(\Fq)}\Psi_{\calU_{\hat{\rho}}}^{-1}(x)f(x'\lambda(t)) \]
\[= \frac{1}{|\calU_{\hat{\rho}}(\Fq)|} \sum_{x\in \calU_{\hat{\rho}}(\Fq)}\Psi_{\calU_{\hat{\rho}}}^{-1}(x)\Psi_U(x')f(\lambda(t)) = f(\lambda(t)).\]
\xpf

Finally, we can show that the specialization of our exponential module \(M_{\exp}(\mathbf{q})\) recovers the usual Whittaker module, after the choice of a nontrivial additive character \(\psi\).

\prop
There is a natural isomorphism (only depending on the choice of \(\psi\))
\[M_{\exp}(\mathbf{q}) \otimes_{\Z[\mathbf{q}], \mathbf{q}\mapsto q} \Lambda \cong C_c\left((U(F),\Psi_U )\backslash G(F)/G(\mathcal{O}),\Lambda\right).\]
\xprop
\pf
Consider the map
\[C_c\left((\calU_{\hat{\rho}}(\Fq),\Psi_{\calU_{\hat{\rho}}}) \backslash G(F)/G(\mathcal{O}),\Lambda\right) \xrightarrow{\gamma} \frac{C_c\left(\calU_{\hat{\rho}}^{\exp}(\Fq)\backslash G(F)/G(\mathcal{O})\right)}{C_c\left(\calU_{\hat{\rho}}(\Fq)\backslash G(F)/G(\mathcal{O})\right)}\]
given by \(\Gm\)-averaging, i.e., which sends a function \(f\) to the function \(x\mapsto \sum_{g\in \Gm(\Fq)} f(gx)\).
It is clearly equivariant for the module structure under the spherical Hecke algebra, and we claim it is an isomorphism.

Since the groups \(\calU_{\hat{\rho}}\) and \(\calU_{\hat{\rho}}^{\exp}\) preserve the \(\calI_{\hat{\rho}}\)-orbits in \(\Gr\), the map \(\gamma\) preserves the set of \(\calI_{\hat{\rho}}\)-orbits on which functions are supported.
Thus, we may replace \(G(F)/G(\mathcal{O}) = \Gr(\Fq)\) by the \(\Fq\)-points of an \(\calI_{\hat{\rho}}\)-orbit \(X\) in \(\Gr\), which is an affine space by \thref{lemm-dimensions of twisted exponential orbits}.
Let \(x\in X(\Fq)\) be the basepoint.

Now, if \(X\) is a \((\calU^{>0}\rtimes U_0)_{\hat{\rho}}\)-orbit, both \(C_c\left((\calU_{\hat{\rho}}(\Fq),\Psi_{\calU_{\hat{\rho}}}) \backslash X(\Fq),\Lambda\right)\) and \(\frac{C_c\left(\calU_{\hat{\rho}}^{\exp}(\Fq)\backslash X(\Fq)\right)}{C_c\left(\calU_{\hat{\rho}}(\Fq)\backslash X(\Fq)\right)}\) are trivial by \thref{twisted exponential orbits in Gr} and \thref{combinatorics of Whittaker} (this uses that \(\psi\) is nontrivial).
So let us assume \(X\) splits up into two \(\calU_{\rho}^{\exp}\)-orbits; denote by \(Y\) the closed orbit (containing \(x\)), and by \(Z\) the open orbit.
Then \(C_c\left((\calU_{\hat{\rho}}(\Fq),\Psi_{\calU_{\hat{\rho}}}) \backslash X(\Fq),\Lambda\right)\) and \(\frac{C_c\left(\calU_{\hat{\rho}}^{\exp}(\Fq)\backslash X(\Fq)\right)}{C_c\left(\calU_{\hat{\rho}}(\Fq)\backslash X(\Fq)\right)}\) are both free \(\Lambda\)-modules of rank 1.
Consider the basis function \(f_1\) of the former, determined by the value \(1\) at \(x\).
Its image under \(\gamma\) sends any \(y\in Y(\Fq)\) to \(\sum_{g\in \Gm(\Fq)} f_1(x) = q-1\), since \(\Gm\) fixes \(x\), and any \(z\in Z(\Fq)\) to \(\sum_{g\in \Gm(\Fq)} \psi(g) f_1(x) = -1\), since \(\psi\) is a nontrivial additive character.
In the quotient, this function is equivalent to the function supported on \(Y\), with constant value \(q\).
Since \(q\) is invertible in \(\Lambda\), this is clearly a basis function, so that \(\gamma\) is an isomorphism.

Combining this with \thref{prop exponential module} and \thref{Whittaker vs baby} then yields the desired isomorphism.
\xpf

\bibliographystyle{alphaurl}
\bibliography{bib}
					
\end{document}